\documentclass[12pt]{article} 
\usepackage{fullpage,amssymb}

\makeatletter\@addtoreset{equation}{section}\makeatother
\makeatletter\@addtoreset{figure}{section}\makeatother
\makeatletter\@addtoreset{table}{section}\makeatother

\newtheorem{theorem}{Theorem}[section]
\newtheorem{proposition}[theorem]{Proposition}
\newtheorem{lemma}[theorem]{Lemma}

\newtheorem{corollary}[theorem]{Corollary}

\newcommand{\R}{{\bf R}}
\newcommand{\C}{{\bf C}}
\newcommand{\Z}{{\bf Z}}

\newcommand{\op}[1]{\!\!\mathop{\rm ~#1}\nolimits}
\newcommand{\fop}[1]{\!\!\mathop{\mbox{\rm \footnotesize ~#1}}\nolimits}
\newcommand{\scriptop}[1]{\!\!\mathop{\mbox{\rm \scriptsize ~#1}}\nolimits}
\newcommand{\tinyop}[1]{\!\!\mathop{\mbox{\rm \tiny ~#1}}\nolimits}
\newcommand{\dd}[2]{\frac{\fop{d}\! #1}{\fop{d}\! #2}}

\newenvironment{proof}{\par\medskip\noindent{\bf Proof}~~}{
\unskip\nobreak\hfill\hbox{ $\Box$}\par \bigskip}        

\newcounter{exerc}[section]
\renewcommand{\theexerc}{\thesection.\arabic{exerc}}

\newenvironment{remark}{\refstepcounter{theorem}\par\medskip\noindent{\bf Remark~\thetheorem~~}}{\unskip\nobreak\hfill\hbox{ $\oslash$}\par\bigskip}

\newenvironment{question}{\refstepcounter{theorem}\par\medskip\noindent{\bf Question~\thetheorem~~}}{\unskip\nobreak\hfill\hbox{ $\oslash$}\par\bigskip}

\newfont{\gothic}{eufm10 scaled\magstep0}
\title{Chaplygin's Sphere}
\author{J.J. Duistermaat\thanks{A large part of this work has been 
done during a sabbatical leave in Berkeley, in the fall of 1994, 
partially supported by AFOSR Contract AFO F 49629-92. A more recent version 
was prepared in July 2000, as a chapter in the planned book 
``The geometry of Nonholonomically Constrained Systems'', 
together with  R.H. Cushman and J. \'Sniaticky.}}
\date{}

\begin{document}
\maketitle

\begin{center}
{\bf\large Abstract}
\end{center}

Chaplygin \cite{chaplsphere} proved the integrability by quadratures 
of a round sphere, rolling without slipping on a horizontal plane, 
with center of mass at the center of the sphere, but with arbitrary 
moments of inertia. Although the system is integrable in every 
sense of the word, it neither arises as a Hamiltonian system, 
nor is the integrability an immediate consequence of the symmetries. 
On the other hand, the constants of motion are obtained as a consequence 
of Noether's principle, cf. Section \ref{noethersec} and \ref{noethercs}. 
The system also turns out to be related to a Hamiltonian system, 
the geodesic flow on the Euclidean motion group for a left invariant metric, 
cf. Subsection \ref{geodsubsec}. 

In this paper we analyse the global dynamics of Chaplygin's sphere. 
In the process we will explain almost all of Chaplygin's results. 
Readers who are mainly interested in these may consult Sections 
\ref{noethersec}, \ref{noethercs}, \ref{eqmotsec}, \ref{volsec}, 
Subsections \ref{redhorsec}, \ref{ellcoordsubsec}, \ref{velellsubsec}, 
\ref{phorsubsec}, and Section \ref{geomsec}. These can be read 
independently from the rest of the paper with some exceptions, 
such as Subsection \ref{redhorsec} in which Subsection \ref{phi0subsec} 
has been used. At the end of each section we describe in a 
subsection "Chaplygin" the relation between our text and Chaplygin's. 

We also obtain some new results, such as the proof in Section 
\ref{levelsec} that the level sets of the constants of motion 
in the reduced phase space are 
two\--dimensional tori. 
In Section \ref{commvvsec} we prove that, after a 
suitable time reparametrization, the rotational motion is 
quasi\--periodic on two\--dimensional tori. After suitable 
completion of the level surfaces, this is also true for the 
complexified system, cf. Section \ref{ccsec}.  This shows that the 
rotational motion is algebraically integrable according to the 
definition of Adler and van Moerbeke. In Subsection \ref{jacsubsec} 
it is explained how this also follows, in a quite different way, 
from Chaplygin's integration in terms of hyperelliptic integrals.

\section{Noether's Principle for Nonholomic Systems}
\label{noethersec}
We use the equations of motion for a system with nonholonomic 
constraints as given by d'Alembert's principle 
\begin{equation}
\langle [L]^{\gamma}(t),\, v\rangle =0\;\mbox{\rm for every}\; 
v\in C_{\gamma (t)},
\label{dalembert}
\end{equation}
as described in \cite[Ch. 1, Sec. 2.5]{Arnold}. We may assume that $Q_0=Q$. 
Let $s\mapsto\gamma_s$ be a smooth family of smooth 
curves in $Q$, with $\gamma =\gamma _0$, 
$\delta (t):=\partial\gamma _s(t)/\partial s|_{s=0}$, 
for which we this time do not assume that 
$\delta (a)=0$ and $\delta (b)=0$. 
Write 
\[
j^{\delta}(t):=\sum_i\,\frac{\partial L(\gamma (t),\, v)}{\partial v^i}
|_{_{v=\gamma'(t)}}\,\delta ^i(t)
\]
for the {\em $\delta$\--component of the momentum}, 
which is a coordinate\--invariant quantity.
The classical variational 
equation, which is obtained by a partial integration, reads 
\begin{equation}
\frac{\fop{d}}{\fop{d}\! s}\,
\int_a^b\, L\left(\gamma _s(t),\, 
\gamma '_s(t)\right)\,\op{d}\! t\; |_{_{s=0}}
=-\int_a^b\, \langle [L]^{\gamma}(t),\,\delta (t)\rangle\,\op{d}\! t 
+j^{\delta}(b)-j^{\delta}(a).
\label{[L]var}
\end{equation}
If one differentiates (\ref{[L]var}) with respect to $b$, then 
one obtains the equivalent form 
\begin{equation}
\frac{\fop{d}}{\fop{d}\! t}j^{\delta}(t)
=\frac{\fop{d}}{\fop{d}\! s}\,
L\left(\gamma _s(t),\, 
\gamma '_s(t)\right) |_{_{s=0}}
+\langle [L]^{\gamma}(t),\,\delta (t)\rangle .
\label{[L]var'}
\end{equation}
The second term in the right hand side of (\ref{[L]var'}) 
can be viewed as the {\em $\delta$\--component of the 
reaction force acting on the system}. It is equal to zero 
if $\gamma$ satisfies d'Alembert's principle (\ref{dalembert}) 
and $\delta (t)\in C_{\gamma (t)}$, which means in words that 
$\delta (t)$ is a virtual displacement. 

The form (\ref{[L]var'}) of the variational equations is due to 
Emmy Noether \cite{Noether}, in a version where the independent 
variable $t$ is replaced by a finite number of real variables 
and $L$ is a smooth function on a jet bundle of arbitrary order. 
In the case that the first term in the right hand side 
of (\ref{[L]var'}) is equal to zero, the formula (\ref{[L]var'}), 
is the ``momentum equation'' of Bloch e.a. \cite[Thm. 4.5]{BKMM}. 

Suppose that $w$ is a smooth vector field on $Q$ such that 
$\frac{\fop{d}}{\fop{d}\! s}\,\gamma _s(t)|_{_{s=0}}=w(\gamma (t))$. 
Let $\widehat{w}$ denote the lift of $w$ to $\op{T}Q$, the 
vector field $\widehat{w}$ on $\op{T}Q$ such that 
\begin{equation}
\op{e}^{s\,\widehat{w}}=\op{T}\left(\op{e}^{s\, w}\right) ,
\label{lift}
\end{equation}
if $\op{e}^{s\, w}$ denotes the flow after time $s$ of the vector 
field $w$. In local coordinates $\widehat{w}$ is given by 
\begin{equation}
\widehat{w}(x,\, v)=\left( w(x),\,\op{D}\! w(x)\cdot v\right) ,
\quad (x,\, v)\in\op{T}Q,
\label{lift'}
\end{equation}
where $\op{D}\! w(x)$ denotes the matrix $\partial w^i(x)/\partial x^j$. 
With this notation, the first term in the right hand side of 
(\ref{[L]var'}) is equal to the derivative of $L$ at 
$(\gamma (t),\gamma '(t))$ in the 
direction of $\widehat{w}$. This leads to 
the following version of Noether's principle for variational 
systems with nonholonomic constraints. 
\begin{lemma}
Let $L$ be a smooth function on $\op{T}Q$, of which $C$ is 
a smooth vector subbundle. Let $w$ be a smooth vector field 
on $Q$ with the following properties
\begin{itemize}
\item[{\em i)}] $w$ is a virtual diplacement, which means 
that $w$ is a section of $C$.
\item[{\em ii)}] At each point of $C$, the derivative of $L$ in the 
direction of $\widehat{w}$ is equal to zero, where $\widehat{w}$ 
is the lift of $w$ to $\op{T}Q$ as defined by {\em (\ref{lift})}.
\end{itemize}
Then the $w$\--component of the momentum is constant along every 
solution of {\em (\ref{dalembert})}. 
\label{noetherlem}
\end{lemma}
If there are no constraints, when $C=\op{T}Q$, then condition i) is 
void and ii) is equivalent to the condition that $L$ is invariant 
under the flow of the vector field $\widehat{w}$ in $\op{T}Q$, 
which is equal to the tangent lift of the flow of $w$ in $Q$. 
In this case Lemma \ref{noetherlem} is due to Emmy Noether 
\cite{Noether}. 

\begin{question}Can all the constants of motion 
in Chaplygin \cite{chaplarea} be obtained as applications of 
Lemma \ref{noetherlem}? 
\end{question}

\begin{remark}
Lemma \ref{liftlem} below leads to the warning that in the 
non\--integrable case the condition ii), under the assumption that 
i) holds, is not a property of only the restriction of $L$ to $C$, 
because at the points of $C$ the vector field $\widehat{w}$ need not 
be tangent to $C$.
\end{remark} 
\begin{lemma}
For each section $w:Q\to C$ of $C$ the vector field 
$\widehat{w}$ is tangent to $C$ if and only if the subbundle 
$C$ of $\op{T}Q$ is integrable.
\label{liftlem}
\end{lemma} 
\begin{proof}
Let $\psi ^s$ be the flow of $w$. Then 
the condition that the lift of $w$ is tangent to 
$C$ is equivalent to the condition that 
the mappings $\op{T}\psi ^s$ leave $C$ invariant, 
or that these mappings send sections of $C$ to 
sections of $C$. This in turn is equivalent to the 
condition that $[w, u]$ is a section of $C$ for 
every section $u$ of $C$. That this holds for every 
section $w$ of $C$ is one of the equivalent Frobenius 
conditions for the integrability of $C$. 
\end{proof} 

\subsection{Chaplygin}
The version of 
Lemma \ref{noetherlem} with nonholonomic constraints 
can be found in Arnol'd \cite[p. 82]{Arnold}, with condition ii) 
replaced by the somewhat stronger condition that $L$ is 
$\widehat{w}$\--invariant. Two applications have been given 
in Arnol'd \cite[p. 83, 84]{Arnold}, the first with a reference to 
Chaplygin \cite{chaplarea} and the second with a reference 
to Chaplygin \cite{chaplsphere}. In Chaplygin \cite{chaplsphere} 
the constants of motion have been described as an application of 
\cite{chaplarea}. 

\section{Noether's Principle for Chaplygin's Sphere}
\label{noethercs}
The position of a rigid body 
is given by a pair $(A,\, a)$, with $A\in\op{SO}(3)$ and 
$a\in\R ^3$. Thus, if $x\in\R ^3$ is the position of 
a material point of the body in its reference position, 
then $y=A\, x+a$ is the position of the corresponding point 
in the moving body. If $s$ is the position on the surface 
$S$ of the body in the reference position, such that 
$p=A\, s+a$ is the point of contact of the moving body 
with the surface $P$ on which the body is rolling, then the 
condition of rolling without slipping means that 
\begin{equation}
\dot{A}\, s+\dot{a}=0,
\label{constraint}
\end{equation}
meaning that the at the point of contact the corresponding material 
point of the body is at rest.  
Correspondingly, $(\widetilde{A},\,\widetilde{a})$ is a virtual 
displacement if and only if 
\begin{equation}
\widetilde{A}\, s+\widetilde{a}=0. 
\label{virtual}
\end{equation}

If $\mu$ denotes the mass distribution of the body in the 
reference position, which is a finite Borel measure on $\R ^3$, 
then the kinetic energy of the moving body is given by 
\begin{equation}
T=\int_{_{\R ^3}}\, \frac{1}{2}\langle \dot{A}\, x+\dot{a},
\,\dot{A}\, x+\dot{a}\rangle\,\mu(\op{d}\! x).
\label{T}
\end{equation}
It follows that the $(\widetilde{A},\,\widetilde{a})$\--component 
of the momentum is equal to 
\begin{equation}
j^{(\widetilde{A},\,\widetilde{a})}=\int_{_{\R ^3}}\, \langle \dot{A}\, x+\dot{a},
\,\widetilde{A}\, x+\widetilde{a}\rangle\,\mu(\op{d}\! x).
\label{I}
\end{equation}
Let $\nu\in\R ^3$ be the unique vector such that 
\begin{equation}
\widetilde{A}\, z=\nu\times (A\, z),\quad z\in\R ^3 
\label{nu}
\end{equation}
--- note that this  
corresponds to the {\em right} trivialization of the tangent bundle. 
If the condition (\ref{virtual}) holds, meaning that 
$(\widetilde{A},\,\widetilde{a})$ is a virtual displacement, then 
\[
j^{(\widetilde{A},\,\widetilde{a})}=\int_{_{\R ^3}}\, \langle \dot{y},
\,\widetilde{A}\, (x-s)\rangle\,\mu(\op{d}\! x)
=\langle j,\, \nu\rangle ,
\]
in which
\begin{equation}
j:=\int_{_{\R ^3}}\,\mu(\op{d}\! x)
\,  (y-p)\times\dot{y}
\label{J}
\end{equation}
is the {\em moment of momentum about the point of contact} $p$. 
Here we have used that $\dot{A}\, x+\dot{a}=\dot{y}$, 
$A\, (x-s)=(y-a)-(p-a)=y-p$, and 
$\langle\dot{y},\,\nu\times (y-p)\rangle 
=\langle (y-p)\times\dot{y},\,\nu\rangle$. 

We now turn to the case of Chaplygin's sphere \cite{chaplsphere}, 
where the surface $S$ of the body in the reference position is a 
sphere, the center of mass is at the center of $S$, and the 
body is rolling without slipping on a horizontal plane $P$. 
We will take the origin of the reference frame at the center of mass 
= the center of $S$. If $r$ denotes the radius of $S$, and we take the 
plane $P$ at height $-r$, then the condition that 
$A(S)+a$ is lying on top of $P$ corresponds to the condition that 
the third (vertical) component of $a$ is equal to zero. 
The point of contact then is equal to $p=a-r\, e_3$ if $e_3$ 
denotes the third standard basis vector, and the corresponding 
point on $S$, in body coordinates, is equal to 
\begin{equation}
s=-r\, A^{-1}\, e_3.
\label{s}
\end{equation}  
The condition (\ref{virtual}) therefore is equivalent to 
\begin{equation}
\widetilde{a}=r\nu\times e_3,
\label{virtchap}
\end{equation}
where we have also used (\ref{nu}). 

From this moment on, we keep the infinitesimal 
rotation vector $\nu\in\R ^3$ constant. Then (\ref{virtchap}) 
implies that $\widetilde{a}$ is a constant horizontal vector. 
If we use the {\em left} trivialization of the tangent bundle 
of $\op{SO}(3)$, corresponding to assigning to 
$\dot{A}\in\op{T}_A\op{SO}(3)$ the infinitesimal rotation 
given by 
\begin{equation}
\dot{A}\, z=A\, (\omega\times z),\quad z\in\R ^3 
\label{dotA}
\end{equation}
for some vector $\omega\in\R ^3$, then the tangent lift 
of the vector field $\widetilde{A}$ such that (\ref{nu}) with a 
constant $\nu$ does not effect $\omega$. Because (\ref{virtchap}) 
implies that the tangent lift of $\widetilde{a}$ does not effect 
$\dot{a}$ either, the conclusion is that condition ii) of 
Lemma \ref{noetherlem} holds if we take $L=T$. Note that 
for Chaplygin's sphere the center of mass remains at the 
same height, which means that the gravitational potential energy 
is constant, and therefore can be disregarded. 
We have arrived at the conclusion that 
\begin{proposition}
For Chaplygin's sphere, the moment of momentum about 
the point of contact is a constant of motion.
\label{jprop}
\end{proposition}

The kinetic energy of the rigid body is 
given by 
\begin{equation}
T=\frac12\,\langle I\,\omega ,\,\omega\rangle +
\frac12m\,\langle\dot{a},\,\dot{a}\rangle ,
\label{Tomega}
\end{equation}
where we have used the left trivialization (\ref{dotA}) of the 
tangent bundle of $\op{SO}(3)$. 
Here $I$ denotes the moment of inertia tensor, which is given 
by a positive definite symmetric matrix, and $m$ denotes the 
total mass of the body. 

In this notation $\widetilde{A}$, given by 
(\ref{nu}), corresponds to 
\begin{equation}
\widetilde{\omega} =A^{-1}\,\nu ,
\label{delomega}
\end{equation}
because 
\[
\widetilde{A}\, z=\nu\times A\, z=A\,\left( A^{-1}\nu\times z\right) .
\]
It follows that 
\begin{eqnarray*}
\langle j,\,\nu\rangle =j^{(\widetilde{A},\,\widetilde{a})}
&=&\langle I\,\omega ,\, A^{-1}\nu\rangle 
+m\,\langle\dot{a},\,\widetilde{a}\rangle 
\\
&=&\langle A\, I\,\omega ,\,\nu\rangle 
+mr^2\,\langle A\,\omega\times e_3,\,\nu\times e_3\rangle\\
&=&\langle A\, I\,\omega +mr^2 e_3\times
\left( A\,\omega\times e_3\right),\,\nu\rangle ,
\end{eqnarray*}
where we have used that (\ref{constraint}), (\ref{s}) and (\ref{dotA}) imply that 
\[
\dot{a}=r\, A\left(\omega\times A^{-1}\, e_3\right)
=r\, A\,\omega\times e_3
\]
and $\delta$ is given by (\ref{virtchap}). 
With the notation 
\begin{equation}
u:=A^{-1}\, e_3
\label{u}
\end{equation} 
this leads to the formula 
\begin{equation}
j=A\,\left( I\,\omega +mr^2 u\times\left(\omega\times u\right)\right) 
\label{Jchap}
\end{equation}
for the moment of momentum about the point of contact. 
Note that $\langle u,\, u\rangle =1$ and therefore 
\begin{equation}
u\times (\omega\times u)=\omega- \langle u,\,\omega\rangle\, u ,
\label{uu}
\end{equation}
which is the orthogonal projection of $\omega$ onto the plane 
which is orthogonal to $u$. Also note that $u$ has the 
concrete interpretation that {\em $-r\, u$ is equal to the 
point of contact on the surface of the sphere, in body coordinates}, 
cf. (\ref{s}).   

In order to simplify the notation somewhat, we write 
\begin{equation}
\rho :=m\, r^2,
\label{rho}
\end{equation}
and define the symmetric linear mapping $I_{\rho ,\, u}:\R ^3\to\R ^3$ by 
\begin{equation}
I_{\rho ,\, u}(\omega ):=(I+\rho )\,\omega -\rho\,\langle u,\,\omega\rangle\, u,  
\label{Iu}
\end{equation}
which is equal to $I$ plus $\rho$ times the orthogonal projection 
to the plane orthogonal to $u$. With these notations, we have that 
\begin{equation}
j=A\, I_{\rho,\, u}\,\omega .
\label{JIu}
\end{equation}

Write
\begin{equation}
J:=(I+\rho )^{-1}.
\label{Jdef}
\end{equation}
If $I_{\rho ,\, u}\,\omega =\nu$ then 
$\omega =J\, (\nu +\theta\, u)$ for some $\theta\in\R$, 
which moreover has to satisfy the equation 
\[
\nu =\nu +\theta\, u-\rho\,\langle u,\, 
J\, (\nu +\theta\, u)\rangle 
\, u,
\]
which holds if and only if 
\[
\theta\, \left( 1-\rho\,\langle u,\, J\, u\rangle\right) 
-\rho\,\langle u,\, J\,\nu\rangle =0.
\]
This leads to the conclusion that the symmetric linear 
mapping $I_{\rho ,\, u}$ is invertible, with inverse given by
\begin{equation}
{I_{\rho ,\, u}}^{-1}(\nu )=
J\,\nu +\frac{\rho\,\langle u,\, 
J\,\nu\rangle}
{1-\rho\,\langle u,\, J\, u\rangle}\, 
J\, u. 
\label{Iuinv}
\end{equation}
If we fix the constant of motion $j$, then we can use the equation \begin{equation}
\omega ={I_{\rho,\, u}}^{-1}\, A^{-1}\, j
\label{omegaJIu}
\end{equation}
in order to express $\omega$ in terms of $A$, where we note that 
$u$ is given in terms of $A$ by means of (\ref{u}). In this way $\omega$ 
can be eliminated from the equations of motion. 

\subsection{Chaplygin}
In \cite[\S 1]{chaplsphere}, Proposition \ref{jprop} is stated 
``$\cdots$ as a consequence of a generalized theorem of areas'', 
with a reference to \cite{chaplarea}. 

In \cite[\S 6]{chaplarea} 
the case of 
\cite{chaplsphere} appears as the limit when the radius of the 
big sphere in which the small sphere rolls tends to infinity. 

\section{The Equations and Constants of Motion}
\label{eqmotsec}
\subsection{The Equations of Motion}
\label{eqmotsubsec}

The equations of motion are
\begin{eqnarray}
&&\dd{p}{t}=
r\, A\, (\omega\times u)=r\, (A\,\omega )\times e_3,  \label{pdot}\\
&&\dd{A}{t}
=A\circ\omega _{\scriptop{op}},\label{Adot}\\
&&\dd{u}{t}=u\times\omega ,\quad u(t):=A(t)^{-1}\, e_3,\label{udot}\\
&&\dd{}{t}I\,\omega -I\,\omega\times\omega 
=m\, r^2\left(\langle u,\,\dd{\omega}{t}\rangle\, u-\dd{\omega}{t}\right) .
\label{omegadot}
\end{eqnarray} 
Here, as in Section \ref{noethercs}, 
$x\mapsto A\, x+a$, $x,\, a\in\R ^3$, $A\in\op{SO}(3)$ is 
the rigid motion which is applied to the body in the 
reference position, with the center of mass at the origin, 
and $p:=A\,\widetilde{s}(u) + a=a-r\, e_3$ denotes the point of 
contact between the sphere and the horizontal plane. 
Furthermore $\omega _{\scriptop{op}}$ denotes 
the antisymmetric linear mapping 
$\nu\mapsto \omega\times\nu :\R ^3\to\R ^3$. 
The equation (\ref{Adot}) expresses that $\omega\in\R ^3$ 
can be viewed as a rotational velocity vector. 

The full system (\ref{pdot}), (\ref{Adot}), (\ref{omegadot}) 
is defined in the eight dimensional $(p,\, A,\,\omega )$\--space 
$\R ^2\times\op{SO}(3)\times\R ^3$. 
The equations (\ref{Adot}), (\ref{omegadot}) for 
\[
(A,\,\omega )\in\op{SO}(3)\times\R ^3\simeq \op{T}\left(\op{SO}(3)\right)
\] 
are the equations for the {\em rotational motion}, 
the system obtained by {\em reduction by the horizontal 
translation group}. The equations (\ref{udot}), (\ref{omegadot}) for  
\[
(u,\,\omega )\in \op{S}^2\times\R ^3
\] 
represent the {\em reduction of the system by the left action 
of the horizontal motion group} $\op{E}(2)$.  

With the notations (\ref{rho}) and (\ref{Iu}), 
the equation (\ref{omegadot}) takes the form 
\begin{equation}
I_{\rho ,\, u}\,\dd{\omega}{t}=(I\,\omega )\times\omega . 
\label{Jomegadot}
\end{equation}
Combining (\ref{Iu}) with (\ref{udot}), we obtain that 
\begin{eqnarray*}
\dd{I_{\rho ,\, u}}{t}\,\omega &=&\, -\rho\,\left( 
\langle\omega ,\,\dd{u}{t}\rangle\, u 
+\langle\omega ,\, u\rangle\,\dd{u}{t}\right)\\ 
&=&\, -\rho\,\left( 
\langle\omega ,\, u\times\omega \rangle\, u 
+\langle\omega ,\, u\rangle\, u\times\omega\right) 
=\, -\rho\,\langle\omega ,\, u\rangle\, u\times\omega ,
\end{eqnarray*}
and therefore 
\[
\dd{}{t}\left( I_{\rho ,\, u}\,\omega\right) 
=\dd{I_{\rho ,\, u}}{t}\,\omega +I_{\rho ,\, u}\,\dd{\omega}{t}
=\, -\rho\,\langle\omega ,\, u\rangle\, u\times\omega 
+(I\,\omega )\times\omega ,
\]
hence 
\begin{equation}
\dd{}{t}\left( I_{\rho ,\, u}\,\omega\right) 
=\left( I_{\rho ,\, u}\,\omega\right)\times\omega ,
\label{dotJomega}
\end{equation}
because $\rho\,\omega\times\omega =0$. 

\subsection{The Constants of Motion}
\label{constmotsubsec}
In general the total energy is a constant of motion, when the 
the equations of motions have been obtained as a consequence 
of the principle of d'Alembert, cf. (\ref{dalembert}). 
Because in our case the potential energy $m\, g\, r$ is a constant, 
it follows that the total kinetic energy 
\begin{equation}
T=\frac12\,\langle I\,\omega ,\,\omega\rangle 
+\frac12\, m\, r^2\langle u\times\omega ,\, u\times\omega\rangle 
=\frac12\,\langle I_{\rho ,\, u}\,\omega ,\,\omega\rangle ,
\label{TI}
\end{equation}
cf. (\ref{Iu}), is a constant of 
motion. This can also be verified directly from 
(\ref{Jomegadot}) and 
(\ref{dotJomega}), because 
\[
\langle I_{\rho ,\, u}\,\omega ,\,\dd{\omega}{t}\rangle 
=\langle \omega ,\, I_{\rho ,\, u}\,\dd{\omega}{t}\rangle 
=\langle\omega ,\, (I\,\omega )\times\omega\rangle =0
\]
and therefore also 
\[
\dd{T}{t}=\langle\dd{}{t}(I_{\rho ,\, u}\,\omega ),\,\omega\rangle 
=\langle\left( I_{\rho ,\, u}\,\omega\right)\times\omega ,\,\omega\rangle =0.
\]

On the other hand, combination of (\ref{Adot}) with (\ref{dotJomega}) 
yields that 
\begin{eqnarray*}
\dd{}{t} \left( A\, I_{\rho ,\, u}\,\omega \right) &=&
\dd{A}{t}\, I_{\rho ,\, u}\,\omega 
+A\,\dd{}{t}\left( I_{\rho ,\, u}\,\omega\right)\\
&=&A\,\left(\omega\times\left( I_{\rho ,\, u}\,\omega\right)\right) 
+A\,\left(\left( I_{\rho ,\, u}\,\omega\right) \times\omega\right) =0.
\end{eqnarray*}
In this way we have verified again that the vector 
$A\, I_{\rho ,\, u}\,\omega$, which according to 
(\ref{JIu}) is equal to the moment $j$ of the 
momentum around the point of contact, 
is a constant of motion. 

\subsection{A Pair of Vectors}
\label{pairsubsec}
If $j$ is not vertical, then the rotation $A$ is determined by 
the pair of vectors 
\begin{equation}
u:=A^{-1}\, e_3\quad\mbox{\rm and}\quad 
v :=A^{-1}\, j.
\label{sigmadef}
\end{equation}
More precisely, in this case the mapping $A\mapsto (u,\, v )$ 
is a diffeomorphism from $\op{SO}(3)$ onto the smooth algebraic submanifold 
of $\R ^6$, which consists of the $(u,\, v )\in\R ^3\times\R ^3$ 
such that 
\begin{equation}
\langle u,\, u\rangle =1,\quad \langle u,\, v\rangle =j_3,\quad 
\langle v ,\, v\rangle =\| j\| ^2. 
\label{SO3eq}
\end{equation}
The equations of motion for the rotational motion are given in these coordinates by 
\begin{equation}
\dd{u}{t}=u\times\omega\quad\mbox{\rm and}\quad
\dd{ v}{t}= v\times\omega ,
\label{sigmadot}
\end{equation}
in which $\omega$ is determined in terms of $u$ and $ v$ by the 
equation 
\begin{equation}
\omega  =\omega (u,\, v)={I_{\rho ,\, u}}^{-1}\, v =
J\, v +\frac{\rho\,\langle u,\, J\, v\rangle}
{1-\rho\,\langle u,\, J\, u\rangle}\, J\, u, 
\label{omegausigma}
\end{equation}
which in view of (\ref{sigmadef}) 
is equivalent to (\ref{JIu}). Here we have used the equation 
(\ref{Iuinv}) in order to write $\omega$  
even more explicitly as a function of $u$ and $ v$. 

In view of (\ref{TI}), the kinetic energy can be expressed in 
terms of $\omega$ and $ v$ as 
\begin{equation}   
T=\frac12\,\langle  v ,\,\omega\rangle ,
\label{Tomegasigma}
\end{equation}
which in view of (\ref{omegausigma}) and (\ref{Iuinv}) 
can be written in the form 
\begin{equation}   
T=\frac12\,\langle v ,\, J\, v\rangle 
+\frac12\,\frac{\rho\,\langle u,\, J\, v\rangle ^2}
{1-\rho\,\langle u,\, J\, u\rangle}. 
\label{Tusigma}
\end{equation}

Later it will turn out to be convenient to write the kinetic 
energy equation in the form 
\begin{equation}
f(u,\, v ):=Y(u,\, v )^2-X(u)\, Z( v )=0,
\label{phipsichi}
\end{equation}
in which 
\begin{eqnarray}
X(u)&:=&\rho ^{-1}-\langle u,\, J\, u\rangle , 
\label{phidef}\\
Y(u,\,  v )&:=&\langle u,\, J\, v\rangle ,
\quad\quad\mbox{\rm and}
\label{chidef}\\
Z( v )&:=&2T
-\langle  v ,\, J\, v\rangle .
\label{psidef}
\end{eqnarray}
Note that $f(u,\, v )$ is a polynomial of degree four, but 
of degree two in each of the variables $u$ and $ v$ separately. 

\subsection{The Left $\op{SO}(2)$ Action}
\label{SO2action}
If $R$ is a rotation about the vertical axis, then 
its action from the left sends 
$A$ and $\dot{A}$ to 
$R\, A$ and $R\,\dot{A}$, respectively. It therefore 
leaves $\omega$ and $u$ invariant and sends $j$ to $R\, j$. 
Note that the action of the group $\op{SO}(2)$ of the 
rotations in $\R ^3$ about the vertical axis is free 
on the set ${\cal J}'$ of $j\in\R ^3$ which are not equal to 
a multiple of $e_3$. In ${\cal J}'$, the $\op{SO}(2)$\--orbits 
are equal to the level curves of the functions 
$F(j)=\langle j,\, e_3\rangle =j_3$ and 
$G(j)=\langle j,\, j\rangle =\| j\| ^2$, 
where $j\in {\cal J}'$ corresponds to the condition that 
$F(j)^2<G(j)$. Substituting (\ref{JIu}) we obtain the 
constants of motion 
\begin{equation}
j_3=\langle j,\, e_3\rangle 
=\langle I_{\rho ,\, u}\,\omega ,\, u\rangle =\langle I\,\omega ,\, u\rangle 
\label{F}
\end{equation}
and 
\begin{equation}
\| j\| ^2=\langle j,\, j\rangle 
=\langle I_{\rho ,\, u}\,\omega ,\, I_{\rho ,\, u}\,\omega\rangle 
\label{G}
\end{equation}
for the left $\op{E}(2)$\--reduced system for $(u,\,\omega )\in 
\op{S}^2\times\R ^3$.  

Let 
\[
\pi :(A,\,\omega )\mapsto (u,\,\omega )=\left( A^{-1}\, e_3,\,\omega\right) 
\] 
denote the projection from the phase space 
$\op{T}\left(\op{SO}(3)\right)\simeq\op{SO}(3)\times\R ^3$ 
of the rotational motion onto the phase space 
$\op{S}^2\times\R ^3$ of the left $\op{E}(2)$\--reduced 
system, which maps each left $\op{SO}(2)$\--orbit to a point. 
The fact that the action of $\op{SO}(2)$ on ${\cal J}'$ is free 
implies that $\pi$ is a diffeomorphism from the 
submanifold of $\op{SO}(3)\times\R ^3$ 
determined by the equation $A\, I_{\rho ,\, u}\,\omega =j$ 
onto the the submanifold of 
$S^2\times\R ^3$ determined by the equations (\ref{F}) and (\ref{G}), 
where each of these 
submanifolds is invariant under motion of the system. 
The left $\op{SO}(2)$\--invariance of the system means 
that $\pi$ intertwines the rotational motion with the 
flow of the $\op{E}(2)$\--reduced system. 

\subsection{Chaplygin}
\label{chaplremeqmot}
In the left column of the the following table we list 
the variables and some formulas which appear in 
Chaplygin \cite[\S 2]{chaplsphere}, with our 
corresponding notations in the right column. 
It is assumed that the moment of inertia tensor $I$  
is in diagonal form, in accordance to the ``principal 
axes of inertia attached to the sphere'' of  
Chaplygin \cite[\S 2]{chaplsphere}. 
\[
\begin{array}{cc}
\mbox{\rm Chaplygin \cite[\S 2]{chaplsphere}} & \mbox{\rm our notation}\\
(p,\, q,\, r) & \omega\\
(u,\, v,\, w) & A^{-1}\,\dd{p}{t}\\
(P,\, Q,\, R) & A^{-1}\, j= v\\
\left(\gamma ,\,\gamma ',\,\gamma ''\right) & A^{-1}\, e_3=u\\
m & m\\
\rho & r\\
D=m\,\rho ^2 & \rho =m\, r^2\\
(L,\, M,\, N) & \mbox{\rm diagonal of}\; I\\
(A,\, B,\, C) & \mbox{\rm diagonal of}\; I+\rho =J^{-1}\\
\mbox{\rm (1)} & \mbox{\rm (\ref{pdot})}\\
\mbox{\rm (2)} & \mbox{\rm (\ref{omegausigma})}\\
\omega\;\mbox{\rm in (3)} & -\langle\omega ,\, u\rangle\\
X & X(u) \;\mbox{\rm in (\ref{phidef})}\\
Y & Y(u,\, v )\;\mbox{\rm in (\ref{chidef})}\\
\mbox{\rm (6)} & 
(I+\rho )\,\omega = v +\frac{Y}{X}\, u,\quad 
\rho\, \langle\omega ,\, u\rangle =\frac{Y}{X}\\
\mbox{\rm (7)} & \mbox{\rm (\ref{sigmadot}) and (\ref{udot})}\\
n & \| j\| ^2=\langle j,\, j\rangle =\langle  v ,\, v\rangle\\
h & j_3=\langle j,\, e_3\rangle =\langle v ,\, u\rangle\\
l & 2T=\langle  v ,\,\omega\rangle ,
\;\mbox{\rm cf. (\ref{Tomegasigma})}\\
Z & Z( v)\;\mbox{\rm in (\ref{psidef})}  \\
\mbox{\rm (10)} & \mbox{\rm (\ref{phipsichi})} 
\end{array}
\]
The only comment of Chaplygin to his formulas (1) and (2) 
consists of the preceding sentence ``We easily find ...''.

No explicit notation has been introduced in 
Chaplygin \cite[\S 2]{chaplsphere} for the rotation $A$.  
However, when Chaplygin said ``(7) are the equations of motions 
of the sphere'', it is clear that he meant our Subsection \ref{pairsubsec}.   

The constants of motion of the $\op{E}(2)$\--reduced system 
are Chaplygin's $n$, $h$ and $l$, which correspond to our 
$\| j\| ^2$, $j_3$ and $2T$, respectively. 

\section{The Level Surfaces of the Constants of Motion}
\label{levelsec}
\subsection{Fixing the Moment}
\label{fixmoment}
The system of equations (\ref{Adot}), (\ref{omegadot}) 
in the $(A,\,\omega )$\--space $\op{SO}(3)\times\R ^3$ 
describes the rotational motion of Chaplygin's ball. 
It is equal to the system which is obtained by ignoring 
the equation (\ref{pdot}) for the motion of the point of 
contact (or the center of gravity), which is 
the same as the $\R ^2$\--reduced system, 
obtained by working modulo the symmetry 
group of the horizontal translations 
$(A,\,\omega ,\, a)\mapsto (A,\,\omega ,\, a+b)$, 
where $b\in\R ^2$ is viewed as a horizontal vector in $\R ^3$. 
The constants of motion, viewed as functions of 
$(A,\,\omega ,\, a)$ in the phase space $\op{SO}(3)\times\R ^3\times\R ^2$, 
do not depend on the horizontal translations $a$, and therefore will be 
considered as functions of $(A,\,\omega)$ in the phase space 
$\op{SO}(3)\times\R ^3$ of the rotational motion.   

As observed at after (\ref{omegaJIu}), the 
constant of motion $j$ (= the moment of the momentum about 
the point of contact) can be used in order to eliminate 
$\omega$ from the equations of motion. 
In other words, $j:\op{SO}(3)\times\R ^3\to\R ^3$ is an  
analytic (rational) fibration, of which each fiber is equal to the graph 
of an analytic (rational) function ($\omega =\omega _j(A)$ a function of $A$), 
such that the projection $(A,\,\omega)\to A$ is an analytic (rational) 
diffeomorphism from the 
level set of $j$ onto $\op{SO}(3)$. 

We note that $j$ is {\em not} invariant under the full 
symmetry group $\op{E}(2)$ (= the horizontal motion group)   
of the system. If $R$ is a rotation around the 
vertical axis then it acts on the phase space by 
sending $(A,\,\omega ,\, a)$ to $(R\, A,\,\omega ,\, a)$. 
It leaves $u=A^{-1}\, e_3$ invariant and we 
read off from (\ref{JIu}) that it sends $j$ to $R\, j$. 
Therefore the level set is $\op{E}(2)$\--invariant if and 
only if the level $j$ is vertical. If $j$ is not vertical 
and $R$ is a non\--trivial rotation around the 
vertical axis, then the action of $R$ on the phase space 
sends the level set at the level $j$ to the different 
(disjoint) level set at the different level $R\, j$. 

The level set of all the constants of motion $j$ and $T$ together 
is diffeomorphic to the 
level set in $\op{SO}(3)$ of the function 
$T_j$ on $\op{SO}(3)$, defined by 
$T_j(A)=T(A,\,\omega )$ when 
$A\, I_{\rho,\, u}\,\omega =j$. 
It follows from (\ref{TI}) that 
\begin{equation}
T_j(A)=\frac12\, \langle I_{\rho ,\, u}\,\omega ,\,\omega\rangle 
=\frac12\,\langle j,\, A\,\omega\rangle 
=\frac12\,\langle j,\, A\, {I_{\rho ,\, u}}^{-1}\, A^{-1}\, j\rangle .
\label{Tj}
\end{equation}

If $j=0$ then, $\omega\equiv 0$, $\dd{A}{t}\equiv 0$, 
$T_j\equiv 0$, and Chaplygin's sphere 
is at rest. We will exclude this rather trivial 
case in the remainder of our discussions. 
The vector field of the motion in $\op{SO}(3)$ is given by 
(\ref{Adot}), where $\omega =\omega _j(A)$ is given by (\ref{JIu}). 
If $j\neq 0$ then $\omega _j(A)\neq 0$ and the 
vector field on $\op{SO}(3)$ has no zeros. 
It may also be observed 
that replacing $j$ by $c\, j$ with a constant $c$ 
leads to replacing $\omega$ by $c\,\omega$, multiplying the vector field 
on $\op{SO}(3)$ by $c$, whereas $T_{c\, j}=c^2\, T_j$. The solutions 
of the equations of motion are changed only by a rescaling of time 
by the constant factor $c$. 

If $T$ is a regular value of $T_j$, then the level set 
$\left\{ A\in\op{SO}(3)\mid T_j(A)=T\right\}$ 
is a smooth (algebraic) closed two\--dimensional 
submanifold of $\op{SO}(3)$, compact because $\op{SO}(3)$ is compact. 
It is oriented by the area form $\Omega/\op{d}\! T_j$, where 
$\Omega$ is a volume form on $\op{SO}(3)$. 
Let $C$ be a connected component of a regular level set.  
We conclude that $C$ is a compact connected oriented two\--dimensional 
smooth (algebraic) manifold which carries a tangent vector field without zeros, 
which implies that the Euler characteristic of $C$ is equal to zero. 
According to the classification of compact oriented surfaces, this 
in turn implies that $C$ is diffeomorphic to the two\--dimensional 
{\em torus} $\R ^2/\Z ^2$. The considerations below will lead to a much more 
detailed description of the level sets, from which the conclusion 
that the regular ones consist of tori can be obtained without using 
the just mentioned facts from differential topology. 

\subsection{The Critical Points of the Energy}
For any vector $\nu\in\R ^3$, let $\op{R}_{\nu}(A)\in\op{T}_A\op{SO}(3)$ 
be the tangent vector 
of $\op{SO}(3)$ at $A\in\op{SO}(3)$ which is given by 
\begin{equation}
\op{R}_{\nu}(A):=A\circ\nu _{\scriptop{op}}.
\label{Lnu}
\end{equation}

In order to determine the derivative $\op{R}_{\nu}\omega _j$ 
of the vector\--valued function $\omega _j$ on $\op{SO}(3)$ 
in the direction of the vector field $\op{R}_{\nu}$, 
we begin with the observation that 
(\ref{JIu}) implies that 
\[
0=\op{R}_{\nu}\left( A\, I_{\rho ,\, u}\,\omega _j\right) 
=A\left(\nu\times I_{\rho ,\, u}\,\omega _j 
+\op{R}_{\nu}I_{\rho ,\, u}\,\omega _j\right) .
\]
Furthermore, 
\[
\op{R}_{\nu}u=\, \op{R}_{\nu}A^{-1}\, e_3
\, -\nu\times A^{-1}\, e_3=\, -\nu\times u, 
\]
and therefore it follows from (\ref{Iu}) that 
\[
\op{R}_{\nu}I_{\rho ,\, u}\,\omega _j
=\rho\,\left(\langle\omega _j,\,\nu\times u\rangle\, u
+\langle\omega _j,\, u\rangle\,\nu\times u\right) 
+I_{\rho ,\, u}\,\op{R}_{\nu}\omega _j. 
\]
Therefore $\op{R}_{\nu}\omega _j$ is determined by the equation
\begin{equation}
0=\nu\times ((I+\rho )\,\omega _j)
+\rho\,\langle\omega _j,\,\nu\times u\rangle\, u
+I_{\rho ,\, u}\,\op{R}_{\nu}\omega _j. 
\label{Rnuomega}
\end{equation}

Because $I_{\rho ,\, u}$ is symmetric, we have that 
\[
\frac12\langle I_{\rho ,\, u}\,\op{R}_{\nu}\omega _j,\,\omega _j\rangle 
+\frac12\langle I_{\rho ,\, u}\,\omega _j,\,\op{R}_{\nu}\omega _j\rangle 
=\langle I_{\rho ,\, u}\,\op{R}_{\nu}\omega _j,\,\omega _j\rangle .
\]
The derivative of $T_j$, cf. (\ref{TI}), 
in the direction of $\op{R}_{\nu}$ therefore is 
equal to 
\begin{eqnarray*}
\op{R}_{\nu}T_j&=&\frac12\,
\langle\left(\op{R}_{\nu}I_{\rho ,\, u}\right)\,\omega _j,\,\omega _j\rangle
+\langle I_{\rho ,\, u}\,\op{R}_{\nu}\omega _j,\,\omega _j\rangle\\
&=&\frac12\langle\langle\omega _j,\,\nu\times u\rangle\, u
+\langle\omega _j,\, u\rangle\,\nu\times u\\
&&-\langle \nu\times ((I+\rho )\,\omega _j)
+\rho\,\langle\omega _j,\,\nu\times u\rangle\, u,\,\omega _j\rangle ,
\end{eqnarray*}
from which we obtain that 
\begin{equation}
\op{R}_{\nu}T_j=\langle\omega _j\times (I+\rho )\,\omega _j,\,\nu\rangle 
=\langle\omega _j\times I\,\omega _j,\,\nu\rangle .  
\label{RnuT}
\end{equation}
 
Let $\Sigma _j$ denote the set of critical points of $T_j$. 
It follows from (\ref{RnuT}) that $A\in\Sigma _j$ 
if and only if $\omega =\omega _j(A)$ satisfies 
$\omega\times I\,\omega =0$. 
In view of (\ref{omegadot}) this is equivalent to 
$\dd{\omega}{t}=0$. Because the function $T_j$ is invariant under the motion,  
the $\op{R}_{\omega _j}$\--flow on $\op{SO}(3)$, 
$\Sigma _j$ is invariant under the 
motion. Therefore, if $A(0)\in\Sigma _j$ then we have for every $t$ that 
$A(t)\in\Sigma _j$, which implies that 
$\dd{\omega _j(A(t))}{t}=0$ for every $t$ 
and therefore $\omega =\omega _j(A(t))$ is a constant. 
It follows then from (\ref{Adot}) that 
$A(t)=A(0)\circ\op{e}^{t\,\omega}$ describes a circle, 
which we will call a {\em critical circle}. 

We have $\omega\times I\,\omega =0$ if and only if 
$I\,\omega =\iota\,\omega$, which means that 
$\omega$ is an eigenvector of the moment of inertia tensor $I$, 
with eigenvalue $\iota$ equal to one of the principal inertial moments 
$I_1$, $I_2$ or $I_3$. If this is the case, 
it then follows from (\ref{JIu}) that 
\begin{equation}
j=(\iota +\rho)\, A\,\omega -\rho\,\langle A\,\omega ,\, e_3\rangle\, e_3, 
\label{jiota}
\end{equation}
Taking the inner product with $e_3$ we obtain that 
\begin{equation}
j_3=\iota\,\langle A\,\omega ,\, e_3\rangle ,
\label{j3}
\end{equation} 
which can be inserted into (\ref{jiota}) in order to yield that 
\begin{equation}
A\,\omega =
\frac{1}{\iota +\rho}\,\left( j+\frac{\rho\, j_3}{\iota}\, e_3\right) . 
\label{Aiota}
\end{equation}
From (\ref{Aiota}) we obtain that 
\begin{equation}
\|\omega \| ^2=
\frac{1}{(\iota +\rho )^2}\, 
\left( \| j\| ^2+2\frac{\rho\, {j_3}^2}{\iota}
+\left(\frac{\rho\, j_3}{\iota}  
\right) ^2\right) , 
\label{lengthomega}
\end{equation}
and combining (\ref{Aiota}) with (\ref{Tj}) we obtain that the 
critical level is equal to 
\begin{equation}
T_{\scriptop{crit}}=\frac12\,\frac{1}{\iota +\rho}\, 
\left( \| j\| ^2+\frac{\rho\, {j_3}^2}{\iota}\right) .  
\label{Tcrit}
\end{equation}
Note that the right hand side of (\ref{Tcrit}) is a 
monotonously decreasing function of $\iota$, which implies that 
if the principal inertial moments are taken in increasing order, 
then the corresponding critical levels of 
the kinetic energy appear in decreasing order. 

It follows from (\ref{Aiota}) and 
(\ref{pdot}) that {\em the point of contact $p(t)$ moves 
along a straight line, with constant velocity equal to} 
\begin{equation}
\dd{p}{t}=
\frac{r}{\iota +\rho}\, j\times e_3. 
\label{pcrit}
\end{equation}
We note that for the uniformly rolling sphere 
the axis of rotation $\R\, A\,\omega$ need not be horizontal 
as one might expect. It follows from (\ref{j3}) that 
it is horizontal if and only if the vector $j$ 
is horizontal. 

The point of contact is at rest if and only if the moment of 
momentum is vertical ($j$ is equal to a multiple of $e_3$), 
which according to (\ref{Aiota}) corresponds to the case that 
the vector $A\,\omega$ is vertical. In this case the sphere 
is spinning around the vertical axis, which then coincides with 
an inertial axis. 

\medskip
Let $E_{j,\, \iota}$ denote the set of $\omega\in\op{ker}(I-\iota )$
such that (\ref{lengthomega}) holds, and  
let $\Sigma _{j,\,\iota}$ denote the set of 
$A\in\op{SO}(3)$ such that (\ref{Aiota}) holds for some 
$\omega\in E_{j,\,\iota}$. Then $\Sigma _j$ is a smooth algebraic 
circle bundle over the smooth algebraic manifold $E_{j,\,\iota}$ 
and therefore is an smooth algebraic submanifold of $\op{SO}(3)$. 
For different values of the eigenvalue $\iota$ of $I$, 
the sets $\Sigma _{j,\,\iota}$ are disjoint and 
$\Sigma _j$ is equal to the union of the 
$\Sigma _{j,\,\iota}$, where $\iota$ runs over the 
principal inertial moments. 

If $\iota$ is a simple eigenvalue of $I$ then 
$E_{j,\,\iota}$ consist of two opposite 
($\omega\leftrightarrow -\omega$) eigenvectors of $I$ 
for the eigenvalue $\iota$, and 
$\Sigma _{j,\,\iota}$ consists of two disjoint critical circles. 
Note that the function $T_j$ (and also the equation of motion) 
is invariant under a transformation $A\mapsto A\, R$, where 
$R\in\op{SO}(3)$ commutes with $I$. There exists such 
$R$ which maps $\omega$ to $-\omega$ and for each such $R$  
the mapping $A\mapsto A\, R$ interchanges the two critical 
circles in $\Sigma _{j,\,\iota}$. 
In the generic case that all the principal inertial moments $I_1$, $I_2$, 
$I_3$ are different, we obtain six critical circles, 
two for every choice of $\iota =I_1,\, I_2,\, I_3$.  

The case of two equal principal 
moments of inertia $I_1=I_2\neq I_3$ is that of a body of revolution 
with surface equal to a sphere and center of mass at 
the center of the sphere. This is the example of {\em Routh's sphere 
with center of mass at the center of the sphere}, 
or {\em Bobylev's sphere}, which will be 
discussed in some more detail in Section \ref{srsec}. 
If $\iota =I_1=I_2\neq I_3$ then $E_{j,\,\iota}$ is a circle in 
$\op{ker}(I-\iota )$ and $\Sigma _{j,\,\iota}$ is a two\--dimensional 
torus in $\op{SO}(3)$. As discussed before, 
$\Sigma _{j,\, I_3}$ consists of two critical circles. 

If all the principal moments of inertia are equal, $\iota =I_1=I_2=I_3$, 
or equivalently $I$ is equal to $\iota$ times the identity, then 
$\Sigma _j=\op{SO}(3)$, the function $T_j$ is constant, 
and all solutions of the equations of motion are of the form $A(t)=A(0)\circ\op{e}^{t\,\omega}$ 
with a constant vector $\omega$. 

\medskip
We will now verify that {\em each $\Sigma _{j,\,\iota}$ 
is a is a nondegenerate critical manifold} of $T_j$ 
in the sense of Bott \cite{Bott}, which means that for each 
$A\in\Sigma _{j,\,\iota}$ the null space of the Hessian 
$T_j''(A)$ of $T_j$ at $A$ is equal to the tangent space 
$\op{T}_A\Sigma _{j,\,\iota}$ of $\Sigma _{j,\,\iota}$ at $A$. 
Because always $\op{T}_A\Sigma _{j,\,\iota}\subset 
\op{ker}T_j''(A)$, we only need to verify that the dimension 
of the null space of the Hessian is at most equal to the dimension of the 
critical submanifold. This follows from

\begin{lemma}
Let $I\,\omega _j(A)=\iota\,\omega _j(A)$ and 
$\op{R}_{\nu}(A)\in\op{ker}T_j''(A)$. 
If $\iota$ is a simple eigenvalue of 
$I$, then $\nu$ is a multiple of $\omega _j(A)$. 
If $\iota$ is a double eigenvalue of $I$ and 
$\alpha$ is a nonzero vector which is orthogonal to 
$\op{ker}(I-\iota )$, then $\nu$ is a linear combination 
of $\omega _j(A)$ and $\alpha$. 
\label{Bottlem}
\end{lemma}

\begin{proof}
It follows from (\ref{RnuT}) that $\op{R}_{\nu}(A)\in\op{ker}T_j''(A)$ 
if and only if, at the point $A\in\op{SO}(3)$,  
\[
0=\op{R}_{\nu}\,\left( I\,\omega _j\times\omega _j\right) 
=I\,\op{R}_{\nu}\,\omega _j\times\omega _j
+I\,\omega _j\times\op{R}_{\nu}\,\omega _j
=\left( I\,\op{R}_{\nu}\,\omega _j-\iota\,\op{R}_{\nu}\,\omega _j\right) 
\times\omega ,
\]
or $(I-\iota )\,\op{R}_{\nu}\,\omega _j$ is equal to a multiple 
of $\omega _j$. This implies that $(I-\iota )^2\,\op{R}_{\nu}\,\omega _j=0$, 
which in turn implies that $(I-\iota )\,\op{R}_{\nu}\,\omega _j=0$, 
because $I-\iota$ is a diagonal matrix. 

Substituting $I\omega _j=\iota\,\omega _j$ and 
$I\,\op{R}_{\nu}\,\omega _j=\iota\,\op{R}_{\nu}\,\omega _j$ in 
(\ref{Rnuomega}), we obtain 
\begin{eqnarray*}
0&=&\nu\times ((\iota +\rho )\,\omega _j)
+\rho\,\langle\omega _j,\,\nu\times u\rangle\, u
+(\iota +\rho )\,\op{R}_{\nu}\omega _j 
-\rho\,\langle\op{R}_{\nu}\,\omega _j,\,\nu\times u\rangle\, u\\
&=&(\iota +\rho )\,\theta -\rho\,\langle\theta,\, u\rangle\, u,
\end{eqnarray*}
in which $\theta :=\nu\times\omega _j+\op{R}_{\nu}\,\omega _j$. 
It follows that $\theta =0$. 

If $\iota$ is a simple eigenvalue of $I$ then 
$I\,\op{R}_{\nu}\,\omega _j=\iota\,\op{R}_{\nu}\,\omega _j$ 
implies that $\op{R}_{\nu}\,\omega _j$ is equal to a multiple 
of $\omega _j$ and it follows from $\theta =0$ that 
$\nu\times\omega _j=0$ and $\op{R}_{\nu}\,\omega _j=0$, 
which in turn implies that $\nu$ is a multiple of $\omega _j$. 

If $\op{ker}(I-\iota )$ is two\--dimensional, 
then $\theta =0$ implies that 
$0=\langle\nu\time\omega _j,\,\alpha\rangle 
=\langle\nu ,\,\omega _j\times\alpha\rangle$. 
Because $\omega _j\in\op{ker}(I-\iota )$, we have that 
$\alpha$ is orthogonal to $\omega _j$ and it follows that  
$\nu$ is a linear combination of $\omega _j(A)$ and $\alpha$. 
\end{proof}

Because $T_j$ is a continuous function on the compact 
set $\op{SO}(3)$, it attains its maximum and its minimum. 
The variational principle says that the points $A$ where $T_j$ 
attains its maximum (minimum) are critical points for 
$T_j$. Therefore, if $I_1\leq I_2\leq I_3$, then 
the maximum (minimum) value of $T_j$ is equal to the 
right hand side of (\ref{Tcrit}), with $\iota =I_1$ ($\iota =I_3$). 
If $I_1<I_2$ ($I_2<I_3$), then $\Sigma _{j,\,\iota}$ consists 
of two critical circles and on a transversal two\--dimensional 
manifold the function $T_j$ has a nondegenerate maximum (minimum), 
the nearby level sets of which are small loops around the critical points. 
It follows that the level sets of $T_j$ near $\Sigma _{j,\,\iota}$ 
consist of narrow tubes around the critical circles. This implies 
that the critical circles corresponding to $\iota =I_1<I_2\leq I_3$ 
or to $\iota =I_3>I_2\geq I_1$ are {\em stable} periodic solutions 
of the system. 

Now assume that $I_1<I_2<I_3$. 
A Morse theoretic argument then yields that $T_j$ has 
cannot have a local maximum or minimum at $\Sigma _{j,\, I_2}$. 
Note that the index of $T_j''$, the number of negative eigenvalues 
of $T_j''$, is constant along each of the 
two critical circles in $\Sigma _{j,\, I_2}$. If 
$R=\op{diag}(-1,\, -1,\, 1)$ or $R=\op{diag}(1,\, -1,\, -1)$, 
then $A\mapsto A\, R$ leaves $T_j$ invariant and interchanges the 
two critical circles in $\Sigma _{j,\, I_2}$, which implies that 
the index of $T_j''$ is constant along $\Sigma _{j,\, I_2}$. 
Suppose that it is equal to two, which means that 
$T_j$ has a local maximum at $\Sigma _{j,\, I_2}$. 
Let $\op{grad}T_j$ denote the gradient vector field of $T_j$ with respect to 
a given Riemannian structure on $\op{SO}(3)$. Let $\gamma ^t$ 
denote the flow of $\op{grad}T_j$. Define $S_1$ and $S_2$ 
as the set of $A\in\op{SO}(3)$ such that, when $t\to\infty$, 
$\gamma ^t(A)$ converges to $\Sigma _{j,\, I_1}$ and $\Sigma _{j,\, I_2}$, 
respectively. $S_1$ and $S_2$ are nonvoid disjoint open subsets 
of $\op{SO}(3)$, with union equal to $\op{SO}(3)\setminus\Sigma _{j,\, I_3}$. 
The set $\op{SO}(3)\setminus\Sigma _{j,\, I_3}$ is connected, 
because $\Sigma _{j,\, I_3}$ is a codimension two submanifold 
of the connected manifold $\op{SO}(3)$. This leads to a contradiction. 
In a similar way the assumption that $T_j$ has a local maximum at 
$\Sigma _{j,\, I_2}$ leads to a contradiction and the conclusion 
is that the index of $T_j''$ is equal to one along $\Sigma _{j,\, I_2}$, 
which means that transversally $T_j$ has a saddle point behaviour. 

Let $P$ denote the linearization of the return map to a transversal 
plane (the Poincar\'e map) of the flow along the critical circles in 
$\Sigma _{j,\, I_2}$. It follows from Corollary \ref{volcor} 
that $\op{det}P=1$.  
This implies that $|\op{trace}P|<2$ if and only if 
$P$ is conjugate to a nontrivial rotation, whereas 
$|\op{trace}P| >2$ if and only if $P$ is a hyperbolic 
map with real eigenvalues $\lambda$, $1/\lambda$ such that 
$\lambda\neq \pm 1$, in which case the critical circles in 
$\Sigma _{j,\, I_2}$ are linearly unstable. 
The saddle point behaviour of $T_j$ near 
$\Sigma _{j,\, I_2}$ (when $I_1<I_2<I_3$) excludes that $P$ is 
conjugate to a nontrivial rotation. 
The number $\op{trace}P$ depends in a real analytic fashion on  
$j\in\R ^3\setminus\{ 0\}$.  As we will see 
at the end of Subsection \ref{vertj} below, if $j$ is 
vertical then the critical circles in $\Sigma _{j,\, I_2}$ 
are linearly unstable, which implies that 
$|\op{trace}P| >2$ when $j$ is vertical. It follows that the set 
$N$ of $j\in\R ^3\setminus\{ 0\}$ such that $\op{trace}P=\pm 2$ is 
a proper closed analytic subvariety of $j\in\R ^3\setminus\{ 0\}$. 
Using the invariance of the equations of motion under the action 
$(A,\,\omega )\mapsto (R\, A,\,\omega )$, $j\mapsto R\, j$ 
of the rotations $R$ around the vertical 
axis and the homogeneity $(A,\,\omega )\mapsto (A,\, c\,\omega )$, 
$j\mapsto c\, j$, it follows that $N$, if not empty, is equal 
to the union of finitely many cones in $\R ^3$ 
which are invariant under the rotations 
around the vertical axis.  
For all $j$ in the complement of $N$, which is an open and dense 
subset of $\R ^3\setminus\{ 0\}$, the critical circles in 
$\Sigma _{j,\, I_2}$ are linearly unstable. This in turn implies that 
{\em every critical circle with $I\,\omega =I_2\,\omega$ is unstable with 
repect to the flow in the full phase space}. 

Conversely, at any critical circle which is linearly unstable the 
invariant function $T_j$ must have transversal saddle point 
behaviour. Because the linearly unstable critical circles 
are dense, this leads to a proof that the index of $T_j$ 
is equal to one along $\Sigma _{j,\, I_2}$ without using Morse theory.  

\begin{question}
When $I_1<I_2<I_3$, is every critical circle with 
$I\,\omega =I_2\,\omega$ linearly unstable? 
\end{question} 

\subsection{The Moment Mapping}
It is also instructive to consider 
$(j,\, T)$ as a mapping from the phase space 
$\op{SO}(3)\times\R ^3$ of the rotational motion 
to $\R ^3\times\R$. The set of {\em singular points}  
of $(j,\, T)$, the set of points where the rank of the 
tangent map is less than four, is equal to the 
union of $\op{SO}(3)\times \{ 0\}$ and the 
set of all $(A,\,\omega )$ where 
$j\in\R ^3\setminus\{ 0\}$, $A\in\Sigma _j$ and 
$\omega =\omega _j(A)$. If $I_1<I_2<I_3$ then the latter 
set is a smooth conic closed submanifold 
of codimension two in $\op{SO}(3)\times\left(\R ^3\setminus\{ 0\}\right)$. 

Because of the scaling  
$(A,\,\omega )\mapsto (A,\, c\,\omega )$, which 
maps $(j,\, T)$ to $\left( c\, j,\, c^2\, T\right)$, 
and because the case $T=0$, when $\omega =0$ and everything is at rest, 
is not very interesting, 
we restrict ourselves to an energy hypersurface 
where $T$ is equal to a positive constant. 
(An isotopy argument as below shows that 
the energy hypersurface is diffeomorphic to the Cartesian 
product of $\op{SO}(3)$ with a two\--dimensional sphere.)  
The moment of the momentum around the point of contact 
then defines a mapping $j_T$ from the energy hypersurface to $\R ^3$. 
The singular points of $j_T$ are the above singular points 
on the energy hypersurface. 

According to 
(\ref{Tcrit}), the set of singular values of 
$j_T$ consists of the points $j\in\R ^3$ such that 
\begin{equation}
\frac{1}{2\left(\iota +\rho\right)}\, 
\left(\| j\| ^2+\frac{\rho}{\iota}\, {j_3}^2\right)=T,
\label{jsing}
\end{equation}
in which $\iota =I_1$, $\iota =I_2$ or $\iota =I_3$. 
If $I_1<I_2<I_3$, then this set is 
equal to the union of three disjoint ellipsoids with center at the origin 
and which are invariant under the rotations around the vertical 
axis. The inner and the outer one correspond to 
$\iota =I_1$ and $\iota =I_3$, respectively. The points in the 
phase space of the rotational motion which by $j_T$ are 
mapped to these ellipsoids correspond to the extremal 
critical circles (the stable ones) of the functions $T_j$, 
and it follows that the inner and outer ellipsoids 
together form the boundary 
of the image of $j_T$. Therefore, the image of $j_T$ 
is equal to the set of all $j\in\R ^3$ such that (\ref{jsing}) holds for 
some $\iota\in\left[ I_1,\, I_3\right]$. 

The unstable critical circles are mapped to the 
intermediate (interior) ellipsoid described by (\ref{jsing}) with 
$\iota =I_2$. For the singular values $j$ in this 
interior ellipsoid, the level sets are two\--dimensional, 
with a singularity of normal crossing type along the two unstable 
critical circles in the level set.

For the regular values of $j_T$, the points $j$ such that (\ref{jsing}) holds for 
some $\iota$ such that $I_1<\iota <I_2$ or $I_2<\iota <I_3$, 
the level sets are smooth (alegbraic) two\--dimensional compact oriented 
submanifolds of the phase space for the rotational motion.

\subsection{Isotopy of the Fibration}
Following Cushman \cite[p. 412]{C83}, a smooth function 
on a compact manifold, for which 
the critical set consists of nondegenerate critical 
manifolds (possibly with varying dimensions), will be  
called a {\em Bott\--Morse function}. 
We will use an isotopy lemma for 
families of Bott\--Morse functions, which 
should be well\--known. However, because we did not 
find a reference in the literature, we include a proof.  
\begin{lemma}
Let $M$ be a compact smooth manifold and $f_{\epsilon}$ a 
familie of smooth functions on $M$, depending smoothly 
on a real parameter $\epsilon$. Furthermore assume that the set 
of critical points of $f_{\epsilon}$ consists of 
finitely many disjoint compact connected smooth submanifolds 
$C_{\epsilon ,\, i}$, $1\leq i\leq N$, depending smoothly on $\epsilon$ 
and such that $C_{\epsilon ,\, i}$ is a nondegenerate critical 
manifold of $f_{\epsilon ,\, i}$. On $C_{\epsilon ,\, i}$ the function $f_{\epsilon}$ is constant, let $F_{\epsilon ,\, i}$ 
be the value of $f_{\epsilon}$ on $C_{\epsilon ,\, i}$. 
We finally assume that the ordering of the real numbers 
$F_{\epsilon ,\, i}$, $1\leq i\leq N$, does not change with 
varying $\epsilon$. Under these assumptions there exist 
smooth diffeomorphisms $\psi _{\epsilon}$ and $\Phi _{\epsilon}$ 
of $\R$ and $M$ respectively, such that 
$\psi _{\epsilon}\circ f_{\epsilon}\circ\Phi _{\epsilon}$ 
does not depend on $\epsilon$. The $\psi _{\epsilon}$ can be 
chosen to be order\--preserving. 
\label{isotoplem}
\end{lemma}

\begin{proof}
The assumption that the ordering of the critical values does not 
change implies that there exists a family of order\--preserving 
smooth diffeomorphisms $\psi _{\epsilon}$ of $\R$ such 
that the real numbers 
$G_i:=\psi _{\epsilon}\left( F_{\epsilon ,\, i}\right)$, 
$1\leq i\leq N$ do not depend on $\epsilon$. 
It is also quite easy to prove that there exists a smooth family of 
diffeomorphisms $\Psi _{\epsilon}$ of $M$, depending smoothly 
on $\epsilon$, such that the manifolds 
$D_i:=\Psi _{\epsilon}^{-1}\left( C_{\epsilon ,\, i}\right)$ 
do not depend on $\epsilon$. 
The functions 
$g_{\epsilon}:=\psi _{\epsilon}\circ f_{\epsilon}\circ\Psi _{\epsilon}$ 
have the same properties as $f_{\epsilon}$, but now with the 
constant nondegenerate critical manifolds $D_i$ on which 
$g_{\epsilon}$ has the constant critical values $G_i$. 

We now follow the idea of the proof of Moser \cite{Moser}. 
The condition for a smooth family of diffeomorphisms $\Xi _{\epsilon}$ 
of $M$, depending smoothly on $\epsilon$, that 
$g_{\epsilon}\circ\Xi _{\epsilon}=\Xi _{\epsilon}^*\, g_{\epsilon}$ 
does not depend on 
$\epsilon$, is equivalent to the condition that 
$
0=\dd{}{\epsilon}\Xi _{\epsilon}^*\, g_{\epsilon}
=\Xi _{\epsilon}^*\,\left( X_{\epsilon}\, g_{\epsilon}
+\frac{\partial g_{\epsilon}}{\partial\epsilon}\right) ,
$
or
\begin{equation}
X_{\epsilon}\, g_{\epsilon}
+\frac{\partial g_{\epsilon}}{\partial\epsilon}=0,
\label{vepseq}
\end{equation}
in which $X_{\epsilon}$ denotes the vector field on $M$ defined by 
\begin{equation}
\frac{\partial\Xi _{\epsilon}(x)}{\partial\epsilon} 
=X_{\epsilon}\left( \Xi _{\epsilon}(x)\right) ,\quad x\in M.
\label{veps}
\end{equation}

If the equation (\ref{vepseq}) for $X_{\epsilon}$ can be solved 
locally near every point of $M$, then a global solution 
can be obtained by means of a smooth 
partition of unity. 

Near a noncritical point of $g_{\epsilon}$ 
we can use $g_{\epsilon}$ as one of the local coordinates and 
(\ref{vepseq}) then amounts to prescribing the 
corresponding component of the vector field $X_{\epsilon}$. 

If $x^{(0)}\in D_i$, then one can introduce a local coordinate 
system near $x^{(0)}$ in which $x^{(0)}=0$ and $x\in D_i$ corresponds to 
$x_j=0$ for $1\leq j\leq c$, if $c$ denotes the codimension 
of $C_i$ in $M$. Writing $x=(y,\, z)$ with $y\in\R ^c$, $z\in\R ^d$, 
we view $g_{\epsilon}(y,\, z)$ as a family of functions 
of $y$, with $\epsilon$ and $z$ as parameters. 
A second order Taylor expansion with respect to $y$ at $y=0$, 
in which the remainder term in intergral form is absorbed into 
the second order term, yields that 
\[
g_{\epsilon}(y,\, z)
=G_i+\frac12\,\langle Q_{\epsilon}(y,\, z)\, y,\, y\rangle ,
\]
where $Q=Q_{\epsilon}(y,\, z)$ is a nondegenerate symmetric matrix, 
depending smoothly on all the variables. (This is also 
one of the steps in the proof of the Morse lemma with parameters of 
H\"ormander \cite[Lemma 3.2.3]{Hor}.) 
If we take $X_{\epsilon}=\left( Y_{\epsilon},\, 0\right)$ with 
$Y_{\epsilon}\in\R ^c$, 
then the equation (\ref{vepseq}) for $X_{\epsilon}$ 
is equivalent to the equation 
\[
\langle Q_{\epsilon}(y,\, z)\, Y_{\epsilon},\, y\rangle 
+\frac12\,\langle 
\left(\frac{\partial Q_{\epsilon}(y,\, z)}{\partial y}\, Y_{\epsilon}\right)
\, y,\, y\rangle 
+\frac12\,\langle \frac{\partial Q_{\epsilon}(y,\, z)}{\partial\epsilon}
\, y,\, y\rangle =0
\]
for $Y_{\epsilon}$, which is satisfied if 
\[
Q_{\epsilon}(y,\, z)\, Y_{\epsilon}
+\frac12\, 
\left(\frac{\partial Q_{\epsilon}(y,\, z)}{\partial y}\, Y_{\epsilon}\right)
\, y+\frac12\,\frac{\partial Q_{\epsilon}(y,\, z)}{\partial\epsilon}
\, y=0.
\]
For sufficiently small $y$ the latter equation has a unique solution 
$Y_{\epsilon}$ which depends smoothly on $y$, $z$, and $\epsilon$. 

Piecing together the local solutions by means of a smooth 
partition of unity, we obtain a smooth vector field $X_{\epsilon}$ 
on $M$, depending smoothly on $\epsilon$, such that (\ref{vepseq}) 
holds. Define $\Xi _{\epsilon}(x)$ as the solution of the 
$\epsilon$\--dependent ordinary differential equation 
(\ref{veps}), with initial condition $\Xi _0(x)=x$. 
Using the compactness of $M$ we obtain that 
the $\Xi _{\epsilon}$ are globally defined smooth diffeomorphisms of $M$, 
depending smoothly on $\epsilon$. Reading the paragraph preceding 
(\ref{veps}) backwards, we obtain that $g_{\epsilon}\circ\Xi _{\epsilon}$ 
does not depend on $\epsilon$. This proves the lemma with 
$\Phi _{\epsilon}=\Psi _{\epsilon}\circ\Xi _{\epsilon}$. 
\end{proof}

In order to emphasize the dependence on $\rho$ of the kinetic 
energy function on $\op{SO}(3)$, we now write $T_{\rho ,\, j}$ 
instead of $T_j$. Applying Lemma \ref{isotoplem} 
to $f_{\epsilon}=T_{\epsilon ,\, j}$, we obtain that there 
exists an order\--preserving smooth diffeomorphism $\psi$ of $\R$ 
and a smooth diffeomorphism $\Phi$ of $\op{SO}(3)$ such that 
$\psi\circ T_{\rho ,\, j}\circ\Phi =T_{0,\, j}$. 

We can make $\psi$ and $\Phi$ to depend smoothly on $j$ when 
$j$ varies over the unit sphere in $\R ^3$. Extending the transformations 
by homogeneity for the scaling $(A,\,\omega )\mapsto (A,\, c\,\omega )$, 
one obtains a diffeomorphism $\Phi$ of 
$\op{SO}(3)\times\left(\R ^3\setminus\{ 0\}\right)$, 
and a diffeomorphism $\Psi$ of $\left(\R ^3\setminus\{ 0\}\right) 
\times\R$ of the form $(j,\, T)\mapsto (j,\,\psi (j,\, T))$, 
such that $\psi\left( c\, j,\, c^2\, T\right) =c^2\,\psi (j,\, T)$ for 
every $c>0$, such that $\Psi\circ\left( j_{\rho},\, T_{\rho}\right)
\circ\Phi =\left( j_0,\, T_0\right)$. This implies that the 
smooth diffeomorphism $\Phi$ maps the whole fibration, together 
with its singularities,  of the phase space defined by the constants 
of motion for $\rho =0$ to the one for our given value of $\rho$.  

If $\rho =0$, then $I_{\rho ,\, u}=I_{0,\, u}=I$. 
The equation of motion (\ref{omegadot}) then turns into Euler's  
equation of motion 
\begin{equation}
\dd{}{t}I\,\omega =I\,\omega\times\omega 
\label{Eulertop}
\end{equation}
for the Euler top, and the constants of motion $j$ of (\ref{JIu}) 
and $T$ of (\ref{Tj}) 
are given by the familiar formulas 
\begin{equation}
j=A\, I\,\omega 
\label{JIu0}
\end{equation}
and 
\begin{equation}
T_{0,\, j}(A)=\frac12\, \langle I\,\omega ,\,\omega\rangle 
=\frac12\,\langle A^{-1}\, j,\, I^{-1}\, A^{-1}\, j\rangle 
\label{Tj0}
\end{equation}
for the moment of momentum around the center of mass and the 
kinetic energy of the Euler top, respectively. 
For details about the Euler top, we refer to 
Cushman and Bates \cite[Ch. III]{CB}. 

It follows from (\ref{Tj0}) that the kinetic energy $T_{0,\, j}$ 
of the Euler top is invariant under the circle action 
$A\mapsto \op{e}^{t\, j_{\tinyop{op}}}\circ A$, the orbits of which are 
the fibers of the mapping $ v _j :A\mapsto A^{-1}\, j$. 
Note that $ v _j:\op{SO}(3)\to\op{S}_{\| j\|}$ is a smooth 
fibration of $\op{SO}(3)$ over the {\em Euler sphere} 
$\op{S}_{\| j\|}$, the sphere in $\R ^3$ with center 
at the origin and radius equal to $\| j\|$. 
On the Euler sphere, the kinetic energy is equal to the 
restriction to $\op{S}_{\| j\|}$ of the 
quadratic form $ v\mapsto 
\frac12\,\langle v ,\, I^{-1}\, v\rangle$, 
defined by the positive definite symmetric matrix $I^{-1}$. 

The critical levels of $T_{0,\, j}$ are equal to 
$\frac12\, \frac{\| j\| ^2}{\iota}$, where $\iota =I_1$, 
$I_2$ or $I_3$, and the regular values are the numbers 
in between the critical levels. On the Euler sphere 
each regular level set has two connected components, 
opposite to each other, each of which is a smooth closed curve, 
diffeomorphic to a circle. The preimages of these 
under the mapping $ v _j$ are circle bundles over these 
circles and therefore each regular level set in $\op{SO}(3)$ 
has two connected components, each of which is diffeomorphic to 
the two\--dimensional torus. 

For each extremal level the level sets consists of two critical 
circles surrounded by narrow tubes. For the intermediate 
critical level, the level set on $\op{S}_{\| j\|}$ 
consists of two opposite critical points. The complement 
of these in the level set has four connected component, 
each of which is a smooth curve running from one of the 
critical points to the opposite one. It follows that 
the level set in $\op{SO}(3)$ of the intermediate 
critical level contains two critical circles, 
the complement of which has four connected components 
each of which is a smooth cylinder running from one of 
the critical circles to the other. 

For more details about the fibration in $\op{SO}(3)\times\R ^3$ for the Euler 
top we refer to Cushman and Bates \cite[Ch. III, Sec. 5]{CB}. 
The point of the isotopy lemma is that there exists a diffeomorphism 
$\Phi$ which sends the whole fibration with singularities 
for the Euler top to the one for Chaplygin's sphere for 
an arbitrary value of $\rho$. In particular all the qualitative 
statements about the level sets remain true. The regular level 
sets have two connected components, each of which is diffeomorphic 
to the two\--dimensional torus. The level set of an intermediate 
critical level contains two critical circles, 
the complement of which has four connected components 
each of which is a smooth cylinder running from one of 
the critical circles to the other. 

\subsection{Chaplygin}
In the beginning of \cite[\S 6]{chaplsphere} 
Chaplygin gave a short description of the critical circles, 
but without relating these solutions to the points where 
the derivatives of the constants of motion are linearly dependent. 
He also stated that the ones corresponding to the extremal 
moments of inertia are stable and the ones corresponding 
to the intermediate moment of inertia are unstable, but 
without any proof. 

The question of the smoothness of the level surface of the 
constant of motion, which is related to the question 
of the linear independence of their derivatives,  
does not occur in Chaplygin \cite{chaplsphere}. 

\section{When the Moment is Vertical}
\label{vertj}

Next to the critical circles, the solutions with {\em vertical moment $j$ 
of the momentum around the point of contact}, $j=j_3\, e_3$, form another interesting special family. 

Note that the condition that $j$ is vertical defines a 
smooth codimension two algebraic 
submanifold of $\op{SO}(3)\times 
\left(\R ^3\setminus\{ 0\}\right)\times\R ^2$. 
If $I_1<I_2<I_3$ then also the set of singular points of the 
constants of motion, corresponding to 
the critical circles, is a smooth codimension two algebraic 
submanifold of $\op{SO}(3)\times 
\left(\R ^3\setminus\{ 0\}\right)\times\R ^2$. The intersection of these 
submanifolds of special motions consists of the rotations of the 
sphere around a vertical axis which is equal to an axis of inertia, 
during which the point of contact is at rest. These motions define 
a submanifold of codimension four in $\op{SO}(3)\times 
\left(\R ^3\setminus\{ 0\}\right)\times\R ^2$. 

\subsection{Invariance under Rotations about the Vertical Axis}
If $j$ is vertical, then $j$ is invariant under the group of rotations 
around the vertical axis, which is the action of the 
$\op{E}(2)$\--symmetry group, the horizontal motion group, 
on $j$. This implies that the level set of $j$ is invariant 
under the transformations $(A,\,\omega )\mapsto (R\, A,\,\omega )$ 
with $R\in\op{SO}(2)$. The orbits of this action are equal to 
the fibers of the projection $(A,\,\omega )\mapsto (u,\,\omega )$, 
$u=A^{-1}\,\ e_3$, where the space $\op{S}\times\R ^3$ 
of the $(u,\,\omega )\in\R ^3\times\R ^3$ such that $\| u\| =1$ 
is viewed as the phase space for the $\op{E}(2)$\--reduced system. 
Because $T$ is $\op{E}(2)$\--invariant, it follows that 
the function $T_j$ on $\op{SO}(3)$, which represents $T$ on the 
$j$\--level set, is $\op{SO}(2)$\--invariant, a fact which 
can also be deduced directly from (\ref{Tj}). In particular the 
two\--dimensional regular level sets of $T_j$ are $\op{SO}(2)$\--invariant, 
and therefore are mapped by the projection 
$A\mapsto (u,\,\omega )$, $u=A^{-1}\, e_3$, $\omega =\omega _j(A)$, 
to smooth compact one\--dimensional algebraic submanifolds of 
$\op{S}\times\R ^3$. Each connected component of the 
regular level set in the $(u,\,\omega )$\--space therefore 
will be a closed curve, which implies that the motion in 
the $\op{E}(2)$\--reduced phase space is {\em periodic}. 

As a consequence, we can apply the reconstruction technique in Hermans  
\cite[Sec. 3.2]{joost} in order to obtain information about the flow 
in the full $(A,\,\omega ,\, a)$\--phase space $\op{SO}(3)\times\R ^3
\times\R ^2$ of the rolling body. There it is assumed that the symmetry 
group is compact, but the only thing which is needed is that 
the centralizer in the group of the shift element is 
a torus. Now the centralizer of the element 
$(B,\, b)\in\op{E}(2)\simeq\op{SO}(2)\times\R ^2$ 
is a circle subgroup of $\op{E}(2)$ when $B\neq I$, 
or $B$ is a nontrivial rotation around 
the vertical axis, whereas it is equal to the translation subgroup 
$\R ^2$ if $B=I$. Because the $E(2)$\--reduced 
phase space $\op{S}^2\times\R ^3$ is equal to the 
$\op{SO}(2)$\--reduced phase space of the $\R ^2$\--reduced 
phase space $\op{SO}(3)\times\R ^3$, we obtain the following conclusion. 
\begin{proposition}
Let $j$ be vertical. Each solution on a regular level set 
in the $(u,\,\omega )$\--space, the $\op{E}(2)$\--reduced 
phase space, is periodic. The corresponding 
rotational motion in the $(A,\,\omega )$\--space is quasiperiodic on 
an analytic two\--dimensional torus, depending analytically 
on the parameters $j$ and $T$.  

For each periodic solution in the 
$(u,\,\omega )$\--space such that the corresponding solution 
in the $(A,\,\omega )$\--space is not periodic with the same period, 
the motion in the full phase space 
$\op{SO}(3)\times\R ^3\times\R ^2$ 
is quasiperiodic on an analytic two\--dimensional torus, depending 
analytically on $j$ and $T$.  
In particular the point of contact $p(t)$ remains in a 
bounded subset of the plane in this case. 

If the rotational motion, the 
motion in the $(A,\,\omega )$\--space, is periodic with the 
same minimal period as the motion in the $(u,\,\omega )$\--space, 
then the translational motion, 
the motion of the point of contact $p(t)$, is equal to the superposition 
of a straight line motion with constant speed and a periodic 
motion with the same period as the motion in the 
$(u,\,\omega )$\--space. 
\label{vertmotion}
\end{proposition}
Some more information about the periodic solutions 
mentioned in Proposition \ref{vertmotion} is given in
Proposition \ref{pervert} below. Proposition 
\ref{perhorvert} implies that for at least half of 
these periodic rotational motions the motion of the point of contact 
$p(t)$ is actually periodic, with the same period. 

For more details on the behaviour of the sphere when the moment 
is vertical, see Kilin \cite[Sec. 3.3]{kilin}, which also contains 
computer pictures of orbits of the point of contact. 

\subsection{Fourier Series}
\label{Fouriersubsec}
In order to appreciate the statements in Proposition \ref{vertmotion}, 
we recall what it means that the motion in a suitable open subset $M$ 
of the phase space is quasiperiodic on analytic $r$\--dimensional tori, 
with analytic dependence on $s$ parameters, varying in some open 
subset $E$ of $\R ^s$. It means that 
there exists an analytic diffeomorphism 
$\Phi$ from $\left(\R ^r/\Z ^r\right)\times E$ to $M$ and an 
analytic function $\nu :E\to\R ^r$, such that the pull\--back 
$w=\Phi ^*\, v$ of the velocity field $v$ on $M$ is of the form 
$w(x,\,\epsilon )=(\nu (\epsilon ),\, 0)$, 
$x\in\R ^r/\Z ^r$, $\epsilon\in E$. This implies that the 
solution curves $\gamma$ in $M$ are of the form 
$\gamma (t)=\Phi\left( x_0+t\,\nu (\epsilon ),\,\epsilon\right)$, 
$t\in\R$. In our case $r=2$. 

If we apply this to the rotational motion, then it follows that 
the right hand side $r\, (A\,\omega )\times e_3$ in (\ref{pdot}) 
is of the form 
$\dot{p}\left( x_0+t\,\nu (\epsilon ),\,\epsilon\right)$, 
where $\dot{p}$ is an analytic mapping from $\left(\R ^2/\Z ^2\right)\times E$ 
to the horizontal plane $\R ^2$. Using Fourier expansion it 
follows that we can write 
\begin{equation}
\dot{p}\left( x_0+t\,\nu (\epsilon ),\,\epsilon\right) 
=\sum_{k\in\Z ^2}\, c_k(\epsilon )\, 
\op{e}^{2\pi\, i\,\langle x_0+t\,\nu (\epsilon ),\, k\rangle}, 
\label{Fourierpdot}
\end{equation}
in which the Fourier coefficients $c_k(\epsilon )$ depend 
analytically on the parameters $\epsilon$. The 
analyticity of the function $\dot{p}$ implies that the 
Fourier coefficients $c_k(\epsilon )$ are rapidly decreasing as 
$\| k\|\to\infty$. 
Formal termwise integration of (\ref{Fourierpdot}) 
would lead to 
\begin{eqnarray}
p(t)&=&p(0)+t\,\sum_{k\in\Z ^2,\,\langle\nu (\epsilon ),\, k\rangle =0}
\, c_k(\epsilon )\,\op{e}^{2\pi\, i\,\langle x_0,\, k\rangle}\\
&&+\sum_{k\in\Z ^2,\,\langle\nu (\epsilon ),\, k\rangle \neq 0}
\,\frac{c_k(\epsilon )}{2\pi\, i\,\langle\nu (\epsilon ),\, k\rangle}
\,\op{e}^{2\pi\, i\,\langle x_0+t\,\nu (\epsilon ),\, k\rangle} .
\label{Fourierp}
\end{eqnarray}

The coefficient of the linear term in $t$, the {\em secular term}, 
contains, apart from the term $c_0(\epsilon )$ (which is equal to the 
average of the function $\dot{p}$ over the torus), terms for nonzero 
$k$ if and only if the components $\nu _1(\epsilon )$ and $\nu _2(\epsilon )$ 
have a rational ratio. If we assume that $\nu _2(\epsilon )\neq 0$ 
and $\nu _1(\epsilon )/\nu _2(\epsilon ) =n_1/n_2$, with $n_1\in\Z$, 
$n_2\in\Z _{>0}$ and $\op{gcd}\left( n_1,\, n_2\right) =1$, then 
\[
t\,\langle\nu (\epsilon ),\, k\rangle =t\,\frac{\nu _2(\epsilon )}{n_2}
\,\left( n_1\, k_1+n_2\, k_2\right) ,
\]
which is equal to an integer if $t$ is equal to an integral 
multiple of $T:=n_2/\nu _2(\epsilon )$. This in turn implies 
that the function 
$t\mapsto\dot{p}\left( x_0+t\,\nu (\epsilon ),\,\epsilon\right)$ is 
periodic with period equal to $T$. For generic analytic functions 
$\nu _1(\epsilon )$, $\nu _2(\epsilon )$ the ratio 
$\nu _1(\epsilon )/\nu _2(\epsilon )$ is rational for a dense 
subset of parameter values, where the denominator $n_2$ and 
therefore also the period $T$ is unbounded in every nonvoid 
open subset of the parameter space. The coefficients of the 
secular term in general would have a correspondingly wild behaviour 
as a function of the parameters. 

In the case that the fraction 
$\nu _1(\epsilon )/\nu _2(\epsilon )$ is irrational, 
when the function 
\[
t\mapsto\dot{p}\left( x_0+t\,\nu (\epsilon ),\,\epsilon\right)
\] 
is not periodic, then the sum in (\ref{Fourierp}) over the $k\in\Z ^2$ 
such that $\langle\nu (\epsilon ),\, k\rangle \neq 0$ need not converge,  
due to the possibility that the denominators 
$\langle\nu (\epsilon ),\, k\rangle$ may become arbitrarily small. 
This problem can arise, despite the rapid decrease 
of the Fourier coefficients $c_k(\epsilon )$ when $\| k\|\to\infty$. 
If the fraction $\nu _1(\epsilon )/\nu _2(\epsilon )$ satisfies suitable diophantine inequalities, 
then the sum is convergent and defines a quasiperiodic function 
of $t$ on a two\--dimensional torus. 

Proposition \ref{vertmotion} implies that for the motion 
of Chaplygin's sphere when $j$ is vertical {\em 
none of the above 
complications occur:} no wild behaviour of the coefficients 
of the secular term and no problem with convergence of the 
Fourier series in (\ref{Fourierp}). Note that the conclusions 
of Proposition \ref{vertmotion} have not been obtained by means of 
an analysis with Fourier series, but by using the 
reconstruction technique in Hermans \cite[Sec. 3.2]{joost} instead. 

\subsection{Euler's Equations}
\label{EEsubsec}
In this subsection we give some explicit formulas, 
which among other show that if $j$ is vertical, then 
the rotational motion is determined by Euler's equations. 
An extensive discussion of the Euler top can be found in   
Cushman and Bates \cite[Ch. III]{CB}. 

If $j=j_3\, e_3$, then it follows from (\ref{Tj}) that 
\begin{equation}
2T=\langle A\,\omega ,\, j\rangle =j_3\,\langle A\,\omega ,\, e_3\rangle 
=j_3\,\langle\omega ,\, u\rangle
\label{omegau}
\end{equation}
Because (\ref{udot}) implies that 
$\langle\omega ,\,\dd{u}{t}\rangle =\langle\omega ,\, u\times\omega\rangle =0$, 
we obtain from (\ref{omegau}) that 
\[
\langle\dd{\omega}{t},\, u\rangle =\dd{}{t}
\langle\omega ,\, u\rangle =0.
\]
But then (\ref{omegadot}) implies that 
\begin{equation}
(I+\rho )\,\omega\times\omega =I\,\omega\times\omega 
=\dd{}{t}(I+\rho )\,\omega ,
\label{omegadotvert}
\end{equation}
which are Euler's equations with $I$ replaced by $I+\rho$. 

The equation $j_3\, e_3=j=A\, I_{\rho ,\, u}\,\omega$ yields 
in view of (\ref{omegau}) that
\[
j_3\, u=(I+\rho )\,\omega -\rho\,\langle\omega ,\, u\rangle\, u
=(I+\rho )\,\omega -\frac{2T\,\rho}{j_3}\, u,
\]
which leads to the formula 
\begin{equation}
u=\frac{j_3}{{j_3}^2+2T\,\rho}\, (I+\rho )\,\omega ,
\label{uomega}
\end{equation}
which expresses $u$ in terms of $\omega$. 
The differential equation $\dd{u}{t}=u\times\omega$ then leads to 
the differential equation 
\begin{equation}
\dd{u}{t}=\frac{{j_3}^2+2T\,\rho}{j_3}\, u\times J\, u,
\label{udotvert}
\end{equation}
for $u$ only, which also is an equation of Euler type.
Here $J=(I+\rho )^{-1}$, cf. (\ref{Jdef}). 

From (\ref{uomega}) and (\ref{omegau}) we obtain that
\begin{equation}
\langle u,\, J\, u\rangle 
=\frac{2T}{{j_3}^2+2T\,\rho},
\label{uellipsoid}
\end{equation}
which shows that the solutions $u(t)$ run on the intersection 
of the sphere $\langle u,\, u\rangle =1$ with an ellipsoid 
defined by the symmetric matrix $J$. 
Using (\ref{uomega}) it follows that also $(I+\rho )\,\omega (t)$ 
runs over the intersection of a sphere with an ellipsoid, explicitly 
given by  
\begin{eqnarray}
\langle (I+\rho )\,\omega ,\, (I+\rho )\,\omega\rangle 
&=&\left(\frac{{j_3}^2+2T\,\rho}{j_3}\right) ^2,
\label{omegasphere}\\
\langle \omega ,\, (I+\rho )\,\omega \rangle 
&=&2T\,\left( 1+2T\,\rho/{j_3}^2\right) . 
\label{omegaellipsoid}
\end{eqnarray}

The kinetic energy $T$ 
is given in terms of $u$ by 
\begin{equation}
2T=\frac{{j_3}^2\,\langle u,\, J\, u\rangle}
{1-\rho\,\langle u,\, J\, u\rangle} ,
\label{Tvert}
\end{equation}
cf. (\ref{uellipsoid}). The critical points of $T$ on the 
unit $u$\--sphere are the unit eigenvectors $u$ of 
$J$, with the eigenvalues $1/(\iota +\rho )$, 
with $\iota =I_1$, $I_2$ or $I_3$. The corresponding 
critical value is equal to  
\begin{equation}
T_{\scriptop{crit}}=\frac{{j_3}^2}{2\iota},
\label{critvert}
\end{equation}
which is in accordance with (\ref{Tcrit}) when $j=j_3\, e_3$. 
It is well\--known that the eigenvectors for the two 
extremal eigenvalues of $I+\rho$ are stable equilibrium 
points of the Euler equation, whereas the eigenvectors 
for the intermediate eigenvalue are linearly unstable 
equilibrium points of Euler's equation (\ref{udotvert}). 
See for instance Cushman and Bates \cite[p. 117]{CB}. 
The latter implies that, when $j$ is vertical, the 
critical circles of $T_j$ for the intermediate critical 
values are linearly unstable. 

\subsection{The Translational Motion}
Let $u=u(t)$ be the solution of (\ref{udotvert}) on the 
intersection curve of the ellipsoid (\ref{uellipsoid}) 
with the unit sphere $U$. The rotational velocity 
vector $\omega =\omega (t)$ is determined in terms 
of $u$ by means of (\ref{uomega}). 
The projection $A\mapsto u=A^{-1}\, e_3$, cf. (\ref{u}), 
exhibits the $T_j$\--level set in $\op{SO}(3)$ as a principal 
$\op{SO}(2)$\--bundle over the aforementioned curve on $U$, 
where $\op{SO}(2)$ is the group of rotations about the vertical 
axis, acting on $\op{SO}(3)$ by means of left multiplications. 
The projection intertwines the flow on the $T_j$\--level set in 
$\op{SO}(3)$ with the motion on the curve on $U$ determined 
by (\ref{udotvert}), where the flow on the $T_j$\--level set in 
$\op{SO}(3)$ is given by the differential equation 
(\ref{Adot}) in which $\omega =\omega (t)$ is determined  
by (\ref{uomega}). 

If $T_j$ is close to the intermediate critical level, then 
the solution $u=u(t)$ of (\ref{udotvert}) will stay for a long time 
near one of the unstable equilibria $u_{\scriptop{crit}}$ 
of (\ref{udotvert}) before it 
moves on to the other one. During this time $\omega =\omega (t)$ 
will stay close to the nonzero vector $\omega _{\scriptop{crit}}$ 
which is determined by (\ref{uomega}) with $u$ replaced by 
$u_{\scriptop{crit}}$. It follows that $A=A(t)$ will make many rotations 
during that time. This leads to the following conclusion. 

\begin{proposition}
Let ${\cal T}_{\scriptop{per}}$ 
denote the levels of $T_j$ such that the motion 
on the corresponding torus in $\op{SO}(3)$ is 
periodic with the same period the motion in the 
$(u,\,\omega )$\--space, cf. Proposition \ref{vertmotion}. 
Then ${\cal T}_{\scriptop{per}}$ is an infinite subset of $\R$ 
with the intermediate critical level $T={j_3}^2/2I_2$ as its 
only accumulation point. This accumulation point is approached by 
${\cal T}_{\scriptop{per}}$ both from above and from below. 
\label{pervert}
\end{proposition}

The translational motion, the motion of the point of contact $p$, 
is obtained by integrating the right hand side 
$r\, (A\,\omega )\times e_3$ of (\ref{pdot}). 
Here the vector $\theta :=A\,\omega$, the rotational velocity vector 
in space coordinates, the {\em herpolhode} in 
Poinsot's description of 
the Euler top, lies in view of (\ref{omegau}) in the fixed horizontal plane 
$\langle\theta ,\, e_3\rangle =2T/j_3$. 
We will use the discussion 
in Cushman and Bates 
\cite[II.7.2]{CB} of the construction of Poinsot. 
(Additional information can be found in Routh 
\cite[Art. 151, p.98 and pp. 471-473]{Routh}. 
An interesting fact is for instance that the herpolhode is 
always concave towards the interior, without inflexion points. 
In his Th\'eorie nouvelle de la
rotation des corps, 1834, Poinsot drew the herpolhode 
like a snake (= herpes in Greek), 
which therefore is misleading.) 

To begin with, the image of the 
$T_j$\--level set in $\op{SO}(3)$ under the projection 
$A\mapsto \theta =A\,\omega$ is invariant under the 
rotations around the vertical axis, and is therefore known if 
we know how $\langle\theta ,\,\theta\rangle =\langle\omega ,\,\omega\rangle$ 
varies as $\omega$ is coupled to $u$ by means of (\ref{omegau})
and $u\in U$ runs over the curve determined by (\ref{uellipsoid}). 

We will restrict ourselves to one of the two opposite connected 
components of the intersection of $U$ with the ellipsoid 
(\ref{uellipsoid}), 
which is the orbit of the motion on the $u$\--sphere. 
On it, the function $\langle\theta ,\,\theta\rangle$ has 
four critical points which all are nondegenerate. Two of these 
correspond to the maximal value $R^2_{\scriptop{max}}$ and two to 
the minimal value $R^2_{\scriptop{min}}$, where 
$R_{\scriptop{max}}>R_{\scriptop{min}}>2T/j_3$. To be more precise, 
if $T_{\scriptop{crit},\, i}={j_3}^2/2I_i$ denote the critical values 
of $T_j$, cf. (\ref{critvert}), then we have the following two cases. 
\begin{itemize}
\item[i)] $T_{\scriptop{crit},\, 3}<T<T_{\scriptop{crit},\, 2}$. 
We have $u_3\neq 0$ on the intersection of $U$ with the 
ellipsoid (\ref{uellipsoid}), 
and therefore $u_3$ has one sign on the $u$\--orbit. 
The critical points corresponding to the maximal value 
are equal to the two intersection points 
$\left(\pm u_1,\, 0,\, u_3\right)$ of the orbit with 
the coordinate plane $u_2=0$. 
The critical points corresponding to the minimal value 
are equal to the two intersection points 
$\left( 0,\,\pm u_2,\, u_3\right)$ of the orbit with 
the coordinate plane $u_1=0$. 
\item[ii)] $T_{\scriptop{crit},\, 2}<T<T_{\scriptop{crit},\, 1}$. 
We have $u_1\neq 0$ on the intersection of $U$ with the 
ellipsoid (\ref{uellipsoid}), 
and therefore $u_1$ has one sign on the $u$\--orbit. 
The critical points corresponding to the maximal value 
are equal to the two intersection points 
$\left( u_1,\, 0,\,\pm u_3\right)$ of the orbit with 
the coordinate plane $u_2=0$. 
The critical points corresponding to the minimal value 
are equal to the two intersection points 
$\left( u_1,\,\pm u_2,\, 0\right)$ of the orbit with 
the coordinate plane $u_3=0$. 
\end{itemize}
It follows that the image of the 
$T_j$\--level set in $\op{SO}(3)$ under the projection 
$A\mapsto \dot{p} =r\, A\,\omega\times e_3$ is equal to a 
circular annulus in the plane, with center at the origin. 
When $u(t)$ runs around the orbit in $U$ then it subsequently 
passes an intersection point with the coordinate plane 
$u_2=0$, then an intersection point with a second coordinate plane, 
then the other intersection point with the coordinate plane 
$u_2=0$, and finally the other intersection point with 
the second coordinate plane, before it closes. 
The corresponding point $\dot{p}$ in the annulus will 
then reach the outer circle, with a second order contact, 
then touch the inner circle, return to 
the outer circle at a point which is rotated over an angle $\alpha$ 
as compared to the first contact point with the outer circle, 
and then touch the inner circle for the second time before 
the curve in $U$ closes. 

If the rotational motion $A(t)$ is periodic with the same period 
as the motion $u(t)$ on $U$, then the third point of contact 
with the outer circle is equal to the first one, which means that 
$2\alpha$ is equal to an integral multiple of $2\pi$. There are 
two cases. 
\begin{itemize}
\item[a)] $\alpha$ itself is not an integral multiple of $2\pi$,  
which correspond to the case that the second point of contact with the outer circle lies opposite to the first one. In this case the  
$\dot{p}$\--orbit is symmetric about the origin, its time average 
is equal to zero and the motion of the point 
of contact $p(t)$ with the horizontal plane is periodic. 
In other words, the speed of the straight line motion 
in Proposition \ref{vertmotion} is equal to zero. 
\item[b)] $\alpha$ is equal to an integral multiple of $2\pi$,  
which correspond to the case that the second point of contact with the outer circle is equal to the first one. Equivalently, 
$\dot{q}(t)$ is periodic with a period equal to half 
the period of the motion $u(t)$ on $U$. 
In this case the $\dot{p}$\--orbit is not symmetric about the origin. 
\end{itemize}
The argument preceding Proposition \ref{pervert} yields that 
$\alpha$ is given by a smooth real\--valued function of $T$ which tends 
to $\pm\infty$ as $T\to T_{\scriptop{crit},\, 2}$. 

Moreover, if $T$ 
converges to the intermediate critical level $T={j_3}^2/2I_2$, 
then the inner boundary circle of the $\dot{p}$\--annulus shrinks to the 
origin. Because $u(t)$ stays for a long time near the 
critical point $\pm e_2$, to which the intersection point 
of the $u$\--orbit with 
the coordinate plane $u_1=0$ in case i) and $u_3=0$ in case ii) 
is close, the conclusion is that $\dot{p}(t)$ stays for a long time 
close to the inner boundary circle, running many times around it 
in the process. It follows that the time average of 
$\dot{p}(t)$, which is equal to the speed of the straight line motion in 
Proposition \ref{vertmotion}, converges to zero when $T$ 
converges to  the intermediate critical level $T={j_3}^2/2I_2$. 

We therefore arrive 
at the following conclusions. 
\begin{proposition}
Let ${\cal T}_{\scriptop{per}}'$ 
denote the levels of $T_j$ such that the motion 
of 
\[
\dot{q(t)}=r\, A(t)\,\omega (t)\times e_3
\] 
is  periodic with half the period of the motion in the 
$(u,\,\omega )$\--space. 
Then ${\cal T}_{\scriptop{per}}'$ is an infinite subset of $\R$ 
with the intermediate critical level $T={j_3}^2/2I_2$ as its 
only accumulation point. This accumulation point is approached by 
${\cal T}_{\scriptop{per}}'$ both from above and from below. 
If the speed of the straight line motion in 
Proposition \ref{vertmotion} is nonzero, then 
necessarily $T\in {\cal T}_{\scriptop{per}}'$. 
If $T\in {\cal T}_{\scriptop{per}}'$ 
converges to  the intermediate critical level 
$T={j_3}^2/2I_2$, then the speed of the corresponding 
straight line motion, cf. Proposition \ref{vertmotion}, 
converges to zero. 
\label{perhorvert}
\end{proposition}
We conjecture that for most  
values of $T\in {\cal T}_{\scriptop{per}}'$ 
the speed of the straight line motion in 
Proposition \ref{vertmotion} is {\em not} equal to zero. 
Here the word ``most'' can mean all except finitely many, or 
for generic values of $I_1,\, I_2,\, I_3$. 

If $T\notin{\cal T}_{\scriptop{per}}'$ approaches a value 
$T_0\in{\cal T}_{\scriptop{per}}'$ for which speed of the straight line motion in 
Proposition \ref{vertmotion} is not equal to zero, then the bounded 
area in which the quasiperiodic motion of $p(t)$ takes place 
``opens up to infinity'', and closes again to a bounded subset when 
$T$ has passed the value $T_0$. If the above conjecture holds, then this 
scenario takes place infinitely often when $T$ approaches the 
intermediate critical level $T={j_3}^2/2I_2$, but with the  
average speed of the point of contact converging to zero.  

\subsection{Chaplygin}
The remainder of Chaplygin \cite[\S 6]{chaplsphere}, 
starting with the 
sentence ``In addition, there is an exceptional case, $\cdots$'',  
consists of a discussion of the case that $j$ is vertical. 
This discussion contains several interesting observations, 
but it does not give the qualitative information 
about the motion as in Proposition \ref{vertmotion}.

\section{Bobylev's Sphere}
\label{srsec}
A rigid body is called a {\em solid of revolution} if it is 
dynamically symmetric with respect to all rotations $R$ 
about a given axis, for which we can take the vertical axis. 
This means that both the 
surface $S$ and the inertial tensor $I$ are invariant 
under such rotations $R$. Because for Chaplygin's sphere 
the surface $S$ already is invariant, it is a solid of 
revolution if and only if $I_1=I_2$. In this case 
Chaplygin's sphere is equal to Routh's sphere with the 
center of mass at the center of the sphere.
% (which in the 
%notation of \cite[Example 5.1]{sr} means that $\sigma _0=0$). 
This case has been studied by 
Bobylev \cite{Bobylev}, where in our case there is no gyroscope 
as mentioned in the title of \cite{Bobylev}. 

The assumption of having a solid of revolution implies that 
the equations of motion are invariant under the 
{\em right $\op{SO}(2)$\--action} 
\[
\left( A,\,\omega ,\, a\right)\mapsto\left( A\, R^{-1},\, 
R\,\omega ,\, R\, a\right) 
\]
of rotations $R$ about the vertical axis. The 
right $\op{SO}(2)$\--action also leaves the kinetic 
energy $T$ and the moment $j$ invariant. The quotients of 
the connected components of the regular $(j,\, T)$\--level 
surfaces by the right $\op{SO}(2)$\--action 
are diffeomorphic to circles, which implies that the 
solutions of the right $\op{SO}(2)$\--reduced system, 
for the regular levels of $(j,\, T)$, are 
periodic. The reconstruction method in \cite[Sec. 3.2]{joost} 
then leads to the following conclusions. Note that 
the vector $A\,\omega$ is invariant under the 
right $\op{SO}(2)$\--action, because 
$\left( A\, R^{-1}\right)\,\left( R\,\omega\right) =A\,\omega$. 

\begin{proposition} 
Suppose that two of the moments of inertia are equal to each other. 
Then the right $\op{SO}(2)$\--reduced rotational motion 
on the regular levels of $(j,\, T)$ is periodic, with a period 
which depends analytically on $j$ and $T$. In particular 
the vector $A\,\omega$, and therefore also 
$\dd{p}{t}=r\, (A\,\omega )\times e_3$ performs a periodic motion. 

It follows that the rotational motion 
on the regular levels of $(j,\, T)$ is quasiperiodic 
on two\--dimensional analytic tori, depending 
analytically on $j$ and $T$. 
Also, the motion of the point of contact $p$ is equal to the superposition 
of a straight line motion with constant speeed and a periodic 
motion with the same period as that of the right 
$\op{SO}(2)$\--reduced motion. 
\label{Bobylevprop}
\end{proposition}
As in the case when $j$ is vertical, cf. Proposition \ref{vertmotion}, 
there are no problems with secular terms or with the convergence 
of Fourier series for the motion of the point of contact $p$. In 
contrast with Proposition \ref{vertmotion}, we obtain here that 
$\dd{p}{t}$ is always periodic. 

\medskip

It follows from (\ref{F}) that  
\[
F=I_1\, x+I_3\, y\, z,\quad x=u_1\,\omega _1+u_2\,\omega _2,
\quad y=\omega _3,\quad z=u_3.
\]
%in the notation of \cite[(5.3), (5.4), (5.1)]{sr}. 
%In other words, 
%Jellet's integral, as mentioned in \cite[Example 5.1]{sr}, is 
%equal to $r\, F$. 
This is known as {\em Jellet's integral}, cf. Routh \cite[Art. 243]{Routh}. 

In the same notation, 
 %of \cite[(5.5), (5.4), (5.3), (5.1)]{sr}
the constant of motion $G$ of (\ref{G}) takes the form 
\[
G=\left( I_1+\rho\right) ^2\, w+\left( I_3+\rho\right) ^2\, y^2
-2\rho\, (x+y\, z)\left(\left( I_1+\rho\right)\, x+
\left( I_3+\rho\right)\, y\, z\right) 
+\rho ^2\, (x+y\, z)^2.
\]

Together with the kinetic energy we thus obtain 
three constants of motion in the four\--dimensional 
$\op{E}(2)\times\op{SO}(2)$\--reduced phase space, 
the ``fully reduced'' phase space.  
%of \cite[Sec. 5]{sr}. 
The regular level sets 
of all the constants of motion are algebraic curves, 
the periodic motion on which can be obtained by means of quadratures. 
We do not go into further details about this here.

\subsection{Chaplygin}
In his Introduction, Chaplygin \cite{chaplsphere} referred to 
the papers of Bobylev \cite{Bobylev} and Zhukovsky \cite{Zhukovsky} 
for the case that two moments of inertia are 
equal. In the beginning of \cite[\S 6]{chaplsphere}, Chaplygin 
wrote ``We will not treat the case when two or all three 
principal moments of inertia are equal, because the motion 
of such a sphere has already been investigated (see the Introduction).''
Apparently Chaplygin did not feel that the articles of Bobylev and Zhukovsky, 
which I have not seen, needed further comments. In particular 
I wonder whether Bobylev and/or Zhukovsky used the moment $j$ of 
the momentum around the point of contact as a constant of motion.  
In the paper \cite{chaplarea}, to which Chaplygin referred for the 
fact that $j$ is a constant of motion, 
there is no reference to Bobylev or Zhukovsky. 

The description of Bobylev's paper \cite{Bobylev} in 
the Fortschritte der Mathematik says: ``After the proposed integration, 
which can be performed with the help of the elliptic 
functions of Weierstrass, the author reaches the conclusion that the 
center of the sphere describes a curve which is enclosed between 
two parallel straight lines and has a periodic character, 
where it successively reaches the one and the other straight 
line with constant distances between the successive 
contact points on each of the straight lines.'' This 
corresponds to the description of the point of contact $p$ in 
Proposition \ref{Bobylevprop}.  

As observed before, 
Chaplygin's sphere with two equal moments of inertia is equal to 
Routh's sphere with the center of mass at the center of the sphere. 
In \cite[Art. 243]{Routh}, where Routh's sphere is treated, no 
special attention is paid to the case that the 
center of mass is at the center of the sphere. 

\section{An Invariant Volume Form}
\label{volsec}
We return to the general case of Chaplygin's sphere, with 
arbitrary moments of inertia, arbitrary total mass and 
radius of the sphere and arbitrary moment of 
the momentum around the point of contact. 
In the following lemma we use the notation of Subsection \ref{pairsubsec}. 
Note that $X(u)$ in (\ref{phidef}) is strictly positive when 
$\| u\| =1$, because the eigenvalues of the symmetric matrix 
$J=(I+\rho )^{-1}$ are equal to $1/\left( I_i+\rho\right)$, 
$i=1,\, 2,\, 3$,  
and therefore strictly smaller than $1/\rho$. 
\begin{lemma}
Consider the equations of motion $\dd{u}{t}=u\times\omega$, 
$\dd{ v}{t}= v\times\omega$ in the 
$(u,\, v )$\--space $\R ^3\times\R ^3$, in which $\omega$ 
is determined in terms of $u$ and $ v$ by {\em (\ref{omegausigma})}. 
Let $\Omega$ be the volume form in $\R ^6$ which is equal to 
$X(u)^{-1/2}$ times the Euclidean volume form, 
in which $X(u)$ is defined by {\em (\ref{phidef})}.  
Then $\Omega$ is invariant under the flow in $\R ^6$. 
\label{vollem}
\end{lemma}
\begin{proof}
The velocity field of the flow is equal to the vector field 
\[
\op{R}_{\omega}:\left( u,\, v\right)\mapsto
\left( u\times\omega ,\, v\times\omega\right) ,
\] 
in which $\omega =\omega (u,\, v )$ is determined by 
(\ref{omegausigma}). The divergence of this vector field is 
equal to the trace of the derivative, and therefore equal to the 
sum of the trace of the derivative $D_u$ of 
$u\mapsto u\times\omega (u,\, v )$ and the trace 
of the derivative $D_{ v}$ of 
$ v\mapsto  v\times\omega (u,\, v )$. 
Because the traces of the linear mappings 
$\delta u\mapsto\delta u\times\omega (u,\, v )$ and 
$\delta v\mapsto\delta v\times\omega (u,\, v )$ 
are equal to zero, the divergence is equal to 
$\op{trace}D^{\omega}_u+\op{trace}D^{\omega}_{ v}$, in which 
$D^{\omega}_u:\delta  u\mapsto u\times
\frac{\partial\omega (u,\, v )}
{\partial u}\,\delta u$ and 
$D^{\omega}_{ v}:\delta  v\mapsto v\times\frac{\partial\omega (u,\, v )}
{\partial v}\,\delta v$. 
From  (\ref{omegausigma}) we obtain that 
\begin{equation}
\omega (u,\, v )=J\, v +\psi (u,\, v )
\, J\, u,
\label{omusig}
\end{equation}
in which 
$\psi (u,\, v ):=Y(u,\, v )/X(u)$ and $Y(u,\, v )$, 
$X(u)$ are given by (\ref{chidef}), (\ref{phidef}), respectively. 

We now use that the trace of $\delta u\mapsto u\times J\,\delta u$ 
is equal to 
\[
\sum_{i=1}^3\,\langle u\times J\left( e_i\right),\, e_i\rangle 
=\sum_{i=1}^3\, 
\left(I_i+\rho\right) ^{-1}\,
\langle u\times e_i,\, e_i\rangle =0.
\]
It follows that $\op{trace}D^{\omega}_u$ is equal to the trace of the rank one 
mapping 
\[
\delta u\mapsto\left(\frac{\partial\psi (u,\, v )}{\partial u}\,\delta u
\right)\, u\times J\, u,
\]
and therefore equal to 
\begin{equation}
\frac{\partial\psi (u,\, v )}{\partial u}\,
\left( u\times J\, u\right) 
=\langle u\times J\, u,\, J\, v\rangle 
/X(u),
\label{div}
\end{equation}
because 
\[
\frac{\partial X(u)}{\partial u}\,
\left( u\times J\, u\right) 
=-2\langle u\times J\, u,\, J\, u\rangle =0.
\]
Similarly the trace of 
$\delta v\mapsto v\times J\,\delta v$ 
is equal to zero, and therefore the trace of $D^{\omega}_{ v}$ is equal to 
the trace of the rank one mapping 
\[
\delta v\mapsto\left(\frac{\partial\psi (u,\, v )}
{\partial  v}\,\delta v
\right)\, u\times J\, u,
\]
which is equal to zero because 
\[
\langle u,\, J\,\left( u\times J\, u\right)\rangle
=\langle J\, u,\, u\times J\, u\rangle =0. 
\]

The conclusion is therefore that the divergence of the vector field 
is equal to (\ref{div}). On the other hand the derivative 
$\op{R}_{\omega}X$ of the function $X$ in the direction of 
the vector field $\op{R}_{\omega}$ is equal to 
\[
-2\langle u\times\omega ,\, J\, u\rangle 
=-2\langle u\times J\, v ,\, J\, u\rangle 
=2\langle u\times J\, u,\, J\, v\rangle ,
\]
and therefore 
\[
\op{div}\,\op{R}_{\omega}=\frac12\, X^{-1}\,\op{R}_{\omega}\, X.
\]
It follows that 
\[
\op{div}\,\left( X^{-1/2}\,\op{R}_{\omega}\right) 
=-\frac12\, X^{-3/2}\,\op{R}_{\omega}X 
+X^{-1/2}\,\op{div}\,\op{R}_{\omega}=0,
\]
which completes the proof of the lemma. 
\end{proof}

Let $M$ be a smooth manifold of dimension $m$, $\Omega$ a smooth volume form on 
$M$ and $f$ a smooth function on $M$ such that $\op{d}\! f\neq 0$ 
at every point of the level set $M_c:=\{ x\in M\mid f(x)=c\}$. Then 
$M_c$ is a smooth $(m-1)$\--dimensional submanifold of $M$, 
and there is a unique volume form $\omega$ on $M_c$ such that 
\begin{equation}
\Omega _x\left( v_1,\,\ldots ,\, v_{m-1},\, v_m\right) 
=\omega _x\left( v_1,\,\ldots ,\, v_{m-1}\right)\, \op{d}\! f_x\left( v_m\right) 
\label{relvoldef}
\end{equation}
whenever $x\in M_c$, $v_i\in\op{T}_xM_c$ every $1\leq i\leq m-1$, 
and $v_m\in\op{T}_xM$. The volume form $\omega$ on $M_c$ is smooth 
and nonzero at every $x\in M$ where $\Omega _x\neq 0$. 
It is called the {\em relative quotient of $\Omega$ and  
$\op{d}\! f$} and denoted by $\omega =\Omega /\op{d}\! f$.  

Jacobi observed in \cite[10--14. Vorlesung]{Jacobi}, that if $v$ 
is a smooth vector field on $M$ such that $v\, f=0$ and its divergence 
${\cal L}_v\Omega /\Omega$ with respect to $\Omega$ is equal to zero, 
then the flow of $v$ leaves $M_c$, $\op{d}\! f$ and $\Omega$ invariant, 
and therefore $\omega$ as well. In other words, if $v_c$ denotes 
the restriction of $v$ to $M_c$, which is tangent to $M_c$, 
then the divergence ${\cal L}_{v_c}\omega /\omega$ of $v_c$ 
with respect to $\omega$ is equal to zero. 

Applying this principle succesively to the functions in the left 
hand sides in (\ref{SO3eq}), which are all invariant under 
the action $(u,\, v )\mapsto (R\, u,\,\ R\, v )$ of 
arbitrary rotations $R$, and using that the set determined 
by (\ref{SO3eq}) can be identified with $\op{SO}(3)$ if 
$j$ is not vertical, we arrive at the following corollary,    
where for vertical $j$ we can apply a continuity argument.  

\begin{corollary}
Let $\op{d}\! A$ be a Haar volume form on $\op{SO}(3)$, 
a volume form 
on $\op{SO}(3)$ which is invariant 
under right (or left) multiplications with elements of $\op{SO}(3)$.  
Then $ X ^{-1/2}\,\op{d}\! A$ is invariant under the 
$\op{R}_{\omega _j}$\--flow on $\op{SO}(3)$. 

On the regular level surfaces for the kinetic energy function 
$T_j$, the area form 
\[
 X ^{-1/2}\,\op{d}\! A/\op{d}\! T_j
\] 
is invariant under the $\op{R}_{\omega _j}$\--flow on $\op{SO}(3)$. 
\label{volcor}
\end{corollary}

If $M$ is a two\--dimensional smooth manifold, $v$ a nowhere 
vanishing smooth vector field on $M$, and $\alpha$ is a nowhere 
vanishing smooth area form on $M$ which is $v$\--invariant, then 
the fact that the three\--form $\op{d}\!\alpha$ is equal to zero 
on $M$ implies that 
\begin{equation}
0={\cal L}_v\,\alpha =\op{d}\left(\op{i}_v\alpha\right) 
+\op{i}_v\left(\op{d}\!\alpha\right) =\op{d}\left(\op{i}_v\alpha\right) ,
\label{calvarea}
\end{equation}
or that the nowhere vanishing one\--form 
\[
\beta :=\op{i}_v\alpha 
\]
is closed. Let $g(x)$ be the function which is obtained by 
integrating $\beta$ along a curve in $M$ starting at some base point 
and ending up at $x$. (This is an allowed procedure in the 
``integration by quadratures'' philosophy.) Then 
$\op{d}\! g=\beta$, hence $v\, g=0$ and the orbits of the 
$v$\--solution curves correspond to the level sets of $g$. 
Note that $\op{d}\! g=\beta$ is nowhere vanishing, which implies 
that the connected components of the level sets are smooth curves. 
Also note that the function $g$ is always globally defined on the 
universal covering space of $M$, but that on $M$ in general it will be 
a multi\--valued function, with the indeterminacy 
that $g(x)$ has to be replaced by 
$g(x)+\langle [\gamma ],\, [\omega ]\rangle$ if 
the curve ending up at $x$ is followed by a loop 
$\gamma$ which starts and ends at $x$. Here 
$[\gamma ]\in\op{H}_1(M)$ and $[\omega ]\in\op{H}^1(M)$ 
denote the homology and (de Rham) cohomology class of 
$\gamma$ and $\omega$, respectively. 

If $\Omega _0$ is a nonzero smooth volume form on a manifold $M$, 
then every smooth volume form $\Omega$ on $M$ is of the form 
$\Omega =\mu\,\Omega _0$, for a unique smooth function $\mu$. 
Therefore the search for an invariant volume form is a matter 
of finding the right factor (multiplier) $\mu$. 
If one has such a multiplier on an $n$\--dimensional manifold 
and one also has $n-2$ independent constants of motion 
$f_1,\,\ldots ,\, f_{n-2}$, then taking successively the 
relative quotient volume forms on the level manifolds 
$f_i=c_i$, one obtains multipliers on the $(n-j)$\--dimensional 
level sets $f_1=c_1,\,\ldots ,\, f_j=c_j$. For $j=n-2$ one 
finally obtains the ``last multiplier'' on the two\--dimensional 
level surface of all the constants of motion, to which one 
then can apply the above integration by quadratures. 
This is method for integration of vector fields by quadratures, 
which has been introduced in \cite[10-14. Vorlesung]{Jacobi}, is 
called {\em Jacobi's last multiplier method}. 

\medskip
If the two\--dimensional $M$ is compact and connected then the fact 
that $v$ has no zeros implies that $M$ is diffeomorphic to a torus, 
as we have observed before at the end of Subsection \ref{fixmoment}. 
Siegel \cite[Lemma 3 and 4]{Siegel} proved that there exists a smooth closed 
loop $C$ in $M$ such that $V$ is everywhere transversal to $C$ and 
that for every such $C$ and every $v$\--solution curve $\gamma$ 
there exists a $t>0$ such that $\gamma (t)\in C$. 
The transversality of $v$ to $C$ implies that if $x\in C$,  
and $t\mapsto\gamma (t,\, x)$ denotes the $v$\--solution curve with 
$\gamma (0,\, x)=x$, and $T(x)$ is the smallest $t>0$ such that 
$\gamma (t,\, x)\in C$, then $T$ depends smoothly on $x\in C$ and 
we obtain a smooth Poincar\'e map $P:x\mapsto\gamma (T(x),\, x):C\to C$. 
The restriction $\beta _C$ of the above one\--form $\beta =\op{i}_v\alpha$ 
to $C$ is a smooth one\--form on $C$ without zeros, and there 
an angle coordinate $\theta$ on $C$ such that $\op{d}\!\theta =\beta _C$, 
which is unique up to an additive constant.  
The $v$\--invariance of the area form $\alpha$ implies that 
$\beta$ is $v$\--invariant. In turn this implies that $\beta _C$ 
is invariant under $P$ and we conclude that $P(\theta )=\theta +c$, 
where $c$ is a constant. In other words, the return map is a rotation. 

By modifying the speed of the solution curves before they 
arrive at $C$, we can arrange that the return time $T(x)$ is a 
constant. In other words, there exists a strictly positive 
smooth function $f$ on $M$ (which we can choose to be non\--constant 
only in a thin strip at one side of $C$), such that $T(x)$ 
is equal to a constant if we replace $v$ by $f\, v$. Let $w_C$ 
be the unique tangent vector field of $C$ such that 
$\op{i}_{w_C}\alpha _C\equiv 1$. Because $w_C$ is invariant under $P$, 
we can carry $w_C$ around with the $v$\--flow and obtain an extension 
$w$ of $w_C$ which is a smooth vector field on $M$ and commutes with 
$f\, v$ by construction. Is is also clear that $w$ and $f\, v$ 
are everywhere linearly independent. 

We now recall the argument of Arnol'd and Avez 
\cite[Appendix 26]{ArnoldAvez} 
that the $f\, v$\--flow is quasiperiodic. 
It follows that if $\op{e}^{t\, v}$ denotes the flow after time 
$t$ of the vector field $v$, then 

\begin{equation}
(t,\, s)\mapsto \op{e}^{t\, f\, v}\circ\op{e}^{s\, w}
\label{vwaction}
\end{equation} 
defines an action of $\R ^2$ on $M$. Because of the linearly 
independence of $f\, v$ and $w$, the orbits are open subsets of $M$. 
because the orbits form a partition of $M$ and $M$ is connected, 
there is only one orbit, equal to $M$. In other words, the 
action is transitive. If $\op{e}^{t\, f\, v}\circ\op{e}^{s\, w}(x)=x$ 
for some $x\in M$ then $\op{e}^{t\, f\, v}\circ\op{e}^{s\, w}(x)=x$ 
for every $x\in M$. The {\em period lattice} 
\[
\Pi :=\{ (t,\, s)\in\R ^2\mid \op{e}^{t\, f\, v}\circ\op{e}^{s\, w}=1\}
\]
is a discrete additive subgroup of $\R ^2$, and because for each $x\in M$ 
the mapping $(t,\, s)\mapsto \op{e}^{t\, f\, v}\circ\op{e}^{s\, w}(x)$ 
induces a diffeomorphism from $\R ^2/\Pi$ onto the compact 
manifold $M$, the conclusion is that the lattice $\Pi$ is 
two\--dimensional. In the coordinates with respect to a
$\Z$\--basis of $\Pi$, $\R ^2/\Pi$ is equal the standard 
torus $\R ^2/\Z ^2$. Because in these coordinates the vector 
field $f\, v$ is constant, we arrive at the conclusion 
that the $f\, v$\--flow is quasiperiodic on a two\--dimensional torus.  

It is clear from its introduction that the function $f$ is far 
from unique. Actually, Kolmogorov \cite{Kol} proved that 
if the rotation number of the Poincar\'e map $P$ 
satisfies suitable diophantine inequalities, ensuring that 
it cannot be approximated too rapidly by means of rational numbers, 
then also the $v$\--flow, without the time\--reparametrizing 
factor $f$, is quasiperiodic. He also showed that  
for this conclusion the diophantine inequalities for the 
rotation number are essential, in the sense that 
in general it is not sufficient to assume that the 
rotation number of $P$ is irrational. 

In Corollary \ref{commvvcor} we will obtain that for the Chaplygin sphere 
the rotational motion is quasiperiodic on two\--dimensional 
tori, if the reparametrization of the time corresponds to 
multiplication of the vector field by the specific 
function $f=X(u)^{1/2}$, 
where $X(u)$ is given by (\ref{phidef}).  

\subsection{Chaplygin}
The last part of Chaplygin \cite[\S 2]{chaplsphere}, starting 
with ``To solve the problem completely, \ldots '', contains the 
proof of Lemma \ref{vollem}, followed by the conclusion, in one line, 
that Jacobi's last multiplier method can be applied 
in order to solve the equations of motion by quadratures. 
Apparently at the time of \cite{chaplsphere} this method 
was so well\--known, that no further explanations or 
references were needed. 

\section{Two Commuting Vector Fields}
\label{commvvsec}
\subsection{The Second Vector Field}
It follows from (\ref{RnuT}) that the tangent spaces of the 
level surfaces in $\op{SO}(3)$ of the function $T_j$ are 
spanned by the vector fields $\op{R}_{\omega _j}$ and 
$\op{R}_{(I+\rho )\,\omega _j}$. It is therefore natural 
to investigate the divergence of the vector field 
$\op{R}_{(I+\rho )\,\omega _j}$ with respect to the area form 
$X^{-1/2}\,\op{d}\! A/\op{d}\! T_j$, in analogy with 
Corollary \ref{volcor}. Note that $I+\rho =J^{-1}$, cf. (\ref{Jdef}). 

\begin{lemma}
Consider the vector field $\op{R}_{\nu}$ 
in the $(u,\, v )$\--space $\R ^3\times\R ^3$, 
defined by $\op{R}_{\nu}\, u:=u\times\nu$, 
$\op{R}_{\nu}\,  v := v\times\nu$, where 
$\nu =\nu (u,\, v ):=(I+\rho )\,\omega (u,\, v )$ and 
$\omega =\omega (u,\, v )$ 
is determined in terms of $u$ and $ v$ by {\em (\ref{omegausigma})}. 
Let $X(u)$ be defined by {\em (\ref{phidef})} and let 
$\Omega$ be the volume form in $\R ^6$ which is equal to 
$X(u)^{-1/2}$ times the Euclidean volume form.
Then the divergence of $\op{R}_{\nu}$ with repect 
to $\Omega$ is equal to zero. 

It follows that if $\op{d}\! A$ is a 
Haar volume form on $\op{SO}(3)$, then the divergence of 
the vector field $\op{R}_{(I+\rho )\,\omega _j}$ on $\op{SO}(3)$ 
with respect to $X^{-1/2}\,\op{d}\! A$ is equal to zero.  
Also, the divergence is equal to zero of 
the vector field $\op{R}_{(I+\rho )\,\omega _j}$ on any 
regular level surfaces of $T_j$, with respect to the 
area form $X^{-1/2}\,\op{d}\! A/\op{d}\! T_j$. 
\label{vollem2}
\end{lemma}
\begin{proof}
The proof follows the same lines as the proof of Lemma \ref{vollem}, 
The calculations are  actually somewhat easier this time, because 
the expression 
\[
\nu (u,\, v )=(I+\rho )\,\omega (u,\, v )= v +\psi (u,\, v )\, u
\]
for $\nu$ is simpler than the formula (\ref{omusig}) 
for $\omega$. 

The divergence of $\op{R}_{\nu}$ is equal to 
$\op{trace}D^{\nu}_u+\op{trace}D^{\nu}_{ v}$, in which 
$D^{\nu}_u$ and $D^{\nu}_{ v}$ is equal to the trace of 
the derivative of $u\times\nu (u,\, v )=u\times v$ and 
$ v\times\nu (u,\, v )=\psi (u,\, v )\cdot 
 v\times u$ with respect to $u$ and $ v$, respectively. 
It follows that the divergence is equal to the trace of the 
rank one linear mapping 
\[
\delta v\mapsto 
\left(\frac{\partial\psi (u,\, v )}{\partial v}\,\delta v\right)  
\, v\times u,
\]
and therefore equal to 
\[
\frac{\partial\psi (u,\, v )}{\partial v}\, ( v\times u)
=\langle u,\, J( v\times u)\rangle /X(u).
\]

On the other hand
\[
\op{R}_{\nu}X =\, -2\langle u\times\nu ,\, J\, u\rangle 
=\, -2\langle u\times v ,\, J\, u\rangle
=2X\,\op{div}\,\op{R}_{\nu},
\]
which in the same way as at the end of the proof of 
Lemma \ref{vollem} implies that 
\[
\op{div}\,\left(X^{-1/2}\,\op{R}_{\nu}\right) 
=-\frac12\, X^{-3/2}\,\op{R}_{\nu}X 
+X^{-1/2}\,\op{div}\,\op{R}_{\nu}=0,
\]
or that the divergence of $\op{R}_{\nu}$ with respect to $\Omega$ 
is equal to zero. 

The statements about the volume form on $\op{SO}(3)$ and the area form 
on the level sets of $T_j$ follow in the same way as Corollary 
\ref{volcor} follows from Lemma \ref{vollem}. 
\end{proof}

\begin{lemma}
Let $\op{d}\! A$ be the Euclidean volume form on $\op{SO}(3)$. 
Then 
\[
\left(\op{d}\! A/\op{d}\! T_j\right) 
\left(\op{R}_{\omega _j},\,\op{R}_{(I+\rho )\,\omega _j}\right) =1. 
\]
\label{areaomeganu}
\end{lemma}
\begin{proof}
It follows from \ref{RnuT} that if $\nu =\omega _j\times (I+\rho )\,\omega _j$, 
then 
\[
\op{R}_{\nu}T_j=\langle\omega _j\times (I+\rho )\,\omega _j,\, 
\omega _j\times (I+\rho )\,\omega _j\rangle =
\op{d}\! A\left(\omega _j,\, (I+\rho )\,\omega _j,\, 
\omega _j\times (I+\rho )\,\omega _j\right) ,
\]
which implies the statement of the lemma in view of the 
defining equation (\ref{relvoldef}) of the relative quotient 
of a volume form and the total derivative of a function. 
\end{proof}

\begin{proposition}
Let the factor $X(u)$ be defined by {\em (\ref{phidef})}. 
Define the vector fields $\xi$ and $\eta$ on $\op{SO}(3)$ by 
$\xi :=X(u)^{1/2}\,\op{R}_{\omega _j}$ and 
$\eta :=X(u)^{1/2}\,\op{R}_{(I+\rho )\,\omega _j}$, respectively. 
Then the vector fields $\xi$ and $\eta$ commute. 

Define the area form $\alpha _j$ and the 
one\--forms $\beta$ and $\gamma$ on the 
regular level surfaces of $T_j$ by 
$\alpha _j:=X(u)^{-1}\,\op{d}\! A/\op{d}\! T_j$, 
$\beta :=\op{i}_{\xi}\alpha _j$ and $\gamma :=\op{i}_{\eta}\alpha _j$, 
where $\op{d}\! A$ is the Euclidean volume form on $\op{SO}(3)$. 
Then $\alpha _j(\xi ,\,\eta )=1$, $\alpha _j=\beta\wedge\gamma$ and 
the one\--forms $\beta$ and $\gamma$ are closed. The area form 
$\alpha _j$ and the one\--forms $\beta$ and $\gamma$ are invariant 
under the flow of both vector fields $\xi$ and $\eta$. 
\label{commvvprop}
\end{proposition}
\begin{proof}
It follows from (\ref{calvarea}) 
and the fact that the 
divergence of $v$ with respect to $\alpha$ is equal to 
zero if and only if the one\--form $\op{i}_v\alpha$ is closed. 
Because for any functions $f$ on $M$ we have that 
$\op{i}_{f\, v}\, g\,\alpha =f\, g\, \op{i}_v\alpha$, it follows that 
for every nowhere vanishing smooth function $f$ we have that 
the divergence of $v$ with respect to $\alpha$ is equal to 
zero if and only if the one\--form 
$\op{i}_{f\, v}\left( f^{-1}\,\alpha\right)$ is closed. 

If we apply this with $v=\op{R}_{\omega _j}$,  
\[
\alpha =X^{-1/2}\,\op{d}\! A/\op{d}\! T_j=X^{1/2}\,\alpha _j,
\]
and $f=X^{1/2}$, then it follows from Corollary 
\ref{volcor} that $\beta =\op{i}_{\xi}\,\alpha _j$ 
is closed, or equivalently that 
${\cal L}_{\xi}\alpha _j=0$. In a similar manner it follows from 
Lemma \ref{vollem2} that $\gamma =\op{i}_{\eta}\,\alpha _j$ is closed, 
or equivalently ${\cal L}_{\eta}\omega _j =0$. 
  
On the other hand Lemma \ref{areaomeganu} implies that 
$\alpha _j(\xi ,\,\eta )=1$, from which it follows in turn that 
\begin{eqnarray}
0&=&{\cal L}_{\xi}\,\alpha _j(\xi ,\,\eta )
=\alpha _j\left(\xi ,\,\left[\xi ,\,\eta\right]\right) ,
\label{LXalphaXY}\\
0&=&{\cal L}_{\eta}\,\alpha _j(\xi ,\,\eta )
=\alpha _j\left(\left[\eta ,\,\xi\right] ,\,\eta\right) =
\alpha _j\left(\eta ,\,\left[\xi ,\,\eta\right]\right) .
\label{LYalphaXY}
\end{eqnarray}
Here we have used in (\ref{LXalphaXY}) that 
${\cal L}_{\xi}\alpha _j=0$, ${\cal L}_{\xi}\xi 
=\left[\xi ,\,\xi\right] =0$, 
and ${\cal L}_{\xi}\eta =\left[\xi ,\,\eta\right] =0$, whereas in 
in (\ref{LYalphaXY}) we have used that 
${\cal L}_{\eta}\alpha _j=0$, ${\cal L}_{\eta}\xi 
=\left[\eta ,\,\xi\right] =
\, -\left[\xi ,\,\eta\right]$, 
${\cal L}_{\eta}\eta =\left[\eta ,\,\eta\right] =0$, 
and the antisymmetry of $\alpha _j$. It follows from 
(\ref{LXalphaXY}) that $\left[\xi ,\,\eta\right]$ is a multiple of $\xi$ 
and from (\ref{LXalphaXY}) that $\left[\xi ,\,\eta\right]$ 
is a multiple of $\eta$. 
Because $\xi$ and $\eta$ are everywhere linearly independent on the 
regular level sets of $T_j$, it follows that 
$\left[\xi ,\,\eta\right] =0$ 
there. Because the regular level sets are dense, it follows by 
continuity that the vector fields $\xi$ and $\eta$ commute on 
all of $\op{SO}(3)$.  
\end{proof}

\begin{corollary}
If the rotational motion is parametrized by a time variable 
$\tau$ which is related to the time $t$ by 
$\dd{\tau}{t}=X(u(t))^{-1/2}$, then the rotational motion 
on the regular level sets is quasi\--periodic on 
two\--dimensional analytic tori, depending analytically 
on the parameters $j$ and $T$. 
\label{commvvcor}
\end{corollary}
\begin{proof}
We have that 
\[
\dd{A}{\tau}=X^{1/2}\,\dd{A}{t}=X^{1/2}\op{R}_{\omega _j}A=\xi\, A.
\]
The conclusion of the corollary follows from the discussion of the 
action (\ref{vwaction}) with $M$, $v$ and $w$ replaced by 
a regular level surface, $\xi$ and $\eta$, respectively. 
\end{proof}

\begin{question}
Is there a proof of Corollary \ref{volcor}, 
Lemma \ref{vollem2} and Lemma \ref{areaomeganu}, 
and therefore also of Proposition \ref{commvvprop}, 
which is based on a general principle, in the same way 
as Proposition \ref{jprop}  
follows from Noether's principle for 
nonholonomic systems in Lemma \ref{noetherlem}? 
\end{question}

\subsection{A Zero Average}
According to (\ref{pdot}), the time derivative of the $j$\--component 
$\langle p,\, j\rangle$ of the point of contact $p$ of the sphere 
with the horizontal plane is equal to $r$ times the quantity  
We begin with 
\begin{eqnarray*}
&&\langle (A\,\omega )\times e_3,\, j\rangle 
=\langle A\,\omega ,\,e_3\times j\rangle 
=\langle\omega ,\,\left( A^{-1}\, e_3\right)\times A^{-1}\, j\rangle\\
&&=\langle\omega ,\, u\times v\rangle 
=\langle\omega ,\, u\times I\,\omega\rangle =\, -\op{det}\left( 
u,\,\omega ,\, I\,\omega\right) . 
\end{eqnarray*}
Here we have used (\ref{u}) and (\ref{sigmadef}) in the third equation, 
and in the fourth equation we have used (\ref{omegausigma}) and 
(\ref{Iu}), together with the facts that $u\times u=0$ and 
$\omega$ is orthogonal to $u\times\omega$. We therefore obtain in view 
of (\ref{pdot}) that 
\begin{equation}
\dd{}{\tau}\langle p(\tau ),\, j\rangle 
=\, -r\, X (u)^{1/2}\,\op{det}(u,\,\omega ,\, I\,\omega ).  
\label{detC}
\end{equation}

\begin{lemma}
Let $M$ be a connected component of a regular 
level set of $T_j$ in $\op{SO}(3)$. 
Let $\op{det}\left( u,\,\omega _j,\, I\,\omega _j\right)$ 
be viewed as a function on $M$. 
Then, for any continuous function $f$ on $U$, the integral 
of $f(u)\,\op{det}\left( u,\,\omega _j,\, I\,\omega _j\right)$ 
over $M$ with respect to the area form $\alpha _j$, 
cf. Proposition \ref{commvvprop},  
is equal to zero.
\label{intzerolem}
\end{lemma}

\begin{proof}
We have 
\begin{eqnarray*}
\xi\, u\times \eta\, u&=&X(u)\, \left( u\times\omega _j\right) 
\times\left( u\times (I+\rho )\,\omega _j\right)\\
&=&-X(u)\,\langle\omega _j,\, u\times (I+\rho )\,\omega _j\rangle\, u 
=X(u)\,\op{det}\left( u,\,\omega _j,\, I\,\omega _j\right)\, u ,
\end{eqnarray*}
because $\langle u,\, u\times (I+\rho )\,\omega\rangle =0$ and 
therefore the $\omega _j$\--term drops out. Because $\alpha _j(\xi ,\,\eta )=1$, 
it follows that $\alpha _j$ is equal to the pull\--back of 
$1/X(u)\,\op{det}\left( u,\,\omega _j,\, I\,\omega _j\right)$ 
times the standard area form $\op{d}\! _2\, u$ on $U$ by means of the 
projection $\pi :A\mapsto u=A^{-1}\, e_3$ from the level set $M$ to $U$. 
Now we have, for any smooth mapping $\pi$ from a compact 
oriented manifold $M$ to a compact oriented manifold $U$, and any 
volume form $\Omega$ on $U$, the formula 
\begin{equation}
\int_M\,\pi ^*\Omega =\op{deg}(\pi )\cdot\int_U\Omega ,
\label{degree}
\end{equation}
see for instance Guillemin and Pollack \cite[p. 188]{GP}. In our case 
\[
f(u)\,\op{det}\left( u,\,\omega _j,\, I\,\omega _j\right)\,\alpha _j 
=\pi ^*\left(\frac{f(u)}{X(u)}\,\op{d}\! _2u\right) .
\]
Moreover, the degree of $\pi$ is equal to zero, because $\pi (M)\neq U$, 
see for instance the description at the end of Remark \ref{branchrem}.   
Therefore the conclusion of the lemma is obtained by applying 
(\ref{degree}) to $\Omega =\frac{f(u)}{X(u)}\,\op{d}\! _2u$. 
\end{proof}
It follows from Lemma \ref{intzerolem} that 
the average of the right and side of (\ref{detC}) over the level 
surface, with respect to the area form $\alpha _j$, is equal to zero. 
In view of Corollary \ref{commvvcor} we can apply (\ref{Fourierp}) with $t$ replaced by $\tau$ and $p(t)$ replaced by $\langle p(\tau ),\, j\rangle$, 
in the case that the rotational motion 
is not periodic. In this case the coefficient of the secular term 
is equal to $c_0(\epsilon )$, which is equal to the average of 
the right hand side of (\ref{detC}) over the level set, with respect 
to the area form $\alpha _j$. This leads to the following conclusion. 
\begin{corollary}
Assume that the rotational motion on the regular level set is 
nonperiodic and that the series in (\ref{Fourierp}), with $t$ replaced by $\tau$ and $p(t)$ replaced by $\langle p(\tau ),\, j\rangle$, is uniformly convergent. 
Then the function $\tau\mapsto\langle p(\tau ),\, j\rangle$ 
is quasiperiodic on a two\--dimensional torus. 
In particular, $\langle p(\tau ),\, j\rangle$ 
remains bounded in this case. 
\label{pjqpcor}
\end{corollary}
Note that the series mentioned in Corollary 
\ref{pjqpcor} converges uniformly when the 
irrational ratio $\nu _1(\epsilon )/\nu _2(\epsilon )$ 
mentioned after (\ref{Fourierp}) 
is sufficiently slowly approximated by rational numbers. 
The set of irrational numbers for which this happens has 
full Lebesque measure on the real axis.  

\begin{question}    
Do the complications with secular terms 
when the rotational motion is periodic, and convergence 
of Fourier series for nonperiodic rotational motions, as 
discussed after (\ref{Fourierp}), really occur? 
\end{question}

For the critical circles we have the simplification that 
$\dd{p}{t}$ is constant, cf. (\ref{pcrit}). When $j$ is 
vertical, and when two of the moments of inertia are equal, we can 
reconstruct the full motion by means of Lie group techniques 
from a periodic motion, which implies that the aforementioned 
complications do not arise. See 
Proposition \ref{vertmotion} and Proposition \ref{Bobylevprop}, 
respectively. 

In the general case I have neither been able to 
find any analogous special features which would eliminate the 
complications with the primitives of the Fourier series, 
nor did I find a proof that the complications really occur. 
See also Question \ref{translthetaq}.

\subsection{Chaplygin}
The themes of Section \ref{commvvsec} do not occur in 
Chaplygin \cite{chaplsphere}. The only exception may be 
that the existence of two closed one\--forms 
$\beta$ and $\gamma$ such that 
$\op{i}_{\xi}\beta =0$ and $\op{i}_{\xi}\gamma =1$, cf. 
Proposition \ref{commvvprop}, is implicitly contained in  
the formulas in Chaplygin \cite[\S 3, (29)]{chaplsphere}, 
cf. Remark \ref{iXbeta'rem}. As discussed after Remark \ref{iXbeta'rem} 
below, 
the vector field $\xi$ corresponds to an explicitly determined constant 
vector field on the Jacobi variety $J(C)$ of a hyperelliptic curve $C$. In 
Remark (\ref{Yjacrem}) we determine the constant vector field 
on $J(C)$ to which $\eta$ corresponds.  

\section{Some Simplifications}
\label{simplesec}

\subsection{A Polynomial System}
\label{phi0subsec}
The square root in the factor $ X (u)^{1/2}$ 
in front of the commuting vector fields $\xi$ and $\eta$ in 
Proposition \ref{commvvprop} becomes double\--valued if 
we extend the vector field holomorphically into the 
complex domain. One can make it single\--valued by passing 
to the space of $(z,\, u,\, v )\in\C\times\C ^3\times\C ^3$ 
where $(u,\, v )$ still satisfy (\ref{SO3eq}) and (\ref{Tusigma}) 
and the new variable $z$ is coupled to $u$ by means of the equation 
\begin{equation}
z^2=X(u)=\rho ^{-1}-\langle u,\, J\, u\rangle .
\label{zu}
\end{equation} 
This means that we pass to the branched (=ramified) double 
covering of the complexification of the level surface, 
which branches along the complex one\--dimensional submanifold 
(complex curve) defined by the equation $ X (u)=0$. 
As observed in front of Lemma \ref{vollem}, we have $ X (u)>0$ 
when $u$ is real, which implies that the curve defined by 
the equation $ X (u)=0$ has no real points. 

We have that $\xi =\op{R}_{z\,\omega}$, 
$\eta =\op{R}_{z\, (I+\rho )\,\omega}$, where  
\begin{eqnarray}
z\,\omega &=&z\, J\, v +
z^{-1}\, Y(u,\, v )\, J\, u,\quad\mbox{\rm and}
\label{zX}\\
z\, (I+\rho )\, \omega &=&z\, v +
z^{-1}\, Y(u,\, v )\, u,
\label{zY}
\end{eqnarray}
in which we have used the abbreviation (\ref{chidef}). 
We recall that $\op{R}_{\nu}\, u=u\times\nu$ and 
$\op{R}_{\nu} v = v\times\nu$. 
The actions on $z$ are given by 
$2z\, \xi\, z=\xi\, z^2=\, -2z\,\langle u\times J\, v ,\, 
J\, u\rangle$ and 
$
\eta\, z=\, -\langle u\times v ,\, J\, u\rangle ,
$
because 
$\langle u\times J\, u ,\, 
J\, u\rangle =0$ and 
$\langle u\times u ,\, 
J\, u\rangle =0$. 
It follows that the vector fields $\xi$ and $\eta$ are rational. 
The functions $\xi\, z$ and $\eta\, z$ are regular,  
but $\xi\, u$ and $\xi\, v$ seem to have poles along $z=0$. 

Combining the kinetic energy equation in the form 
(\ref{phipsichi}) with (\ref{zu}), we see that 
$Y(u,\, v )=\pm\, z\, Z( v )^{1/2}$, 
and it follows that 
$z\,\omega =z\, J\, v 
\pm\, Z( v )^{1/2}\, J\, u$
and 
$z\, (I+\rho )\, \omega =z\, v 
\pm\, Z( v )^{1/2}\, u$. 
Note that the sign choice in $\pm\, Z( v )^{1/2}$ 
is coupled to the choice of the sign of $z$ by means of (\ref{chidef}). 

However, the vector fields are still double\--valued at $z=0$, 
where also the manifold of the solutions $((u,\, z),\, v )$ of the 
equations (\ref{SO3eq}), (\ref{phipsichi}) and (\ref{zu}) 
is singular. These singularities can be resolved by introducing 
one more variable $\zeta$ which is coupled to $v$ by means of the equation 
\begin{equation}
\zeta ^2=Z(v)=2T-\langle v,\, J\, v\rangle .
\label{zetav}
\end{equation} 
The kinetic energy equation $Y(u,\, v)^2=X(u)\, Z(v)=
z^2\,\zeta ^2$ leads to $Y(u,\, v)=\pm z\,\zeta$. 
With the choice of the {\em minus} sign, 
\begin{equation}
Y(u,\, v)=\, -z\,\zeta ,
\label{zzetauv}
\end{equation}
the vector fields $\xi =\op{R}_{z\,\omega}$ and 
$\eta =\op{R}_{z\, (I+\rho)\,\omega}$ are given by 
\begin{eqnarray}
\xi\, u&=&z\, u\times J\, v 
-\zeta\, u\times J\, u,\label{Xuz}\\
\xi\, v&=&z\, v\times J\, v 
-\zeta\, v\times J\, u,\label{Xsigmaz}\\
\xi\, z&=&-\langle u\times J\, v,\, J\, u\rangle 
=\op{det}(u,\, J\, u,\, J\, v),\label{Xz}\\
\xi\,\zeta &=&\langle v\times J\, u,\, J\, v\rangle 
=\op{det}(v,\, J\, u,\, J\, v),\label{Xzeta}
\end{eqnarray}
and 
\begin{eqnarray}
\eta\, u&=&z\, u\times  v ,\label{Yuz}\\
\eta\, v&=&-\zeta\, v\times u,\label{Ysigmaz}\\
\eta\, z&=&-\langle u\times v,\,  J\, u\rangle  
=\op{det}(v,\, u,\, J\, u), \label{Yz}\\
\eta\,\zeta &=&\langle v\times u,\, J\, v\rangle 
=\op{det}(v,\, u,\, J\, v),\label{Yzeta}
\end{eqnarray}
respectively. 

The equations (\ref{Xuz})---(\ref{Yzeta}) 
define polynomial vector fields 
$\xi$ and $\eta$ in $\C ^8$, 
homogeneous of degree three. 
The vector fields $\xi$ and $\eta$ in $\C ^8$ commute with 
each other. Finally both vector fields $\xi$ and $\eta$ 
are divergence\--free, and have the six functions 
\begin{equation}
\begin{array}{ccc}
f_1:=\langle u,\, u\rangle ,&
f_2:=\langle u,\, v\rangle, &
f_3:=\langle v,\, v\rangle ,\\
f_4:=\langle u,\, J\, u\rangle +z^2,&
f_5:=\langle u,\, J\, v\rangle +z\,\zeta ,&
f_6:=\langle v,\, J\, v\rangle +\zeta ^2
\end{array}
\label{fi}
\end{equation}
as constants of motion. (We have chosen the minus sign in 
(\ref{zzetauv}) in order to get a pkus sign in $f_5$.) 
All these statements are true 
when $J$, which appears as the parameter in 
(\ref{Xuz})---(\ref{Yzeta}), 
is a symmetric $3\times 3$\--matrix. 

For the generic values of 
the $f_i$, the equations (\ref{fi}) define a smooth 
two\--dimensional affine algebraic variety $M$ in $\C ^8$, 
on which $\xi$ and $\eta$ are commuting vector fields. 
Note that Chaplygin's sphere corresponds to the case that 
\begin{equation}
f_1=1, \; f_2=j_3,\; f_3=\| j\| ^2,\;  
f_4=1/\rho , \; f_5=0,\; f_6=2T. 
\label{Mdef}
\end{equation}
In contrast with the real case, the complex 
surface is not compact and therefore we cannot conclude that 
the flows of $\xi$ and $\eta$ with complex times lead to 
an identification of $M$ with a complex torus. Also, the 
flows of $\xi$ and $\eta$ on $M$ with complex times 
will not be complete in the sense that these are not 
defined for all complex times and therefore do not define 
an action of $\C ^2$ on $M$. In Proposition \ref{hatMcprop} 
we will obtain a completion $\hat{M}$ of $M$ which is 
isomorphic to a complex torus on which the 
vector fields $\xi$ and $\eta$ are constant (and linearly 
independent). A very different construction, based 
on Chaplygin's integration of the system 
in terms of hyperelliptic integrals, will 
be described in Subsection \ref{jacsubsec}. 

The system (\ref{Xuz})---(\ref{Xzeta}) and even more so the integrals 
(\ref{fi}) resemble the system (2) and the integrals $\langle X,\, X\rangle$, 
(3), (5) and (6) in the article of Adler and van Moerbeke \cite{AM}. 
However, there are also differences: we have the intersection of 
six quadrics in an eight\--dimensional space, whereas in \cite{AM} one has the intersection of ``only'' four quadrics in a six\--dimensional space. 
In Subsection \ref{geodsubsec} we will show that the system 
(\ref{Xuz})---(\ref{Xzeta}) can be mapped to the geodesic flow on 
the Eulidean motion group for a left invariant metric, 
cf. Subsection \ref{geodsubsec}. The latter system resembles the 
ones in the article of Adler and van Moerbeke \cite{AM} even more 
closely, only with the six\--dimensional Lie algebra of $\op{SO}(4)$ 
replaced by the six\--dimensional Lie algebra of the Euclidean motion group 
in the three\--dimensional space. 
Like in \cite{AM}, the vector field in this 
Lie algebra is homogeneous of degree two and 
has four quadratic constants of motion. 

\subsection{Reduction to Horizontal Moment}
\label{redhorsec}
In this subsection we assume that 
the moment $j$ of the momentum around the point of 
contact is not vertical. We will investigate how the vector field $\xi$ 
defined by (\ref{Xuz})---(\ref{Xzeta}) changes if we apply 
a linear substitution of variables of the form 
\begin{equation}
u=a\,\widetilde{u}+b\,\widetilde{v}\quad\mbox{\rm and}\quad 
v=c\,\widetilde{u}+d\,\widetilde{v}.
\label{abcd}
\end{equation}

The equations (\ref{abcd}) are equivalent to 
$d\, u-b\, v=D\,\widetilde{u}$, $-c\, u+a\, v=D\,\widetilde{v}$, 
in which $D:=a\, d-b\, c$.  
Substituting (\ref{Xuz}) and (\ref{Xsigmaz}) in 
$D\,\xi\,\widetilde{u}=d\,\xi\, u-b\,\xi\, v$, 
$D\,\xi\,\widetilde{v}=-c\,\xi\, u+a\,\xi\, v=$, we obtain 
with a straightforward calculation that 
\[
D\,\xi\,\widetilde{u}=D\, (z\, c-\zeta\, a)\,
\widetilde{u}\times J\,\widetilde{u}
+D\, (z\, d-\zeta\, b)\,\widetilde{u}\times J\,\widetilde{v},
\]
and in a similar fashion that 
\[
\xi\,\widetilde{v}=(z\, c-\zeta\, a)\,\widetilde{v}\times J\,\widetilde{u}
+(z\, d-\zeta\, b)\,\widetilde{v}\times J\,\widetilde{v}.
\]
These equations are of the form (\ref{Xuz}), (\ref{Xsigmaz}), 
with $u$, $v$, $z$ and $\zeta$ replaced by 
$\widetilde{u}$, $\widetilde{v}$, $\widetilde{z}$ and $\widetilde{\zeta}$,  respectively, if we take 
\begin{equation}
\widetilde{z}:=z\, d-\zeta\, b ,\quad
-\widetilde{\zeta}:=z\, c-\zeta\, a . 
\label{tildezzeta}
\end{equation}
Applying $\xi$ to (\ref{tildezzeta}) and substituting 
(\ref{Xz}), (\ref{Xzeta}) and (\ref{abcd}), 
we obtain with a straightforward calculation that 
\[
\xi\,\widetilde{z}=D^2\,\op{det}\left(
\widetilde{u},\,\widetilde{J}\,\widetilde{u},
\,\widetilde{J}\,\widetilde{v}\right) ,\quad
-\xi\,\widetilde{\zeta}=D^2\,\op{det}\left(
\widetilde{v},\,\widetilde{J}\,\widetilde{v},
\,\widetilde{J}\,\widetilde{u}\right) .
\]
It follows that the vector field $\xi$ is invariant if we arrange that 
\begin{equation}
(a\, d-b\, c)^2=D^2=1.
\label{thetaD}
\end{equation}

The constants of motion (\ref{fi}) are related to the corresponding ones 
\begin{equation}
\begin{array}{ccc}
\widetilde{f_1}:=\langle\widetilde{u},\, \widetilde{u}\rangle ,&
\widetilde{f_2}:=\langle \widetilde{u},\, \widetilde{v}\rangle, &
\widetilde{f_3}:=\langle \widetilde{v},\, \widetilde{v}\rangle ,\\
\widetilde{f_4}:=
\langle \widetilde{u},\, \widetilde{J}\, \widetilde{u}\rangle 
+\widetilde{z}^2,&
\widetilde{f_5}:=
\langle \widetilde{u},\, \widetilde{J}\, \widetilde{v}\rangle 
+\widetilde{z}\,\widetilde{\zeta} ,&
\widetilde{f_6}:=
\langle \widetilde{v},\, \widetilde{J}\, \widetilde{v}\rangle
+\widetilde{\zeta} ^2
\end{array}
\label{tildefi}
\end{equation}
with tildes over all the variables, by means of the formulas
\begin{eqnarray}
f_1&=&a^2\,\widetilde{f_1}+2a\, b\,\widetilde{f_2}+b^2\,\widetilde{f_3},
\label{f1tilde}\\
f_2&=&a\, c\,\widetilde{f_1}+(a\, d+b\, c)\,\widetilde{f_2}
+b\, d\,\widetilde{f_3},
\label{f2tilde}\\
f_3&=&c^2\,\widetilde{f_1}+2c\, d\,\widetilde{f_2}+d^2\,\widetilde{f_3},
\label{f3tilde}\\
f_4&=&a^2\,\widetilde{f_4}+2a\, b\,\widetilde{f_5}+b^2\,\widetilde{f_6},
\label{f4tilde}\\
f_5&=&a\, c\,\widetilde{f_4}+(a\, d+b\, c)\,\widetilde{f_5}
+b\, d\,\widetilde{f_6},
\label{f5tilde}\\
f_6&=&c^2\,\widetilde{f_4}+2c\, d\,\widetilde{f_5}
+d^2\,\widetilde{f_6}.
\label{f6tilde}
\end{eqnarray}
Here we have used in (\ref{f4tilde}), (\ref{f5tilde}), (\ref{f6tilde}) 
that (\ref{tildezzeta}) imply that 
$a\,\widetilde{z}+b\,\widetilde{\zeta}=D\, z$ and 
$c\,\widetilde{z}+d\,\widetilde{\zeta}=D\,\zeta$, 
whereas (\ref{thetaD}) 
implies that $D^{-2}=1$. 

Consider the matrices $M=\left(\begin{array}{cc}a&b\\c&d\end{array}\right)$,  
$F=\left(\begin{array}{cc}f_1&f_2\\f_2&f_3\end{array}\right)$, 
$G=\left(\begin{array}{cc}f_4&f_5\\f_5&f_6\end{array}\right)$. 
If $\widetilde{F}$ is equal to the matrix $F$ 
with tildes over the coefficients, then the equations 
(\ref{f1tilde}), (\ref{f2tilde}), (\ref{f3tilde}) are equivalent 
to the matrix equation $F=M\,\widetilde{F}\, M^*$. 
Here $M^*$ denotes the transposed of $M$. 
Similarly, if $\widetilde{G}$ is equal to the matrix $G$ 
with tildes over the coefficients, then the equations 
(\ref{f4tilde}), (\ref{f5tilde}), (\ref{f6tilde}) are equivalent 
to the matrix equation $G=M\,\widetilde{G}\, M^*$. 
An obvious consequence is that, for any $\lambda\in\C$, 
\begin{equation}
p(\lambda ):=\op{det}(G-\lambda\, F)=\op{det}(\widetilde{G}-\lambda\,\widetilde{F}).
\label{HPeq}
\end{equation}
Note that $p(\lambda )=\alpha\,\lambda ^2-\beta\,\lambda +\gamma$, 
in which 
\[
\alpha =\op{det}F=f_1\, f_3-{f_2}^2, \quad
\beta =f_1\, f_6+f_3\, f_4-2f_2\, f_5, \quad\mbox{\rm 
and}\quad \gamma =\op{det}G=f_4\, f_6-{f_5}^2.
\] 
Therefore the invariance of the polynomial $p$ is equivalent to 
the three equations $\op{det}F=\op{det}\widetilde{F}$, 
$\op{det}G=\op{det}\widetilde{G}$ and 
$f_1\, f_6+f_3\, f_4-2f_2\, f_5=
\widetilde{f_1}\,\widetilde{f_6}
+\widetilde{f_3}\,\widetilde{f_4}
-2\widetilde{f_2}\widetilde{f_5}$. 
In the case of Chaplygin's sphere, we have 
$p(\lambda )=\alpha\,\lambda ^2-\beta\,\lambda +\gamma$ 
with 
\[
\alpha =\| j\| ^2-{j_3}^2={j_1}^2+{j_2}^2,\quad 
\beta = 2T+\| j\| ^2/\rho,\quad\mbox{\rm and}\quad \gamma =2T/\rho ,
\]
cf. (\ref{Mdef}). 

Suppose that all the coefficients are real and that 
$F$ is positive definite. Then, first diagonalizing 
$F$ by means of an orthogonal transformation, one 
susbsequently can obtain a real $2\times 2$\--matrix $A$ 
such that $F= A\, A^*$. It follows that 
$(\op{det}A) ^2=\op{det}F>0$, which implies that $A$ is invertible. 
There exists an orthogonal transformation $O$ such that 
$D:=O^{-1}\, G'\, O=O^{-1}\, G'\, (O^{-1})^*$ is diagonal, where 
$G'$ is equal to the symmetric matrix $A^{-1}\, G\, (A^{-1})^*$. 
In other words, $G=A\, G'\, A^*=A\, O\, D\, (A\, O)^*$.  
For an arbitrary invertible diagonal matrix $B$ we have now 
arranged that $G=M\,\widetilde{G}\, M^*$, with $M=A\, O\, B^{-1}$ 
and $\widetilde{G}:=B\, D\, B^*$ diagonal. 
In order to arrange that $(\op{det}M)^2=1$, it is sufficient to 
take $1=(\op{det}A)^2\, (\op{det}B)^{-2}=\op{det}F\, (\op{det}B)^{-2}$, 
or $(\op{det}B)^2=\op{det}F$. We now have 
$F=M\,\widetilde{F}\, M^*$ with 
$\widetilde{F}=B\, B^*$. If we choose 
$B=\op{diag}\left( 1,\, (\op{det}F)^{1/2}\right)$, then we have arrived 
at the situation that $\widetilde{f_1}=1$, $\widetilde{f_2}=0$, 
$\widetilde{f_3}=\op{det}F$, $\widetilde{f_5}=0$. 
Note also that in $O$ we still have the freedom to precede it 
by the matrix which switches the two basis vectors, which 
means that we still can switch $\widetilde{f_4}$ and 
$\widetilde{f_6}$. 

\begin{remark}
A diagonal matrix remains unchanged if it is multiplied from 
the left and the right by the matrix $\op{diag}(1,\, -1)$, 
which has determinant equal to $-1$. Therefore we can arrive 
at the same diagonal matrices $\widetilde{F}$, $\widetilde{G}$ 
with the help of a matrix $M$ which satisfies the stronger 
condition $\op{det}M=1$ instead of $\op{det}M=\pm 1$ as 
required in (\ref{thetaD}). 

In other words, $\op{SL}(2,\,\C )$ acts as a symmetry group 
for the vector fields $\xi$ and $\eta$. The functions $f_i$ 
are not invariant under the action of $\op{SL}(2,\,\C )$ and 
actually the $\op{SL}(2,\,\C )$\--action can be used to 
change to levels $\widetilde{f_i}$ such that 
$\widetilde{f_1}=1$, $\widetilde{f_2}=0$, 
$\widetilde{f_3}=\op{det}F$, $\widetilde{f_5}=0$. 
On the other hand, the coefficients of the polynomial 
$p$ in (\ref{HPeq}) are invariant, which fact can be used 
in order to determine $\widetilde{f_3}$, $\widetilde{f_4}$ 
and $\widetilde{f_6}$ in terms of the $f_i$. 

In Subsection \ref{geodsubsec} we will give a description of the 
quotient space under the action of $\op{SL}(2,\,\C )$
\label{symmgrouprem}
\end{remark}

In the case of Chaplygin's sphere, when we have (\ref{Mdef}),      
we have that $\op{det}F=\| j\| ^2-{j_3}^2={j_1}^2+{j_2}^2>0$ because 
$j$ is not vertical. Therefore $F$ is positive definite 
and we have a reduction to the situation that $\widetilde{f_1}=1$, $\widetilde{f_2}=0$ and $\widetilde{f_5}=0$. It follows then 
from (\ref{HPeq}) that $\widetilde{f_3}=\op{det}F={j_1}^2+{j_2}^2$, 
whereas $\widetilde{f_4}$, $\widetilde{f_6}$ satisfy the 
equations 
$\widetilde{f_4}\,\widetilde{f_6}=\op{det}G=2T/\rho$ and 
\[
2T+\| j\| ^2/\rho =\widetilde{f_6}+\left({j_1}^2+{j_2}^2\right)
\,\widetilde{f_4}.
\]
These $\widetilde{f_i}$ are again the levels of a Chaplygin's sphere, but 
with $j$, $\rho$, $T$ replaced by $\widetilde{j}$, $\widetilde{\rho}$, 
$\widetilde{T}$, respectively. If $\widetilde{j}_3=0$, 
which means that {\em the new moment $\widetilde{j}$ is horizontal}, 
the length of $\widetilde{j}$ is equal to the length of the 
horizontal projection of $j$, and $\widetilde{\rho}$ and 
$\widetilde{T}$ satisfy the equations 
$\widetilde{T}/\widetilde{\rho}=T/\rho$ and 
\[
2\widetilde{T}+\left({j_1}^2+{j_2}^2\right)
/\widetilde{\rho}=2T+\| j\| ^2/\rho .
\]
The solution which depends continuously on the parameters 
and satisfies $\widetilde{T}=T$, $\widetilde{\rho}=\rho$ 
when $j_3=0$, is given by 
\begin{equation}
\widetilde{\rho}=\rho\,\widetilde{T}/T
\quad\mbox{\rm and}\quad 
2\widetilde{T}=T+\| j\| ^2/2\rho +
\left[\left( T-\| j\| ^2/2\rho\right) ^2+2{j_3}^2\, T/\rho\right] ^{1/2}.
\label{tildeT}
\end{equation}
Because $\widetilde{T}\neq T$ when $j$ is not already 
horizontal, we have $\widetilde{\rho}\neq\rho$. 
{\em Because we do not change the matrix $J=(I+\rho )^{-1}$, 
this means that we have to change the moment of inertia 
tensor $I$ to the new one $\widetilde{I}=I+\rho -\widetilde{\rho}$.} 

If we allow complex coefficients, then we can arrive at any 
$\widetilde{F}$, $\widetilde{G}$ such that 
(\ref{HPeq}) holds, provided that the polynomial $p$ is of 
degree two and has two distinct zeros, 
cf. Hodge and Pedoe \cite[Vol. II, p. 278]{HP}. 
These conditions are equivalent to the conditions 
that $f_1\, f_3-{f_2}^2\neq 0$ and the discriminant 
\begin{equation}
\Delta =\left( f_1\, f_6+f_3\, f_4-2f_2\, f_5\right) ^2 
-4\left( f_1\, f_3-{f_2}^2\right)\,
\left( f_4\, f_6-{f_5}^2\right)
\label{discrf}
\end{equation}
of $p$ is not equal to zero. 
Again we can arrange that $\widetilde{f_1}=1$, 
$\widetilde{f_3}=f_1\, f_3-{f_2}^2$. The new 
values $\widetilde{f_4}$ and $\widetilde{f_6}$ are 
determined from the condition that $\widetilde{f_4}$ 
and $\widetilde{f_6}/\widetilde{f_3}$ are equal to the zeros 
of $p$, which are unique up to their ordering. 

\begin{remark}
In the case of Chaplygin's sphere, when we have (\ref{Mdef}).    
and similar equations with tildes over all the symbols,  
the equations (\ref{f1tilde}), (\ref{f2tilde}), (\ref{f3tilde}) 
are equivalent to the statement that there exists a rotation 
$C\in\op{SO}(3)$ such that $C\, e_3=a\, e_3+b\,\widetilde{j}$ and 
$C\, j=c\, e_3+d\,\widetilde{j}$. These equations are equivalent to 
(\ref{abcd}) if $u=A^{-1}\, e_3$, 
$v=A^{-1}\, j$, $\widetilde{u}=\widetilde{A}^{-1}\, e_3$ and  
$\widetilde{v}=\widetilde{A}^{-1}\,\widetilde{j}$,   
in which $\widetilde{A}=C\, A$, $A\in\op{SO}(3)$. 
The invariance of the vector 
field $\xi$ under the substitutions (\ref{abcd}) then means that 
if $A(\tau )$ denotes the rotational motion as a function of 
the reparametrized time $\tau$, then $\widetilde{A}(\tau )=C\, 
A(\tau )$ satisfies a differential equation of the same form 
as $A(\tau )$, but with $j$, $\rho$, $I$ and $T$ replaced by 
$\widetilde{j}$, $\widetilde{\rho}$, $\widetilde{I}$ and 
$\widetilde{T}$, respectively. 
\end{remark}

\begin{remark}
If $j$ is not vertical, then it is impossible to make 
$\widetilde{j}$ vertical. 
Actually the rotational system with nonvertical $j$ can 
not be transformed 
in any algebraic fashion to the system with vertical moment of 
Section \ref{vertj}. 

Indeed, if $j$ is vertical, then the completion of the 
complexification of the phase space of the $\op{SO}(2)$\--reduced system 
is equal to the elliptic curve 
defined by (\ref{omegasphere}), (\ref{omegaellipsoid}). 
On the other hand the complexification of the 
$\op{SO}(2)$\--action is a free action of $\C ^{\times}$, 
where $\C ^{\times}$ denotes the multiplicative 
group of the nonzero complex numbers. 
(The mapping $t\mapsto \op{e}^t$ is an isomorphism from the 
additive group $\C /2\pi\op{i}\Z$ onto the multiplicative 
group $\C ^{\times}$.)  In this way the completion of the complexification 
of the phase space of the system for vertical $j$ is a 
$\C ^{\times}$\-- bundle over an elliptic curve, 
which is not isomorphic to the 
Jacobi variety $\op{Jac}(C)$ of the hyperelliptic 
curve $C$ which we obtain when the moment is not vertical. 
In particular the $\C ^{\times}$\--bundle is not compact 
because $\C ^{\times}$ is not compact, whereas $\op{Jac}(C)$ is compact. 

The passage from nonvertical to vertical $j$ is an example of  
``a degenerate limit of an abelian variety, as an extension of 
a power of $\C ^{\times}$ by an abelian variety'', mentioned 
by Mumford \cite[p. 3.53]{TthII}. 
\end{remark}

For the motion of the point of contact $p$, we observe that 
our substitutions imply that, with $\omega$ as in 
(\ref{omegausigma}), $z\,\omega =z\, J\, v-\zeta\, J\, u
=\widetilde{z}\,\widetilde{J}\,\widetilde{v}-
\widetilde{\zeta}\,\widetilde{J}\,\widetilde{u}
=\widetilde{z}\,\widetilde{\omega}$. Therefore 
(\ref{pdot}) yields that 
\[
\dd{\langle p,\, j\rangle}{\tau}
=z\, r\,\langle (A\,\omega )\times e_3,\, j\rangle 
=z\, r\,\langle \omega\times u,\, v\rangle =
\widetilde{z}\, r\, D\,\langle  
\widetilde{\omega}\times\widetilde{u},\,\widetilde{v}\rangle ,
\]
which is equal to a constant times the same function for 
the new system with horizontal moment, and therefore equal 
to a constant times the rational function on the double 
covering of $\op{Jac}(\widetilde{C})$, described after (\ref{undetermined}).  
Here $\widetilde{C}$ denotes the hyperelliptic curve 
corresponding to the new system with the horizontal moment. 

A similar calculations yields for the $j\times e_3$\--component of $p$ that 
\[
\dd{\langle p,\, j\times e_3\rangle}{\tau}
=D\,\widetilde{z}\,\left( a\, 2T+b\,\widetilde{\zeta}/\widetilde{\rho}\right) .
\]
The function $\widetilde{z}$ has a similar description as a 
rational function on a double covering of 
$\op{Jac}(\widetilde{C})$, with 
$\widetilde{\gamma}$ replaced by $-1/\widetilde{\rho}$. 
However, the function 
$\widetilde{z}\,\widetilde{\zeta}=\, -\langle\widetilde{u},
\,\widetilde{J}\,\widetilde{v}\rangle$ 
does not seem to have an equally straightforward description 
in terms of $\op{Jac}(\widetilde{C})$. 

\subsection{Geodesic Flow on the Euclidean Motion Group}
\label{geodsubsec}
The form of the equations 
(\ref{Xuz})---(\ref{Yzeta}) suggests to introduce the vectors 
\begin{equation}
q:=u\times v\quad\mbox{\rm and}\quad r:=z\, v-\zeta\, u.
\label{qr}
\end{equation}
The vector $(q,\, r)\in\C ^3\times\C ^3=\C ^6$ represents 
the exterior product of the vectors 
$(u,\, z)$ and $(v,\,\zeta )$ 
in $\C ^4$, and therefore the mapping 
$((u,\, z),\, (v,\,\zeta ))\mapsto (q,\, r)$ from 
$\C ^8$ to $\C ^6$ 
will be denoted by $\wedge$. 

For Chaplygin's sphere, it follows from (\ref{sigmadef}) that 
\[
q=u\times v=A^{-1}\, \left( e_3\right)\times A^{-1}\, j
=A^{-1}\,\left( e_3\times j\right) ,
\]
which means that $q$ is equal to the vector $e_3\times j$ in 
body coordinates. Furthermore, it follows from 
(\ref{zX}) and (\ref{zzetauv}) that $J\, r=z\,\omega$, 
which is equal to the rotational velocity with respect to the 
time variable $\tau$ which is related to $t$ by 
$\op{d}\! t/\op{d}\!\tau =X^{1/2}=z$, cf. Corollary \ref{commvvcor} 
and (\ref{zu}).  

We have   
\[
\xi q=(\xi u)\times v+u\times (\xi v)=
(u\times J\, r)\times v+u\times (v\times J\, r)=(u\times v)\times J\, r
=q\times J\, r, 
\]
in which the third identity follows from the Jacobi identity 
$(u\times r)\times v+(r\times v)\times u+(v\times u)\times r$ 
in $\mbox{\gothic so}(3)$. Similarly 
\begin{eqnarray*}
\xi r&=&(\xi z)\, v-(\xi\zeta )\, u+z\, (\xi v)-\zeta\, (\xi u)\\
&=&\langle u,\, J\, u\times J\, v\rangle\, v
-\langle v,\, J\, u\times J\, v\rangle\, v
+z\, v\times J\, r-\zeta\, u\times J\, r\\
&=&(u\times v)\times (J\, u\times J\, v)+r\times J\, r.
\end{eqnarray*}
For any vector $w$ we have that 
\[
\langle J\, u\times J\, v,\, J\, w\rangle =
\op{det}(J\, u,\, J\, v,\, J\, w)=\op{det}J\,\op{det}(u,\, v,\, w)
=\op{det}J\,\langle u\times v,\, w\rangle ,
\]
which in view of the symmetry of $J$ implies that 
$J\, (J\, u\times J\, v) =\op{det}J\, u\times v$, 
or $J\, u\times J\, v=\op{det}J\, J^{-1}(u\times v)$ when $J$ 
is invertible. It follows that 
\begin{equation}
J\, u\times J\, v=
J^{\scriptop{co}}\, (u\times v),\quad\mbox{\rm in which}
\quad J^{\scriptop{co}}=\op{diag}\left( 
J_2\, J_3,\, J_3\, J_1,\, J_1\, J_2\right) . 
\label{coJ}
\end{equation}
Note that $J^{\scriptop{co}}=(\op{det}J)\, J^{-1}$ when $J$ is 
invertible, but (\ref{coJ}) also holds for noninvertible $J$.  
This means that the mapping $\wedge$ intertwines the 
vector field $\xi$ in $\C ^8$ with the vector field 
$\xi$ in $\C ^6$ defined by 
\begin{equation}
\xi q=q\times J\, r\quad\mbox{\rm and}\quad 
\xi r=q\times J^{\scriptop{co}}\, q+r\times J\, r.
\label{xiqr}
\end{equation}
Similarly, we have 
\[
\eta q=(\eta u)\times v+u\times (\eta v)=
z\, q\times v+u\times \zeta\, q=q\times (z\, v-\zeta u)=q\times r,
\]
and 
\begin{eqnarray*}
\eta r&=&(\eta z)\, v-(\eta\zeta )\, u+z\, (\eta v)-\zeta\, (\eta u)
=\, -\langle J\, u,\, q\rangle\, v
+\langle J\, v,\, q\rangle\, u
+z\, v\times r-\zeta\, u\times r\\
&=&\, -\langle u,\, J\, q\rangle\, v
+\langle v,\, J\, q\rangle\, u+r\times r
=\, -(u\times v)\times J\, q=\, -q\times J\, q. 
\end{eqnarray*}
Therefore the mapping $\wedge$ intertwines the 
vector field $\eta$ in $\C ^8$ with the vector field 
$\eta$ in $\C ^6$ defined by 
\begin{equation}
\eta q=q\times r\quad\mbox{\rm and}\quad 
\eta r=\, -q\times J\, q.
\label{etaqr}
\end{equation}

If $q$ and $r$ is interpreted as a rotational and translational 
velocity vector, then the $(q,\, r)$\--space can be identified 
with the (complexified) Lie algebra $\mbox{\gothic e}(3)$ 
of the Euclidean motion 
group $\op{E}(3)$ in the three dimensional Euclidean space, 
with the Lie brackets defined by 
\begin{equation}
\left[ (q,\, r),\, (q',\, r')\right] =(q\times q',\, 
q\times r'+r\times q'). 
\label{e3}
\end{equation}
Therefore the vector field $\xi$ defined by (\ref{xiqr}) has the form of 
a Lax pair $\op{d}\! X/\op{d}\!\tau =\left[ X,\, L_{\xi}(X)\right]$, 
where $L_{\xi}$ is the linear transformation in $\mbox{\gothic e}(3)$ 
defined by 
\begin{equation}
L_{\xi}(q,\, r)=\left( J\, r,\, J^{\scriptop{co}}\, q\right) .
\label{Lxi}
\end{equation}
Similarly the vector field $\eta$ defined by (\ref{etaqr}) has the form of 
a Lax pair $\op{d}\! X/\op{d}\!\tau =\left[ X,\, L_{\eta}(X)\right]$, 
where $L_{\eta}$ is the linear transformation in $\mbox{\gothic e}(3)$ 
defined by 
\begin{equation}
L_{\eta}(q,\, r)=\left( r,\, -J\, q\right) .
\label{Leta}
\end{equation}
The Lax pair form of $\xi$ and $\eta$ implies that 
the flows of both vector fields leave the conjugacy 
classes in $\mbox{\gothic e}(3)$ invariant. Because the functions 
\begin{equation}
h_1(q,\, r):=\langle q,\, q\rangle\quad\mbox{\rm and}\quad
h_2(q,\, r):=\langle q,\, r\rangle
\label{h1h2}
\end{equation}
are constant on the conjugacy classes, it follows 
that {\em both functions $h_1$ and $h_2$ 
are constants of motion for $\xi$ and $\eta$}. 

Let $G$ be a Lie group and $f$ a function on the cotangent bundle 
$\op{T}^*G$ of 
$G$ which is invariant under all left multiplications by 
elements of $G$. The quotient space of $\op{T}^*G$ by means of 
the left action of $G$ on $\op{T}^*G$ is naturally identified 
with the dual space $\mbox{\gothic g}^*$ of the 
Lie algebra $\mbox{\gothic g}$ of $G$, and we denote the 
restriction of $f$ to $\mbox{\gothic g}^*$ with the same letter. 
The canonical Possion structure on $\op{T}^*G$ induces 
a Poisson structure on $\mbox{\gothic g}^*$ in such a 
way that the Hamiltonian vector field $\op{H}_f$ 
of $f$ on $\mbox{\gothic g}^*$ 
is given by 
\begin{equation}
-\langle X,\,\op{H}_f(l)\rangle =\langle\left[ X,\,\op{d}\! f(l)\right] ,
\, l\rangle ,\quad X\in\mbox{\gothic g},\; l\in\mbox{\gothic g}^*.  
\label{Poisson}
\end{equation}
The vector field $\op{H}_f$ is tangent to the coadjoint orbits 
in $\mbox{\gothic g}^*$, on which the Poisson structure 
is given by a symplectic structure. This means that on each 
coadjoint orbit the vector field $\op{H}_f$ is Hamiltonian 
with respect to this symplectic structure. This construction 
has been introduced already by Lie in \cite[Kap. 19]{Lie}  
under the name ``Die dualistische der adjungierte Gruppe''. 
It has been rediscovered independently by Kostant, Kirillov, 
and Souriau. The coadjoint orbits with their symplectic structure 
are the Marsden\--Weinstein reduced phase spaces of 
$\op{T}^*G$ for the left action of $G$ on $\op{T}^*G$. 
See also Abraham and Marsden \cite[Sections 4.3, 4.4]{AbM}. 

When $\mbox{\gothic g}=\mbox{\gothic e}(3)$, cf. (\ref{e3}, then 
the right hand side in (\ref{Poisson}) takes the form 
\begin{eqnarray*}
&&\langle\left[ (q',\, r'),\,\left(\frac{\partial f}{\partial a}, 
\,\frac{\partial f}{\partial b}\right)\right] ,
\, (a,\, b)\rangle 
=\langle q'\times\frac{\partial f}{\partial a},\, a\rangle 
+\langle q'\times\frac{\partial f}{\partial b}+
r'\times\frac{\partial f}{\partial a},\, b\rangle \\
&&=\langle q',\,\frac{\partial f}{\partial a}\times a 
+\frac{\partial f}{\partial b}\times b\rangle 
+\langle r',\,\frac{\partial f}{\partial a}\times b\rangle .
\end{eqnarray*}
Therefore, if we identify the point $(a,\, b)\in\mbox{\gothic g}^*$ 
with the point $(q,\, r)\in\mbox{\gothic g}$ with $q=b$ and $r=a$, 
then we recoginze from (\ref{xiqr}) that $\xi$ is equal 
to the Hamiltonian vector field of the function $f=h_3/2$, in which 
\begin{equation}
h_3(q,\, r):=\langle q,\, J^{\scriptop{co}}\, q\rangle 
+\langle r,\, J\, r\rangle .
\label{h3}
\end{equation}
With the same identifications the vector field $\eta$ is equal to the 
Hamiltonian vector field defined by the function 
\begin{equation}
h_4(q,\, r):=\, -\langle q,\, J\, q\rangle 
+\langle r,\, r\rangle .
\label{h4}
\end{equation}

\begin{remark}
The Lie algebra $\mbox{\gothic e}(3)$ has not been identified with 
its dual $\mbox{\gothic e}(3)^*$ by means of the Killing 
form $-4\langle q,\, q\rangle$, which is degenerate, 
but by means of the nondegenerate invariant quadratic 
form $h_2=\langle q,\, r\rangle$. 

The function $h_3$ defines a left invariant metric on 
$\op{T}^*\op{E}(3)$, which can be used to identify 
$\op{T}^*\op{E}(3)$ with the tangent bundle of $\op{E}(3)$. 
Under this identification the flow of the Hamiltonian system of $h_3/2$ 
corresponds to the geodesic flow on  
$\op{T}\op{E}(3)$ defined by the dual metric on $\op{T}\op{E}(3)$. 

We therefore obtain the following, somewhat roundabout correspondence 
between this geodesic flow and our vector field $\xi$. 
First pass from the geodesic flow of the left invariant metric 
on the tangent bundle to the Hamiltonian system of the function $h_3/2$ 
on the cotangent bundle, using the 
metric in order to identify the tangent bundle with the 
cotangent bundle. Then pass to the reduced system 
on $\mbox{\gothic e}(3)^*$ by means of the left action of $\op{E}(3)$ 
on $\op{T}^*\op{E}(3)$. In the next step, pass to the 
vector field $\xi$ on $\mbox{\gothic e}(3)$ using the 
identification of $\mbox{\gothic e}(3)$ with $\mbox{\gothic e}(3)^*$ 
by means of the quadratic form $\langle q,\, r\rangle$. 
Finally the mapping $\wedge$ intertwines the vector 
field $\xi$ in the $((u,\, z),\, (v,\,\zeta ))$\--space with 
the vector field $\xi$ in the $(q,\, r)$\--space $\mbox{\gothic e}(3)$. 
\label{geodrem}
\end{remark}

The Poisson brackets $\{ f,\, g\} :=\op{H}_fg$ define a Lie algebra 
structure on the space of functions. In particular it is antisymmetric, 
which implies that $\op{H}_ff=0$ and $\op{H}_fg=0$ if and only if 
$\op{H}_gf=0$. It follows from $\op{H}_ff=0$ that 
$h_3$ and $h_4$ is a constant of motion 
for $\xi$ and $\eta$, respectively. Moreover, 
\[
\eta h_3/2=\langle\eta q,\, J^{\scriptop{co}}\, q\rangle 
+\langle\eta r,\, J\, r\rangle 
=\langle q\times r,\, J^{\scriptop{co}}\, q\rangle 
-\langle q\times J\, q,\, J\, r\rangle =0,
\]
because 
\[
\langle q\times J\, q,\, J\, r\rangle 
=\langle q,\, J\, q\times J\, r\rangle 
=\langle q,\, J^{\scriptop{co}}\, (q\times r)\rangle 
=\langle J^{\scriptop{co}}\, q,\, q\times r\rangle ,
\]
cf. (\ref{coJ}). This implies that {\em $h_3$ and $h_4$ are constants 
of motion for both vector fields $\xi$ and $\eta$}. 
The Jacobi identity of the Poisson structure implies that 
$\left[\op{H}_f,\,\op{H}_g\right] =\op{H}_{\{ f,\, g\}}$. 
Because we just proved that $\{ h_3,\, h_4\} =0$, it follows 
also that {\em the vector fields $\xi$ and $\eta$ commute}. 

In particular the system in $\mbox{\gothic e}(3)\simeq\C ^6$ 
defined by the vector field $\xi$ is completely integrable, 
as a Hamiltonian system on the four\--dimensional coadjoint 
orbits determined by fixing the values of the functions 
$h_1$ and $h_2$, defined in (\ref{h1h2}). The Hamiltonian 
function is the function $h_3/2$ with $h_3$ defined in 
(\ref{h3}) and the function $h_4$ defined in (\ref{h4}) 
is the additional constant of motion which yields the 
complete integrability of the Hamiltonian system. 

\begin{remark}
The system in the $(q,\, r)$\--space resembles the one in the 
paper of Adler and van Moerbeke \cite{AM} very closely: it 
is defined on a six\--dimensional (dual of a) Lie algebra 
($\mbox{\gothic e}(3)$ instead of the Lie algebra 
$\mbox{\gothic so}(4)$ of \cite{AM}) and has four 
quadratic constants of motion. Furthermore it 
is Hamiltonian on coadjoint orbits, and we have 
two polynomial vector fields which are homogeneous of 
degree two. In this respect the vector fields 
are simpler than the vector fields $\xi$ and $\eta$ 
defined in (\ref{Xuz})---(\ref{Yzeta}).  

Also the behaviour at infinity the level surface of the $h_i$ 
(with $h_2=0$) in the complex projective space 
is very similar to the behaviour of the level 
surface in $\mbox{\gothic so}(4)$ as described 
in Mumford's appendix to \cite{AM}. 
See Subsection \ref{e3inftysec}. 
 
On the other hand it turns out that this behaviour 
is more singular than that of the projective 
closure of the level surface $M$ of the functions 
(\ref{fi}), see Subsection \ref{e3inftysec}. 
For this reason we start Section \ref{ccsec} 
with a discussion of the system 
in $\C ^8$, because this seemingly more complicated
system has a simpler behaviour at infinity. 
\label{AM2rem}
\end{remark}

We now turn to a closer examination of the mapping 
$\wedge$ from $\C ^8$ to $\C ^6$ defined by (\ref{qr}). 
To begin with, if $(q,\, r)$ belongs to the image of 
$\wedge$, then $\langle q,\, r\rangle =0$, which means that 
$\wedge$ is a mapping from $\C ^8$ to the hypersurface 
$h_2=0$ in $\C ^6$, cf. (\ref{h1h2}). Conversely, 
any $q\in\C ^3$ can be written as $u\times v$ for 
some $u,\, v\in\C ^3$. If $q\neq 0$, then $u$ and $v$ 
are linearly independent and span the orthogonal complement 
of $q$, which implies that for every $r\in\C ^3$ such 
that $\langle q,\, r\rangle =0$ there exist $z,\,\zeta\in\C$ such 
that $r=z\, v-\zeta\, u$. If $q=0$, then we have 
(\ref{h1h2}) for $z=1$, $v=r$, $\zeta =0$, $u=0$. 
It follows that $\wedge$ is surjective from $\C ^8$ onto  
the hypersurface $h_2=0$ in $\C ^6$. 

\begin{remark}
If $q\neq 0$, then $\langle q,\, r\rangle =0$ if and only 
if there exists a vector $x_0$ such that $q\times x_0+r=0$, 
and every vector $x$ such that $q\times x+r=0$ is of the form 
$x=x_0+c\, q$ for some scalar $c$. Therefore the condition 
that $\langle q,\, r\rangle$ means that the infinitesimal 
motion $(q,\, r)\in\mbox{\gothic e}(3)$ is either 
equal to an infinitesimal translation ($q=0$), or 
to an infinitesimal rotation about some axis in the three\--dimensional 
case: ``no spiralling''.
\end{remark}  

If we interpret $(q,\, r)$ as an element of $\bigwedge ^2\C ^4$, 
then the condition $\langle q,\, r\rangle =0$ means that the rank 
of $(q,\, r)$ is smaller than four. If $(q,\, r)\neq (0,\, 0)$, 
then its null space is two\--dimensional and is spanned by 
$(u,\, z)$ and $(v,\,\zeta )$ if 
$(q,\, r)=(u,\, z)\wedge (v,\,\zeta )=(u\times v,\, z\, v-\zeta\, u)$. 
Note that $(u,\, z)$ and $(v,\,\zeta )$ are linearly independent 
in this case. If $(q,\, r)=(u',\, z')\wedge (v',\,\zeta ')$, 
then $(u',\, z')$ and $(v',\,\zeta ')$ are also 
contained in the null space of $(q,\, r)$, and therefore 
there are unique $a,\, b,\, c,\, d\in\C$ such that 
\begin{equation}
(u',\, z')=a\, (u,\, z)+b\, (v,\,\zeta )\quad\mbox{\rm and}\quad 
(v',\,\zeta ')=c\, (u,\, z)+d\, (v,\,\zeta ), 
\label{sl2act}
\end{equation}
and for such vectors $(u',\, z')$ and 
$(v',\,\zeta ')$ we have that $(q,\, r)=(u',\, z')\wedge (v',\,\zeta ')$ 
if and only if $a\, d-b\, c=1$. In other words, if 
$\langle q,\, r\rangle =0$ and $(q,\, r)\neq (0,\, 0)$, 
then the fiber of $(q,\, r)$ for the mapping $\wedge$ is equal 
to the orbit in the $((u,\, z),\, (v,\,\zeta ))$\--space $\C ^8$ 
of the action of $\op{SL}(2,\,\C )$ defined by (\ref{sl2act}). 

Let $U$ denote the set of the $((u,\, z),\, (v,\,\zeta ))\in\C ^8$ 
such that the vectors $(u,\, z)$ and $(z,\,\zeta )$ in $\C ^4$ 
are linearly independent, and let $V$ be the set of 
$(q,\, r)\in\C ^6$ such that $\langle q,\, r\rangle =0$ 
and $(q,\, r)\neq (0,\, 0)$. Then $U$ is an open subset 
of $\C ^8$ and $V$ is an open subset, equal to the smooth part, 
of the 5\--dimensional hypersurface $\langle q,\, r\rangle =0$ 
in $\C ^6$. The action of $\op{SL}(2,\,\C )$ on $U$ is free 
and the mapping $\wedge$ identifies $V$ with the orbit space 
of the $\op{SL}(2,\,\C )$\--action on $U$. 

In Subsection \ref{redhorsec} we had observed that the 
action of $\op{SL}(2,\,\C )$ leaves the vector fields 
$\xi$ and $\eta$ invariant, cf. Remark \ref{symmgrouprem}.  
Therefore the projection $\wedge :U\to V$ intertwines 
the vector fields $\xi$ and $\eta$ on $U$ with uniquely 
defined vector fields on $V$, which we denoted with the 
same letters. The fact that the vector fields 
$\xi$ and $\eta$ in $U$ commute implies that 
their push\--forwards under $\wedge$, the vector fields 
$\xi$ and $\eta$ in $V$, commute as well. 
In the beginning of this subsection we 
showed that the vector fields $\xi$ and $\eta$ in 
$\C ^6$ defined by (\ref{xiqr}) and (\ref{etaqr}) 
extend the vector fields $\xi$ and $\eta$ in $V$. 
The fact that the vector fields $\xi$ and $\eta$ 
in $\C ^6$ commute is stronger than the fact that 
their restrictions to $V$ commute. 

As observed in Remark \ref{symmgrouprem}, the functions 
$f_i$ are not invariant under the $\op{SL}(2,\,\C )$\--action.  
However, the coefficients of the polynomial $p$ in 
(\ref{HPeq}) are invariant under the $\op{SL}(2,\,\C )$\--action, 
which means that these coefficents 
can be written as functions of $(q,\, r)\in V$. 
Actually, we have:
\begin{eqnarray}
h_1\circ\wedge &=&f_1\, f_3-{f_2}^2,
\label{h1fi}\\
h_2\circ\wedge &=&0,
\label{h2fi}\\
h_3\circ\wedge &=&f_4\, f_6-{f_5}^2,
\label{h3fi}\\
(\op{trace}J)\, h_1\circ\wedge +h_4\circ\wedge 
&=&f_1\, f_6+f_3\, f_4-2f_2\, f_5, 
\label{h4fi}
\end{eqnarray}
from which the compositions $h_i\circ\wedge$ of the 
functions $h_i$ with the mapping $\wedge$ can be determined 
in terms of the functions $f_i$ defined in (\ref{fi}). 
Here the functions $h_i$ are defined by (\ref{h1h2}), 
(\ref{h3}) and (\ref{h4}). 
\begin{proof}
The equation (\ref{h1fi}) follows from 
\[
\op{det}F=f_1\, f_3-{f_2}^2
=\langle u,\, u\rangle\cdot\langle v,\, v\rangle 
-\langle u,\, v\rangle ^2=\langle u\times v,\, u\times v\rangle 
=\langle q,\, q\rangle =h_1. 
\]
For (\ref{h3fi}) we write 
\begin{eqnarray*}
\op{det}G&=&f_4\, f_6-{f_5}^2
=\left(\langle u,\, J\, u\rangle +z^2\right)\, 
\left(\langle v,\, J\, v\rangle +\zeta ^2\right) 
-\left(\langle u,\, J\, v\rangle +z\,\zeta\right) ^2 \\
&=&\langle u,\, J\, u\rangle\cdot\langle v,\, J\, v\rangle 
-\langle u,\, J\, v\rangle ^2
+z^2\,\langle v,\, J\, v\rangle +\zeta ^2\,\langle u,\, J\, u\rangle 
-2z\,\zeta\,\langle u,\, J\, v\rangle \\
&=&\langle u,\, \langle v,\, J\, v\rangle\, J\, u 
-\langle v,\, J\, u\rangle\, J\, v\rangle 
+\langle z\, v-\zeta\, u,\, J\, (z\, v-\zeta\, u)\rangle \\
&=&\langle u,\, v\times (J\, u\times J\, v)\rangle +\langle r,\, J\, r\rangle
=\langle q,\, J^{\scriptop{co}}\, q\rangle +\langle r,\, J\, r\rangle 
=h_3,
\end{eqnarray*}
where we used (\ref{coJ}) in the fifth identity. 

Finally, we have that $f_1\, f_6+f_3\, f_4-2f_2\, f_5$ is equal to 
\[
\langle u,\, u\rangle\,\left(\langle v,\, J\, v\rangle +\zeta ^2\right) 
+\langle v,\, v\rangle\,\left(\langle u,\, J\, u\rangle +z^2\right) 
-2\langle u,\, v\rangle\,\left(\langle u,\, J\, v\rangle 
+z\,\zeta\right) ,
\]
for the computation of which we write 
\[
\langle u,\, u\rangle\cdot\langle v,\, J\, v\rangle 
+\langle v,\, v\rangle\cdot\langle u,\, J\, u\rangle 
-2\langle u,\, v\rangle\cdot\langle u,\, J\, v\rangle =\langle u,\, a\rangle ,
\]
in which 
\begin{eqnarray*}
a&=&\langle v,\, J\, v\rangle\, u+\langle v,\, v\rangle\, J\, u
-2\langle u,\, J\, v\rangle\, v\\
&=&\langle v,\, J\, v\rangle\, u-\langle u,\, J\, v\rangle\, v
+\langle v,\, v\rangle\, J\, u-\langle v,\, J\, u\rangle\, v\\
&=&J\, v\times (u\times v)+v\times (J\, u\times v) .
\end{eqnarray*}
It follows that $\langle u,\, a\rangle =\langle u\times v,\, b\rangle$, 
in which $b=u\times J\, v+J\, u\times v$. Now we have, for 
any vector $c$, 
\begin{eqnarray*}
\langle b,\, c\rangle &=&\op{det}(u,\, J\, v,\, c)+\op{det}(J\, u,\, v,\, c)
=(\op{trace}J)\op{det}(u,\, v,\, c)-\op{det}(u,\, v\, J\, c) \\
&=&(\op{trace}J)\langle u\times v,\, c\rangle -\langle u\times v,\, J\, c\rangle 
=\langle (\op{trace}J-J)(u\times v),\, c\rangle,
\end{eqnarray*}
which implies that $b=(\op{trace}J-J)(u\times v)$. 
Collecting all the results, we arrive at 
\[
f_1\, f_6+f_3\, f_4-2f_2\, f_5 = \langle q,\, (\op{trace}J-J)(q)\rangle 
+\langle r,\, r\rangle =(\op{trace}J)\, h_1+h_4,  
\]
from which (\ref{h4fi}) follows.
\end{proof}
In the case of Chaplygin's sphere, we can insert the values 
(\ref{Mdef}) of the functions $f_i$, which leads to the values 
\begin{equation}
h_1={j_1}^2+{j_2}^2,\; h_2=0,\; 
h_3=2T/\rho ,\; (\op{trace}J)\, h_1+h_4=2T+\| j\| ^2/\rho 
\label{chaplhi}
\end{equation} 
for the constants of motion $h_i$ of the system in the six\--dimensional 
$(q,\, r)$\--space. 

The action of the matrix $M\in\op{SL}(2,\,\C )$ with the coefficients 
$a,\, b,\, c,\, d$ leaves a given level surface of the functions 
$f_i$ invariant if and only if, in the notation of Subsection 
\ref{redhorsec}, we have that $M\, F\, M^*=F$ and 
$M\, G\, M^*=G$. Assume that $\op{det}F\neq 0$ and that 
$F$ and $G$ are not proportional. Then we obtained in 
Subsection \ref{redhorsec} that there exists an 
$A\in\op{SL}(2,\,\C )$ such that the matrices 
$\widetilde{F}=A\, F\, A^*$ and $\widetilde{G}=A\, G\, A^*$ are 
diagonal, and it follws that $\widetilde{F}$ and 
$\widetilde{G}$ are not proportional. With the 
notation $\widetilde{M}=A\, M\, A^{-1}$, we now have 
$\widetilde{M}\,\widetilde{F}\,\widetilde{M}^*=\widetilde{F}$ 
and $\widetilde{M}\,\widetilde{G}\,\widetilde{M}^*=\widetilde{G}$. 
A straighforward calculation, in which we use that the diagonal matrices 
$\widetilde{F}$ and $\widetilde{G}$ are not proportional, 
leads to the conclusion that $\widetilde{M}=\pm 1$, 
which in turn implies that $M=\pm 1$. If $M=\, -1$ then it 
acts on the $((u,\, z),\, (v,\,\zeta ))$\--space as the antipodal 
map $x\mapsto -x$. 

It follows that the restriction of 
the mapping $\wedge$ to the level surface $M$ of the function 
$f_i$, where we assume that $\op{det}F\neq 0$ and 
$F$ and $G$ are not proportional, defines a twofold 
unbranched covering from $M$ onto the level surface of the 
functions $h_i$, for the levels given by (\ref{h1fi})---(\ref{h4fi}). 
The fibers of $\wedge |_M$ are pairs of antipodal points, 
and therefore the mapping $\wedge$ leads to an identification of 
the level surface $\wedge (M)$ on $h_2=0$ of the functions 
$h_1$, $h_3$, $h_4$ with the quotient 
$M/{\pm 1}$ of the surface $M$ by means of the antipodal 
mapping. 

Because the commutation of vector fields is a local property, 
the fact that the vector fields $\xi$ and $\eta$ on $\wedge (M)$ 
commute implies, together with the fact that $\wedge : M\to\wedge (M)$ 
is a covering, that the vector fields $\xi$ and $\eta$ on $M$ 
commute. This leads to a proof of 
Corollary \ref{commvvcor} which is based on 
the facts that $\xi$ is a Hamiltonian vector field 
on a coadjoint orbit and has the functions $f_i$ as constants 
of motion.   

\subsection{Symmetric Matrices}
\label{matsec}
The equations of motion (\ref{Adot}), in which $\omega =\omega _j(A)$ 
is given in terms of $A$ by (\ref{omegaJIu}), (\ref{Iuinv}), (\ref{u}), 
can be entirely expressed in terms of the (positive definite) 
symmetric matrix 
\begin{equation}
B:=A\, (I+\rho )^{-1}\, A^{-1}. 
\label{B}
\end{equation}
We have 
\begin{equation}
\dd{B}{t}=\left[\xi (B)_{\fop{op}},\, B\right] :=
\xi (B)_{\fop{op}}\circ B-B\circ\xi (B)_{\fop{op}}.
\label{Bdot}
\end{equation} 
Here $\xi (B)_{\fop{op}}$ denotes the antisymmetric 
linear mapping $\nu\mapsto \xi (B)\times\nu$ and 
\begin{equation}
\xi (B):=A\,\omega _j(A)=B\, j+\frac{\rho\,\langle B\, j,\, e_3\rangle}
{1-\rho\,\langle B\, e_3,\, e_3\rangle}\, B\, e_3. 
\label{xiB}
\end{equation}
The velocity (\ref{pdot}) of the point of contact $p$ also 
is a function of $B$:
\begin{equation}
\dd{p}{t}=r\,\xi (B)\times e_3.
\label{pdotB}
\end{equation}

The equation (\ref{Bdot}) is a Lax system, cf. \cite{Lax}, 
and therefore the eigenvalues of $B$ are constants of motion. 
More explicitly, it follows directly from (\ref{B}) that we have the 
the following three constants of motion in the six\--dimensional 
vector space of the symmetric $3\times 3$\--matrices:
\begin{equation}
\op{trace} B^i=\op{trace}\, (I+\rho )^{-i},\quad i=1,\, 2,\, 3.
\label{traceBi}
\end{equation}
These equations are homogeneous of degree one, two and three, 
respectively. The kinetic energy equation (\ref{phipsichi}) 
is a polynomial equation of degree two in $B$:
\begin{equation}
\langle B\, j,\, e_3\rangle ^2
-\left( \rho ^{-1}-\langle B\, e_3,\, e_3\rangle\right) \,
\left( 2T-\langle B\, j,\, j\rangle\right) =0. 
\label{TB}
\end{equation}

When the moment of inertia tensor $I$, assumed to be diagonal, 
has three different 
eigenvalues, then the mapping $A\mapsto B=A\, (I+\rho )^{-1}\, A^{-1}$ 
is a fourfold covering from $\op{SO}(3)$ onto the manifold 
determined by (\ref{traceBi}), where $A$ and $A'$ are mapped 
to the same symmetric matrix if and only $A'=A\circ R$ in which 
$R$ is one of the four diagonal rotations. Apparently the reduction 
of this symmetry group leads to quite a reduction of the degrees 
of the constants of motion. 
After passing to a covering on which the vector field 
$X$ is single valued and regular, Proposition 
\ref{mapprop} would lead to a mapping to a 
fourfold covering of the Jacobi variety of a hyperelliptic curve 
of genus two, on which $X$ corresponds 
to a constant vector field. 

The system (\ref{Bdot}) has a strong resemblance to 
the equations of van Moerbeke \cite[formula (17)]{Moerbeke} in the 
lowest dimensional case $N=3$. In his case the kinetic 
energy equation (\ref{TB}) is replaced by the condition that 
the ``modulus'', the product of the upper triangular elements of $B$ 
is kept constant.  
His equations of motion also lead to a constant 
vector field on the Jacobi variety of a hyperelliptic curve of 
genus two.

\subsection{Chaplygin}
Subsection \ref{redhorsec} reflects our understanding of 
Chaplygin's \cite[\S 4]{chaplsphere}. His 
$\lambda ,\,\lambda ',\,\mu ,\,\mu '$ in (35) 
correspond to our $a,\, b,\, c,\, d$ in (\ref{abcd}). 
The equations in his (37) correspond to our
(\ref{f1tilde}), (\ref{f2tilde}), (\ref{f3tilde}) 
in which the $f_i$ and $\widetilde{f_i}$ are given 
by (\ref{Mdef}) and the same formulas with tildes over all 
the variables. The equation 
$lD\,\lambda\,\lambda '+\mu\,\,\mu '=0$ after 
Chaplygin's formula (40) 
corresponds to our (\ref{f5tilde}) with $f_5=\widetilde{f_5}=0$. 

The sentence ``The sphere rolls in a direction perpendicular ...'' 
in front of \cite[(47)]{chaplsphere} has not been formulated 
very accurately. The velocity 
of the point of contact is neither exactly orthogonal to 
the moment, nor is it a periodic function of (the reparametrized) time. 

The formula \cite[(47)]{chaplsphere} corresponds to our description 
of $\op{d}\!\langle p,\, j\rangle /\op{d}\!\tau$. The variable 
$\xi$ in \cite[\S 4]{chaplsphere} corresponds to our 
$\langle p,\, j\times e_3\rangle$, but our equations 
for it differ from the equations which Chaplygin obtained 
for it at the end of \cite[\S 4]{chaplsphere}. In the very last 
formula in \cite[\S 4]{chaplsphere}, we believe that the 
right hand side has to be replaced by its primitive with 
respect to the time $t$, an expression which is not much more transparent 
than our constant times the primitive with respect to $\tau$ 
of $\widetilde{z}\,\widetilde{\zeta}$. 

\section{Complexification and Completion}
\label{ccsec}

\subsection{A Smooth Complex Surface}
\label{Msubsec} 
Let $M=M(c)$ denote the set of the solutions 
$((u,\, z),\, (v ,\,\zeta ))\in\C ^8$ 
of the equations $f_i=c_i$, in which $f_i$ are the quadratic 
forms defined in (\ref{fi}). In this subsection we will 
assume that $J$ 
is a diagonal matrix with three different eigenvalues $J_1$, $J_2$, 
$J_3$ on the diagonal. We will furthermore assume that 
the constants $c_i$ satisfy 
\begin{equation}
c_1\, c_3\neq 0,\; c_2=0,\; c_5=0,\;  
c_1\, c_6\neq c_3\, c_4,\;\left( c_4-J_i\, c_1\right)\,
\left( c_6-J_i\, c_3\right)\neq 0
\;\mbox{\rm for}\; i=1,\, 2\;\mbox{\rm and}\; 3.
\label{cicond}
\end{equation}  
As we have seen in Subsection \ref{redhorsec}, we can arrive at the 
first three conditions $c_1\, c_3\neq 0$, $c_2=0$, $c_5=0$ and  
$c_1\, c_6\neq c_3\, c_4$, if and only if, in the original system, 
$\op{det}F=f_1\, f_3-{f_2}^2\neq 0$ and the discriminant 
$\Delta$ of the polynomial $p:\lambda\mapsto\op{det}(G-\lambda\, F)$, 
given by (\ref{discrf}), is not equal to zero. Because 
the polynomial $p$ is invariant under the transformations 
in Subsection \ref{redhorsec}, the last condition means that 
none of the $J_i$ is a zero of $p$. Summarizing, the conditions 
mean for the original system that the polynomial $p$ is of 
second order and has two distinct zeros, none of these equal to 
one of the $J_i$'s. In Subsection \ref{redhorsec} we actually 
arranged also that $c_1=1$, as one always has for 
Chaplygin's sphere. 
 
In the case of Chaplygin's sphere, where we have (\ref{Mdef}) 
and $J_i=1/\left( I_i+\rho\right)$, these assumptions  
mean the following. 
\begin{itemize} 
\item[i)] The $J_i$ are different: 
there are three different principal moments of inertia $I_i$. 
\item[ii)] $f_1\, f_3-{f_2}^2=\| j\| ^2-{j_3}^2=
{j_1}^2+{j_2}^2\neq 0$: the moment vector $j$ is not vertical. 
\item[iii)] $\Delta =\left( 2T-\| j\| ^2/\rho\right) ^2
+8{j_3}^2\, T/\rho\neq 0$: this is automatically true when  
$j_3\neq 0$, or $j$ is not horizontal, because this 
implies that $j\neq 0$, which in turn implies that $T>0$. 
When $j$ is horizontal, then the critical energy levels are equal to 
$T_{\scriptop{crit},\, i}=\|j\| ^2/2\left( I_i+\rho\right) <\| j\| ^2/2\rho$,  
cf. (\ref{Tcrit}). Because $T$ is less than or equal to the 
largest critical energy level, we have that $2T<\| j\| ^2/2\rho$ 
when $j_3=0$, and therefore we have always that $\Delta >0$. 
\item[iv)] For Chaplygin's sphere, we have that 
\[
p(\lambda )=\left(\| j\| ^2-{j_3}^2\right)\,\lambda ^2 
-\left( 2T+\| j\| ^2/\rho\right)\,\lambda +2T/\rho .
\]
The equation (\ref{Tcrit}) for $\iota =I_i$, 
$J_i=1/\left( I_i+\rho\right)$, turns out to be equivalent to 
$p\left( J_i\right) =0$ when $T=T_{\scriptop{crit},\, i}$. 
Therefore the last condition in (\ref{cicond}) follows from the 
condition that $T$ is not equal to one of the critical 
energy levels, or that the real part 
of the level set is a smooth two\--dimensional manifold. 
\end{itemize}
After the reduction to the situation that the moment $j$ 
is a nonzero horizontal vector, we also have that 
$c_2=0$, whereas $c_5=0$ always holds for Chaplygin's sphere. 
We conclude that for Chaplygin's sphere the conditions are 
satisfied if the moment vector $j$ is not vertical 
and $T$ is in between the critical energy energy levels. 

\begin{proposition} 
The derivatives $\op{d}\! f_i$ of the functions $f_i$ defined in 
{\em (\ref{fi})} 
are linearly independent at each point of $M$, and therefore 
$M$ is a smooth complex two\--dimensional affine algebraic variety. 
\label{Msmooth}
\end{proposition}
\begin{proof}
We have to prove that 
if $((u,\, z),\, (v ,\,\zeta ))\in M$ and $\alpha _i$, $1\leq i\leq 6$ 
are constants such that $\sum_{i=1}^6\,\alpha _i\,\op{d}\! f_i=0$, 
then all the $\alpha _i$ are equal to zero. 
The equations for the $\alpha _i$ amount to 
\begin{eqnarray}
2\alpha _1\, u+\alpha _2\, v +2\alpha _4\, J\, u+\alpha _5\, J\, v &=&0,
\label{(2)}\\
2\alpha _4\, z+\alpha _5\,\zeta &=&0,\label{(1)}\\
\alpha _2\, u+2\alpha _3\, v +\alpha _5\, J\, u +2\alpha _6\, J\, v&=&0
\quad\mbox{\rm and}
\label{(3)}\\
\alpha _5\,\ z +2\alpha _6\,\zeta &=&0, 
\label{(4)}
\end{eqnarray}
corresponding to the derivatives with respect to 
$u$, $z$, $v$ and $\zeta$, respectively.  

If we take the inner product of (\ref{(2)}) with $v$, we obtain that  
\[
0=\alpha _2\, c_3-2\alpha _4\, z\,\zeta +\alpha _5\,\left( c_6-\zeta ^2\right)
=\alpha _2\, c_3+\alpha _5\, c_6,
\]
where we have used $\langle v,\, u\rangle =0$, 
$\langle v,\, v\rangle =c_3$, 
$\langle v,\, J\, v\rangle =\, -z\,\zeta$, 
$\langle v,\, J\, v\rangle =c_6-\zeta ^2$, and 
(\ref{(1)}). If we take the inner product of 
(\ref{(3)}) with $u$, we obtain that 
\[
0=\alpha _2\, c_1+\alpha _5\,\left( c_4-z^2\right) 
-2\alpha _6\, z\,\zeta =\alpha _2\, c_1+\alpha _5\, c_4,
\]
where we have used $\langle u,\, u\rangle =1$, 
$\langle u ,\, v\rangle =0$, 
$\langle u,\, J\, v\rangle =\, -z\,\zeta$, 
$\langle u,\, J\, u\rangle =c_4-z^2$, 
and (\ref{(4)}). 
These two equations for $\alpha _2$ and $\alpha _5$ lead in combination 
with $c_3\, c_4\neq c_1\, c_6$ to the conclusion that  
$\alpha _2=\alpha _5=0$. 

If we substitute this in (\ref{(2)}) then it follows that, 
unless $\alpha _1=\alpha _4=0$, the vectors $J\, u$ and $u$ are linearly 
dependent, which in turn implies that $J\, u=J_i\, u$ for some 
$i=1,\, 2,\, 3$. It follows that $\langle u,\, J\, u\rangle =
J_i\,\langle u,\, u\rangle =J_i\, c_1$, and therefore 
\begin{equation}
z^2=c_4-J_i\, c_1\neq 0.
\label{z2neq0}
\end{equation}
From (\ref{(1)}) with $\alpha _5=0$ we obtain that $\alpha _4\, z=0$, which in view of (\ref{z2neq0}) implies that $\alpha _4=0$. Now 
(\ref{(2)}) is equivalent to $\alpha _1\, u=0$, which implies that 
$\alpha _1 =0$ because $\langle u,\, u\rangle =c_1\neq 0$ implies that 
$u\neq 0$.   

Unless $\alpha _3=\alpha _6=0$, it follows from (\ref{(3)}), 
in which $\alpha _2=\alpha _5=0$, that $J\, v$ and $ v$ 
are linearly dependent, which implies that 
$J\, v =J_i\, v$ for some 
$i=1,\, 2,\, 3$. It follows that $\langle v,\, J\, v\rangle =
J_i\,\langle v,\, v\rangle =J_i\, c_3$, and therefore 
 \begin{equation}
\zeta ^2=c_6-J_i\, c_3\neq 0.
\label{zeta2neq0}
\end{equation}
From (\ref{(4)}) with $\alpha _5=0$ we obtain that 
$\alpha _6\,\zeta=0$, which in view of (\ref{zeta2neq0}) 
implies that $\alpha _6=0$. Now 
(\ref{(3)}) is equivalent to $\alpha _3\, v=0$, which implies that 
$\alpha _3 =0$ because $\langle v,\, v\rangle =c_3\neq 0$ implies that 
$v\neq 0$.   
\end{proof}

\begin{proposition}
The polynomial vector fields $\xi$ and $\eta$ are 
linearly independent at every point of $M$.
\label{XYindep}
\end{proposition} 
\begin{proof}
If $u\times J\, u=0$, which means that $u=u_i\, e_i$ for 
$i=1$, 2 or 3. This implies (\ref{z2neq0}) and therefore $z\neq 0$. 
Because $\langle u,\, v\rangle =0$, it follows 
that the vectors $\xi\, u=z\, u\times J\, v$ and  
$\eta\, u=z\, u\times v$ can only be linearly dependent if 
$v$ and $J\, v$ are linearly dependent, which implies that 
$v=v_h\, e_h$ for some $h=1$, 2 or 3. 
This implies (\ref{zeta2neq0}) with $i$ replaced by $h$ 
and therefore $\zeta\neq 0$. 
On the other hand we have that ${u_i}^2=\langle u,\, u\rangle =c_1\neq 0$, 
$0=\langle u,\, v\rangle =u_i\, v_i$ hence $v_i=0$ and therefore 
$-z\,\zeta =\langle u,\, J\, v\rangle =u_i\, J_i\, v_i=0$, which leads 
to a contradiction. In a similar way we obtain that 
$\xi$ and $\eta$ are linearly independent when 
$ v\times J\, v =0$. 

In the sequel of the proof we therefore may assume that 
$u$ and $J\, u$ are linearly independent and 
that $v$ and $J\, v$ are linearly independent. 
Assume that $z\,\zeta =0$, which in turn implies that 
$\langle u,\, J\, v\rangle =\, -z\,\zeta =0$, and therefore 
$(u\times J\, v)\times u=c_1\, J\, v$, whereas 
$(u\times v)\times v=c_1\, v$ because $\langle u,\, v\rangle =0$. 
It follows that $u\times J\, v$ and $u\times v$ 
are linearly independent. 
If $z\neq 0$ then $\zeta =0$ and we 
obtain that $\xi\, u=z\, u\times J\, v$ 
and $\eta\, u =z\, u\times  v$ are linearly independent. 
If $\zeta\neq 0$ then $z=0$ and we obtain that 
$\xi\, v =\, -\zeta\, v\times J\, u$ 
and $\eta\, v =\, -\zeta\, v\times u$ are 
linearly independent. 

If $z=\zeta =0$, then it follows from 
$0=\, -z\,\zeta 
=\langle u,\, J\, v\rangle 
=\langle J\, u,\, v\rangle $
and $\langle u, v\rangle =0$ that 
there are nonzero $a,\, b\in \C$ such that 
$u=a\, v\times J\, v$ and 
$ v =b\, u\times J\, u$. 
Inserting this in (\ref{Xz}), (\ref{Xzeta}), (\ref{Yz}), 
(\ref{Yzeta}), we obtain that 
\begin{eqnarray*}
b\, \xi\, z&=&\langle v ,\, J\, v \rangle
=c_6-\zeta ^2=c_6,\\  
-a\, \xi\,\zeta &=&\langle u,\, J\, u\rangle 
=c_4-z\,\zeta =c_4, \\
b\, \eta\, z&=&\langle v ,\, v\rangle =c_3,
\quad\mbox{\rm and}\\
-a\, \eta\,\zeta &=&\langle u,\, u\rangle =c_1 .
\end{eqnarray*}
It follows that 
\[
-a\, b\, \left[ (\xi\, z)\, (\eta\,\zeta )-(\xi\,\zeta )\, (\eta\, z)\right] 
=c_1\, c_6-c_3\, c_4\neq 0,
\]
cf. (\ref{cicond}), which implies that $\xi$ and $\eta$ are linearly independent. 

Finally, suppose that $u\times J\, u\neq 0$, 
$ v\times J\neq 0$, $z\neq 0$, $\zeta\neq 0$, 
and $\alpha\,\xi +\beta\,\eta =0$ for some $\alpha ,\,\beta\in\C$. 
It follows that there are $\gamma ,\,\delta\in \C$ such that 
\begin{eqnarray}
\alpha\,\left( z\, J\, v -\zeta\, J\, u\right) 
+\beta\, z\, v &=&\gamma\, u\quad\mbox{\rm and}
\label{abc}\\
\alpha\,\left( z\, J\, v -\zeta\, J\, u\right) 
-\beta\,\zeta \, u&=&\delta\, v .
\label{abd}
\end{eqnarray} 
If we take the inner product of (\ref{abc}) with $ v$, then we 
obtain that 
\[
0=\alpha\,\left( z\,\left[ c_6 -\zeta ^2\right] 
+\zeta\, z\,\zeta\right) +\beta\, z\, c_3
=z\, \left(\alpha\, c_6+\beta\, c_3\right) ,
\]
where we have used that $\langle v,\, u\rangle =0$, 
$\langle v,\, J\, v\rangle =c_6-\zeta ^2$,
$\langle v,\, J\, u\rangle =\, -z\,\zeta$,  
and $\langle v,\, v\rangle =c_3$. Similarly 
the inner product of (\ref{abd}) with $ u$ yields 
\[
0=\alpha\,\left( -z^2\,\zeta 
-\zeta\,\left[ c_4-z^2\right]\right) 
-\beta\,\zeta\, c_1 =\, -\zeta\,\left(\alpha\, c_4 +\beta\, c_1\right) .
\]
Because $c_6\, c_1\neq c_4\, c_3$, $z\neq 0$ and $\zeta\neq 0$, these two equations 
for $\alpha$ and $\beta$ imply that $\alpha =\beta =0$.  
\end{proof}

Let $L$ denote the manifold of the 
solutions $(u,\, v )\in\C ^3\times\C ^3$ of (\ref{SO3eq}) and (\ref{phipsichi}). 
The projection $((u,\, z),\, (v ,\,\zeta ))\mapsto (u,\, v )$ 
exhibits $M$ as a two\--fold covering of $L$  
and intertwines the vector fields $\eta$ and $\eta$ defined by 
(\ref{Xuz})---(\ref{Yzeta}) with the previously defined 
unoriented vector fields $\xi =X(u)^{1/2}\,\op{R}_{\omega}$ 
and $\eta =X(u)^{1/2}\,\op{R}_{(I+\rho )\,\omega }$ on $M_0$, where 
the word ``unoriented'' refers to the fact that the latter vector fields 
are only determined up to their signs.  

\subsection{At Infinity}
\label{infsubsec}
The (homogeneous) polynomial nature of the vector fields $\xi$ and $\eta$, 
given by (\ref{Xuz})---(\ref{Xzeta}) and (\ref{Yuz})---(\ref{Yzeta}), 
respectively, and of the constants of motion (\ref{fi}) makes it natural 
to investigate the system in the eight\--dimensional 
complex projective space $\C{\bf P} ^8$. 
The complex projective space is obtained by adding one more variable, 
which we denote 
by $\epsilon$, and then taking the quotient of $\C ^9\setminus\{ 0\}$ 
with respect to the actions $x\mapsto c\, x$ of the multiplicative 
group $\C ^{\times}$ of the nonzero complex numbers $c$. The standard 
coordinate charts correspond to the quotients of the sets of $x$ 
for which one of the coordinates, say $x_i$, is nonzero, and then 
the coordinates for this chart are obtained by putting $x_i=1$ and 
using the the $x_j$ with $j\neq i$ as the coordinates. The changes 
of coordinates are obtained by using the identification of 
$x$ with $c\, x$. Although it would be clearer to do so, we 
will not introduce different notations for the coordinates in the 
various charts, in order to avoid heavy notations. 
If we put $\epsilon =1$, then we obtain the affine space 
$\C ^8$ as a subset, equal to one of the standard coordnate 
charts, of $\C{\bf P} ^8$. The complement 
$\C{\bf P}^8_{\infty}:=\C{\bf P} ^8\setminus\C ^8$ 
corresponds to taking $\epsilon =0$ in the other coordinate charts, in which 
one of the coordinates of $((u,\, z),\, (v,\,\zeta ))$ is taken 
equal to 1. In this way $\C{\bf P} ^8\setminus\C ^8$ is identified 
with $\C{\bf P} ^7$. We will refer to $\C ^8\subset\C{\bf P} ^8$ and 
$\C{\bf P}^8_{\infty}\simeq\C{\bf P}^7$ as the affine (or finite) part of 
$\C{\bf P} ^8$ and the projective space at infinity, respectively. 

For any choice of the constants $c_i$, we denote by 
$M(c)$ the set of the solutions in $\C ^8$ of the equations 
$f_i=c_i$, the level set of the constants of motion. 
As a subset of $\C{\bf P} ^8$, the set $M(c)$ is obtained 
by homogenizing the equations. In view of the homogeneity 
of the $f_i$ of degree two, this corresponds to replacing the equations 
$f_i=c_i$ by 
\begin{equation}
g_i((u,\, z),\, (v,\,\zeta ),\,\epsilon ):=
f_i((u,\, z),\, (v,\,\zeta ))-c_i\,\epsilon ^2=0.   
\label{gieq}
\end{equation}
Let $N(c)$ denote the set of solutions of (\ref{gieq}) 
in $\C{\bf P} ^8$. Note that $M(c)=N(c)\cap\C ^8=N(c)\setminus 
\C{\bf P}^8_{\infty}$ is equal to the affine part of 
$N(c)$. 

Our goal in this section is to study the closure 
$\overline{M(c)}$ of $M(c)$ in $\C{\bf P}^8$, especially 
in the case that $M(c)=M$ with $c$ and $M$ as in Subsection \ref{Msubsec}. 
Here the closure is taken with respect to the ordinary 
topology, but it is known that $\overline{M(c)}$ is 
equal to a projective algebraic variety, and therefore 
also closed in the Zariski topology. See {\L}ojasiewicz 
\cite[p. 383]{Loj}.  

The solutions of the equations (\ref{gieq}) at $\C{\bf P}^8_{\infty}$ 
are obtained by putting $\epsilon =0$ in (\ref{gieq}), 
in which case we obtain the equations 
\begin{equation}
\begin{array}{ccc}
\langle u,\, u\rangle =0,&\langle u,\, v\rangle =0,& 
\langle v,\, v\rangle =0,\\
\langle u,\, J\, u\rangle +z^2=0,
&\langle u,\, J\, v\rangle +z\,\zeta =0,&
\langle v,\, J\, v\rangle +\zeta ^2=0.
\end{array}   
\label{gieq0}
\end{equation}
The solutions of (\ref{gieq0}) in $\C ^8$ form the conic 
affine algebraic variety $M(0)$, and the corresponding 
projective variety in $\C{\bf P}^7$, which we denote by 
$M(0)_{\infty}$, is equal to the set of 
the solutions in $\C{\bf P}^8_{\infty}$ of the equations 
(\ref{gieq}). Note that $M(0)_{\infty}=N(c)\cap\C{\bf P}^8_{\infty}$ 
does not depend on the choice of the $c_i$. 

\begin{lemma}
The set $M(0)$ is a three\--dimensional 
conic subvariety of $\C ^8$, consisting of the 
$((u,\, z),\, (v,\,\zeta ))$ such that 
$\langle u,\, u\rangle =0$, $\langle u,\, J\, u\rangle +z^2=0$ and  
the vectors $(u,\, z)$ and $(v,\,\zeta )$ in $\C ^4$ are linearly 
dependent. If $u=0$ and $z=0$, then we have to add the conditions 
that $\langle v,\, v\rangle =0$ and $\langle v,\, J\, v\rangle 
+\zeta ^2=0$.  
The corresponding projective variety  
$M(0)_{\infty}=N(c)\cap\C{\bf P}^8_{\infty}$ at infinity is a 
smooth two\--dimensional subvariety of 
$\C{\bf P}^8_{\infty}$. 
\label{M0inftylem}
\end{lemma}

\begin{proof}
If $u\neq 0$, then $\langle u,\, v\rangle =0$ implies that 
$v=u\times w$ for some $w\in\C ^3$. Using $\langle u,\, u\rangle =0$, 
we obtain that 
\[
0=\langle v,\, v\rangle =\langle u\times w,\, u\times w
=\langle u,\, u\rangle\cdot\langle w,\, w\rangle 
-\langle u,\, w\rangle ^2=\, -\langle u,\, w\rangle ^2,
\]
or $\langle u,\, w\rangle =0$, which in turn implies that 
$w=u\times a$ for some $a\in\C ^3$. But then 
\[
v=u\times (u\times a)=\langle u,\, a\rangle\, u-\langle u,\, u\rangle\, a
=\langle u,\, a\rangle\, u,
\]
which shows that $v=\lambda\, u$ for some $\lambda\in\C$. 

Now assume conversely that $u$, $z$ are solutions of $\langle u,\, u\rangle =0$, 
$\langle u,\, J\, u\rangle +z^2=0$ and that $v=\lambda\, u$ for 
some $\lambda\in C$. Then we have automatically that 
$\langle u,\, v\rangle =\lambda\,\langle u,\, u\rangle =0$ 
and $\langle v,\, v\rangle =\lambda ^2\,\langle u,\, u\rangle =0$,  
whereas the equations 
\begin{eqnarray*}
0&=&\langle u,\, J\, v\rangle +z\,\zeta 
=\lambda\,\langle u,\, J\, u\rangle +z\,\zeta 
=\, -\lambda\, z^2+z\,\zeta =z\, (\zeta -\lambda\, z),\\
0&=&\langle v,\, J\, v\rangle +\zeta ^2
=\lambda ^2\langle u,\, J\, u\rangle +\zeta ^2
=\, -\lambda ^2\, z^2+\zeta ^2=(\zeta +\lambda\, z)\, (\zeta -\lambda\, z)
\end{eqnarray*}
are equivalent to $\zeta =\lambda\, z$ or 
$z=\zeta +\lambda\, z=0$. In the second case $z=\zeta =0$, 
and therefore the conclusion is that the equations 
(\ref{gieq0}) hold if and only if $\zeta =\lambda\, z$. 

In a similar way we obtain that if $v$, $\zeta$ are solutions of 
$\langle v,\, v\rangle =0$, 
$\langle v,\, J\, v\rangle +\zeta ^2=0$ and $u=\mu\, v$, 
then the equations 
(\ref{gieq0}) hold if and only if $z=\mu\,\zeta$. 
For $\mu\neq 0$ this corresponds to the solutions in the previous 
paragraph with $\lambda =1/\mu$, whereas for $\mu =0$ 
we obtain the missing solutions with $u=0$, which implies that $z=0$ 
in view of $0=\langle u,\, u\rangle +z^2=z^2$. 
\end{proof}

Suppose that $M(c)=M$ with $c$ and $M$ as in Subsection \ref{Msubsec}. 
Then $M$ is a smooth two\--dimensional affine algebraic variety, 
and its closure $\overline{M}$ in $\C{\bf P}^8$ with respect 
to the ordinary topology is a 
projective algebraic variety, cf. {\L}ojasiewicz 
\cite[p. 383]{Loj}. It follows that the intersection 
$M_{\infty}:=\overline{M}\cap\C{\bf P}^8_{\infty}$ 
of $\overline{M}$ with the 
projective space at infinity, the set of the limit points of 
$M$ at infinity, is an algebraic variety 
in $\C{\bf P}^8_{\infty}$. It is known that in general 
$\op{dim}M_{\infty}=\op{dim}M-1$, cf. {\L}ojasiewicz 
\cite[p. 388]{Loj}, and therefore $M_{\infty}$ 
is aan algebraic curve in the projective space at infinity. 
(Actually, the explicit computations below lead to 
an independent verification of this, 
see Proposition \ref{barMprop}.) 

It follows from Lemma \ref{M0inftylem} that $N(c)=M\cup M(0)_{\infty}$, 
which implies that $\overline{M}\subset N(c)$, or 
$M_{\infty}=\overline{M}\cap M(0)_{\infty}$. 
$N(c)$ is not irreducible,
because it has the two\--dimensional 
varieties $\overline{M}$ and $M(0)_{\infty}$, which intersect 
along the curve $M_{\infty}$, as proper components. 
As we will see below, the curve $M_{\infty}=M(c)_{\infty}$ 
depends on the choice of the constants $c_i$, and actually 
the surface $M(0)_{\infty}$ is equal to the union of the 
curves $M(c)_{\infty}$ for the various $c$'s such that 
$\op{dim}M(c)=2$.  

The fact that $M(0)_{\infty}$ is higher\--dimensional than 
$M(c)_{\infty}$ is surprising, because for generic polynomials $g_i$ 
the codimension of $M(0)_{\infty}$ 
in the projective space at infinity is equal to the number of the equations. 
Because the codimension of $M(c)_{\infty}\subset M(0)_{\infty}$ 
cannot be larger, it follows that, for generic $g_i$, $M(c)_{\infty}$ 
is equal to the union of some of the 
components of $M(0)_{\infty}$. 
Compared to this, our set 
of polynomials $g_i$ is quite degenerate. We still have to determine, 
in the case that $\op{dim}M(c)=2<\op{dim}M(0)$, 
which curve in the projective surface $M(0)_{\infty}$ is equal to the 
limit curve $M(c)_{\infty}$ of $M$ at infinity.   

\begin{remark}
At the subset $M(0)$, both vector fields $\xi$ and $\eta$ are 
equal to zero. Actually, the set where both $\xi$ and $\eta$ 
are equal to zero is much larger. One component consists of the 
$((u,\, z),\, (v,\,\zeta ))$ for which the vectors 
$(u,\, z)$ and $(v,\,\zeta )$ in $\C ^4$ are linearly 
dependent, this component is five\--dimensional. 
At this component the values $c_i$ of the functions $f_i$ have the 
property that the matrices 
\begin{equation}
F:=\left(\begin{array}{cc}c_1&c_2\\c_2&c_3\end{array}\right)\quad 
\mbox{\rm and}\quad 
G:=\left(\begin{array}{cc}c_4&c_5\\c_5&c_6\end{array}\right)
\label{FG}
\end{equation}
both are singular and one is a multiple of the other. 

Further components consist, for each $i=1$, $2$ or $3$, 
of the $((u,\, z),\, (v,\,\zeta ))$ such that 
$u$ and $v$ are multiples of $e_i$ and $z$ and $\zeta$ are 
arbitrary. These components are four\--dimensional.  
At this component the matrices $F$ and $G-J_i\, F$ are singular. 

It follows that the condition that $\op{det}F\neq 0$, which allowed us in 
Subection \ref{redhorsec} to make a reduction to the 
case that $c_2=0$ and $c_1\, c_3\neq 0$, implies that 
we avoid the variety where both vector fields vanish. 

Note that if $u,\, z,\, v,\,\zeta$ are real, then the equations 
$\langle u,\, u\rangle =0$ and $\langle v,\, v\rangle =0$ imply that 
$u=0$ and $v=0$. Subsequently the equations 
$\langle u,\, J\, u\rangle +z^2=0$ and 
$\langle v,\, J\, v\rangle +\zeta ^2=0$ imply that 
$z=0$ and $\zeta =0$. Therefore $M(0)_{\infty}$ has no real 
points, which implies that $M_{\infty}$ has no real points either. 
This corresponds to the fact that the real points of 
$M$ form a compact subset of $\R ^8$.  
\end{remark}

Suppose that we are at a point 
$p=((u,\, z),\, (v,\,\zeta ),\,\epsilon )$ of $M(0)_{\infty}$, 
where $\epsilon =0$, $\langle u,\, u\rangle =0$, 
$\langle u,\, J\, u\rangle +z^2=0$, $(u,\, z)\neq (0,\, 0)$, 
and, for some 
$\lambda\in\C$, $v=\lambda\, u$ and $\zeta =\lambda\, z$. 
Recall that in the standard charts we have to put  
one of the coordinates of $u,\, z,\, v,\,\zeta $ 
identically equal to 1.
The case $u=0$, $z=0$ is covered by interchanging the role 
of $(u,\, z)$ and $(v,\,\zeta )$. 

At such a point $p$ 
the equation $\sum_{i=1}^6\,\alpha _i\,\op{d}\! g_i(p)=0$ 
for the constants $\alpha _i$, $1\leq i\leq 6$, amounts to the equations 
\begin{eqnarray*}
0&=&\left( 2\alpha _1+\lambda\,\alpha _2\right)\, u
+\left( 2\alpha _4+\lambda\,\alpha _5\right)\, J\, u, 
\quad 0=\left( 2\alpha _4+\lambda\,\alpha _5\right)\, z,\\
0&=&\left(\alpha _2+2\lambda\,\alpha _3\right)\, u
+\left(\alpha _5+2\lambda\,\alpha _6\right)\, J\, u,\quad 
0=\left(\alpha _5+2\lambda\,\alpha _6\right)\, z, 
\end{eqnarray*}
cf. (\ref{(1)})---(\ref{(4)}). If $u\neq 0$, then $J\, u$ is 
linearly independent of $u$, because otherwise $u$ would be 
a nonzero multiple of one of the basis vectors $e_i$, in 
contradiction with $\langle u,\, u\rangle =0$. 
In that case the equations are equivalent to the four 
equations 
\begin{equation}
2\alpha _1+\lambda\,\alpha _2=0,\; 2\alpha _4+\lambda\,\alpha _5=0, 
\;\alpha _2+2\lambda\,\alpha _3=0,\; \alpha _5+2\lambda\,\alpha _6=0. 
\label{alphainfty}
\end{equation}
On the other hand, if $u=0$ then $z\neq 0$ and the same conclusion holds. 
It follows that at all points of $M_{\infty}$ the rank of the matrix of the 
$\op{d}\! g_i$'s is equal to four, instead of the expected five. 

The equations (\ref{alphainfty}) are equivalent to 
$\alpha _2=\, -2\lambda\,\alpha _3$, $\alpha _1=\lambda ^2\,\alpha _3$ 
and $\alpha _5=\, -\lambda\,\alpha _6$, $\alpha _4=\lambda ^2\,\alpha _6$, 
in which $\alpha _3$ and $\alpha _6$ are free. In other words, 
at the aforementioned points of $M(0)_{\infty}$ we have that 
the derivatives at $p$ of $\lambda ^2\, g_1-2\lambda\, g_2+g_3$ and 
$\lambda ^2\, g_4-2\lambda\, g_5+g_6$ are equal to zero. 
This implies that $\op{d}\! g_3(p)$ and $\op{d}\! g_6(p)$ 
are equal to linear combinations of the $\op{d}\! g_j(p)$ with 
$j=1,\, 2,\, 4,\, 5$. 

Let $B$ denote the common zeroset in $\C{\bf P} ^8$ of the $g_j$ with 
$j=1,\; 2,\; 4,\; 5$. Because the $\op{d}\! g_j(p)$ are linearly 
independent, we have that near $p$ the set $B$ is a smooth 
four\--dimensional complex projective subvariety of $\C{\bf P} ^8$. 
The tangent space of $B$ at $p$ is equal to the common 
null space of the $\op{d}\! g_j(p)$, 
$j=1,\, 2,\, 3,\, 4$, which in turn is equal to 
the common null space of the $\op{d}\! g_i(p)$, $1\leq i\leq 6$. 
Note that in the following lemma one of the coordinates of $\widehat{p}$ 
is kept equal to zero, corresponding to the 
projective coordinate chart in which we are working

\begin{lemma}
The tangent space of $B$ at $p$ is equal to the common null space 
of the $\op{d}\! g_i(p)$, $1\leq i\leq 6$. It  
consists of the vectors 
$\widehat{p}=(\widehat{z},\,\widehat{u},\,\widehat{v},
\,\widehat{\zeta},\,\widehat{\epsilon} )$, such that 
$\langle u,\,\widehat{u}\rangle =0$, 
$\langle u,\,\widehat{v}\rangle =0$, 
$\langle J\, u,\,\widehat{u}\rangle +z\,\widehat{z}=0$ and 
$\langle J\, u,\,\widehat{v}\rangle +z\,\widehat{\zeta}=0$.
\label{hatp0}
\end{lemma}
\begin{proof}
The first equation is equivalent to $\op{d}\! g_1(p)\,\widehat p=0$. 
Assuming $v=\lambda\, u$ and the first equation, 
the second equation is equivalent 
to $\op{d}\! g_2(p)\,\widehat p=0$. The equation 
$\op{d}\! g_3(p)\,\widehat p=0$ follows from the combination of 
$v=\lambda\, u$, the first and the second equation. 

The third equation is equivalent to $\op{d}\! g_4(p)\,\widehat p=0$. 
Assuming $v=\lambda\, u$ and the third equation, 
the fourth equation is equivalent 
to $\op{d}\! g_5(p)\,\widehat p=0$. 
The equation $\op{d}\! g_6(p)\,\widehat p=0$ follows from the 
combination of 
$v=\lambda\, u$, the third and the fourth equation. 
\end{proof}

Assume that $z\neq 0$, which means that we can 
work in the projective coordinate system for which $z\equiv 1$. 
In this case $\zeta =\lambda\, z=\lambda$. Therefore, if we 
define the polynomials $g$ and $h$ by 
\begin{eqnarray}
h&:=&\zeta ^2\, g_1-2\zeta\, g_2+g_3
=\langle\zeta\, u-v,\,\zeta\, u-v\rangle -\left( c_1\,\zeta ^2-2c_2\,\zeta +c_3\right)\,\epsilon ^2\quad\mbox{\rm and}
\label{hzeta}\\
k&:=&\zeta ^2\, g_4-2\zeta\, g_5+g_6
=\langle\zeta\, u-v,\, J\, (\zeta\, u-v)\rangle 
-\left( c_4\,\zeta ^2-2c_5\,\zeta +c_6\right)\,\epsilon ^2,
\label{kzeta}
\end{eqnarray}
respectively, then $\op{d}\! h(p)=0$ and $\op{d}\! k(p)=0$ 
for all $p\in M(0)_{\infty}$. Moreover, 
the set $N(c)$, the common zeroset in $\C{\bf P} ^8$ of {\em all} the $g_i$, $1\leq i\leq 6$, is equal to 
the zeroset in $B$ of two functions $h$ and $k$. 
Note that $M(0)_{\infty}\subset N(c)\subset B$. 
Also recall 
that $M(0)_{\infty}=N(c)\cap\C{\bf P}^8_{\infty}$, 
that $M=M(c)=N(c)\setminus\C{\bf P}^8_{\infty}$, and that 
$M_{\infty}=\overline{M}\cap\C{\bf P}^8_{\infty}=$ the set of 
limit points for $\epsilon\to 0$ of solutions of (\ref{gieq}) 
with $\epsilon\neq 0$. 

Because $h$ and $k$ vanish up to second order at $p\in M(0)_{\infty}$, 
their second order Taylor expansions at $p$ are canonically 
defined quadratic forms on $\op{T}_pB$, given by 
\begin{eqnarray}
h^{(2)}(\widehat{p})&=&\langle \widehat{w},\,\widehat{w}\rangle 
-\left( c_1\,\zeta ^2-2c_2\,\zeta +c_3\right)\,\widehat{\epsilon}^2,
\label{hinfty}\\
k^{(2)}(\widehat{p})&=&\langle\widehat{w},\, J\, 
\widehat{w}\rangle 
-\left( c_4\,\zeta ^2-2c_5\,\zeta +c_6\right)\,\widehat{\epsilon}^2,
\label{kinfty}
\end{eqnarray}
in which $\widehat{w}:=\zeta\,\widehat{u}+\widehat{\zeta}\, u-
\widehat{v}$. These formulas are obtained by replacing 
$\zeta\, u-v$ and $\epsilon$ in (\ref{hzeta}) and (\ref{kzeta}) 
by their first order approximations 
$\widehat{w}$ and $\widehat{\epsilon}$. 

The equations 
$\langle u,\, u\rangle =0$ and $\langle u,\, J\, u\rangle +1=0$ 
imply that  
\begin{equation}
\op{det}(u,\, J\, u,\, u\times J\, u)=
\langle u\times J\, u,\, u\times J\, u\rangle 
=\, -\langle u,\, J\, u\rangle ^2=\, -1,
\label{detuJJ}
\end{equation}
and therefore the vectors $u$, $J\, u$ and $u\times J\, u$ 
form a basis of $\C ^3$. 
Furthermore $z\equiv 1$ implies that $\widehat{z}=0$, and the equations 
in Lemma \ref{hatp0} for $\widehat{p}\in\op{T}_pB$ imply that 
$\langle u,\widehat{w}\rangle =0$ and $\langle J\, u,\widehat{w}\rangle =0$, 
which in turn imply that $\widehat{w}=\widehat{\delta}\, u\times J\, u$ 
for some $\widehat{\delta}\in\C$. Also note that the 
condition that $\widehat{p}\in\op{T}_pM(0)_{\infty}$ is 
equivalent to $\widehat{w}=0$ and $\widehat{\epsilon}=0$, 
or $\widehat{\delta}=\widehat{\epsilon}=0$. 
 
Let $a$, $b$, $\delta$ be functions of $z,\, u,\, v,\,\zeta$, 
which together with $\epsilon$ form a regular system of coordinates for $B$ 
near $p$, in such a way that $a=b=\delta =\epsilon =0$ 
corresponds to the point $p$ and, near $p$, 
the equations $\delta =\epsilon =0$ define $M(0)_{\infty}$.
Then the tangent vector $\partial /\partial\delta$ at the origin 
corresponds to a tangent vector $\widehat{p}$ 
of $N(c)$ at $p$, such that $\widehat{\epsilon}=0$ and 
$\widehat{p}$ is not tangent to $M(0)_{\infty}$. 
We can arrange this such that $\widehat{w}=u\times J\, u$. 
At $\delta =\epsilon =0$ the functions $h$ and $k$ and their 
first order derivatives with respect to $\delta$ and $\epsilon$ 
are equal to zero. Their Taylor expansions with respect 
to $\delta$ and $\epsilon$ start with quadratic terms of the 
following special form
\begin{equation}
h^{(2)}=h_1(a,\, b)\,\delta ^2-h_2(a,\, b)\,\epsilon ^2,\quad 
k^{(2)}=k_1(a,\, b)\,\delta ^2-k_2(a,\, b)\,\epsilon ^2.
\label{h(2)k(2)}
\end{equation}
The structure of the common zeroset of $h$ and $k$ in $B$ will 
now be clarified in the following lemma. 
 
\begin{lemma}
Let $h=h(a,\, b,\,\delta ,\,\epsilon )$ and 
$k=k(a,\, b,\,\delta ,\,\epsilon )$ be two holomorphic functions 
defined in an open neighborhood of the origin in $\C ^4$. 
Assume that their Taylor expansion at $\delta =\epsilon =0$ with 
respect to $\delta$ and $\epsilon$ start with quadratic 
terms as in {\em (\ref{h(2)k(2)})}. 
Write $\Delta (a,\, b):=h_1(a,\, b)\, k_2(a,\, b)-h_2(a,\, b)\, k_1(a,\, b)$, 
so that the equation $\Delta =0$ means that the quadratic 
forms $h^{(2)}$ and $k^{(2)}$ are proportional. 
If not both $h_1(0,\, 0)$ and $h_2(0,\, 0)$ 
are equal to zero and not both $k_1(0,\, 0)$ and $k_2(0,\, 0)$ 
are equal to zero, then the origin can only be approached 
by points in the common zeroset of $h$ and $k$ for which 
$(\delta ,\,\epsilon )\neq (0,\, 0)$ if 
$\Delta (0,\, 0)=0$. 

If conversely $\Delta (0,\, 0)=0$, $h_1(0,\, 0)\neq 0$, 
$h_2(0,\, 0)\neq 0$ and the derivative at $(0,\, 0)$ of $\Delta (a,\, b)$ 
with respect to $(a,\, b)$ is not equal to zero, then 
near $(0,\, 0)$ the common zeroset of $h$ and $k$ is equal to the 
union of two smooth complex analytic surfaces  
which intersect cleanly along the smooth curve through the 
origin which is determined by the equations 
$\delta =\epsilon =0$, $\Delta (a,\, b)=0$. 
\label{hklem}
\end{lemma}

\begin{proof}
Suppose that 
$h_1(0,\, 0)\neq 0$ and $h_2(0,\, 0)\neq 0$,  
which conditions are equivalent to the condition that 
$h^{(2)}$  is a nondegenerate quadratic form in $\delta$ and $\epsilon$.  
(If $k_1(0,\, 0)\neq 0$ and 
$k_2(0,\, 0)\neq 0$, then we can interchange the roles of $h$ and $k$.)  
Let $\theta =\theta (a,\, b)$ be a square root 
of $h_2(a,\, b)/h_1(a,\, b)$ which depends holomorphically 
on $(a,\, b)$ in a neighborhood of $(0,\, 0)$. 
The Morse lemma with parameters, cf. H\"ormander \cite[Lemma 3.2.3]{Hor}, 
yields that there is a holomorphic change of the coordinates 
$(\delta ,\,\epsilon )$ to coordinates $(x,\, y)$, 
depending holomorphically on $(a,\, b)$, such that, 
near the origin, $h=h_1(a,\, b)\, x\, y$. 
We can moreover arrange that in first order approximation 
at $\delta =\epsilon =0$ we have that 
$x=\delta +\theta (a,\, b)\,\epsilon$ 
and $y=\delta -\theta (a,\, b)\,\epsilon$. 

The Taylor expansion 
of $x\mapsto k(a,\, b,\, x,\, 0)$ at $x=0$ now starts with 
a quadratic term, which implies that we can write 
$k(a,\, b,\, x,\, 0)=K(a,\, b,\, x)\, x^2$, 
in which $K(a,\, b,\, x)$ is a holomorphic function 
of $(a,\, b,\, x)$ near the origin, and 
\[
K(a,\, b,\, 0):=k_1(a,\, b)/4-k_2(a,\, b)/4\theta (a,\, b)^2
=\left( h_2\, k_1-h_1\, k_2\right) 
/4h_2. 
\]
For $x\neq 0$ the equation $k(a,\, b,\, x)=0$ is equivalent 
to the equation $K(a,\, b,\, x)=0$, and we conclude that 
the point $p$ cannot be a limit point of $M$ when 
$K(0,\, 0,\, 0)\neq 0$, or $\Delta (0,\, 0)\neq 0$. 

Assume conversely that $K(0,\, 0,\, 0)=0$, 
which means that $k_2/k_1=h_2/h_1$ at $(0,\, 0)$,  
and that the derivative at $(0,\, 0)$ of $(a,\, b)\mapsto K(a,\, b,\, 0)$ 
is not equal to zero, which is equivalent to the 
condition that the derivative at $(0,\, 0)$ of $\Delta$  
is not equal to zero. 
For instance, assume that 
$\partial K(0,\, b,\, 0)/\partial b\neq 0$ when $b=0$.  
Then the implicit function theorem yields that 
there exists a holomorphic function $B(a,\, x)$ of 
$(a,\, x)$ near $(0,\, 0)$ with $B(0,\, 0)=0$, such that,  
for $(a,\, b,\, x)$ near $(0,\, 0,\, 0)$ the equation 
$K(a,\, b,\, x)=0$ is equivalent to $b=B(a,\, x)$. 
This describes a smooth complex analytic surface, and 
we obtain the description of the common zeroset 
of $h$ and $k$ near the origin as in the lemma. 

The only case which we have not discussed yet is 
that $\Delta (0,\, 0)\neq 0$ but not $h_1(0,\, 0)\neq 0$ 
and $h_2(0,\, 0)\neq 0$ and not $k_1(0,\, 0)\neq 0$ 
and $k_2(0,\, 0)\neq 0$, for instance when 
$h_1(0,\, 0)\neq 0$, $h_2(0,\, 0)=0$, $k_1(0,\, 0)=0$, 
and $k_2(0,\, 0)\neq 0$. However, in this case we 
obtain, for an arbitrarily small positive constant $c$, 
that the points $(a,\, b,\,\delta ,\epsilon )$ 
near the origin in the zeroset of $h$ satisfy 
an estimate of the form $|\delta |\leq c\, |\epsilon |$ 
and those in the zeroset of $k$ satisfy 
$|\epsilon |\leq c\, |\delta |$, and the conclusion 
is that $\delta =\epsilon =0$ for the 
points $(a,\, b,\,\delta ,\epsilon )$ 
near the origin in the common zeroset of $h$ and $k$. 
\end{proof}

In our case the coefficients in (\ref{h(2)k(2)}) are given by
\begin{equation}
\begin{array}{cc}
h_1(a,\, b)=\langle u\times J\, u,\, u\times J\, u\rangle =\, -1,&
h_2(a,\, b)=c_1\,\zeta ^2-2c_2\,\zeta +c_3,\\
k_1(a,\, b)=\langle u\times J\, u,\, J\, (u\times J\, u)\rangle ,&
k_2(a,\, b)=c_4\,\zeta ^2-2c_5\,\zeta +c_6.
\end{array}
\label{h2k2}
\end{equation}
Until now we did not really use that 
$M(c)=M$ with $c$ and $M$ as in Subsection \ref{Msubsec}, 
but from now on this assumption will be essential.  
Then $h_2=c_1\,\zeta ^2+c_3$ and $k_2=c_4\,\zeta ^2+c_6$. 
If $h_2=0$ then $k_2\neq 0$ because of the assumption that 
$c_1\, c_6\neq c_3\, c_4$. It follows from Lemma \ref{hklem} 
that the points of $M(0)_{\infty}$ 
where $h_2=k_1=0$ do not belong to $\overline{M}$. 

Inserting $h_2=c_1\,\zeta ^2+c_3$ and $k_2=c_4\,\zeta ^2+c_6$ 
in (\ref{h2k2}) we obtain that
\begin{equation}
\Delta :=h_2\, k_1-h_1\, k_2=\left( c_1\, k_1+c_4\right)\,\zeta ^2
+c_3\, k_1+c_6, 
\label{Deltazeta}
\end{equation}
with $k_1=\langle u\times J\, u,\, J\, (u\times J\, u)\rangle$.  
Here the vectors $u$ run over the elliptic curve $E$ given by 
\begin{equation}
\langle u,\, u\rangle =0,\quad 
\langle u,\, J\, u\rangle +1=0. 
\label{Edef}
\end{equation}
The $u\in E$ together with 
the free $\zeta$ are parametrizing $M(0)_{\infty}$. 
The equation $\Delta =0$ determines a curve in the 
$(u,\,\zeta )$\--space. We have $\op{d}\!\Delta =0$ 
at a zero of $\Delta$ if and only if $\Delta =0$, 
$\left( c_1\, k_1+c_4\right)\,\zeta =0$ and 
$\left(c_1\,\zeta ^2+c_3\right)\,\op{d}\! k_1=0$. 

If $\zeta\neq 0$, then $c_1\, k_1+c_4=0$ and 
$\Delta =0$ yields that $c_3\, k_1+c_6=0$. This leads to a
contradiction with the assumption that $c_1\, c_6\neq c_3\, c_4$. 

If $\zeta =0$ then $\op{d}\! k_1=0$ because $c_1\,\zeta ^2+c_3
=c_3\neq 0$ 
and $\Delta =0$ yields that $c_3\, k_1+c_6=0$. 
The tangent space of $E$ is spanned by the vector $u\times J\, u$, 
on which $\op{d}\! k_1=0$ if and only if 
\begin{eqnarray*}
0&=&\langle (u\times J\, u)\times J\, u,\, J\, (u\times J\, u)\rangle 
+\langle u\times J\, (u\times J\, u),\, J\, (u\times J\, u)\rangle\\ 
&=&\langle u,\, J\, u\rangle\cdot\langle J\, u,\, J\, (u\times J\, u)\rangle 
-\langle J\, u,\, J\, u\rangle\cdot\langle u,\, J\, (u\times J\, u)\rangle 
=\, -\langle J\, u,\, J\, (u\times J\, u)\rangle ,
\end{eqnarray*} 
where in the first equality we used that $J$ is symmetric, 
in the second that $(u\times J\, u)\times J\, u=
\langle u,\, J\, u\rangle\, J\, u-\langle J\, u,\, J\, u\rangle\, u$, 
and in the last that $\langle u,\, J\, u\rangle =\, -1$ 
and once more that $J$ is symmetric. A straightforward 
calculation shows that 
\[
\langle J\, u,\, J\, (u\times J\, u)\rangle 
=\left( J_1-J_2\right)\,\left( J_2-J_3\right)\,\left( J_3-J_1\right)\,
u_1\, u_2\, u_3,
\]
which, in view of the assumption that the $J_i$ are different from each other, 
is equal to zero if and only if $u_i=0$ for some $i=1$, 2 or 3. 
Writing the indices of the coordinates of $u$ modulo 3, the condition 
that $u\in E$ now amounts to ${u_{i+1}}^2+{u_{i+2}}^2=0$, 
$J_{i+1}\,{u_{i+1}}^2+J_{i+2}\,{u_{i+2}}^2+1=0$, or 
${u_{i+1}}^2=1/\left( J_{i+2}-J_{i+1}\right)$, 
${u_{i+2}}^2=1/\left( J_{i+1}-J_{i+2}\right)$. 
This implies that 
\[
k_1=\langle u\times J\, u,\, J\, (u\times J\, u)\rangle 
=J_i\,\left( J_{i+2}-J_{i+1}\right) ^2\,{u_{i+1}}^2\,{u_{i+2}}^2 
=\, -J_i,
\]
and we obtain a contradiction with the equation $c_3\, k_1+c_6=0$, 
in view of the assumption that $c_6-c_3\, J_i\neq 0$. 

Applying Lemma \ref{hklem}, we obtain the conclusion 
that in the domain where $z\neq 0$ the 
curve $M_{\infty}$ coincides with the subset of $M(0)_{\infty}$ 
determined by the equation $\Delta =0$. It is smooth and near it 
$\overline{M}$ is equal to the union 
of two smooth complex analytic surfaces which 
intersect cleanly along $M_{\infty}$. 

If $z=0$ and $u\neq 0$, then we work in a chart where, 
for some $i=1$, 2 or 3, $u_i\equiv 1$. Then 
$v_i=\lambda\, u_i=\lambda$ and it becomes expedient to 
replace the functions $h$ and $k$ in (\ref{hzeta}) and 
(\ref{kzeta}) near such a point $p$ by 
\begin{equation}
\begin{array}{c}
h:={v_i}^2\, g_1-2v_i\, g_2+g_3
=\langle w,\, w\rangle -\left( c_1\,{v_i}^2 +c_3\right)\,\epsilon ^2
\quad\mbox{\rm and}
\\
k:={v_i}^2\, g_4-2v_i\, g_5+g_6
=\langle w,\, J\, w\rangle 
+(v_i\, z-\zeta )^2
-\left( c_4\, {v_i}^2+c_6\right)\,\epsilon ^2,
\end{array}
\label{hki}
\end{equation}
respectively, where $w:=v_i\, u-v$. 
The manifold $M(0)_{\infty}$ near $p$ is now 
parametrized with the curve of the $(u,\, z)$ with 
$u_i=1$, $\langle u,\, u\rangle =0$, $\langle u,\, J\, u\rangle +z^2=0$, 
and the coordinate $v_i$.  

At $z=0$ we have that $\langle u,\, u\rangle =0$ and 
$\langle u,\, J\, u\rangle =0$, which imply that 
$\langle u\times J\, u,\, u\times J\, u\rangle =0$. 
On the other hand it follows from 
Lemma \ref{hatp0} that $\langle u,\,\widehat{u}\rangle =0$, 
$\langle u,\,\widehat{v}\rangle =0$, 
$\langle J\, u,\,\widehat{u}\rangle =0$, 
$\langle J\, u,\,\widehat{v}\rangle =0$, whereas 
$\widehat{z}$ and 
$\widehat{\zeta}$ are free. It follows that 
$u$, $\widehat{u}$, $\widehat{v}$ and therefore also 
$\widehat{w}$ are multiples of 
$u\times J\, u$, and we conclude that 
$\langle\widehat{w},\,\widehat{w}\rangle =0$. 
If we let the vector $\partial /\partial\delta$ 
in the paragraph preceding (\ref{h(2)k(2)}) correspond 
this time to the vector $\widehat{p}$ such that 
$\widehat{z}=0$, $\widehat{u}=0$, $\widehat{v}=0$, 
$\widehat{\zeta}=1$ and $\widehat{\epsilon}=0$, then 
we obtain (\ref{h(2)k(2)}) with 
$h_1=0$, $h_2=c_1\, {v_i}^2+c_3$, $k_1=1$, $k_2=c_4\, {v_i}^2+c_6$. 
It follows from Lemma \ref{hklem} that, even when $k_2=0$, the point $p$ can 
only be approached by $M$ if $h_2=0$, which implies that 
$k_2\neq 0$ because of the assumption that $c_1\, c_6\neq c_3\, c_4$. 
Moreover, because $c_1\neq 0$, 
the derivative of $\Delta =h_2\, k_1=c_1\, {v_i}^2+c_3$ 
with respect to $v_i$ is nonzero when $\Delta =0$. 
Again we can apply Lemma \ref{hklem} and conclude that the curve 
$M_{\infty}$ is smooth 
at $p$, and that $\overline{M}$ near $p$ 
is equal to the union of two smooth complex analytic surfaces 
which intersect cleanly along $M_{\infty}$.  

The case that $z=0$ and $u=0$ is treated by interchanging the role 
of the vectors $(u,\, z)$ and $(v,\,\zeta )$. In this case we use the 
assumption that $c_4-c_1\, J_i\neq 0$ for every $i$. 
Again we can apply Lemma \ref{hklem} and conclude that the curve 
$M_{\infty}$ is smooth 
at $p$, and that $\overline{M}$ near $p$ 
is equal to the union of two smooth complex analytic surfaces 
which intersect cleanly along $M_{\infty}$.  

\medskip
The projection $(u,\,\zeta )\mapsto u$ exhibits 
$M_{\infty}$ as a branched covering over the 
elliptic curve $E$, where $E$ is defined by (\ref{Edef}) 
and $\zeta\in\C\cup\{\infty\}$ are the solutions of 
$\Delta =0$, with $\Delta$ as in (\ref{Deltazeta}). 
Here $\zeta =\infty$ corresponds to $z=0$, in which case 
we interchange the role of the vectors $(u,\, z)$ and 
$(v,\,\zeta )$. 
The branching occurs when $c_3\, k_1+c_6=0$ or $c_1\, k_1+c_4=0$, 
and all these branch points are simple.  
A straightforward calculation shows that the equations 
$u\in E$ and $k_1+c=0$ are equivalent to 
\begin{equation}
{u_i}^2=\left( J_i-c\right)/
\left( J_{i-1}-J_i\right)\,\left( J_i-J_{i+1}\right) ,
\quad i\in\Z /3\Z .
\label{Ebranch}
\end{equation}
It follows that there are $2\cdot 2^3=16$ branch points, 
all of which are simple, because $J_i-c\neq 0$ for 
$c=c_4/c_1$ and for $c=c_6/c_3$.  

The Riemann\--Hurwitz formula says that if one has 
an $n$\--fold branched covering from a curve $\Gamma$ onto a curve $C$ 
and $B$ is the set of branch points in $\Gamma$, then 
\begin{equation}
\op{genus}(\Gamma )-1=n\,\left(\op{genus}(C)-1\right) +
\sum_{b\in B}\,\op{order}(b)/2,
\label{RH}
\end{equation}
cf. Farkas and Kra, \cite[p. 18]{FK}. Here the order 
of the branch point $b$ is equal to $m$ if the first $m$ 
derivatives of the mapping at $b$ are equal to zero. 
Because the genus of 
an elliptic curve is equal to one, it follows that the 
genus of $M_{\infty}$ minus one is equal to $16/2=8$, 
or the genus of $M_{\infty}$ is equal to 9. 
We have proved:

\begin{proposition} 
Suppose that 
$M(c)=M$, with $c$ and $M$ as in Subsection {\em \ref{Msubsec}}. 
Then $M_{\infty}:=\overline{M}\cap\C{\bf P}^8_{\infty}$ 
is determined by the condition that the quadratic forms 
$h^{(2)}$ and $k^{(2)}$ are proportional. 
$M_{\infty}$ is a 
smooth closed algebraic curve in $\C{\bf P}^8_{\infty}$ 
of genus equal to $9$.  
Near $M_{\infty}$, the variety $\overline{M}$ is equal to the union 
of two smooth complex analytic surfaces which 
intersect cleanly along $M_{\infty}$. 
\label{barMprop}
\end{proposition}

The singularities of $\overline{M}$ can be resolved by 
considering the bundle $G$ over $\C{\bf P}^8$, of which 
the fiber $G_p$ at $p\in\C{\bf P}^8$ consists of the space 
of all two\--dimensional linear subspaces of the tangent 
space at $p$ of $\C{\bf P}^8$. Let $G_M$ denote the restriction 
of $G$ to $M$ and let $\tau _M$ be the 
section of $G_M$ which is obtained by assigning to 
$p\in M$ the tangent space $\op{T}_pM$ of $M$ at $p$, 
which is regarded as an element of $G_p$. The projection 
$\pi :\tau _M\to M$ is an isomorphism from 
$\tau _M$ onto $M$. 
Let $\widehat{M}$ denote the closure of 
$\tau _M$ in the projective variety $G$. 
Then $\widehat{M}$ is a closed smooth two\--dimensional 
subvariety of $G$. Define $\widehat{M}_{\infty}:=
\widehat{M}\setminus\op{T}M$.  
The projection $\pi :\widehat{M}\to\overline{M}$ is an 
isomorphism from the complement $\tau _M$ of 
$\widehat{M}_{\infty}$ in $\widehat{M}$, onto 
the complement $M$ of $M_{\infty}$ in $\overline{M}$. 
On the other hand $\widehat{M}_{\infty}$ is a smooth closed 
curve in $\widehat{M}$ and the projection 
$\pi :\widehat{M}_{\infty}\to M_{\infty}$ is an 
unbranched two\--fold covering. 

The mapping $\pi :\widehat{M}\to\overline{M}$ is a 
so\--called {\em normalization} of $\overline{M}$, 
a regular mapping from an irreducible normal variety 
(every smooth variety is normal) onto $\overline{M}$, 
which is a birational mapping and finite\--to\--one over 
every point of $\overline{M}$, cf. \cite[II.5.2]{Shaf}. 
Because normalizations are unique up to isomorphisms, 
one talks about {\em the} normalization of $\overline{M}$.  
Our $\pi :\widehat{M}\to\overline{M}$ is a 
simple, explicit one. 

The Riemann\--Hurwitz formula 
(\ref{RH}) yields that the genus of $\widehat{M}_{\infty}$ 
minus one is equal to $2\cdot (9-1)=16$, or that 
the genus of $\widehat{M}_{\infty}$ is equal to 17. 
In this way we obtain a smooth completion 
$\widehat{M}$ of $M$ which is obtained by adding a smooth curve 
of genus 17 at infinity. Here the word ``completion'' 
is used in the algebraic sense. Proposition 
\ref{hatMcprop} below says that it can also be used in the sense that the 
flows of $\xi$ and $\eta$, with complex times, are complete on 
$\widehat{M}$ in the sense that they 
define a transitive action on 
$\widehat{M/\Sigma}$ of the additive group $\C ^2$.

\begin{proposition}
Suppose that 
$M(c)=M$, with $c$ and $M$ as in Subsection {\em \ref{Msubsec}}. 
Let $\widehat{M}$ be the smooth completion of $M$ described above, 
the normalization of the projective clusure of $M$, 
which is obtained by adding to $M$ a smooth curve of genus $17$ at 
infinity. Then the rational vector fields $\xi$ and $\eta$ on 
$\widehat{M}$ are everywhere finite and linearly 
independent. Their respective flows $\op{e}^{t\,\xi}$ and 
$\op{e}^{s\,\eta}$ with complex times $t$ and $s$ define 
a transitive action of the additive group 
$\C ^2$ on $\widehat{M}$, and for each 
$p\in\widehat{M}$ the mapping 
$(t,\, s)\mapsto\op{e}^{t\,\xi}\circ\op{e}^{s\,\eta}(p)$ 
defines an isomorphism from the complex torus 
$\C ^2/\Lambda$ onto $\widehat{M}$. Here 
\[
\Lambda :=\{ (s,\, t)\mid \op{e}^{t\,\xi}\circ\op{e}^{s\,\eta}(p)=p\}
\] 
denotes the period lattice. It does not depend on the choice 
of $p$ and has a $\Z$\--basis consisting of four elements 
of $\C ^2\simeq\R ^4$ which are linearly independent over $\R$.  
\label{hatMcprop}
\end{proposition}

\begin{proof}
We first investigate the vector fields 
(\ref{Xuz})---(\ref{Yzeta}) near infinity when 
$z\neq 0$, where we use projective coordinates with $z\equiv 1$. 
With $\epsilon$ as the last coordinate, this means that 
we identify the affine coordinates 
$((\widetilde{u},\,\widetilde{z}),\, (\widetilde{v},\,\widetilde{\zeta}),\, 1)$ 
with $p=((u,\, 1),\, (v,\,\zeta ),\,\epsilon )$, 
where
\begin{eqnarray}
u=\widetilde{z}^{-1}\, \widetilde{u},\; 
v=\widetilde{z}^{-1}\,\widetilde{v},\;
\zeta =\widetilde{z}^{-1}\,\widetilde{\zeta},\;
\epsilon =\widetilde{z}^{-1},\quad
\mbox{\rm or}
\label{utildeu}\\
\widetilde{u}=\epsilon ^{-1}\, u,\; 
\widetilde{v}=\epsilon ^{-1}\, v,\; 
\widetilde{\zeta}=\epsilon ^{-1}\,\zeta ,\; 
\widetilde{z}=\epsilon ^{-1}.
\label{tildeuu}
\end{eqnarray}
Note that our notation means that we have to put tilde's 
over all the coordinates in the formulas (\ref{Xuz})---(\ref{Yzeta})
for the vector fields $\xi$ and $\eta$ in the affine coordinate system. 

We will write the point $p$, at which we consider the vector fields 
$\xi$ and $\eta$, as an analytic function of $\epsilon$ and 
a base point $p_0$ which varies in the curve at infinity. 
It follows that we have a convergent power series expansion  
$p=\sum_{j\geq 0}\,\epsilon ^j\, p_j$, in which 
the coefficients $p_j$ for $j\geq 1$ depend analytically 
on the point $p_0$ in the curve at infinity. We may also assume 
that the vector $p_1$ is not tangent to the curve at infinity, 
which means that it can be identified with the vector 
$\widehat{p}=\partial/\partial\delta$ in the paragraph 
preceding (\ref{h(2)k(2)}). Recall also that 
\begin{equation}
\zeta\,\widehat{u}+\widehat{\zeta}\, u-\widehat{v}=\widehat{w}=\theta\, 
u\times J\, u,   
\label{hatw}
\end{equation}
in which the nonzero factor $\theta$ is equal to a square root of 
$-\zeta ^2-1$. Note that there the coordinates 
of the base point $p_0$ are denoted by $((u,\, 1),\, (v,\,\zeta ),\, 0)$, 
instead of the $((u_0,\, 1),\, (v_0,\,\zeta _0),\, 0)$ which we will use here. 

With these notations, we have that 
\[
\xi\epsilon =\, -\widetilde{z}^{-2}\,\xi\widetilde{z} 
=\, -\epsilon ^2\,\langle\widetilde{u}\times J\,\widetilde{u},\, 
J\,\widetilde{v}\rangle
=\, -\epsilon ^{-1}\,\langle u\times J\, u,\, J\, v\rangle .
\]
The constant term in the expression following $\epsilon ^{-1}$ 
is equal to zero, because $v_0=\zeta _0\, u_0$. 
The first order term in its Taylor expansion with respect to $\epsilon$ 
is equal to 
\[
\langle u_1\times J\, u_0,\, J\, v_0\rangle 
+\langle u_0\times J\, u_1,\, J\, v_0\rangle 
+\langle u_0\times J\, u_0,\, J\, v_1\rangle 
=\langle u_0\times J\, u_0,\, J\, (v_1-\zeta _0\, u_0)\rangle 
\]
where we again have used that $v_0=\zeta _0\, u_0$. 
Using (\ref{hatw}), we obtain that $\xi\epsilon$ attains the 
finite value
\begin{equation}
\xi\epsilon =\langle u\times J\, u,\, J\,\widehat{w}\rangle 
\quad\mbox{\rm at infinity}.
\label{xieps}
\end{equation}

Using that $u=\epsilon\,\widetilde{u}$, we subsequently obtain that 
\[
\xi\, u=(\xi\epsilon )\,\widetilde{u}+\epsilon\,\xi\widetilde{u}
=\epsilon ^{-2}\,\left(\epsilon\, (\xi\epsilon )\, u 
+u\times J\, (v-\zeta\, u)\right) ,
\]
where we have used the homogeneity of $\xi$ of degree 3 and $z\equiv 1$. 
The constant term in the expression following $\epsilon ^{-2}$ 
is equal to zero, because $v_0=\zeta _0\, u_0$. 
Using (\ref{xieps}), we obtain that 
the first order term in its Taylor expansion with respect to $\epsilon$ 
is equal to 
\[
\langle u\times J\, u,\, J\,\widehat{w}\rangle \, u
+u_1\times J\, (v-\zeta \, u) 
+u\times J\, (v_1-\zeta _1\, u-\zeta\, u_1) 
=\langle u\times J\, u,\, J\,\widehat{w}\rangle \, u
-u\times J\,\widehat{w} ,
\]
where we have dropped all the subscripts 0 in the notation. The inner product 
of this expression with $u$ is equal to zero. Using that 
$\langle u,\, J\, u\rangle =1$, we obtain that the inner product 
with $J\, u$ is equal to zero as well. Finally the inner 
product with $u\times J\, u$ is equal to 
\[
-\langle u\times J\,\widehat{w},\, u\times J\, u\rangle 
=\langle u\times (u\times J\, u),\, J\,\widehat{w}\rangle 
=\, -\langle u,\, J\,\widehat{w}\rangle 
\]
because $u\times (u\times J\, u)=\langle u,\, J\, u\rangle\, u
-\langle u,\, u\rangle\, J\, u$, $\langle u,\, J\, u\rangle =\, -1$ and 
$\langle u,\, u\rangle =0$. Now it follows from 
Lemma \ref{hatp0} with $z=1$, $\widehat{z}=0$, and 
$\langle J\, u,\, u\rangle =\, -1$ that 
\[
\langle J\, u,\,\widehat{w}\rangle 
=\zeta\,\langle J\, u,\,\widehat{u}\rangle 
+\widehat{\zeta}\,\langle J\, u,\, u\rangle 
-\langle J\, u,\,\widehat{v}\rangle =0.
\]
Because $u$, $J\, u$ and $u\times J\, u$ form a basis of $\C ^3$, 
the conclusion is that the $\epsilon ^{-1}$\--term in 
the expansion of $\xi u$ in powers of $\epsilon$  
is equal to zero as well, or that $\xi u$ is finite. 

Using that $\xi$ is tangent to the surface $\widehat{M}$, 
we have obtained sufficient evidence to conclude that $\xi$ is 
finite in the complement of at most finitely many points 
of the curve $\widehat{M}_{\infty}$. In combination 
with the rationality of $\xi$ this implies that $\xi$ is finite 
on $\widehat{M}$. 

For the vector field $\eta$ we begin with 
\[
\eta\epsilon =\, -\epsilon ^{-1}\,\langle u\times J\, u,\, v\rangle .
\]
The constant term in the expression following $\epsilon ^{-1}$ 
is equal to zero, because $v_0=\zeta _0\, u_0$. 
The first order term in its Taylor expansion with respect to $\epsilon$ 
is equal to 
\[
\langle u_1\times J\, u_0,\, v_0\rangle 
+\langle u_0\times J\, u_1,\, v_0\rangle 
+\langle u_0\times J\, u_0,\, v_1\rangle 
=\langle u_0\times J\, u_0,\, v_1-\zeta _0\, u_0\rangle 
\]
where we again have used that $v_0=\zeta _0\, u_0$. 
Using (\ref{hatw}), we obtain that $\eta\epsilon$ attains the 
finite value
\begin{equation}
\eta\epsilon =\langle u\times J\, u,\, \widehat{w}\rangle 
\quad\mbox{\rm at infinity}.
\label{etaeps}
\end{equation}
Note that $\widehat{w}=\theta\, u\times J\, u$ for a nonzero
factor $\theta$, and that 
$\langle u\times J\, u,\, u\times J\, u\rangle =\, -1$, 
cf. (\ref{detuJJ}). Therefore $\eta\epsilon =\, -\theta\neq 0$ 
at every point on the curve at infinity where $z\neq 0$. 

Using that $u=\epsilon\,\widetilde{u}$, we subsequently obtain that 
\[
\eta\, u=(\eta\epsilon )\,\widetilde{u}+\epsilon\,\eta\widetilde{u}
=\epsilon ^{-2}\,\left(\epsilon\, (\eta\epsilon )\, u 
+u\times v\right) ,
\]
where we have used the homogeneity of $\eta$ of degree 3 and $z\equiv 1$. 
The constant term in the expression following $\epsilon ^{-2}$ 
is equal to zero, because $v_0=\zeta _0\, u_0$. 
Using (\ref{etaeps}), we obtain that 
the first order term in its Taylor expansion with respect to $\epsilon$ 
is equal to 
\[
\langle u\times J\, u,\,\widehat{w}\rangle \, u
+u_1\times v 
+u\times v_1
=\langle u\times J\, u,\,\widehat{w}\rangle \, u
-u\times\widehat{w} ,
\]
where we have dropped all the subscripts 0 in the notation. The inner product 
of this expression with $u$ is equal to zero. Using that 
$\langle u,\, J\, u\rangle =1$, we obtain that the inner product 
with $J\, u$ is equal to zero as well. Finally the inner 
product with $u\times J\, u$ is equal to 
$
-\langle u\times\widehat{w},\, u\times J\, u\rangle =0
$, 
because $\widehat{w}=\theta\, u\times J\, u$. 
Again using that $u$, $J\, u$ and $u\times J\, u$ form a basis of $\C ^3$, 
we obtain that the $\epsilon ^{-1}$\--term in 
the expansion of $\eta u$ in powers of $\epsilon$  
is equal to zero as well, or that $\eta u$ is finite. 
In the same way as for $\xi$, we conclude that the vector field 
$\eta$ is finite on $\widehat{M}$.  

Let $S$ denote the set of points in $\widehat{M}$ where 
$\xi$ and $\eta$ are linearly dependent. 
Proposition \ref{XYindep} implies that 
$S\cap M=\emptyset$, which means that 
$S$ is contained in the curve 
$\widehat{M}_{\infty}$ at infinity. Because $\xi$ and 
$\eta$ commute, the set $S$ 
is invariant under the flow of both vector fields, 
and it follows that at every point of $S$ both vector fields 
must be tangent to the curve at infinity. Because 
$\eta\epsilon \neq 0$ at every point of the curve at infinity 
where $z\neq 0$, we are left with the points at infinity where $z=0$. 

If $z=0$, then we have $\langle u,\, u\rangle =0$ and  
$\langle u,\, J\, u\rangle =0$ and it would follow 
that $u=0$ if $u_i=0$ for some $i$. Therefore, assuming 
that $u\neq 0$, we have for every $i$ that $u_i\neq 0$.  
In the projective coordinate chart where $u_i\equiv 1$ 
we have that $\epsilon ={\widetilde{u}_i}^{-1}$. 
It turns out that then $\eta\epsilon =0$ at $\epsilon =0$, 
which means that $\eta$ is tangent to the curve at infinity 
when $z=0$. For this reason we turn to the computation 
of $\xi\epsilon$, which is equal to the 
$i$\--th coordinate of 
\[
-\epsilon ^2\,\xi\widetilde{u}=\, -\epsilon ^{-1}\, 
u\times (z\, J\, v-\zeta\, J\, u) .
\]
Note that $v_0=\lambda\, u_0$ and $\zeta _0=\lambda\, z_0$ 
for the same factor $\lambda$, and $z_0=0$, which implies that 
$\zeta _0=0$ as well. Therefore the constant term in the expression 
after $\epsilon ^{-1}$ is equal to zero and the first order 
term in its Taylor expansion with respect to $\epsilon$ is 
equal to 
\[
(\widehat{z}\,\lambda -\widehat{\zeta})\, u\times J\, u,
\]
where we have dropped all the subscripts 0 in the notation. 
Because $\lambda =v_i$ when $u_i\equiv 1$, we conclude from 
(\ref{hki}) that the factor $\widehat{z}\,\lambda -\widehat{\zeta}$ 
is not equal to zero, where we also use that 
the equations in Lemma \ref{hatp0} with $z=0$ imply that 
$\langle\widehat{w},\, J\,\widehat{w}\rangle =0$ when 
$\widehat{w}=v_i\,\widehat{u}+\widehat{v}_i\, u-\widehat{v}$. 
Because $u\times J\, u\neq 0$, it follows that for 
at least one choice of $i$ we obtain that $\xi\epsilon\neq 0$, 
which proves that $\xi$ is not tangent to the curve at infinity 
when $z=0$. 

The case that $z=0$ and $u=0$ is treated by interchanging 
the role of the vectors $(u,\, z)$ and $(v,\,\zeta )$. 
Collecting all results, we have proved that 
$S=\emptyset$, or that $\xi$ and $\eta$ are linearly 
independent at every point of $\overline{M}$. 

Using the branched covering over $U_{\C}$ in Subsection 
\ref{sphereprojsubsec}, one obtains that the complex 
level surface $M$ is 
connected (in contrast to the real one), 
and therefore $\widehat{M}$ is connected 
as well. The remaining conclusions of the proposition 
now follow by applying the argument of Arnol'd and Avez 
\cite[Appendix 26]{ArnoldAvez} as at the end of Section \ref{volsec}. 
\end{proof}

\begin{remark}
In the complex time coordinates 
on $\widehat{M}\simeq\C ^2/\Lambda$, the 
vector fields $\xi$ and $\eta$ are constant (and 
linearly independent). 
Proposition \ref{hatMcprop} implies that the rotational motion of 
Chaplygin's sphere with horizontal moment is 
{\em algebraically integrable} according to the 
definition of Adler and van Moerbeke \cite[p. 297]{AM}.  
In view of Subsection \ref{redhorsec}, this result remains true 
for arbitrary non\--vertical moment. 

A very different proof of the algebraic integrability 
can be given by means of Chaplygin's integration of the system 
in terms of hyperelliptic integrals as 
described in Subsections \ref{ellcoordsubsec} and \ref{velellsubsec}.  
See Subsection \ref{jacsubsec}. 
\label{AMrem2}
\end{remark}

\begin{remark}
The surface $M$ is invariant under the antipodal 
mapping $x \mapsto -x$. In projective coordinates 
near infinity, where we take one of the affine 
coordinates equal to 1, this mapping is given 
by $\epsilon \mapsto -\epsilon$, keeping the 
affine coordinates fixed. The set where 
$\epsilon =0$ is the projective space at infinity, 
which belongs to the fixed point set of the 
antipodal mapping. The coordinates for $\C{\bf P}^8/{\pm 1}$ 
near infinity are obtained by replacing $\epsilon$ by $\epsilon ^2$. 
The antipodal mapping interchanges the two sheets along 
$M_{\infty}$, and it follows that 
$\overline{M}/{\pm 1}$ is a smooth variety. Its 
curve at infinity, 
$(\overline{M}/{\pm 1}) \setminus (M/{\pm 1})$, 
is isomorphic to $M_{\infty}$. 

The antipodal mapping extends to an involution 
in $\widehat{M}$ without fixed points, which 
leaves the vector fields $\xi$ and $\eta$ invariant. 
It follows that the projection from 
$\widehat{M}$ to $\overline{M}/{\pm 1}$ is a 
twofold unbranched covering, which intertwines 
$\xi$ and $\eta$ with two vector fields on 
$\overline{M}/{\pm 1}$, which we also denote by $\xi$ and 
$\eta$, which at every point are regular and linearly independent. 
Therefore the complex times of the flows of $\xi$ and $\eta$ 
lead to an identification of $\overline{M}/{\pm 1}$ with 
a complex torus, on which the vector fields $\xi$ and 
$\eta$ are constant. 

Because the antipodal mapping belongs to the group 
$\Sigma$ in (\ref{Sdef}), we obtain an eightfold 
unbranched covering $\pi :\overline{M}/{\pm 1}\to 
\widehat{M}/\Sigma$ such that the projection from $\widehat{M}$ onto 
$\widehat{M}/\Sigma$ in Proposition \ref{MSigmasmooth} below is 
equal to the composition of the twofold covering 
from $\widehat{M}$ onto $\overline{M}/{\pm 1}$, followed by 
$\pi$. In this way the curve $\widehat{M}_{\infty}/\Sigma$ of genus two 
is isomorphic to $M_{\infty}/(\Sigma /{\pm 1})$. 
\label{Mpm1rem}
\end{remark}

\subsection{A Discrete Symmetry Group}
\label{discrsymss}
Let $\Sigma$ denote the group of the 16 transformations $S$ 
in $\C ^8$ of the form 
\begin{equation}
S((u,\, z),\, (v,\,\zeta ))=\left(\left(\epsilon _1\, R\, u,\, 
\epsilon _2\, z\right) ,
\,\left(\epsilon _2\, R\, v ,\,\epsilon _1\,\zeta\right)\right) , 
\label{Sdef}
\end{equation}
in which $\epsilon _i=\pm 1$ and $R\in\op{SO}(3)$ is a 
diagonal matrix, with $\pm 1$ on the diagonal, two of them 
equal to $-1$ if $R\neq 1$. A straighforward computation show that 
every $S\in\Sigma$ leaves the functions $f_i$ in (\ref{fi}) 
invariant, and therefore leaves $M$ invariant as well, 
for any choice of the constants $c_i$. Moreover, every $S\in\Sigma$ 
alos leaves both vector fields $\xi$ and $\eta$ invariant. 
Each linear transformations $S$ has a natural extension to 
a projective linear transformation of $\C{\bf P}^8$, 
which leaves $\overline{M}$ invariant. It also has a 
natural extension to the bundle $G$ mentioned after 
Proposition \ref{hatMcprop}, and this extension leaves 
the smooth variety $\widehat{M}$ and the vector fields 
$\xi$ and $\eta$ on it 
invariant. 

\begin{proposition}
Suppose that 
$M(c)=M$, with $c$ and $M$ as in Subsection {\em \ref{Msubsec}}. 
If $S\in\Sigma$ and $S\neq 1$ then $S$ has no fixed points in 
$\widehat{M}$.  

As a consequence, the quotient $\widehat{M}/\Sigma$ 
is a smooth complex projective algebraic 
surface. The projection from $\widehat{M}$ onto 
$\widehat{M}/\Sigma$ intertwines the vector fields $\xi$ and $\eta$ 
with vector fields on $\widehat{M}/\Sigma$ which we denote by 
the same symbols. The vector fields $\xi$ and $\eta$ 
on $\widehat{M}/\Sigma$ are regular and linearly independent 
at every point, and therefore $\widehat{M}/\Sigma$ 
is isomorphic to a complex torus as well. 

Under the projection from $\widehat{M}$ onto 
$\widehat{M}/\Sigma$, the curve $\widehat{M}_{\infty}$ 
of genus $17$ is mapped onto a smooth curve 
$\widehat{M}_{\infty}/\Sigma$ of 
genus equal to $2$. 
\label{MSigmasmooth}
\end{proposition} 

\begin{proof}
Let $F$ denote the set of fixed points of $S$. 
Because the vector fields $\xi$ and $\eta$ are invariant under $S$, 
$F$ is invariant under the flows of $\xi$ and $\eta$ with complex times. 
Because these flows define a transitive action of $\C ^2$ 
on $\widehat{M}$, it follows that $F$ is either void or equal to 
$\widehat{M}$. Because it is easily 
verified that $M$ is not contained in $F$, 
the conclusion is that $S$ has no fixed points in 
$\widehat{M}$. 

The restriction to $\Gamma:=\widehat{M}_{\infty}$ 
of the projection from $\widehat{M}$ onto 
$\widehat{M}/\Sigma$ defines a 16\--fold unbranched 
covering map from $\Gamma$ onto $C:=\widehat{M}_{\infty}/\Sigma$. 
The Riemann\--Hurwitz formula \ref{RH} 
therefore yields that 
$16\,\left(\op{genus}(C)-1\right)=\op{genus}(\Gamma )-1=17-1=16$, 
which implies that $\op{genus}(C)-1=1$, or 
the genus of $C$ is equal to 2.  
\end{proof}

\begin{remark}
Because every curve of genus 2 is hyperelliptic, cf. 
Farkas and Kra \cite[Prop. III.7.2]{FK}, we conclude 
that by adding a hyperelliptic curve of genus 2 at infinity, 
the manifold $M/\Sigma$ can be completed to a 
complex torus, on which $\xi$ and $\eta$ are linearly 
independent and constant vector fields. 

The torus 
$\widehat{M}/\Sigma$ is isomorphic to the Jacobi variety 
of the hyperelliptic curve $C$ which appears in 
Chaplygin's integration by means of hyperelliptic 
integrals. See Remark \ref{hatMJacrem}. 
\label{AMrem3}
\end{remark}

\begin{remark}
Let $\widehat{M}_{\infty}/\Sigma$ denote the 
hyperelliptic curve of genus 2 which is added 
to $M/\Sigma$ at infinity in order to obtain the torus 
$\widehat{M}/\Sigma$ as the completion of $M/\Sigma$, 
cf. Remark \ref{AMrem3}. 
Let $J^{\scriptop{co}}=\op{diag}\left( 
J_2\, J_3,\, J_3\, J_1,\, J_1\, J_2\right)$ be the comatrix of $J$ 
as defined in (\ref{coJ}). The rational function 
\[
M_{\infty}\ni ((u,\, z),\, (v,\,\zeta ))\mapsto 
-\langle u,\, J^{\scriptop{co}}\, u\rangle /
\langle u,\, J\, u\rangle 
\]
induces a twofold branched covering from $\widehat{M}_{\infty}/\Sigma$ 
onto $\C {\bf P}^1$, which branches over the points 
$\lambda =J_i$ (corresponding to $u_i=0$) for $i=1,\, 2,\, 3$, 
$\lambda =\infty$ (corresponding to $z=0$), and the two 
zeros of the polynomial $p(\lambda )$ given by (\ref{HPeq}). 
For the role of $p(\lambda )$, see also iv) in Subsection \ref{Msubsec}, 
or (\ref{bicond2}) where the values $b_i$ of the functions 
$h_i$ are given in terms of the 
values $c_i$ of the functions $f_i$ by means of 
(\ref{h1fi})---(\ref{h4fi}). 
It follows that $\widehat{M}_{\infty}/\Sigma$ is isomorphic to the 
hyperelliptic curve which is defined by the equation 
\[
\mu ^2=p(\lambda )\,\prod_{i=1}^3\,\left( J_i-\lambda\right) 
\]
between the projective coordinates $(\lambda ,\,\mu )$ in 
$\C {\bf P}^2$. 

Of the six fixed points of the hyperelliptic involution 
$(\lambda ,\,\mu ) \mapsto (\lambda ,\, -\mu )$, 
the four corresponding to $\lambda =J_1,\, J_2,\, J_3,\,\infty$ 
do not depend on the values $c_i$ of the functions $f_i$, 
whereas the other two, the zeros of $p(\lambda )$, 
move freely with the $c_i$, even with the constants 
of motion $T$ and $j$ of Chaplygin's sphere. 
This means that the curves $\widehat{M}_{\infty}/\Sigma$ are 
non\--isomorphic for the generic variation of the constants 
of motion, and describe a two\--dimensional subvariety of the 
three\--dimensional moduli space of curves of genus two. 
If we are also vary the constants $J_i$ freely, then there is 
no restriction on the isomorphism class of 
the curve $\widehat{M}_{\infty}/\Sigma$. 

Remark \ref{Cmodulirem} contains an explicit verification that 
the curve $\widehat{M}_{\infty}/\Sigma$ is isomorphic to the 
hyperelliptic curve $C$ introduced in (\ref{Cdef}).  
\label{modulirem}
\end{remark}

\begin{question}
As observed in Remark \ref{AMrem3}, the torus 
$\widehat{M}/\Sigma$ is isomorphic to the Jacobi variety 
of the hyperelliptic curve $C$. According to Remark 
\ref{modulirem}, $C$ is isomorphic to 
the curve $\widehat{M}_{\infty}/\Sigma$ which is added 
at infinity to the affine algebraic surface $M/\Sigma$ 
in order to obtain the toral completion $\widehat{M}/\Sigma$. 
It follows from Matsusaka \cite{Mats} 
that  $\widehat{M}/\Sigma$ is isomorphic to the Jacobi variety 
of the curve $\widehat{M}_{\infty}/\Sigma$, and that 
$\widehat{M}_{\infty}/\Sigma$ is canonically embedded 
in its Jacobi variety $\widehat{M}/\Sigma$, if and 
only if the self\--intersection number of 
the curve $\widehat{M}_{\infty}/\Sigma$ in $\widehat{M}/\Sigma$ 
is equal to two. (I owe this reference to Ben Moonen.) 
Is it possible to verify directly that 
the self\--intersection number of 
the curve $\widehat{M}_{\infty}/\Sigma$ in $\widehat{M}/\Sigma$ 
is equal to two?
\label{compljacq}
\end{question}

\begin{remark}
In terms of the parametrization of $\widehat{M}$ by means of 
the complex times of the flows of the vector fields $\xi$ and 
$\eta$, cf. Proposition \ref{hatMcprop}, the condition that 
$S$ commutes with these flows implies that $S$ is a translation. 
For each $S\in\Sigma$ we have that $S^2=1$, cf. (\ref{Sdef}).  
Therefore, if we provide $\C ^2\simeq\R ^4$ with a real 
basis with respect to which the period lattice $\Lambda$ 
is equal to $\Z ^4$, we obtain that $S$ is equal to a translation 
over a vector $v\in \left(\frac12\Z\right) ^4/\Z ^4$. 
Because there are $2^4=16$ such vectors $v$, it follows that 
{\em $\Sigma$ is equal to the group of all translations 
of order two in the torus $\widehat{M}$}. 
In other words, the covering of $\widehat{M}\to\widehat{M}/\Sigma$ 
of the complex torus $\widehat{M}/\Sigma$ is obtained by 
replacing the period lattice $\Lambda _0$ of $\widehat{M}/\Sigma$ 
by $\Lambda =2\Lambda _0$. 

The group $\Sigma$ is closely related to the 
{\em set of theta characteristics}, as discussed in 
Mumford \cite[p. 163]{TthI}.  

I owe this remark, and the encouragement that Proposition 
\ref{hatMcprop} and Proposition \ref{MSigmasmooth} 
might be the right picture of the projective completion of $M$, 
to Frans Oort.
\label{Oortrem} 
\end{remark}

\subsection{Jordan Rizov's answer to Question \ref{compljacq}}
\label{rizovss}
The following answer to Question \ref{compljacq} has been kindly 
provided to me by Jordan Rizov. 

Let us collect in i) -- vi) below 
the abstract data we shall be working with. 
\begin{itemize}
\item[i)] Following Chapter \ref{simplesec}, we consider the 
affine variety $M\subset \C ^8$ defined by the equations 
(\ref{fi}), in which $f_i$, $1\leq i\leq 6$ are constants. 
According to Proposition \ref{Msmooth}, 
$M$ is a nonsingular two\--dimensional 
complex variety. 
\item[ii)] Consider the ``standard'' embedding $\C ^8\subset 
\C {\bf P}^8$ and let $\overline{M}$ be the projective closure 
of $M$ in $\C {\bf P}^8$ with respect to the complex topology. Then
$\overline{M}$ is also the projective closure of $M$ in the 
Zariski topolgy, because $M$ is defined as the zeroset of polynomials. 
According to Proposition \ref{barMprop}, 
$M_{\infty}:=\overline{M}\setminus M$ is a smooth algebraic 
curve of genus 9 and $\overline{M}$ is singular along $M_{\infty}$. 
\item[iii)] Consider the normalization $\pi :\widehat{M}\to\overline{M}$ 
of $\overline{M}$ as constructed after Proposition \ref{barMprop}, 
where $\widehat{M}$ is a nonsingular two\--dimensional projective 
variety. The preimage $\widehat{M}_{\infty}$ of $M_{\infty}$ 
in $\overline{M}$ is a nonsingular projective curve and 
\[
\pi :\widehat{M}_{\infty}\to M_{\infty}
\]
is an unramified two\--fold covering. Hence, by the Riemann\--Hurwitz 
theorem, the genus of $\widehat{M}_{\infty}$ is 17. 
\item[iv)] Proposition \ref{hatMcprop} says that 
$\widehat{M}$ is a two\--dimensional complex torus, and hence an Abelian 
variety (it is projective). 
\item[v)] As described in the beginning of Subsection \ref{discrsymss}, 
there is a group $\Sigma$ of order 16 acting on 
$\C ^8$ such that its action extends to an action on 
$\overline{M}$ and on $\widehat{M}$. Furthermore Proposition 
\ref{MSigmasmooth} says that $\Sigma$ acts freely on 
$\widehat{M}$. Hence, by \cite[Ch. II, \S 7, Thm. 1]{Mum74}, 
the quotient map 
\[
\pi_{\Sigma}:\widehat{M}\to\widehat{M}/\Sigma
\]
is \'etale (= an unramified covering map). 
Moreover, because the action is free, $\Sigma$ 
acts a a finite group of translations on the Abelian surface 
$\widehat{M}$ and by \cite[Ch. II, \S 7, Thm. 4]{Mum74} 
the quotient $\widehat{M}/\Sigma$ is an Abelian variety. 
\item[vi)] Since $\pi _{\Sigma}$ is \'etale, the nonsingular curve 
$\widehat{M}_{\infty}$ is mapped onto a nonsignular complete 
curve $\widehat{M}_{\infty}/\Sigma$ of genus 2. 
\end{itemize}
Before going on with any computations, let us simplify 
the notations a little bit by putting
\begin{eqnarray*}
\Gamma &:=&\widehat{M}_{\infty}/\Sigma ,\\
S&:=&\widehat{M}/\Sigma ,\\
K_S&:=&\mbox{\rm the canonical class of }S,
\end{eqnarray*}
i.e. the divisor class of a top degree differential form. 
The question posed at the end of Question \ref{compljacq} 
is whether one can compute directly the self\--instersection of 
$\Gamma$ on $S$. We will do this using the 

\medskip\noindent{\bf Adjunction formula}~~
{\em Let $\Gamma$ be a nonsingular curve of genus $g_{_{\Gamma}}$ on a 
nonsingular surface $S$ with canonical class $K_S$. Then 
the following relation holds}
\[
2g_{_{\Gamma}}-2=\Gamma\cdot (\Gamma +K_S).
\]
\begin{proof}
The proof and the construction of the intersection pairing on a 
nonsingular surface, in an ``algebraic'' way, can be found in 
\cite[Ch. 5, \S 1]{Har77}, where Proposition 1.5 is the adjunction formula. 
An ``analytic'' proof is given in \cite[Ch. 4, \S 1]{GH78}. 
\end{proof}

Therefore, in order to compute the self\--intersection 
$\Gamma\cdot\Gamma$ one has to enquire a little bit about the 
canicial class $K_S$ of $S$. As we already saw, $S$ is an Abelian 
surface, and the next result gives all we need. 

\medskip\noindent{\bf Fact}~~
{\em If $A$ is an Abelian variety of dimension $g$, then}
\[
\Omega^g_A\simeq {\cal O}_A,
\]
{\em or equivalently, the canonical class $K_A$ of $A$ is trivial.}
\begin{proof}
See for instance \cite[Ch. 1, Prop. 1.5]{vdGM} or 
\cite[Ch. 1, (5)]{Mum74}, especially (*) on page 4, and 
\cite[Ch. 2, \S 4, (4) on p. 42]{Mum74}.  
\end{proof}

\medskip\noindent{\bf The computation}~~
Applying the adjunction formula to $S$ and $\Gamma$ with 
$g_{_{|Gamma}}=2$ and $K_S$ trivial, one gets
\[
\Gamma\cdot\Gamma =2\times 2-2=2,
\]
which yields that the self\--intersection number of $\Gamma$ 
in $S$ is equal to two.  

\subsection{The system in $\C ^6=\bigwedge ^2\,\C ^4$}
\label{e3inftysec}
Let $h_i$ be the functions on $\C ^6$ given by 
(\ref{h1h2})---(\ref{h4}). 
In this subsection we will investigate the level set $L=L(b)$ 
defined by the equations 
$h_1=b_1$, $h_2=0$, $h_3=b_3$, $h_4=b_4$. 
We will assume that the constants $b_1$, 
$b_3$, $b_4$ satisfy the following conditions: 
\begin{eqnarray}
&&b_1\neq 0,\quad\left( (\op{trace}J)\, b_1+b_4\right) ^2\neq b_1\, b_3
\quad\mbox{\rm and}\label{bicond1}\\
&&p\left( J_i\right) =b_1\, {J_i}^2-\left( (\op{trace}J)\, b_1+b_4\right)\, J_i+b_3\neq 0
\quad\mbox{\rm for}\; i=1,\, 2,\, 3.\label{bicond2}
\end{eqnarray}
  
Let $U$ denote the set of $((u,\, z),\, (v,\,\zeta ))\in\C ^8$ 
such that $(u,\, z)$ and $(v,\,\zeta )$ are linearly independent, 
$U$ is an open subset of $\C ^8$.  
Let $f$ denote the mapping 
from $U$ to $\C ^6$ defined by the 
functions $f_i$ in (\ref{fi}). Let $V$ be the smooth 
hypersurface in $\C ^6\setminus\{ 0\}$ 
defined by the equation $h_2=\langle q,\, r\rangle =0$, 
and let $H:V\to\C ^3$ denote the mapping defined by 
$H_1=h_1$, $H_2=h_3$ and $H_3=(\op{trace}J)\, h_1+h_4$, 
in which the functions $h_i$ are given by 
(\ref{h1h2})---(\ref{h4}). Let $K$ be the mapping 
from $\C ^6$ to $\C ^3$ defined by 
\[
K\left( c\right) 
=\left( c_1\, c_3-{c_2}^2,\, 
c_4\, c_6-{c_5}^2,\, 
c_1\, c_6+c_3\, c_4-2c_2\, c_5\right) .
\]
Then the equations (\ref{h1fi})---(\ref{h4fi}) mean that 
$\wedge \left(\C ^8\right)\subset V$ and 
$H\circ\wedge =K\circ f$. 

Let $c\in\C ^6$ be such that 
\begin{equation}
\left( b_1,\, b_3,\, (\op{trace}J)\, b_1+b_4\right) =K(c). 
\label{bc}
\end{equation} 
Write
$F=\left(\begin{array}{cc}c_1&c_2\\c_2&c_3\end{array}\right)$ and   
$G=\left(\begin{array}{cc}c_4&c_5\\c_5&c_6\end{array}\right)$.  
Then the assumptions (\ref{bicond1}), (\ref{bicond2}) just 
mean that the polynomial $p:\lambda\mapsto\op{det}(G-\lambda\, F)$  
given by (\ref{discrf}) is of 
second order and has two distinct zeros, none of these equal to 
one of the $J_i$'s. In other words, $c$ satisfies the assumptions 
in Subsection \ref{Msubsec}. 
Note that these conditions 
imply that $F$ and $G$ are linearly independent, which in turn 
imply that $M(c)\subset U$ and that 
the Jacobi matrix $\op{T}_cK$ at the point $c$ of 
$K$ is surjective. Let $x\in M(c)$, which means that 
$f(x)=c$. Proposition \ref{Msmooth} implies that 
the Jacobi matrix $\op{T}_xf$ at the point $x$ of 
$f$ is surjective. Write $y=\wedge (x)$. Then $H(y)=H(\wedge (x))=
K(f(x))=K(c)$. The chain rule implies that 
$\op{T}_yH\circ\op{T}_x\wedge =\op{T}_cK\circ\op{T}_xf$, 
which is surjective, and therefore $\op{T}_yH$ is surjective as well. 
Using also that $\wedge$ intertwines the vector 
fields $\xi$ and $\eta$ in $\C ^8$ with the 
vector fields 
$\xi$ and $\eta$ in $\mbox{\gothic e}(3)\simeq\C ^6$, 
as we have seen in Subsection \ref{geodsubsec}, we have proved:

\begin{proposition}
$L=L(b)$ is a smooth two\--dimensional affine subvariety of $\C ^6$. 
If {\em (\ref{bc})} holds then $\wedge |_{M(c)}$ defines 
a twofold unbranched covering from $M(c)$ onto $L(b)$. 
It intertwines the vector fields $\xi$ and $\eta$ on 
$M(c)$ with the vector fields $\xi$ and $\eta$ in 
$\C ^6$. The latter vector fields are 
tangent to $L(b)$ and linearly independent at every 
point of $L(b)$. The mapping $\wedge |_{M(c)}$ induces 
a birational isomorphism from $M(c)/{\pm 1}$ onto $L(b)$, which we will 
also denote by $\wedge$.  
\label{MLprop}
\end{proposition}

The statements that $L(b)$ is smooth and 
$\xi$ and $\eta$ are linearly independent at every point 
of $L(b)$ can also be checked directly, but we found 
the proof which uses the system in $\C ^8$ simpler.  

That the inverse of the rational map 
$\wedge :M(c)/{\pm 1}\to L(b)$ is rational follows from the 
general fact that if $f:X\to Y$ is a rational map 
between irreducible varieties $X$ and $Y$ of the same dimension, 
$f(X)$ is dense in $Y$ and $f$ is injective over the preimage 
of a dense subset of $Y$, then $f$ is a birational isomorphism. 
Indeed, the homomorphism $f^*$ from the field $\C (Y)$ 
of rational functions on $Y$ to the field $\C (X)$ 
of rational functions on $X$ is injective because 
$f(X)$ is dense in $Y$. Furthermore the degree of the 
field extension of $f^*\C (Y)$ by 
$\C (X)$ is equal to the number of the 
elements of the generic fiber of $f$, cf. 
\cite[II.5.2.Thm. 7]{Shaf} (in characteristic zero 
every field extension is separable). In our case 
this implies that 
$f^*$ is surjective. Clearly $f$ has a rational inverse if 
and only $f^*:\C (y)\to\C (X)$ is an isomorphism. 
 
According to Remark \ref{Mpm1rem}, the manifold $\overline{M(c)}/{\pm 1}$ 
is a complex torus to which the vector fields $\xi$ and $\eta$ on 
$M(c)/{\pm 1}$ extend as constant vector fields. Therefore 
$\overline{M(c)}/{\pm 1}$ can be viewed as a toral completion of $L(b)$ 
to which the vector fields $\xi$ and $\eta$ on 
$L(b)$ extend as constant vector fields. This completion is obtained 
by adding a curve to infinity which is isomorphic to 
the smooth curve $M(c)_{\infty}$ of genus 9.  

\begin{proposition}
Let $\overline{M(c)}$ and $\overline{L(b)}$ denote the closure 
of $M(c)$ and $L(b)$ in $\C{\bf P}^8$ and $\C{\bf P}^6$, respectively. 
Assume that {\em (\ref{bc})} holds. Then $\wedge$ extends by 
continuity to a finite morphism from $\overline{M(c)}$ onto 
$\overline{L(b)}$, which factorizes through a morphism 
from $\overline{M(c)}/{\pm 1}$ onto $\overline{L(b)}$, 
which we also denote by $\wedge$. $\overline{M(c)}/{\pm 1}$ 
is the complex torus of Remark {\em \ref{Mpm1rem}} and 
$\wedge :\overline{M(c)}/{\pm 1}\to\overline{L(b)}$ 
is a normalization of $\overline{L(b)}$. 

The restriction of $\wedge$ to $(\overline{M(c)}/{\pm 1})\setminus 
(M(c)/{\pm 1})=\overline{M(c)}\setminus M(c)=M(c)_{\infty}$ 
maps $M(c)_{\infty}$ onto the curve 
$L(b)_{\infty}:=\overline{L(b)}\setminus L(b)$ of $L(b)$ at infinity. 
It assigns to the one\--dimensional 
linear subspace 
\[
\C\, ((u,\, z),\, (v,\,\zeta ),\, 0)\in M(c)_{\infty}
\]
of $\C ^9$ the one\--dimensional linear subspace 
\[
\C\, (z\, u,\, J\, u\times u,\, 0) 
\]
of $\C ^7$, where $(z\, u,\, J\, u\times u,\, 0)$ 
has to be replaced by $(\zeta\, v,\, J\, v\times v,\, 0)$ 
when $u=0$ and $z=0$. The image $L(b)_{\infty}$ 
is a smooth elliptic curve in $\C{\bf P}^6_{\infty}$. 
The mapping $\wedge :M(c)_{\infty}\to L(b)_{\infty}$ 
is a twofold branched covering of the curve $M(c)_{\infty}$ of genus $9$ 
over the elliptic curve $L(b)_{\infty}$, where the $16$ branch 
points in $M(c)_{\infty}$ coincide with the branch points 
mentioned in the text preceding Proposition {\em \ref{barMprop}}. 
\label{barMLprop}
\end{proposition}

\begin{proof}
The mapping $\wedge :\C ^8\to\C ^6$ extends to a homogeneous 
polynomial mapping $\wedge :\C ^9\to\C ^7$ of degree two by means 
of the formula 
\[
\wedge ((u,\, z),\, (v,\,\zeta ),\,\epsilon )=
(u\times v,\, z\, v-\zeta\, u,\,\epsilon ^2).
\]
Near infinity, where $\epsilon =0$, we can, as in the proof 
of Proposition \ref{hatMcprop}, 
write $p=((u,\, z),\, (v,\,\zeta ))$ as a convergent power series  
$p=\sum_{k\geq 0}\,\epsilon ^k\, p_k$, in which $\C\, \left( p_0,\, 0\right)$ 
varies over the curve $M(c)_{\infty}$ and the coefficients 
$p_k$ with $k\geq 1$ depend analytically on $p_0$. 
Because $\left( u_0,\, z_0\right)$ and $\left( v_0,\,\zeta _0\right)$ 
are linearly dependent when $\C\, \left( p_0,\, 0\right)\in M(c)_{\infty}$, 
we have that $\wedge\left( p_0\right) =0$. If the $i$\--th 
coordinate of 
\begin{equation}
\wedge ':=\op{D}\wedge\left( p_0\right)\, p_1=
\left( u_1\times v_0+u_0\wedge v_1,\, 
z_1\, v_0+z_0\, v_1-\zeta _1\, u_0-\zeta _0\, v_1\right) 
\label{wedge'}
\end{equation}
is nonzero, then a division of all the other coordinates of $\wedge (p)$ 
by $\wedge (p)_i$ yields the coordinates of $\wedge (p)$ in the 
standard projective coordinates in which the $i$\--th coordinates 
is kept equal to 1. Because $\wedge (p)_i=a\,\epsilon$ in which 
$a$ has a nonzero limit as $\epsilon\to 0$, we obtain that 
$\wedge (p(\epsilon ),\,\epsilon)$ converges in these coordinates 
as $\epsilon\to 0$.  
Moreover, its last coordinate $\epsilon ^2/\wedge (p)_i$ 
converges to zero as $\epsilon\to 0$, which means that the 
limit point belongs to the projective space 
$\C {\bf P}^6_{\infty}=\C {\bf P}^6\setminus\C ^6$ at infinity. 

The proof of Proposition \ref{hatMcprop} yields that the 
$r$\--component of the vector 
$\wedge '$ in (\ref{wedge'}) is nonzero whenever 
$\left( p_0,\, 0\right)\in M(c)_{\infty}$. 
When $z_0\neq 0$, we can work 
in the projective coordinate system where $z\equiv 1$, 
hence $z_0=1$, $z_1=0$, in which case the second component of 
$\wedge'$ is equal to the vector $r:=v_1-\zeta _1\, u_0-\zeta _0\, v_1$.  
In the proof of Proposition \ref{hatMcprop} we obtained that 
$r=\theta\, u_0\times J\, u_0$ for 
some nonzero factor $\theta$. Using that $v_0=\zeta _0\, u_0$, 
we find that the $p$\--component of $\wedge '$ then is 
equal to 
\[
u_0\times\left( v_1-\zeta _0\, v_1\right) =u_0\times r
=\theta\, u_0\times\left( u_0\times J\, u_0\right) =\, -\theta\, u_0,
\]
because $\langle u_0,\, u_0\rangle =0$ and 
$\langle u_0,\, J\, u_0\rangle +1=0$. 

If $z_0=0$ and $u_0\neq 0$, then we can work in the 
projective coordinate system where one of the coordinates 
$u_{0,\, i}$ of $u_0$ is identically equal to 1. 
We have $\zeta _0=0$ and $v_0=v_{0,\, i}\, u_0$, and therefore 
$r=\left( z_1\, v_{0,\, i}-\zeta _1\right)\, u_0$, 
where in the proof of Porposition \ref{hatMcprop} we obtained that 
$z_1\, v_{0,\, i}-\zeta _1\neq 0$. If $z_0$ and $u_0=0$, 
then we interchange the roles of the vectors $\left( u_0,\, z_0\right)$ 
and $\left( v_0,\,\zeta _0\right)$. 

This concludes the proof that $\wedge$ has a continuous extension 
$\wedge$ to $\overline{M(c)}$ which maps $M(c)_{\infty}$ into the 
projective space $\C{\bf P}^6_{\infty}\simeq\C{\bf P}^5$ 
at infinity. Furthermore, on the dense subset of 
$M(c)_{\infty}$ where $z\neq 0$ it matches the description 
in Proposition \ref{barMLprop}, which therefore is valid at 
every point of $M(c)_{\infty}$. We also obtain that for $z\equiv 1$ 
the restriction of $\wedge$ to $M(c)_{\infty}$ is equal to 
the composition of the projection onto $u$, followed by 
the embedding $u\mapsto \C\, (u,\, J\, u\times u)$, where 
$u$ runs over the elliptic curve $E$ defined by $\langle u,\, u\rangle =0$, 
$\langle u,\, J\, u\rangle +1=0$. This shows that 
$\wedge \left( M(c)_{\infty}\right)$ 
is isomorphic to $E$ and that the restriction of $\wedge$ 
to $M(c)_{\infty}$ is a twofold branched covering with 
the branch points as mentioned in the text preceding 
Proposition \ref{barMprop}. (It is easy to verify that the 
points on $M(c)_{\infty}$ with $z=0$ are no branch points.) 

The continuity of $\wedge$, together with $\wedge (M(c))= L(b)$ 
implies that $L(b)\subset\wedge 
\left(\overline{M(c)}\right)\subset\overline{L(b)}$. 
On the other hand, because $\overline{M(c)}$ is compact, 
the continuity of $\wedge$ 
also implies that $\wedge \left(\overline{M(c)}\right)$ is a compact, hence 
closed subset of $\C{\bf P}^6$. We conclude 
that $\wedge \left(\overline{M(c)}\right)$ is equal to the 
closure $\overline{L(b)}$ of $L(b)$ in $\C{\bf P}^6$. 

The graph of $\wedge :\overline{M(c)}/{\pm 1}\to\overline{L(b)}$ 
is equal to the projective closure of the graph of 
$\wedge :M(c)/{\pm 1}\to L(b)$, where the latter graph is 
an affine algebraic variety. It follows that the 
graph of $\wedge :\overline{M(c)}/{\pm 1}\to\overline{L(b)}$ 
is an algebraic variety, cf. {\L}ojasiewicz 
\cite[p. 383]{Loj}, which implies that 
$\wedge :\overline{M(c)}/{\pm 1}\to\overline{L(b)}$ 
is an algebraic morphism. Because it is everywhere finte 
and a birational isomorphism from $M(c)$ onto $L(b)$, 
it is a normalization of $\overline{L(b)}$. Note that 
$\overline{M(c)}/{\pm 1}$ is normal, because it is smooth. 
\end{proof}

\begin{proposition}
We have $L(b)_{\infty}=L(0)_{\infty}$ for every choice of 
the constants $b_i$ in the equations $h_i=b_i$ which define $L(b)$. 
$L(0)_{\infty}$ is a smooth elliptic curve in the 
projective space $\C{\bf P}^6_{\infty}=\C{\bf P}^6\setminus \C ^6\simeq 
\C{\bf P}^5$ at infinity. 
\label{L0prop}
\end{proposition}

\begin{proof}
Let $\delta$ be the additional projective coordinate for 
$\C {\bf P}^6$ such that $\delta =0$ corresponds to the 
projective space at infinity. (We had $\delta =\epsilon ^2$ in the 
proof of Proposition \ref{barMLprop}.) 
Then the equations $h_i=b_i$ 
for $L(b)$ correspond to the homogenized equations 
$h_i-b_i\,\delta ^2=0$. It follows that $L(b)_{\infty}$ 
is contained in the subvariety $L(0)_{\infty}$ 
of $\C{\bf P}^6_{\infty}\simeq\C{\bf P}^5$,  
which is determined by the equations $\delta =0$ and $h_i=0$. 

It follows from 
(\ref{h1h2}), (\ref{h3}) and (\ref{h4}) that the 
equations $h_1=0$, $h_3=0$, $h_4=0$ are three 
independent linear equations for ${q_1}^2$, ${q_2}^2$, 
${q_3}^2$, which have the solutions 
\begin{equation}
{q_i}^2={r_{i+1}}^2/\left( J_i-J_{i-1}\right) 
+{r_{i-1}}^2/\left( J_i-J_{i+1}\right) .
\label{q^2r}
\end{equation}
Here the index $i$ is counted modulo 3 (cyclic notation). 
This determines the $q_i$ up to their signs in terms of $r$. 
Also note that $r=0$ implies that $q=0$. because $(q,\, r)=(0,\, 0)$ 
is excluded for the projective space at infinity, we have 
always that $q\neq 0$. 

Let $\psi$ be equal to minus the product of $q_1\, r_1+q_2\, r_2+q_3\, r_3$, 
$-q_1\, r_1+q_2\, r_2+q_3\, r_3$, $q_1\, r_1-q_2\, r_2+q_3\, r_3$, 
and $q_1\, r_1+q_2\, r_2-q_3\, r_3$. Then, for given $r$, 
there exists a solution $q$ of (\ref{q^2r}) and $h_2=0$ 
if and only if there exists a solution $q$ of (\ref{q^2r}) 
and $\psi =0$. Two of the eight sign choices for 
the coordinates $q_i$, then lead to two solution $q$ and $-q$ of 
$\langle q,\, r\rangle =0$. 

On the other hand 
\begin{equation}
\psi =\sum_{i\in\Z /3\Z}\, {q_i}^4\, {r_i}^4
-2\sum_{i\in\Z/ 3\Z}\, {q_i}^2\, {q_{i+1}}^2\, {r_i}^2\, {r_{i+1}}^2.
\label{psiq}
\end{equation}
Substituting (\ref{q^2r}) in (\ref{psiq}), we obtain after a 
straightforward calculation that 
\[
\psi =\left[\sum_{i\in\Z /3\Z}\,\left( J_i-J_{i+1}\right) ^2\, 
{r_i}^2\, {r_{i+1}}^2\right] ^2/\prod_{i\in\Z /3\Z}\, 
\left( J_i-J_{i+1}\right) ^2,
\]
and therefore the equation $\psi =0$ is equivalent to 
\begin{equation}
\sum_{i\in\Z /3\Z}\,\left( J_i-J_{i+1}\right) ^2\, 
{r_i}^2\, {r_{i+1}}^2=0. 
\label{psir0}
\end{equation}

Let $F$ be the curve in $\C {\bf P}^2$ defined by (\ref{psir0}). 
The mapping $r\mapsto\left( {r_1}^2,\, {r_2}^2,\, {r_3}^2\right)$ 
defines a 16\--fold branched covering of $F$ over a nondegenerate 
quadric in $\C{\bf P}^2$, which is isomorphic to $\C{\bf P}^1$. 
The branching occurs when one of the coordinates of $r$ is 
equal to zero, in which case another coordinate of $r$ has to 
be equal to zero as well. Therefore the branching occurs at the 
three coordinate axes, where the sheets of the covering are 
connected to each other. At these points, for instance 
$r=e_3$, $q_1=1/\left( J_1-J_2\right)$, $q_2=-q_1$, $q_3=0$, 
a straightforward check shows that the derivatives of the 
functions $h_i$ defined by (\ref{h1h2}), (\ref{h3}) and (\ref{h4}) 
are linearly independent. Therefore, 
although $F$ is singular (has ordinary 
double points) at the coordinate axes, the curve $M(0)_{\infty}$ 
is smooth at the corresponding points $(q,\, r)$. 

For each $r\in F$ we have two opposite 
$q$'s which satisfy (\ref{q^2r}) and $\langle q,\, r\rangle =0$. 
The equations (\ref{q^2r}) have the solution $q=0$ if and only if 
\begin{equation}
{r_{i+1}}^2=\frac{J_i-J_{i-1}}{J_{i+1}-J_i}\, {r_{i-1}}^2,
\quad i\in\Z /3\Z, 
\label{q0req}
\end{equation}
and it is easily verified that these equations imply 
(\ref{psir0}). Therefore the projection $(q,\, r)\mapsto r$ 
defines a twofold covering from $L(0)_{\infty}$ onto $F$, 
which is branched at the four points in $F$ defined by 
(\ref{q0req}). These are smooth points of $F$. If $y$ is 
a local analytic coordinate of $F$ near such a point, and 
we substitute $y=z^2$, then we obtain that the corresponding 
points $(q,\, r)\in L(0)_{\infty}$ 
can be written as $q=z\, u\left( z^2\right) ^2$, 
$r=r\left( z^2\right)$, in which $u(y)$ and 
$r(y)$ are analytic functions of $y$ 
and $u(0)\neq 0$. It follows that $L(0)_{\infty}$ is 
smooth at $(0,\, r(0))$, and that this point is a simple 
branch point for the covering $(q,\, r)\mapsto r$. 
We conclude that $L(0)_{\infty}$ is smooth and connected, and 
therefore irreducible. Because the curve $L(b)_{\infty}$ 
is a component of $L(0)_{\infty}$, it follows that 
$L(b)_{\infty}=L(0)_{\infty}$. We know already    
from Proposition \ref{barMLprop} that $L(b)_{\infty}$ is 
an elliptic curve, but the above description can be 
also be used to verify directly that the curve $L(0)_{\infty}$ 
is elliptic.  
\end{proof}

\begin{remark}
At first sight the fact that the curve $L(b)_{\infty}$ in 
$\overline{L(b)}$ at infinity does not depend on the values of 
$b_1$, $b_3$, $b_4$ is quite disturbing. According to 
Remark \ref{Mpm1rem}, the quotient of $M(c)_{\infty}$ 
by the group $\Sigma /{\pm 1}$ is isomorphic to 
$\widehat{M}_{\infty}/\Sigma$. According to 
Remark \ref{modulirem} its isomorphism class 
varies in a two\--dimensional subvariety of the 
three\--dimensional moduli space of curves of genus two. 
As a consequence the curves $M(c)_{\infty}$ 
in general will not be isomorphic either as we vary the constants 
of motion. The question is where these moduli appear 
in the completion of $L(b)$, if the curve 
at infinity of the projective closure $\overline{L(b)}$ 
is the same for all $b$. 

The answer is that $L(b)_{\infty}=L(0)_{\infty}$ is an 
ordinary double curve of $\overline{L(b)}$ at all points 
except the branch points of the twofold covering 
$\wedge :M(c)_{\infty}\to L(0)_{\infty}$. At these 
branch points, the variety $\overline{L(b)}$ has worse 
singularities. (We conjecture that, 
as in Mumford's appendix to \cite{AM}, these are ordinary 
pinch points, where $\overline{L(b)}$ has local analytic equations 
$x^2=y\, z^2$.) According to Proposition \ref{barMLprop}, 
the branch points in $L(0)_{\infty}$ 
of $\wedge :M(c)_{\infty}\to L(0)_{\infty}$ correspond to points 
$(q,\, r)\in\C{\bf P}^5$ for which (\ref{Ebranch}) holds 
with $u$ replaced by $q$. Here $c$ runs over the two 
solutions of the equation 
\[
b_1\, c^2-\left( (\op{trace}J)\, b_1+b_4\right)\, c+b_3=0,
\]
cf. (\ref{bicond2}), and therefore the branch points 
move as a function of the moduli. 

The situation is very much similar to the description 
of $\overline{F_c}$ in Mumford's appendix to \cite{AM}, 
in the text starting with ``Thus $C$ is an ordinary 
double curve of $\overline{F_c}$ ...'' and ending with 
``... , hence $\widetilde{C}$ has genus 9'', on p. 330 and 331. 
One difference is that our normalization 
$\wedge :\overline{M(c)}/{\pm 1}\to\overline{L(b)}$ is a quite simple, 
concrete one, whereas Mumford's normalization $\pi :\widetilde{F_c}
\to\overline{F_c}$ is abstract. 

Another difference is that our normalization is equal to 
an 8\--fold unbranched covering of the Jacobi variety of 
a hyperelliptic curve of genus 2, a characterization 
which does not appear in \cite{AM}. 
The symmetry group is the group $\Sigma /{\pm 1}$, 
with $\Sigma$ as in (\ref{Sdef}), which on the $(q,\, r)$\--space acts 
by means of the transformations 
\[ 
S(q,\, r)=(\epsilon\, R\, q,\, R\, r) , 
\]
in which $\epsilon =\pm 1$ and $R$  is a 
diagonal rotation as in (\ref{Sdef}). 

If we take $L(b)$ as defined by 
$h_i=b_i$ with $b_1\neq 0$, then $S^*h_i=h_i$ when $i\neq 2$, 
but $h_2=\langle q,\, r\rangle$ satisfies 
$S^*h_2=\, -h_2$. If $h_2\neq 0$, we therefore 
can only divide out the subgroup of four elements $S$ 
for which $\epsilon =1$. According to the Riemann\--Hurwitz 
formula (\ref{RH}), the quotient of the curve 
$M(c)_{\infty}$ of genus 9 by this group of four elements 
has genus equal to 3. The possibility of arriving at a curve 
of genus 2 may therefore be related to the fact that we 
restricted ourselves to the hypersurface $\langle q,\, r\rangle =0$. 
\label{AMrem4} 
\end{remark}

\begin{question}
What happens with the sytem on the surface $h_i=b_i$ when 
$b_2\neq 0$? Is it still algebraically integrable? 
\end{question}

\subsection{Chaplygin}
The themes of Section \ref{ccsec} do not occur in 
Chaplygin \cite{chaplsphere}. Because theta functions 
are defined in terms of complex coordinates, one might argue 
that the sentence ``From (29) we see that $u$ and $v$ can be 
expressed in terms of theta functions of the two arguments 
$\alpha$ and $\tau$'' in Chaplygin 
\cite[after (30)]{chaplsphere} yields implicit evidence 
that Chaplygin did think of complex variables, as does the   
sentence ``Solving equation (41) gives two real values for the 
quantity $f$'' in Chaplygin \cite[after (41)]{chaplsphere}. 
However, the inequalities between (27) and (28) in 
Chaplygin \cite[\S 3]{chaplsphere} indicate that Chaplygin 
mainly focussed on the real system, whereas he also emphasizes that 
(41) has real solutions. Completion of the complexified system and 
tori (real or complex) definitely do not occur at all in Chaplygin 
\cite{chaplsphere}.

\section{Hyperelliptic Integrals}
\label{horjsec}
In this section we assume that the moments of inertia $I_i$ 
are different from each other, that $\rho\neq 0$, 
that the constants of motion $(j,\, T)$ are 
at a nonsingular level. 
We also assume in this section that the {\em the moment $j$ of 
the momentum around the point of contact is nonzero and horizontal}, 
which means that $\| j\| ^2\neq 0$ and 
$j_3=0$ in (\ref{SO3eq}). As shown in Subsection \ref{redhorsec}, 
the rotational motion with arbitrary nonvertical $j$ 
can be reduced to this case.   

In order to obtain a smooth level surface, we will assume that the kinetic 
energy $T$ is not equal to any of the critical levels 
\begin{equation}
T_{\scriptop{crit},\, i}:=\| j\| ^2/2\left( I_i+\rho\right) ,
\quad i=1,\, 2,\, 3, 
\label{Tcrithor}
\end{equation}
of the function $T_j$ on $\op{SO}(3)$, cf. (\ref{Tcrit}) with $j_3=0$.  
In order to obtain that the complex level surface is smooth, we 
will need furthermore that  
\begin{equation}
2T\,\rho\neq \| j\| ^2.
\label{Tneqj}
\end{equation}
If $M$ has real points, then (\ref{Tneqj}) is a consequence of 
the assumption that $T\neq T_{\scriptop{crit},\, i}$, because 
$T<T_{\scriptop{crit},\, 1}=\| j\| ^2/2\left( I_1+\rho\right)
<\| j\| ^2/2\rho $, 
cf. (\ref{Tcrithor}), and therefore $2T\,\rho <\| j\| ^2$. 
In the case that we allow arbitrary complex values for 
the parameters $I_i$, $\rho$ and $\| j\| ^2$, we have to add 
(\ref{Tneqj}) to the list of conditions. In other words, we 
make the same assumptions as in Subsection \ref{Msubsec}.

\subsection{The Projection onto the First Vector}
\label{sphereprojsubsec}
Our next goal is to simplify the vector field $\xi$ by means of a 
suitable substitution of variables in the $u$\--space. 
We recall that $u$ has the concrete interpretation that $-r\, u$ is 
equal to the point of contact on the surface of the sphere in body 
coordinates, cf. (\ref{s}). 

We have the two complex surfaces $L$ and $U_{\C}$ in which 
$L$ is the surface in $\C ^3\times\C ^3$ determined by the equations (\ref{SO3eq}) and (\ref{phipsichi}) and $U{_\C}$ is the quadric 
\begin{equation}
U_{\C}:=\{ u\in\C ^3\mid\langle u,\, u\rangle =1\} 
\label{UC}
\end{equation}
in $\C ^3$. The projection $(u,\, v)\mapsto u$ is a 
branched covering from $L$ onto $U_{\C}$, 
branching over the set of $u\in U_{\C}$ for which there exists a solution 
$ v\in\C ^3$ of $\langle u,\, v\rangle =0$, 
$\langle v ,\, v\rangle =\| j\| ^2$, and (\ref{phipsichi}), where the 
derivatives of $\langle u,\, v\rangle$, $\langle v ,\, v\rangle$, 
and $f(u,\, v )$ with respect to $ v$ are linearly dependent.  
If $X (u)=0$, then $ v$ is a solution of (\ref{phipsichi}) 
if and only if $Y (u,\, v )=0$, in which case 
\[
\partial f(u,\, v )/\partial v =
2Y (u,\, v )\,\partial Y (u,\, v )/\partial v
-X (u)\, 
\partial Z ( v )/\partial v =0. 
\]
Therefore the zeroset of $X$, which does not contain any real 
points in view of the remark preceding Lemma \ref{vollem}, 
is contained in the branch locus. 

When $u\in U_{\C}$ and $X (u)\neq 0$, then $u$ is a branch point if 
and only if there exists a solution $ v$ of (\ref{SO3eq}) and (\ref{phipsichi}) such that the vectors $u$, $ v$ and 
\begin{equation}
w:=X(u)\,\omega =Y(u,\, v)\, J\, u+X(u)\, J\, v
\label{partialf3}
\end{equation} 
are linearly dependent. Here we have used the formula 
(\ref{omegausigma}) for $\omega$. Assuming that $j$ is not vertical, 
it follows from (\ref{sigmadef}) that $u$ and $ v$ are linearly 
independent, and we obtain that $w =\alpha\, u+\beta\, v$ 
for suitable constants $\alpha$ and $\beta$. 
With the abbreviations $X=X(u)$, $Y=Y(u,\, v)$, we have 
\[
\alpha =\langle u,\, w\rangle =
Y\,\langle u,\, J\, u\rangle +X\,\langle u,\, J\, v\rangle 
=Y/\rho ,
\]
cf. (\ref{phidef}) and (\ref{chidef}), and 
\[
\beta\,\| j\| ^2 =\langle v,\, w\rangle =Y\,\langle v,\, J\, u\rangle 
+X\,\langle v,\, J\, v\rangle =2T\, X,
\]
cf. (\ref{chidef}), (\ref{phipsichi}), and (\ref{psidef}). 
With the notation  
\begin{equation}
\tau :=2T/\| j\| ^2,   
\label{taudef}
\end{equation}
the equation $w =\alpha\, u+\beta\, v$ takes the form 
$X\, \left( J-\tau\right)\, v =Y\, \left( 1/\rho -J\right)\, u$. 
Because $J^{-1}=I+\rho$, cf. (\ref{Jdef}), we have that 
$J^{-1}\, (1-\rho\, J)=I$ and it follows that $v$ is 
equal to a nonzero multiple of $K\, u$, in which 
\begin{equation}
K:=\left[ 1-\tau\, (I+\rho )\right]^{-1}\, I. 
\label{Kdef}
\end{equation}
Because $\langle u,\, v\rangle =0$, we arrive 
at the conclusion that 
\begin{equation}
\langle u,\, K\, u \rangle =0. 
\label{sigmabranch}
\end{equation}
In other words, we have proved that away from $ X (u)=0$ 
the branch locus is contained in the 
intersection of the quadric $U_{\C}$ with the quadratic cone 
defined by (\ref{sigmabranch}). Because this intersection 
is irreducible, the conclusion is that the branch locus 
away from $ X (u)=0$ is equal to the intersection 
of $U_{\C}$ with the quadratic cone defined by (\ref{sigmabranch}). 

\begin{remark}
The matrix $K$ in (\ref{Kdef}) 
is a diagonal matrix with eigenvalues equal to 
\begin{equation}
K_i=I_i/\left( 1-\tau\, (I_i+\rho )\right)
=I_i/\left( 1-T/T_{\scriptop{crit},\, i}\right),
\label{Bi}
\end{equation}
in which the $T_{\scriptop{crit},\, i}$  
are the critical levels of $T_j$, cf. (\ref{taudef}) and (\ref{Tcrithor}). 
The assumption that we are on a regular level set is equivalent 
to the condition that $T\neq T_{\scriptop{crit},\, i}$ for every 
$i=1,\, 2,\, 3$. If the level set contains real points, then 
we have that $T_{\scriptop{crit},\, 3}<T<T_{\scriptop{crit},\, 2}<
T_{\scriptop{crit},\, 1}$ or 
$T_{\scriptop{crit},\, 3}<T_{\scriptop{crit},\, 2}<
T<T_{\scriptop{crit},\, 1}$. In the first case $K$ has one negative and 
two positive eigenvalues and in the second case $K$ has 
two negative and one positive eigenvalue. Therefore, in both cases 
the cone defined by (\ref{sigmabranch}) has nonzero real points. 
The real cone intersects the real unit sphere $U$ in two closed curves, 
diametrically opposite to each other. 

Because each of the two connected components of the real level surface 
of $(j,\, T)$ is a torus, it can only be mapped onto 
the annular region between these two curves. 
(This can also be verified by means of explicit calculations, 
cf. the last part of Remark \ref{sigmaremark}.)  
The images of the quasiperiodic 
solution curves on the level set, cf. Corollary \ref{commvvcor}, 
while running around the sphere in this band, run from one of the 
bounding curves to the other, every time with a contact of order 
two at the bounding curve. A similar behaviour occurs when 
$j$ is neither horizontal nor vertical, but in that case the 
bounding curves are determined by more complicated equations. 
\label{uremark}
\end{remark}

\begin{remark}
If $j_3=0$ then, for given $u\in U_{\C}$, the solutions $ v$ 
of the equations (\ref{SO3eq}) and (\ref{phipsichi}) can be explicitly 
computed in the following way. Using the equations 
$\langle v ,\, v\rangle =\| j\| ^2$  
and (\ref{taudef}), 
we can write the kinetic energy equation (\ref{phipsichi}) in the form 
\begin{equation}
X(u)\,
\langle v ,\,\left( J-\tau\right) 
\, v\rangle +Y(u,\, v )^2 =0, 
\label{Thom}
\end{equation} 
which is homogenous (of degree two) in the variable $v$. 
The condition that $j$ is horizontal 
corresponds to the equation 
\begin{equation}
\langle u,\, v\rangle =0. 
\label{usigmazero}
\end{equation}
In solving the homogenous equations (\ref{Thom}) and 
(\ref{usigmazero}) for $v$ with $ v _3\neq 0$ we 
may put $ v _3=1$. If $u_2\neq 0$, then we can 
solve $ v _2$ in terms of $ v _1$ from (\ref{usigmazero}). 
This leads to a quadratic equation for $ v _1$ of the form 
\[
a_3(u)\,{ v _1} ^2+2 b_2(u)\, v _1+a_2(u)=0,
\]
in which 
\begin{eqnarray}
a_3(u)&:=&X(u)\,\left( c_1\, {u_2}^2+c_2\, {u_1} ^2\right) 
+\left( c_1-c_2\right) ^2\,{u_1} ^2\, {u_2} ^2,
\label{a3def}\\
b_2(u)&:=&\left[X(u)\, c_2+\left( c_1-c_2\right) 
\,\left( c_3-c_2\right)\, {u_2} ^2\right]\, u_1\, u_3
\quad\mbox{\rm and}
\label{b2def}\\
a_1(u)&:=&X(u)\,\left( c_2\, {u_3}^2+c_3\, {u_2} ^2\right) 
+\left( c_2-c_3\right) ^2\,{u_2} ^2\, {u_3} ^2.
\label{a1def}
\end{eqnarray}
Here we have used the abbreviation 
\begin{equation}
c_i:=1/\left( I_i+\rho\right) -\tau .
\label{tauci}
\end{equation}
With these notations, we obtain, using repeatedly that  
$\langle u,\, u\rangle =1$, that 
\begin{equation}
 v _1=\left( -b_2(u)\pm\, \Delta (u)^{1/2}\right) /a_3(u), 
\label{sigma1}
\end{equation}
in which the discriminant turns out to be given by 
\begin{equation}
\Delta _2(u):=b_2(u)^2-a_3(u)\, a_1(u)
=\, - c_1\, c_2\, c_3\,{u_2} ^2\, X(u)\,   
\langle u,\, K\, u\rangle /\rho ,
\label{Delta2}
\end{equation}
with $K$ as in (\ref{Kdef}). 

The solution vector $v$ is determined by taking  
$v_2:=\, -\left( u_1\, v_1+u_3\right) /u_2$. 
In order to obtain a vector of length $\| j\|$, we 
have to replace $v$ by $c\, v$, in which  
$c= \pm\,\| j\|/\| v\|$. 
We obtain four solutions, consisting of two oppositie pairs. 

As expected, 
the discriminant is equal to a multiple of $X (u)\, 
\langle u,\, K\, u\rangle$. Note that for real $u$ we have 
real solutions $ v$ if and only if 
$c_1\, c_2\, c_3\,\langle u,\, K\, u\rangle\leq 0$.  
It follows from (\ref{Bi}), (\ref{taudef}) and (\ref{tauci}) that 
$c_i\, B_i=I_i/\left( I_i+\rho\right) >0$, 
and therefore the determinant of $c_1\, c_2\, c_3\, K$ 
is positive, which implies in view of 
Remark \ref{uremark} that $c_1\, c_2\, c_3\, K$ 
has two negative eigenvalues. It follows that the 
part on the unit sphere where 
$c_1\, c_2\, c_3\,\langle u,\, K\, u\rangle\leq 0$ 
is the connected annular region, bounded by the two smooth curves 
determined by the equations $\langle u,\, u\rangle =1$, 
$\langle u,\, K\, u\rangle =0$. 

For $u$ in the domain $c_1\, c_2\, c_3\,\langle u,\, K\, u\rangle  <0$, 
the interior of the annulus \label{sigmaremark}, we obtain four 
solutions $ v$, which correspond to four possibilities for 
the velocity vector $\dd{u}{t}$. One pair 
of these velocity vectors correspond to one of the 
connected components of the level set, and their opposites 
correspond to the other connected component, cf. the discussion 
at the end of Section \ref{levelsec}. These connected components 
are mapped to each other by means of a rotation over $\pi$ which maps 
$j$ to its opposite, cf. the discussion after Lemma \ref{Bottlem}. 
\label{branchrem}
\end{remark}

The intersection of the quadric $\langle u,\, u\rangle =1$ 
with the cone $\langle u,\, K\, u\rangle =0$ is equal to the 
intersection of the quadric $\langle u,\, u\rangle =1$ with 
the quadric $\langle u, (1+\beta\, K)\, u\rangle =1$, 
in which $\beta$ is any nonzero constant. If we choose 
$\beta =\tau$, then we obtain that the branch locus is equal to 
the intersection with $U_{\C}$ of the quadric defined by the equation 
\begin{equation}
\langle u,\, (1-\gamma\, I)^{-1}\, u\rangle =1,
\quad\mbox{\rm in which}\quad \gamma :=\tau /(1-\tau\,\rho ).
\label{gamma}
\end{equation}
The equation $X (u)=0$ is equivalent to 
$\langle u,\,\rho\, (I+\rho )^{-1}u\rangle 
=1$, cf. (\ref{phidef}) and (\ref{Jdef}).  
The three quadrics $\langle u,\, u\rangle =1$, $X (u)=0$ 
and (\ref{gamma}) therefore belong to the one\--parameter family 
of quadrics 
$\langle u,\, (1-\lambda\, I)^{-1}\, u\rangle =1$, 
in which the parameter $\lambda$ takes the values $\lambda =0$, 
$\lambda =\, -1/\rho$ and $\lambda =\gamma$, respectively. 
The substitution of variables $u=I^{1/2}\, x$ leads to 
\[
\langle u,\, (1-\lambda\, I)^{-1}\, u\rangle 
=\langle x,\, I\, (1-\lambda\, I)^{-1}\, x\rangle
=\langle x,\, \left( I^{-1}-\lambda\right) ^{-1}\, x\rangle ,
\]
and our family of quadrics is turned into the 
{\em one\--parameter family (pencil) of confocal quadrics} 
\begin{equation}
\sum_{i=1}^3\frac{{x_i} ^2}{a_i-\lambda}=1, 
\label{confocal}
\end{equation}  
in which
\begin{equation}
a_i:=1/I_i,\quad\mbox{\rm hence}\quad u_i=x_i/{a_i}^{1/2}. 
\label{aiIi}
\end{equation}

\subsection{Jacobi's Elliptic Coordinates}
\label{ellcoordsubsec}
Given $x\in\C ^3$, the equation (\ref{confocal}) corresponds to 
a polynomial equation of degree three for $\lambda$, which 
has three solutions $\left(\lambda _1,\,\lambda _2,\,\lambda _3\right)$, 
which are called {\em Jacobi's elliptic coordinates}, 
cf. Jacobi \cite[26. Vorlesung]{Jacobi}. 
We briefly recall some of Jacobi's   
observations.  

Let $a_i$, $1\leq i\leq n$, be $n$ different numbers, and let 
$s\in\Z _{\geq 0}$.  
Then the function $z\mapsto z^s/\prod_k\,\left( z-a_k\right)$ 
has simple poles at the points $z=a_i$, with residue equal to 
${a_i}^s/\prod_{k\mid k\neq i}\,\left( a_i-a_m\right)$. 
It follows that 
\[
z\mapsto \frac{z^s}{\prod_k\,\left( z-a_k\right)} 
-\sum_{i=1}^n\,\frac{{a_i}^s}{\prod_{k\mid k\neq i}\,\left( a_i-a_m\right)}
\,\frac{1}{z-a_i}
\]
is equal to a polynomial of degree $s-n$ when $s\geq n$ and equal to 
zero when $s<n$. Therefore, if $s<n$ we have that 
\[
\frac{z^s}{\prod_k\,\left( z-a_k\right)} 
=\sum_{i=1}^n\,\frac{{a_i}^s}{\prod_{k\mid k\neq i}\,\left( a_i-a_m\right)}
\,\frac{1}{z-a_i}.
\]
If we compare the expansions in powers of $1/z$ for $z\to\infty$ 
in the left and right hand side, we obtain that 
\begin{equation}
\sum_{i=1}^n\,\frac{{a_i}^s}{\prod_{k\mid k\neq i}\,\left( a_i-a_m\right)} 
=\left\{\begin{array}{ccc}1&\mbox{\rm when}&s=n-1,\\
0&\mbox{\rm when}&0\leq s<n-1 ,\end{array}\right.
\label{partfrac}
\end{equation}

If $\lambda _i$, $1\leq i\leq n$, are arbitrary numbers 
and we write 
\begin{equation}
{x_i} ^2=y_i=\frac{\prod_k\,\left( a_i-\lambda _k\right)}
{\prod_{k\mid k\neq i}\,\left( a_i-a_k\right)} ,
\label{xlambda}
\end{equation} 
then an expansion of the numerators in the right hand side 
of  
\[
\sum_{i=1}^n\,\frac{y_i}{a_i-\lambda _l}=
\sum_{i=1}^n\,\frac{\prod_{k\mid k\neq l}\,\left( a_i-\lambda _k\right)}
{\prod_{k\mid k\neq i}\,\left( a_i-a_k\right)}
\]
in powers of $a_i$ yields in combination with (\ref{partfrac}) that 
\begin{equation}
\sum_{i=1}^n\,\frac{y_i}{a_i-\lambda _l}=1,\quad 1\leq l\leq n. 
\label{yfrac}
\end{equation}
In other words, if the $y_i$ are defined by (\ref{xlambda}), 
then each of the $\lambda _i$'s is a solution of the equation 
$\sum_i\, y_i/\left( a_i-\lambda\right) =1$ for $\lambda$. 
Note that (\ref{xlambda}) defines a polynomial mapping 
$\left(\lambda _1,\,\cdots ,\,\lambda _n\right)\mapsto 
\left( y_1,\,\cdots ,\, y_n\right)$ 
from $\C ^n$ to $\C ^n$. It is a branched covering because 
for every permutation $\pi$ of $\{ 1,\,\cdots ,\, n\}$ the 
vector $\left(\lambda _{\pi (1)},
\,\cdots ,\,\lambda _{\pi (n)}\right)$ 
has the same image as $\left(\lambda _1,\,\cdots ,\,\lambda _n\right)$. 
If the $x_i$ satisfy (\ref{xlambda}), then 
each $\lambda _i$ is a solution of the equation (\ref{confocal}). 
However, (\ref{xlambda}) does not define a mapping 
$\left(\lambda _1,\,\cdots ,\,\lambda _n\right)\mapsto 
\left( x_1,\,\cdots ,\, x_n\right)$ , because the mapping 
$\left( x_1,\,\cdots ,\, x_n\right)\mapsto 
\left( y_1,\,\cdots ,\, y_n\right)$ is a $2^n$\--fold branched covering --- 
for any choice of signs $\epsilon _i$ the vector 
$\left(\epsilon _1\, x_1,\,\ldots ,\,\epsilon _n\, x_n\right)$ 
has the same image as 
$\left( x_1,\,\ldots ,\, x_n\right)$. 
The branching occurs at the set where one of the coordinates $x_i$ is equal to zero, 
which corresponds to 
the condition that one of the coordinates of $\lambda$ 
is equal to $a_i$. 

In our case we have $n=3$ and (\ref{aiIi}). 
The sign changes in the $x_i$ correspond to 
sign changes in the $u_i$. For the equations of motion in the elliptic 
coordinates these will not cause too much trouble, because 
if $R$ is a diagonal matrix with $\pm 1$'s on the diagonal, 
then the transformation $(u,\, v )\mapsto (R\, u,\, R\, v )$ 
will leave the vector field invariant when $\op{det}R=1$ and 
turns the vector field into its opposite (corresponding with a 
time reversal) when $\op{det}R=-1$. 

We will keep $\lambda _1=0$, which implies that the equation 
$\langle u,\, u\rangle =1$ is fulfilled, and regard the remaining two 
elliptic coordinates $\left(\lambda _2,\,\lambda _3\right)$ 
as coordinates on $U_{\C}$. Apart from the problem that the 
velocity field for $u$ is multi\--valued, it will have singularities 
at all the branch loci, corresponding to the condition that 
$\lambda _2$ or $\lambda _3$ attains any of the five values 
$-1/\rho$, $\gamma$, $1/I_1$, $1/I_2$ or $1/I_3$. 

\medskip
The relation (\ref{xlambda}) induces a relation between tangent 
vectors $X$ and $\Lambda$ in the $x$\--space and the $\lambda$\--space, 
respectively.  
If we take the logarithm of the left and right hand side and 
differentiate,then we obtain that 
\begin{equation}
2\frac{X_i}{x_i}=\sum_k\,\frac{-\Lambda _k}{a_i-\lambda _k}.
\label{Xi/xi}
\end{equation}
Squaring this and inserting the formula (\ref{xlambda}) for ${x_i} ^2$, 
we obtain that 
\begin{equation}
4{X_i} ^2=
\sum_{k,\, l\mid k\neq l}\,\frac{\prod_{m\mid m\neq k,\, m\neq l}
\,\left( a_i-\lambda _m\right)}{\prod_{m\mid m\neq i}\left( a_i-a_m\right)} 
\,\Lambda _k\,\Lambda _l
+\sum_k\,\frac{\prod_{m\mid m\neq k}\,\left( a_i-\lambda _m\right)}
{\prod_{m\mid m\neq i}\,\left( a_i-a_m\right)}
\,\frac{{\Lambda _k} ^2}{a_i-\lambda _k}.
\label{Xi2}
\end{equation}
The sum over $i$ of the first sum in the right hand side, 
over the $k$ and $l$ with $k\neq l$, vanishes in view of 
(\ref{partfrac}) with $s<n-1$. On the other hand 
\[
\sum_{i=1}^n\,\frac{\prod_{m\mid m\neq k}\,\left( a_i-\lambda _m\right)}
{\prod_{m\mid m\neq i}\,\left( a_i-a_m\right)}\cdot\frac{1}{a_i-\lambda _k}
=\frac{\prod_{l\mid l\neq k}\,\left(\lambda _l-\lambda _k\right)}
{\prod_l\,\left( a_l-\lambda _k\right)},
\]
because as a function of $\lambda _k$ the left and the right 
hand side both vanish at infinity and have the same poles and residues. 
This leads to Jacobi's conclusion, cf. \cite[26. Vorlesung]{Jacobi}, that 
\begin{eqnarray}
4\sum_{i=1}^n\,{X_i} ^2&=&\sum_k\, 
\prod_{l\mid l\neq k}\,\left(\lambda _l-\lambda _k\right) 
\,\frac{{\Lambda _k} ^2}{S\left(\lambda _k\right)},
\quad\mbox{\rm in which}
\label{deltaxsquared}\\
S\left(\lambda _k\right)&:=&\prod_l\,\left( a_l-\lambda _k\right) . 
\label{S}
\end{eqnarray}

From now on we use that $\lambda _1=0$, $\Lambda _1=0$, and $n=3$.  
Then the quotient of the 
first sum in the right hand side of (\ref{Xi2}) by $a_i$ vanishes 
in view of (\ref{partfrac}) with $s\leq n-3$, where we have used 
that no terms appear with $k=1$ or $l=1$. 
On the other hand, for every $k\neq 1$ we have that  
\[
\sum_{i=1}^3\,\frac{\prod_{m\mid m\neq k,\, m\neq 1}\,\left( a_i-\lambda _m\right)}
{\prod_{m\mid m\neq i}\,\left( a_i-a_m\right)}\cdot\frac{1}{a_i-\lambda _k}
=\, -\frac{\prod_{l\mid l\neq k,\, l\neq 1}\,\left(\lambda _l-\lambda _k\right)}
{\prod_l\,\left( a_l-\lambda _k\right)},
\]
because as a function of $\lambda _k$ the left and the right 
hand side both vanish at infinity and have the same poles and residues. 
With the notation (\ref{S}) we therefore obtain in our case $n=3$ that 
\begin{equation}
4\sum_{i=1}^3\,\frac{{X_i} ^2}{a_i}=
\left(\lambda _2-\lambda _3\right)\,
\left(\frac{{\Lambda _2} ^2}{S\left(\lambda _2\right)} 
-\frac{{\Lambda _3} ^2}{S\left(\lambda _3\right)}\right) .
\label{Xi2ai}
\end{equation}

For the next equation we will use the index $i$ in a cyclic 
manner, by taking $i\in\Z /3\Z$. Then (\ref{Xi/xi}) implies that 
\[
2\left(\frac{a_i\, X_i}{x_i}-\frac{a_{i+1}\, X_{i+1}}{x_{i+1}}\right) 
=\, -\sum_k\,\left(\frac{a_i}{a_i-\lambda _k}-
\frac{a_{i+1}}{a_{i+1}-\lambda _k}\right)\,\Lambda _k
=
\sum_k\,\frac{\left( a_i-a_{i+1}\right)\,\lambda _k\,\Lambda _k}
{\left( a_i-\lambda _k\right)
\,\left( a_{i+1}-\lambda _k\right)} .
\]
On the other hand it follows from (\ref{xlambda}), 
$\lambda _1=0$ and $n=3$ that 
\[
\frac{{x_i} ^2\, {x_{i+1}} ^2}{a_i\, a_{i+1}}
=\frac{\left( a_i-\lambda _2\right)\, \left( a_i-\lambda _3\right) 
\,\left( a_{i+1}-\lambda _2\right)\,
\left( a_{i+1}-\lambda _3\right)}
{\left( a_i-a_{i+1}\right)\,\left( a_i-a_{i+2}\right)\, 
\left( a_{i+1}-a_{i+2}\right)\,\left( a_{i+1}-a_i\right)} .
\]
It follows that 
\begin{eqnarray*}
&& 4\frac{{x_i} ^2\, {x_{i+1}} ^2}{a_i\, a_{i+1}}
\,\left(\frac{a_i\, X_i}{x_i}-\frac{a_{i+1}\, X_{i+1}}{x_{i+1}}\right) ^2\\
&=&\frac{\left( a_i-\lambda _2\right)\, \left( a_i-\lambda _3\right) 
\,\left( a_{i+1}-\lambda _2\right)\,
\left( a_{i+1}-\lambda _3\right)}
{\left( a_{i+2}-a_i\right)\, 
\left( a_{i+1}-a_{i+2}\right)}
\,\left(\sum_k\,\frac{\lambda _k\,\Lambda _k}{\left( a_i-\lambda _k\right)
\,\left( a_{i+1}-\lambda _k\right)}\right) ^2, 
\end{eqnarray*}
which is equal to $1/\left( a_{i+2}-a_i\right)\, 
\left( a_{i+1}-a_{i+2}\right)$ times 
\[
\frac{\left( a_i-\lambda _3\right)\,\left( a_{i+1}-\lambda _3\right)}
{\left( a_i-\lambda _2\right)\,\left( a_{i+1}-\lambda _2\right)}
{\lambda _2} ^2\, {\Lambda _2} ^2
+2\lambda _2\,\Lambda _2\,\lambda _3\,\Lambda _3
+\frac{\left( a_i-\lambda _2\right)\,\left( a_{i+1}-\lambda _2\right)}
{\left( a_i-\lambda _3\right)\,\left( a_{i+1}-\lambda _3\right)}
{\lambda _3} ^2\, {\Lambda _3} ^2.
\]
Now we have 
\[
\frac{1}{\left( a_{i+2}-a_i\right)\, 
\left( a_{i+1}-a_{i+2}\right)}=\frac{a_i-a_{i+1}}{p},
\quad\mbox{\rm with}\quad p:=\prod_{j\in\Z /3\Z}\,\left( a_j-a_{j+1}\right) ,
\]
\[
\frac{\left( a_i-\lambda _3\right)\,\left( a_{i+1}-\lambda _3\right)}
{\left( a_i-\lambda _2\right)\,\left( a_{i+1}-\lambda _2\right)}
=\frac{S\left(\lambda _3\right)}{S\left(\lambda _2\right)}
\,\frac{a_{i+2}-\lambda _2}{a_{i+2}-\lambda _3} 
=\frac{S\left(\lambda _3\right)}{S\left(\lambda _2\right)}
\,\left( 1+\frac{\lambda _3-\lambda _2}{a_{i+2}-\lambda _3}\right) ,
\]
the cyclic sum over $i$ of $a_i-a_{i+1}$ is 
equal to zero, and finally
\[
\sum_{i\in\Z /3\Z}\frac{a_i-a_{i+1}}{a_{i+2}-\lambda _3}
=\frac{1}{S\left(\lambda _3\right)}\sum_{i\in\Z /3\Z}
\,\left( a_i-a_{i+1}\right)\,\left( a_i-\lambda _3\right) 
\,\left( a_{i+1}-\lambda _3\right) =\, -\frac{p}{S\left(\lambda _3\right)}.
\]
It follows that 
\begin{equation}
4\sum_{i\in\Z /3\Z}\,
\frac{{x_i} ^2\, {x_{i+1}} ^2}{a_i\, a_{i+1}}
\,\left(\frac{a_i\, X_i}{x_i}-\frac{a_{i+1}\, X_{i+1}}{x_{i+1}}\right) ^2
=\left(\lambda _2-\lambda _3\right)\, \left(
\frac{{\lambda _2} ^2\, {\Lambda _2} ^2}{S\left(\lambda _2\right)}
-\frac{{\lambda _3} ^2\, {\Lambda _3} ^2}{S\left(\lambda _3\right)}\right) . 
\label{chaplsum}
\end{equation}

\subsection{The Motion in Elliptic Coordinates}
\label{velellsubsec}
In order to obtain the time derivative $\Lambda =\dot{\lambda}$ of $\lambda$ 
corresponding to the velocity vector $X=\dot{x}$ which in turn corresponds to the 
time derivative 
$\dot{u}:=u\times\omega$ of $u$, 
we start with the observation that the horizontality of $j$ implies that 
\begin{equation}
0=\langle  v ,\, u\rangle =\langle I_{\rho ,\, u}\,\omega ,\, u\rangle 
=\langle I\,\omega ,\, u\rangle =\langle\omega ,\, I\, u\rangle ,
 \label{Iomegau}
\end{equation}
cf. (\ref{usigmazero}), (\ref{omegausigma}), and (\ref{Iu}). 
It follows that 
\begin{equation}
\dot{u} \times I\, u=(u\times\omega )\times I\, u
=\langle u,\, I\, u\rangle\,\omega , 
\label{udotIomega}
\end{equation}
which equation  will be used in order to express $\omega$ in terms of $\dot{u}$. 

From (\ref{TI}) we obtain that 
\begin{equation}
2T=\langle I\,\omega ,\,\omega\rangle +\rho\langle\dot{u},\,\dot{u}\rangle .
\label{2Tomegau}
\end{equation}
Furthermore,
\begin{eqnarray*} 
&&\langle I\,\left(\dot{u}\times I\, u\right) ,\, \dot{u}\times I\, u\rangle 
=\sum_{i\in\Z /3\Z}\, I_i\,\left(
\dot{u}_{i+1}\, I_{i+2}\, u_{i+2}-\dot{u}_{i+2}\, I_{i+1}\, u_{i+1}\right) ^2\\
&=&\sum_{i\in\Z /3\Z}\, \left( 
I_{i-1}\,{I_{i+1}} ^2\, {u_{i+1}} ^2+I_{i-2}\, {I_{i-1}} ^2\, {u_{i-1}} ^2\right) \, 
{\dot{u}_i} ^2 -2I_1\, I_2\, I_3\, \dot{u}_{i+1}\, u_{i+2}\,\dot{u}_{i+2}\, u_{i+1}\\
&=&I_1\, I_2\, I_3\,\sum_{i\in\Z /3\Z}\,\left( 
\frac{I_{i+1}}{I_i}\, {u_{i+1}} ^2+\frac{I_{i-1}}{I_i}\, {u_{i-1}} ^2+{u_i} ^2\right) 
\, {\dot{u}_i} ^2=I_1\, I_2\, I_3\,\langle u,\, I\, u\rangle\cdot
\langle I^{-1}\,\dot{u},\,\dot{u}\rangle ,
\end{eqnarray*}
where in the third equation we have used that 
\[
0=\langle\dot{u},\, u\rangle ^2=\sum_{i\in\Z /3\Z}\, {\dot{u}_i} ^2\, {u_i} ^2 
+2\dot{u}_{i+1}\, u_{i+1}\,\dot{u}_{i+2}\, u_{i+2}.
\]
It follows that the first term in the equation for the kinetic energy can be written in the form 
\begin{equation}
\langle I\,\omega ,\,\omega\rangle =I_1\, I_2\, I_3\,\langle I^{-1}\,\dot{u},\,\dot{u}\rangle/
\langle u,\, I\, u\rangle . 
\label{Tudot}
\end{equation}

With the substitutions (\ref{aiIi}), we have 
\begin{equation}
\langle I^{-1}\dot{u},\,\dot{u}\rangle =\sum_{i=1}^3\, {\dot{x}_i} ^2
\label{I-1dotu}
\end{equation}
and 
\begin{equation}
\langle \dot{u},\,\dot{u}\rangle 
=\langle I\, \dot{x},\,\dot{x}\rangle 
=\sum_{i=1}^3\,{\dot{x}_i} ^2/a_i,
\label{dotudotu}
\end{equation}
which by means of  (\ref{deltaxsquared}) and (\ref{Xi2ai}), respectively, 
can be expressed in terms 
of the velocities of the elliptic coordinates. 

In order to express the denominator 
\begin{equation}
\langle u,\, I\, u\rangle =\langle x,\, I^2\, x\rangle 
=\sum_{i=1}^3\,{x_i} ^2/{a_i} ^2
\label{uIu}
\end{equation}
in (\ref{Tudot}) 
in terms of  the elliptic coordinates, we observe that 
\begin{equation}
1=\langle u,\, u\rangle =\sum_{i=1}^3\, {x_i} ^2/a_i
\label{xi2/ai}
\end{equation}
implies that the third degree polynomial equation, which is obtained 
from (\ref{confocal}) by multiplication with 
$S(\lambda )=\prod_{i=1}^3\,\left( a_i-\lambda\right)$, 
has $\lambda =\lambda _1=0$ as a solution. The second order equation for 
the two remaining solutions $\lambda _2$, $\lambda _3$ takes the form 
\[
\lambda ^2+\left[\sum_{i=1}^3\,\left( {x_i} ^2-a_i\right)\right]\,\lambda 
-\sum_{i\in\Z /3\Z}\,  {x_i} ^2\,\left( a_{i+1}+a_{i+2}\right)  
+\sum_{h\in\Z /3\Z}\, a_{h}\, a_{h+1}=0,
\]
in which the constant term can be simplified to 
\[
\sum_{i\in\Z /3\Z}\,  {x_i} ^2\,\left[  
-a_{i+1}-a_{i+2}+\sum_{h\in\Z /3\Z}\,\frac{a_{h}\, a_{h+1}}{a_i}\right] 
= \sum_{i\in\Z /3\Z}\,  {x_i} ^2\,\frac{a_{i+1}\, a_{i+2}}{a_i}=a_1\, a_2\, a_3\,
\sum_{i\in\Z /3\Z}\, \frac{{x_i} ^2}{{a_i} ^2}.
\]
It follows that 
\begin{equation}
\lambda _2+\lambda _3=\sum_{i=1}^3\,\left( a_i-{x_i} ^2\right) 
\label{sumlambda}
\end{equation}
and 
\begin{equation}
\lambda _2\,\lambda _3=a_1\, a_2\, a_3\,
\sum_{i=1}^3\, \frac{{x_i} ^2}{{a_i} ^2}.
\label{prodlambda}
\end{equation}
Combining (\ref{Tudot}), (\ref{I-1dotu}, (\ref{deltaxsquared}) for $X=\dot{x}$, 
(\ref{uIu}) and (\ref{prodlambda}), we conclude that  
\begin{equation}
4\langle I\,\omega ,\,\omega\rangle =\frac{\lambda _2-\lambda _3}
{\lambda _2\,\lambda _3}\, 
\left(\frac{\lambda _2\, {\dot{\lambda}_2} ^2}{S\left(\lambda _2\right)}
-\frac{\lambda _3\, {\dot{\lambda}_3} ^2}
{S\left(\lambda _3\right)}\right) .
\label{omegaIomega}
\end{equation}

Combining (\ref{2Tomegau}) and (\ref{omegaIomega}) 
with (\ref{dotudotu}) and (\ref{Xi2ai}) for $X=\dot{x}$, we arrive at 
\begin{equation}
8T=\frac{\lambda _2-\lambda _3}{\lambda _2\,\lambda _3}\,\left(
\frac{\lambda _2\, {\dot{\lambda}_2} ^2}
{S\left(\lambda _2\right)}
-\frac{\lambda _3{\dot{\lambda}_3} ^2}
{S\left(\lambda _3\right)}\right) +\rho\,\left(\lambda _2-\lambda _3\right)\,\left(
\frac{{\dot{\lambda}_2} ^2}{S\left(\lambda _2\right)}
-\frac{{\dot{\lambda}_2} ^2}{S\left(\lambda _2\right)}\right) .
\label{Tlambdadot}
\end{equation} 

In order to obtain the second equation for the two unknowns 
$\dot{\lambda}_2$, $\dot{\lambda}_3$, we start with the equation 
\begin{eqnarray*}
\|  j\| ^2-\rho\, 2T&=&\langle I_{\rho ,\, u}\,\omega ,\, I_{\rho ,\, u}\,\omega\rangle 
-\rho\,\langle I_{\rho ,\, u}\,\omega ,\,\omega\rangle \\
&=&\langle I\,\omega +\rho\,\omega -\rho\,\langle\omega ,\, u\rangle\, u,\, 
I\,\omega -\rho\,\langle\omega ,\, u\rangle\, u\rangle 
=\langle I\,\omega ,\, I\,\omega\rangle +\rho\,\langle\omega ,\, I\,\omega\rangle .
\end{eqnarray*}
Here we have used (\ref{SO3eq}) and (\ref{TI}) in the first equation, 
(\ref{Iu}) in the second equation, and (\ref{Iomegau}) in the third one. 
The first term in the right hand side is equal to 
\begin{eqnarray*}
&&\langle I\,\omega ,\, I\,\omega\rangle =\langle u,\, I\, u\rangle ^{-2}\cdot 
\langle I\, (\dot{u}\times I\, u),\, I\, (\dot{u}\times I\, u)\rangle \\
&=&\langle u,\, I\, u\rangle ^{-2}\,\sum_{i\in\Z /3\Z}\, {a_i}^{-2}\,\left( 
{a_{i+1}}^{-1/2}\,\dot{x}_{i+1}\, {a_{i+2}}^{-3/2}\, x_{i+2}
-{a_{i+2}}^{-1/2}\,\dot{x}_{i+2}\, {a_{i+1}}^{-3/2}\, x_{i+1}\right)\\
&=&\langle u,\, I\, u\rangle ^{-2}\, 
\left( a_1\, a_2\, a_3\right) ^{-2}\,
\sum_{i\in\Z /3\Z}\, \frac{{x_{i+1}} ^2\, {x_{i+2}} ^2}
{a_{i+1}\, a_{i+2}}\left( 
\frac{a_{i+1}\,\dot{x}_{i+1}}{x_{i+1}}-\frac{a_{i+2}\,\dot{x}_{i+2}}{x_{i+2}}
\right) ^2 ,
\end{eqnarray*}
where we have used (\ref{udotIomega}) in the first equation and 
(\ref{aiIi}) in the second one. 
Combining this with (\ref{uIu}),  (\ref{prodlambda}), 
and (\ref{chaplsum}) for $X=\dot{x}$, we obtain that 
\begin{equation}
4\langle I\,\omega ,\, I\,\omega\rangle =
\frac{\lambda _2-\lambda _3}{{\lambda _2} ^2\, {\lambda _3} ^2}
\,\left(\frac{{\lambda _2} ^2\, {\dot{\lambda}_2} ^2}{S\left(\lambda _2\right)}
-\frac{{\lambda _3} ^2\, {\dot{\lambda}_3} ^2}{S\left(\lambda _3\right)}\right) . \label{4IomIom}
\end{equation}
Also using (\ref{Tlambdadot}), we therefore arrive at 
\begin{equation}
4\| j\| ^2-8\rho\, T=\frac{\lambda _2-\lambda _3}{{\lambda _2} ^2\, {\lambda _3} ^2}
\,\left(\frac{{\lambda _2} ^2\, {\dot{\lambda}_2} ^2}{S\left(\lambda _2\right)}
-\frac{{\lambda _3} ^2\, {\dot{\lambda}_3} ^2}{S\left(\lambda _3\right)}\right) 
+\rho\,\frac{\lambda _2-\lambda _3}{\lambda _2\, \lambda _3}
\,\left(\frac{\lambda _2\,{\dot{\lambda}_2} ^2}{S\left(\lambda _2\right)}
-\frac{\lambda _3\,{\dot{\lambda}_3} ^2}{S\left(\lambda _3\right)}\right) .
\label{j2-2rhoT}
\end{equation}

For the 
unknowns 
\[
\xi _2:=\left(\lambda _2-\lambda_3\right)\,\left(\frac{1}{\lambda _3}+\rho\right) 
\,\frac{{\dot{\lambda}_2} ^2}{S\left(\lambda _2\right)}, 
\quad
\xi _3:=\left(\lambda _2-\lambda_3\right)\,\left(\frac{1}{\lambda _2}+\rho\right) 
\,\frac{{\dot{\lambda}_3} ^2}{S\left(\lambda _3\right)}
\]
the equations (\ref{Tlambdadot}) and (\ref{j2-2rhoT}) take the form 
\[
8T=\xi _2-\xi _3\quad\mbox{\rm and}\quad 
4\| j\| ^2-8\rho\, T=\frac{1}{\lambda _3}\,\xi _2-\frac{1}{\lambda _2}\,\xi _3,
\]
respectively. Substracting $\lambda _2$ times the second equation from the first one, 
we obtain that 
\begin{equation}
\left( 1-\frac{\lambda _2}{\lambda _3}\right) \, 
\left(\lambda _2-\lambda_3\right)\,\left(\frac{1}{\lambda _3}+\rho\right) 
\,\frac{{\dot{\lambda}_2} ^2}{S\left(\lambda _2\right)}
=\left( 4\| j\| ^2-8\rho\, T\right)\, \left(\gamma -\lambda _2\right) ,
\label{lambda2dot}
\end{equation}
where we have used (\ref{gamma}) in the second equation. 

At this stage we recall the change of the time parametrization defined by 
\[
\dd{\tau}{t}=X (u)^{-1/2},
\]
cf. Corollary \ref{commvvcor} and (\ref{phidef}). In order to express 
$X (u)$ in terms of the elliptic coordinates, we recall that 
\[
\rho\, X (u) = 1-\sum_{i=1}^3\,\frac{{x_i} ^2}{a_i-\lambda}\quad 
\mbox{\rm with}\quad \lambda =\, -1/\rho .
\]
Now 
\begin{equation}
S(\lambda )\,\left( 1-\sum_{i=1}^3\,\frac{{x_i} ^2}{a_i-\lambda}\right) 
=\, -\lambda\,\left(\lambda _2-\lambda\right) 
\,\left(\lambda _3-\lambda\right) ,
\label{Squadr}
\end{equation}
because as a polynomial in $\lambda$ the left and the right hand side have the 
same zeros $\lambda _1=0$, $\lambda _2$ and $\lambda _3$ and have the same 
leading coefficient. It follows that 
\begin{equation}
 X (u)=\rho ^{-2}\,\left(\lambda _2+1/\rho\right)\,
\left(\lambda _3+1/\rho\right) /\prod_{i=1}^3\,\left( a_i+1/\rho\right) . 
\label{philambda}
\end{equation}
Substituting ${\dot{\lambda}_2} ^2
=\left(\dd{\lambda _2}{\tau}\right) ^2/ X (u)$ and 
(\ref{philambda}) in (\ref{lambda2dot}), we arrive at 
\begin{equation}
\left(\dd{\lambda _2}{\tau}\right) ^2/P\left(\lambda _2\right) =\left(\frac{-c\,\lambda _3}
{\lambda _2-\lambda _3}\right) ^2,
\label{dlambda2dtau}
\end{equation}
in which the polynomial $P$ is defined by 
\begin{equation}
P(\lambda ):=\left( -1/\rho -\lambda\right)\, 
\left(\gamma -\lambda\right)\,\prod_{i=1}^3\,\left( a_i-\lambda\right)  
\label{Pdef}
\end{equation}
and the constant $c$ is determined by 
\begin{equation}
c^2:=\left( 4\| j\| ^2-8\rho\, T\right) /\prod_{i=1}^3\,
\left(\rho\,  a_i+1\right) .
\label{c2def}
\end{equation}
In a similar way we obtain the equation 
\begin{equation}
\left(\dd{\lambda _3}{\tau}\right) ^2/P\left(\lambda _3\right) =\left(\frac{-c\,\lambda _2}
{\lambda _2-\lambda _3}\right) ^2
\label{dlambda3dtau}
\end{equation}
for $\op{d}\!\lambda _3/\op{d}\!\tau$. 

At this point it becomes appropriate to introduce the {\em hyperelliptic curve} 
$C$ defined by the 
polynomial $P$ as the one point completion (compactification) of 
the affine curve 
\begin{equation}
\{ (\lambda ,\, z)\in\C ^2\mid z^2=P(\lambda )\} . 
\label{Cdef}
\end{equation}
The genus of the hyperelliptic curve $C$ is equal to $g$ if the degree of the polynomial 
$P$ is equal to $2g+1$ or $2g+2$, cf. Shafarevich \cite[Ch. III, \S 5]{Shaf}, 
or Farkas and Kra \cite[III.7.4]{FK}. Because in our case 
the degree of $P$ is equal to five, we have $g=2$.
 
Consider the motion on $C\times C$ which is determined by the choice 
\begin{equation}
\frac{1}{z_2}\,\dd{\lambda _2}{\tau} 
=\frac{-c\,\lambda _3}{\lambda _2-\lambda _3}
\quad\mbox{\rm and}\quad 
\frac{1}{z_3}\,\dd{\lambda _3}{\tau} 
=\frac{-c\,\lambda _2}{\lambda _2-\lambda _3}
\label{CxCmotion}
\end{equation}
of the square roots in (\ref{dlambda2dtau}) and (\ref{dlambda3dtau}), respectively, where 
$\left( \lambda _2,\, z_2\right)$ and $\left( \lambda _3,\, z_3\right)$ denote 
the points in $C$ which, under the two\--fold branched covering 
$C\ni (\lambda ,\, z)\mapsto\lambda$, project to $\lambda _2$ and $\lambda _3$, 
respectively. 

The equations (\ref{CxCmotion}) imply that 
\begin{equation}
\frac{\lambda _2}{z_2}\,\dd{\lambda _2}{\tau}
-\frac{\lambda _3}{z_3}\,\dd{\lambda _3}{\tau}=0,\quad
\frac{1}{z_2}\,\dd{\lambda _2}{\tau}
-\frac{1}{z_3}\,\dd{\lambda _3}{\tau}=c,
\label{formCxC}
\end{equation}
or
\begin{equation}
\int_{\left(\lambda _3(\tau ),\, z_3(\tau)\right)}
^{\left(\lambda _2(\tau ),\, z_2(\tau )\right)}
\,\frac{\lambda}{z}\,\op{d}\! {\lambda}
=\mbox{\rm constant},\quad
\int_{\left(\lambda _3(\tau ),\, z_3(\tau)\right)}
^{\left(\lambda _2(\tau ),\, z_2(\tau )\right)}\,
z^{-1}\,\op{d}\! {\lambda}
=c\,\tau +\mbox{\rm constant}.
\label{formCxCint}
\end{equation}
Here the integration is over a curve in $C$, running from 
$\left(\lambda _3(\tau ),\, z_3(\tau )\right)$ to $\left(\lambda _2(\tau ),\, z_2(\tau )\right)$, 
and depending smoothly on $\tau$. The primitives of the differential forms 
$z^{-1}\,\op{d}\!\lambda$ and $\frac{\lambda}{z}\,\op{d}\!\lambda$ 
are the {\em hyperelliptic integrals} of the hyperelliptic curve $C$, 
and for this reason one says that (\ref{formCxCint}) implies that 
{\em the problem is solved by quadratures in terms of the hyperelliptic 
integrals} corresponding to the hyperelliptic curve $C$. 
The quantities $\lambda _2-\lambda _3$ and $\lambda _2\,\lambda _3$ 
are given by {\em Jacobi's theta functions} of the integrals 
in the left hand sides of (\ref{formCxCint}), 
cf. Shafarevich \cite[p. 419]{Shaf}, 
and therefore one also talks about {\em solving by quadratures 
in terms of theta functions}. 

\subsection{The Jacobi Variety of the Hyperelliptic Curve}
\label{jacsubsec}
The equations in (\ref{formCxCint}) lead to a 
beautiful interpretation of the system in terms of 
the Jacobi variety $\op{Jac}(C)$ 
of the hyperelliptic curve $C$. 
The differential forms 
\begin{equation}
\beta _i:=\frac{\lambda ^i}{z}\,\op{d}\!\lambda ,\quad i=0,\, 1,
\label{betai}
\end{equation}
extend to holomorphic differential forms on $C$ and actually form a basis of the 
two\--dimensional complex vector space ${\cal H}^1(C)$ of all holomorphic differential 
forms of degree one on $C$, cf. Shafarevich \cite[Ch. III, \S 5]{Shaf} 
or Farkas and Kra \cite[III.7.5]{FK}. 

\begin{remark}
It is clear that the differential 
forms 
\begin{equation}
\beta _i':=
\frac{{\lambda _2}^i}{z_2}\,\op{d}\!\lambda _2
-\frac{{\lambda _3}^i}{z_3}\,\op{d}\!\lambda _3,\quad i=0,\, 1
\label{betai'}
\end{equation}
of degree one on $C\times C$ are holomorphic, and closed because of the 
separation of variables. In this notation, the equations in (\ref{formCxC}) 
mean that the vector field 
\begin{equation}
\xi :=\left(\dd{\lambda _2}{\tau},\,
\dd{\lambda _3}{\tau}\right) 
\label{XCxC}
\end{equation}
on $C\times C$ satisfies 
\begin{equation}
\op{i}_{\xi}\,\beta _1'=0,\quad\op{i}_{\xi}\,\beta _0'=c,
\label{iXbeta}
\end{equation}
if $\op{i}_{\xi}\beta$ denotes the inner product of the differential form 
$\beta$ with the vector field $\xi$. 
This is a complex version of Proposition \ref{commvvprop}, where we had  
the differential forms $\beta$ and $\gamma$ instead of 
$\beta _1'$, $\beta _0'$, and $\op{i}_{\xi}\beta =0$ and 
$\op{i}_{\xi}\gamma =\, -1$. 
\label{iXbeta'rem}
\end{remark}

For each smooth curve $\gamma$ in $C$, we have the complex linear form 
\begin{equation}
{\cal H}^1(C)\ni\beta\mapsto\int_{\gamma}\,\beta 
\label{[gamma]}
\end{equation}
on ${\cal H}^1(C)$ --- this defines an element $\int_{\gamma}$ of the 
dual space ${\cal H}^1(C)^*$ of ${\cal H}^1(C)$. The integral depends 
only on the homology class $[\gamma ]\in\op{H}_1(C,\,\Z )$ 
of  $\gamma$. If  we restrict ourselves to 
closed loops $\gamma$, then this leads to a homomorphism $\int$ from the 
$\op{H}_1(C,\,\Z )$ to ${\cal H}^1(C)^*$. The image 
\begin{equation}
\Lambda (C) :=\int\left(\op{H}_1(C,\,\Z )\right)\subset {\cal H}^1(C)^*
\label{LambdaC}
\end{equation}
is an additive subgroup of ${\cal H}^1(C)^*$, which is called the 
{\em period lattice} of $C$. 
For any Riemann surface (complete algebraic curves over $\C$) of genus $g$ 
we have that the complex dimension of  ${\cal H}^1(C)$ is equal to $g$
and therefore the real dimension of ${\cal H}^1(C)^*$ is equal to $2g$. 
Furthermore, 
the period lattice $\Lambda (C)$ has a $\Z$\--basis which at the same  
time is an $\R$\--basis of  ${\cal H}^1(C)^*$, cf. Farkas and Kra 
\cite[III.2.8]{FK}. 
Therefore the quotient space 
\begin{equation}
\op{Jac}(C):={\cal H}^1(C)^*/\Lambda (C),
\label{JC}
\end{equation}
is compact, a {\em torus of 
real dimension equal to $2g$}. Definition 
(\ref{JC}) is the analytic definition of 
the {\em Jacobi variety of the curve $C$}, 
cf. \cite[p. 87]{FK} or \cite[p. 143]{TthI}, whereas 
the algebraic definition is formulated in terms of divisors, 
cf. \cite[p. 155]{Shaf} or cite[p. 3.28]{TthII}. 

If we fix $p,\, q\in C$, then the difference of the homology classes of two 
curves in $C$ which run from $p$ to $q$ is an element of $\op{H}_1(M,\,\Z )$, 
and it follows that, for each $\beta\in {\cal H}^1(C)$,  the element 
\begin{equation}
\int_p^q:=(\int_{\gamma})+\Lambda (C)\in \op{Jac}(C)
\label{Adef}
\end{equation}
does not depend on the choice of the curve $\gamma$ from $p$ to $q$. 
This defines a smooth mapping $\int :(p,\, q)\mapsto\int_p^q$ from $C\times C$ to $\op{Jac}(C)$. 
Note that $\int_q^p=\, -\int_p^q$ for every $(p,\, q)\in C\times C$, and $\int_p^p=0$. 
If $\iota$ denotes the involution $(\lambda ,\, z)\mapsto (\lambda ,\, -z)$ in $C$, then 
$\iota ^*\beta =\, -\beta$ for every $\beta\in{\cal H}^1(C)$, and it follows that 
$\int$ maps $(p,\, q)$ and $(\iota (q),\,\iota (p))$ to the same element of $\op{Jac}(C)$. 
Therefore the mapping $\int$ can also be viewed as a mapping to 
$\op{Jac}(C)$ from 
the quotient $C\times_{\iota '}C$ of $C\times C$ by the involution 
$\iota ':(p,\, q)\mapsto(\iota (q),\,\iota (p))$ in $C\times C$. 
The diagonal $\{ (p,\, q)\in C\times C\mid p=q\}$  
in $C\times C$, which is isomorphic to $C$, is mapped by the quotient map to 
a curve $D$ in $C\times_{\iota '}C$ which is isomorphic to 
$C/\iota\simeq {\bf  P}_1(\C)$. 
Because the diagonal is mapped to the origin, we have that $\int (D)=\{ 0\}$ as well. 

The mapping  $A$ is closely related to the {\em Abel\--Jacobi map} 
$A:C\times C\to \op{Jac}(C)$, which is defined as follows. 
Choose a base point $p_0\in C$. 
Then $A(p,\, q)$ is equal to the linear form on ${\cal H}^1(C)$ modulo $\Lambda (C)$ 
which is defined by 
\[
A(p,\, q)(\omega )=\int_{p_0}^p\,\omega+\int_{p_0}^q\,\omega .
\]
Because $A(p,\, q)=A(q,\, p)$, $A$ can be viewed as a mapping to 
$\op{Jac}(C)$ from 
the {\em symmetric power} $C^{(2)}$ of $C$, which is defined as the quotient 
of $C\times C$ by the involution $(p,\, q)\mapsto (q,\, p)$. 
It is a classical theorem of Abel and Jacobi that for genus two curves $C$ the 
Abel\--Jacobi map is surjective, cf. Farkas and Kra 
\cite[III.6.6]{FK}. More precisely, it maps a 
genus zero curve $D_0$ in $C^{(2)}$ to a point $j_0$ in $\op{Jac}(C)$ and 
is a diffeomorphism from $C^{(2)}\setminus D_0$ onto 
$\op{Jac}(C)\setminus\{ j_0\}$, 
cf. Farkas and Kra \cite[III.11.8 and III.11.11]{FK}. 
This is an example of a {\em blowing down}, also called a 
{\em sigma\--process}, 
cf. Shafarevich \cite[Ch. II, \S 4 and Ch. IV, \S 3]{Shaf}. 
Now, if $\iota\left( p_0\right) =p_0$, then 
\[
A(p,\, q)=\int_{p_0}^q+\int_{p_0}^p=\int_{p_0}^q-\int_{\iota\left( p_0\right)}
^{\iota (p)}=\int_{p_0}^q-\int_{p_0}^{\iota (p)}=\int_{\iota (p)}^q.
\]
The condition that $\iota\left( p_0\right) =p_0$ means that $p_0$ 
corresponds to one of the five zeros of $P$ or to the point on $C$ at infinity.
It follows that the mapping  
$\int$ also maps a genus zero curve to a point and is a diffeomorphism from the 
complement of the curve to the complement in $\op{Jac}(C)$ of the point. Because 
$\int (D)=\{ 0\}$, the curve must be equal to $D$ and $\int$ defines 
a diffeomorphism from $C\times _{\iota '}C\setminus D$ onto 
$\op{Jac}(C)\setminus\{ 0\}$. 
 
In these terms the equations in (\ref{formCxCint})  
imply that the tangent map of  $\int$ maps the vector field $\xi$ in 
$C\times C$ to a {\em constant} vector field on $\op{Jac}(C)$. 
More precisely, if we identify ${\cal H}^1(C)$ with $\C ^2$ by means 
of the basis $\beta _0$, $\beta _1$ of  (\ref{betai}), then 
(\ref{formCxCint}) yields that 
$\xi$ corresponds to the constant vector field with coordinates $(c,\, 0)$. 
We therefore have verified that the rotational motion of 
Chaplygin's sphere with horizontal moment is 
{\em algebraically integrable} according to the 
definition of Adler and van Moerbeke \cite{AM}, and in view of Subsection 
\ref{redhorsec} this result remains true if we only assume that 
the moment is not vertical. This verification is very 
different from the one via proposition \ref{hatMcprop}.  

Because every holomorphic vector field on $\op{Jac}(C)$ lifts to a bounded 
holomorphic vector field on ${\cal H}^1(C)^*$, and every bounded 
holomorphic function on $\C ^2$ is equal to a constant, we have that 
every holomorphic vector field on $\op{Jac}(C)$ is constant. Therefore 
our conclusion that $\xi$ is mapped to a constant vector field on $\op{Jac}(C)$ 
is equivalent to the conclusion that 
$\xi$ is mapped to a holomorphic vector field 
on $\op{Jac}(C)$. According to (\ref{CxCmotion}), the vector field $\xi$ on 
$C\times C$ is rational with poles along the diagonal. The blowing down of the 
diagonal by the map $\int$ to a point (the origin) in $\op{Jac}(C)$ apparently 
has the effect that it regularizes the vector field $\xi$. 

\begin{remark}
If the rotational motion after time parametrization as in Corollary 
\ref{commvvcor} is truly quasi\--periodic in the sense that the orbit on the $(j,\, T)$\--level set is dense, 
then each continuous vector field which commutes with $\xi$ is equal to 
$a\,\xi +b\, \eta$, in which $a=a(j,\, T)$ and $b=b(j,\, T)$ are constants which 
may depend on $j$ and $T$. Because the levels $(j,\, T)$, for which the 
$\xi$\--orbits are dense in the level set, form a dense subset of the 
set of all $(j,\, T)$, 
it follows that if a continuous vector 
field is defined on all the regular level surfaces and also depends continuously on $(j,\, T)$, then it is of the form $a\,\xi +b\,\eta$ with 
$a=a(j,\, T)$ and $b=b(j,\, T)$ depending continuously on $(j,\, T)$. 
If $b_0$ and $b_1$ are arbitrary constants, then 
the vector field $\zeta$, which in the elliptic 
coordinates corresponds to 
\begin{equation}
\op{i}_{\zeta}\,\beta _1'=b_1,\quad\op{i}_{\zeta}\,\beta _0'=b_0,
\label{Z}
\end{equation}
compare (\ref{iXbeta}), or on $\op{Jac}(C)$ corresponds to the 
constant vector field with coordinates $\left( b_0,\, b_1\right)$, 
commutes with $\xi$ and depends continuously on $(j,\, T)$. 
The conclusion is that there is a bijective linear correspondence 
between the vector fields which commute with $\xi$ 
and have the aforementioned continuity properties, and the 
constant vector fields on $\op{Jac}(C)$.  

In order to determine the coefficients $b_0$ and $b_1$ for 
which $\zeta =\eta$, where $\eta$ is as in Proposition \ref{commvvprop}, 
we start with the time derivative $\dot{x}$ which 
corresponds to
\[
\dot{u}:=u\times (I+\rho )\,\omega =u\times v ,
\]
where in the second identity we have used (\ref{omegausigma}), 
(\ref{Iu}), and $u\times u=0$. It follows that 
\[
\| j\| ^2=\langle v ,\, v\rangle 
=\langle u\times v ,\,\ u\times v\rangle 
=\langle\dot{u},\,\dot{u}\rangle 
=\sum_{i=1}^3\,\frac{{\dot{x}_i}^2}{a_i}.
\]
Here we have used (\ref{SO3eq}) in the first equation, 
$\langle u,\, v\rangle =-j_3=0$ in the second one 
and $u_i={I_i}^{1/2}\, x_i={a_i}^{-1/2}\, x_i$ in the last one, 
cf. (\ref{aiIi}). In view of (\ref{Xi2ai}), 
we obtain the relation 
\begin{equation}
4\| j\| ^2=\left(\lambda _2-\lambda _3\right)\,
\left(\frac{{\dot{\lambda}_2}^2}{S\left(\lambda _2\right)}
-\frac{{\dot{\lambda}_3}^2}{S\left(\lambda _3\right)}\right) 
\label{Yj2}
\end{equation}
for the corresponding time derivative $\Lambda =\dot{\lambda}$ 
in elliptic coordinates. We have that 
$\dd{\lambda}{\tau}= X (u)^{1/2}\,\dot{\lambda}$, 
with $ X (u)$ given by (\ref{philambda}), 
satisfies (\ref{Z}) if and only if 
\[
\frac{{\lambda _2}^i}{z_2}\,\dd{\lambda _2}{\tau}
-\frac{{\lambda _3}^i}{z_3}\,\dd{\lambda _3}{\tau}=b_i,\quad i=0,\, 1,
\]
or, equivalently,  
\[
\dd{\lambda _2}{\tau}/z_2
=\frac{b_1-\lambda _3\, b_0}{\lambda _2-\lambda _3},\quad 
\dd{\lambda _3}{\tau}/z_3
=\frac{b_1-\lambda _2\, b_0}{\lambda _2-\lambda _3}.
\]
Squaring each of the expressions and 
using (\ref{Cdef}), (\ref{Pdef}), we can bring the equation (\ref{Yj2})  
into the form 
\[
\frac{4\| j\| ^2}{\rho ^2\,\prod_i\left( a_i+1/\rho\right)}=
\frac{\left(\lambda _2+1/\rho\right)\, 
\left(\lambda _2-\gamma\right)\left( b_1-\lambda _3\, b_0\right) ^2
-\left(\lambda _3+1/\rho\right)\, 
\left(\lambda _3-\gamma\right)\left( b_1-\lambda _2\, b_0\right) ^2}
{\left(\lambda _2+1/\rho\right)\,
\left(\lambda _2+1/\rho\right)\,
\left(\lambda _2-\lambda _3\right)}. 
\]
The numerator in the right hand side has to be equal to zero when 
$\lambda _2=\, -1/\rho$, which implies that 
\begin{equation}
b_1=\, -b_0/\rho .
\label{b1b0}
\end{equation}
Substitution of (\ref{b1b0}) in the previous formula yields that 
\[
\frac{4\| j\| ^2}{\rho ^2\,\prod_{i=1}^3\,\left( a_i+1/\rho\right)}={b_0}^2
\,\frac{\left(\lambda _2-\gamma\right)\,\left(\lambda _2+1/\rho\right)
-\left(\lambda _2-\gamma\right)\,\left(\lambda _3+1/\rho\right)}
{\lambda _2-\lambda _3}
={b_0}^2\,\left( 1/\rho +\gamma\right) ,
\]
or 
\begin{equation}
{b_0}^2=4\| j\| ^2/\rho ^2\,\left( 1/\rho +\gamma\right)\,
\prod_{i=1}^3\,\left( a_i+1/\rho\right) , 
\label{b0}
\end{equation}
a formula similar to (\ref{c2def}). The formulas 
(\ref{b0}) and (\ref{b1b0}) determine the constant vector field
$\left( b_0,\, b_1\right)$ on $\op{Jac}(C)$ corresponding to $\eta$, up to its sign. 
\label{Yjacrem}
\end{remark}

Until now we have not paid much attention to the fact that the 
substitutions which we have 
used are not bijective but, except for the mapping 
$\int: C\times _{\iota '}C\to \op{Jac}(C)$, are branched coverings. 
Recall that in order to obtain a single\--valued 
vector field $\xi$, we passed from the manifold $L$ of solutions 
$(u,\, v )$ of the equations (\ref{SO3eq}) and (\ref{phipsichi}) 
to its two\--fold covering $M$, defined by (\ref{Mdef}).  

The projection $(u,\, v )\mapsto u$ exhibits $L$  
as a fourfold branched covering over the quadric 
\[
U_{\C}:=\{ u\in\C ^3\mid \langle u,\, u\rangle =1\}
\]
in $\C ^3$, with branch locus at $ X (u)\,\langle u,\, K\, u\rangle =0$, 
cf. Subsection \ref{sphereprojsubsec} and (\ref{sigmabranch}). 
At the generic points of $U_{\C}$ we have four possibilities 
for $\dd{u}{\tau}$, of the form $\pm \xi '$, $\pm \xi ''$, where 
$\xi '$ and $\xi ''$ become equal at the branch locus. 

The branched covering $L\to U_{\C}$ is less than fourfold 
over the set $D_{\infty}$ of the limit points 
of the $u\in U_{\C}$ for which at least one of the 
solutions $ v$ of the equations (\ref{SO3eq}), (\ref{phipsichi}) 
runs off to infinity. 
This is the set of $u\in U_{\C}$ for which 
there exists a nonzero solution $ v$ of the 
homogeneous equations
\begin{equation}
\langle u,\, v\rangle =0,\quad\langle v ,\, v\rangle =0,
\quad\langle u,\, J\, v\rangle ^2
+ X (u)\,\langle v ,\, J\, v\rangle =0.
\label{Dinfty}
\end{equation}
This exhibits $D_{\infty}$ as the resultant set of the 
three polynomials in $ v$ which appear in (\ref{Dinfty}), 
and therefore $D_{\infty}$ is a closed algebraic curve in $U_{\C}$. 

The set $D_{\infty}$ contains the set $D$ of the points 
which correspond to points on 
the diagonal $\lambda _2=\lambda _3$ in the elliptic coordinates, 
because according to (\ref{CxCmotion}) the vector field  
in the $\left(\lambda _2,\,\lambda _3\right)$\--space is 
infinite along the diagonal. 
At the points of $D$ 
both vectors $v$ and or $w$ become infinite, which means that 
all solutions $ v$ of the equations (\ref{SO3eq}), (\ref{phipsichi}) 
run off to infinity when $u$ approaches a point in $D$. 
This implies that the image of $L$ under the projection 
$(u,\, v )\mapsto u$ is contained in $U_{\C}\setminus D$. 
According to (\ref{sumlambda}) and (\ref{prodlambda}), 
the point $u\in U_{\C}$ corresponds to $\lambda _2=\lambda _3$ in the elliptic coordinates, if and only if 
\begin{equation}
\Delta (u):=\left(\langle I^{-1}\, u,\, u\rangle -\op{trace}I^{-1}\right) ^2 
-4\langle I\, u,\, u\rangle/\op{det}I =0. 
\label{diageq}
\end{equation}

\begin{lemma}
We have that $D_{\infty}=D$, which implies that 
the projection $L\to U_{\C}$ 
is a fourfold branched covering from $L$ onto 
$U_{\C}\setminus D$.
\label{Dlem}
\end{lemma} 
\begin{proof}~~ (Sketch.) 
For $\langle u,\, u\rangle =1$ and ${u_1}^2+{u_2}^2\neq 0$ 
the nonzero multiples of the $ v\in\C ^3$ such that 
\[
 v _1=\, -u_1\, u_3\pm\sqrt{-1}\, u_2,\quad
 v _2=\, -u_3\, u_2\mp\sqrt{-1}\, u_1,\quad
 v _3={u_2}^2+{u_1}^2
\]
parametrize the set of nonzero solutions $ v$ of the 
equations $\langle u,\, v\rangle =0$ and 
$\langle v ,\, v\rangle =0$. Inserting this $ v$ 
in the last equation in (\ref{Dinfty}), where we use that 
$J=(I=\rho )^{-1}$, cf. (\ref{Jdef}), we obtain the 
polynomial equation 
\begin{eqnarray*}
&&\left[\left( I_1-I_3\right)\, I_2\, {u_1}^2
+\left( I_2-I_3\right)\, I_1\, {u_2}^2\right]\,\left[ {u_1}^2+{u_2}^2\right] 
+\left(I _2-I_1\right)\, I_3\,\left( {u_1}^2-{u_2}^2\right) \\
&&=2\pm\sqrt{-1}\,\left( I_2-I_1\right)\, I_3\, u_1\, u_2\, u_3.
\end{eqnarray*}
Squaring both sides we obtain a polynomial equation of the 
form 
\[
\Delta (u)\,\left( {u_1}^2+{u_2}^2\right) ^2=0,
\] 
where we have used repeatedly that $\langle u,\, u\rangle =1$. 
A straightforward calculation shows that $D_{\infty}$ 
does not contain the points $u$ 
such that $\langle u,\, u\rangle =1$ 
and ${u_1}^2+{u_2}^2=0$, and therefore 
the equation for $D_{\infty}$ is equivalent to the equation 
(\ref{Dinfty}) for $D$. 
\end{proof}

The mapping $u\mapsto y$, $y_i ={x_i} ^2=a_i\, {u_i} ^2$ exhibits 
$U_{\C}$ as an eightfold branched covering of the plane 
\[
P:=\{ y\in\C ^3\mid \sum_{i=1}^3\,\frac{y_i}{a_i}=1\} ,
\]
with branch locus at the three coordinate planes, at  $y_1\, y_2\, y_3=0$. 
We arrive at a 64\--fold branched covering $M\to P\setminus D$, 
where we denoted the image of $D\subset U_{\C}$ in $P$ 
under the mapping from $U_{\C}$ to $P$ with the same symbol $D$.  
 
In order to reduce the order of the covering, we will 
use the group $\Sigma$ introduced in (\ref{Sdef}). 
The group $\Sigma$ leaves the fibers 
of the 64\--fold branched covering 
$M\to P\setminus D$ invariant, and therefore the latter is equal to the composition of the projection $M\to M/\Sigma$ 
and a uniquely determined mapping 
\begin{equation}
\pi _{_{M/\Sigma}}:M/\Sigma\to P\setminus D,
\label{piMSigma}
\end{equation} 
which is a four\--fold branched covering.

\medskip
In the other direction we have the mapping 
$\left(\lambda _2,\,\lambda _3\right)\mapsto y$ defined by Jacobi's elliptic 
coordinates, this is a twofold branched covering with branch locus equal to the image of the diagonal $\lambda _2=\lambda _3$, cf. Subsection \ref{ellcoordsubsec}. If we replace the 
$\left(\lambda _2,\,\lambda _3\right)$\--space $\C ^2$ by the symmetric 
power $\C ^{(2)}$, then we actually obtain a birational isomorphism 
with the inverse of (\ref{xlambda}) given by (\ref{sumlambda}), 
(\ref{prodlambda}). (Don Zagier pointed this out to me.) 

Let $C_{\scriptop{aff}}$ denote the affine (finite) part of the 
hyperelliptic curve $C$. Then the projection 
$(\lambda ,\, z)\mapsto\lambda$ exhibits $C_{\scriptop{aff}}$ as a twofold 
branched covering over $\C$. This leads to the following 
commuting diagram of branched coverings:
\[
\begin{array}{ccc}
C_{\scriptop{aff}}\times C_{\scriptop{aff}}&\stackrel{2:1}{\rightarrow}
&C_{\scriptop{aff}}\times _{\iota '} C_{\scriptop{aff}}\\
~~~\downarrow\mbox{\rm\scriptsize 4:1}&&~~~
\downarrow\mbox{\rm\scriptsize 4:1}\\
\C\times\C&\stackrel{2:1}{\rightarrow}&\C ^{(2)}
\end{array}
\]

Finally we have the mapping $\int: C\times _{\iota '}C\to \op{Jac}(C)$ 
which blows down the diagonal $D$ to the origin $0$ in $\op{Jac}(C)$ and 
defines an isomorphism from $\left( C\times _{\iota '}C\right)\setminus D$ 
onto $\op{Jac}(C)\setminus\{ 0\}$. The image of 
$C_{\scriptop{aff}}\times _{\iota '} C_{\scriptop{aff}}$ in $\op{Jac}(C)$ is 
equal to the complement in $\op{Jac}(C)$ of the image of 
$C$ under the 
mapping $\int_{\infty}:q\mapsto\int_{\infty}^q$ from $C$ to $\op{Jac}(C)$. 
(On the real axis, we 
have that $\infty =\, -\infty$.) 
This mapping is an embedding of $C$ into $\op{Jac}(C)$, cf. 
Farkas and Kra \cite[III.6.4]{FK}. We denote the image of $C$ in 
$\op{Jac}(C)$ under this mapping by $\int_{\infty}^C$. Note that $\int_{\infty}^{\infty}=0$, 
hence $0\in\int_{\infty}^C$. It follows that $\int$ defines an 
isomorphism from 
$\left( C_{\scriptop{aff}}\times _{\iota '}C_{\scriptop{aff}}\right)\setminus D$ 
onto $\op{Jac}(C)\setminus \int_{\infty}^C$. 

The inverse of $\int$, followed by the projection from 
$\left( C_{\scriptop{aff}}\times _{\iota '}C_{\scriptop{aff}}\right)\setminus D$ 
to $\C ^{(2)}\simeq P$, defines a fourfold 
branched covering 
\begin{equation}
\pi_{\op{Jac}(C)}:\op{Jac}(C)\setminus\int_{\infty}^C\to P\setminus D.  
\label{piJC}
\end{equation}

In combination with the fourfold branched covering 
(\ref{piMSigma}), this leads to the following identification 
of $M/\Sigma$ with $\op{Jac}(C)\setminus\int_{\infty}^C$. 

\begin{proposition}
Let $\xi$ be the velocity vector field on 
$M/\Sigma$ and denote by $\xi ^{\scriptop{Jac}(C)}$ the constant vector 
field on $\op{Jac}(C)$ defined by {\em (\ref{formCxC})}. 
Then there is a unique isomorphism $\psi$ from 
$M/\Sigma$ onto $\op{Jac}(C)\setminus\int_{\infty}^C$, such 
that $\pi _{_{\tinyop{Jac}(C)}}\circ\psi =\pi _{_{M/\Sigma}}$, and 
\[
\op{T}_m\pi _{_{M/\Sigma}}\left( \xi _m\right) 
=\op{T}_{\psi (m)}\pi _{_{\tinyop{Jac}(C)}}
\left( \xi ^{\scriptop{Jac}(C)}_{\psi (m)}\right) 
\]
for every $m\in M/\Sigma$. In other words, the diagram 
\[
\begin{array}{ccccc}
M/\Sigma&&\stackrel{\psi}{\longrightarrow}&&
\op{Jac}(C)\setminus\int_{\infty}^C\\
&\pi _{_{M/\Sigma}}\searrow&&~~~\swarrow\pi _{_{\tinyop{Jac}(C)}}&\\
&&P\setminus D&&
\end{array}
\]
commutes and $\psi$ intertwines the vector field $\xi$ on 
$M/\Sigma$ with the constant vector field $\xi ^{\scriptop{Jac}(C)}$ 
on $\op{Jac}(C)$. The mapping $\psi$ intertwines the vector field $\eta$ 
defined by {\em (\ref{Yuz}) --- (\ref{Yzeta})} with a constant 
vector field $\eta ^{\scriptop{Jac}(C)}$ on $\op{Jac}(C)$, 
which up to a sign choice is 
determined by {\em (\ref{b1b0}), (\ref{b0})}.
\label{mapprop}
\end{proposition}
\begin{proof}
For each $p\in P\setminus D$  
there are four (not necessarily distinct) points $m(i)\in M/\Sigma$, 
$1\leq i\leq 4$ and 
$j(k)\in \op{Jac}(C)\setminus\int_{\infty}^C$, $1\leq k\leq 4$, 
such that $\pi _{_{M/\Sigma}}\left( m(i)\right) =p
=\pi _{_{\tinyop{Jac}(C)}}\left( j(k)\right)$. 
The relation between the velocity field in the $(u,\, v )$\--space  
and in Jacobi's elliptic coordinates, as described in 
Subsection \ref{velellsubsec}, yields that the set $V(p)$ of the 
$\op{T}_{m(i)}\pi _{_{M/\Sigma}}\left( \xi _{m(i)}\right)$ is 
equal to the set of the 
$\op{T}_{j(k)}\pi _{_{\tinyop{Jac}(C)}}
\left( \xi ^{\scriptop{Jac}(C)}_{j(k)}\right)$. 
The set $P'$ of all $p\in P\setminus D$ such that $V(p)$ 
consists of four distinct elements, of the form $\pm\xi '$, 
$\pm\xi ''$, is equal to the complement 
of a closed curve in $P\setminus D$. 
The set $M':=\left(\pi _{_{M/\Sigma}}\right) ^{-1}(P')$ and 
$J':=\left(\pi _{_{\tinyop{Jac}(C)}}\right) ^{-1}(P')$ is equal to the 
complement of a closed curve in $M/\Sigma$ and 
$\op{Jac}(C)\setminus\int_{\infty}^C$, 
respectively. The mapping which assigns to 
$m\in M/\Sigma$ the pair 
\[
\left(\pi _{_{M/\Sigma}}(m),\,\op{T}_m\pi _{_{M/\Sigma}}\left( 
\xi _m\right)\right)
\]
is holomorphic from $M/\Sigma$ to the tangent bundle $\op{T}P$ of $P$. 
Its restriction 
$\Pi _{_{M'}}$ to $M'$ is equal to an analytic diffeomorphism 
from $M'$ onto a smooth two\--dimensional submanifold 
$X_P'$ of $\op{T}P$. Similarly 
the mapping which assigns to 
$j\in \op{Jac}(C)\setminus\int_{\infty}^C$ the pair 
\[
\left(\pi _{_{\tinyop{Jac}(C)}}(j),\,
\op{T}_j\pi _{_{\tinyop{Jac}(C)}}\left( 
\xi ^{\scriptop{Jac}(C)}_j\right)\right)
\]
is holomorphic from $\op{Jac}(C)\setminus\int_{\infty}^C$ 
to $\op{T}P$, and its restriction 
$\Pi _{_{J'}}$ to $J'$ is equal to an analytic diffeomorphism 
from $J'$ onto the same $X_P'$. It follows that 
\[
\psi ':=\left(\Pi _{_{J'}}\right) ^{-1}\circ\Pi _{_{M'}}
\]
is an analytic diffeomorphism from $M'$ onto $J'$. 

The branching properties imply that $\psi '$ has a 
continuous, hence analytic extension 
$\psi :M/\Sigma\to \op{Jac}(C)\setminus\int_{\infty}^C$. 
Similarly $\left(\psi '\right) ^{-1}:J'\to M'$ has a 
continuous, hence analytic extension 
$\chi :\op{Jac}(C)\setminus\int_{\infty}^C\to M/\Sigma$. 
Because $\chi\circ\psi$ is equal to the identity on $M'$ 
and $\psi$ and $\chi$ are continuous, we obtain that 
$\chi\circ\psi$ is equal to the identity on $M/\Sigma$. 
Similarly we obtain that $\psi\circ\chi$ is equal to the 
identity on $\op{Jac}(C)\setminus\int_{\infty}^C$. 
\end{proof}

Proposition \ref{mapprop} implies that there exists a completion $\widehat{M/\Sigma}$ of $M/\Sigma$, obtained by adding 
a curve at infinity which is isomorphic to the hyperelliptic curve $C$, 
such that the mapping $\psi$ in Proposition \ref{mapprop} extends to an isomorphism from $\widehat{M/\Sigma}$ onto $\op{Jac}(C)$. The vector fields 
$\xi$ and $\eta$ extend to algebraic vector fields on $\widehat{M/\Sigma}$, 
which commute and are linearly independent at every point of 
$\widehat{M/\Sigma}$. The word ``completion'' is meant in the 
algebraic sense, but it can also be used in the sense that the 
flows of $\xi$ and $\eta$, with complex times, are complete on 
$\widehat{M/\Sigma}$ in the sense that they 
define a transitive action on 
$\widehat{M/\Sigma}$ of the additive group $\C ^2$. 
A completion with this property is unique up to isomorphism. 

\begin{remark}
In Proposition \ref{MSigmasmooth} we had obtained such 
a completion by adding the genus 2 hyperelliptic curve 
$\widehat{M}_{\infty}/\Sigma$ at infinity to 
$M/\Sigma =M/\Sigma$. There the completion is the complex 
torus $\widehat{M}/\Sigma$, in which $\widehat{M}$ 
is a normalization of the projective closure 
$\overline{M}$ of $M$. The isomorphism 
of Proposition \ref{mapprop} leads to an isomorphism   
between $\op{Jac}(C)$ and $\widehat{M}/\Sigma$ 
and an isomorphism between $C$ and the curve 
$\widehat{M}_{\infty}/\Sigma$. 

The complex torus $\widehat{M}$ is a 16\--fold covering 
of the torus $\widehat{M}/\Sigma\simeq\op{Jac}(C)$, which 
can also be obtained by replacing the 
period lattice $\Lambda$ of $\op{Jac}(C)$ by $2\Lambda$. 
See Remark \ref{Oortrem}. 
\label{hatMJacrem}
\end{remark}

\begin{remark}
The six fixed points of the hyperelliptic involution 
$(\lambda ,\, z)\mapsto (\lambda ,\, -z)$ of $C$ correspond to 
$\lambda =1/I_1,\, 1/I_2,\, 1/I_3,\, -1/\rho ,\,\gamma ,\, \infty$, 
where $\gamma$ is given by (\ref{gamma}), (\ref{taudef}). 
Of these, only $\gamma$ varies when we vary the constants of 
motion $T$ and $j$. However, we are working on the assumption that 
$j_3=0$. The transformation of Subsection \ref{redhorsec} to this 
case involved that we had to change the values of both 
$T$ and $\rho$ according to (\ref{tildeT}), and it follows that 
the isomorphism class of the hyperelliptic curve $C$ varies 
in a two\--dimensional subvariety of the three\--dimensional 
moduli space of curves of genus two, as we vary the constants of 
motion $T$ and $j$. If we also vary the moments of inertia 
$I_i$ freely, then there is 
no restriction on the isomorphism class of 
the curve $C$.

The fractional linear transformation $\lambda\mapsto 
\lambda /(1-\rho\,\lambda )$ maps $J_i=1/\left( I_i+\rho\right)$ to 
$1/I_i$, $\infty$ to $-1/\rho$, the zero $\tau =2T/\| j\| ^2$ 
of $p(\lambda )$ to $\gamma$, and the zero $1/\rho$ 
of $p(\lambda )$ to $\infty$. In this description of the zeros of 
$p(\lambda )$ we have used that $j_3=0$. This leads to an explicit verification that the curve $\widehat{M}_{\infty}/\Sigma$ in Remark \ref{modulirem} 
is isomorphic to $C$.  
\label{Cmodulirem}
\end{remark}

\begin{question}
Can the constant vector fields on 
$\op{Jac}(C)$ be determined more easily than in Subsection 
\ref{velellsubsec} by means of calculations at some special points, 
for instance at points corresponding to one of the fixed points of 
the hyperelliptic involution $\iota$? 
\label{AMq}
\end{question}

\begin{question} (Richard Cushman)  
Could the isomorphism in Proposition \ref{mapprop}
be obtained in a similar way as in Mumford 
\cite[p. 3.57, 3.58]{TthII} for Neumann's system? 
This question is suggested by the form 
(\ref{phipsichi}) of the kinetic energy equation. 
\label{cushmanq}
\end{question}

\subsection{The Translational Motion}
\label{phorsubsec}
In order to obtain the motion of the point of contact $p(t)$, 
we have to integrate the right hand side of (\ref{pdot}). 
Because this vector is horizontal, and the moment of momentum $j$ 
around the point of contact is assumed to be horizontal as well, 
it is sufficient to determine the inner product with $j$ and 
$j\times e_3$. For the latter one we have that 
\[
\langle (A\,\omega )\times e_3,\, j\times e_3\rangle 
=\langle A\,\omega ,\, j\rangle =\langle\omega ,\, A^{-1}\, j\rangle 
=\langle\omega ,\, v\rangle =2T,
\]
where we have used (\ref{sigmadef}) and (\ref{Tomegasigma}). 
It follows therefore that {\em the component of $p(t)$ 
orthogonal to $j$ grows linearly in time:}
\begin{equation}
\langle p(t),\, j\times e_3\rangle =\langle p(0),\, j\times e_3\rangle 
+2r\, T\, t. 
\label{pe3j}
\end{equation}

In order to determine the inner product of $\dd{p}{\tau}$ with $j$ 
in terms of Jacobi's elliptic coordinates, we begin with 
recalling the formula (\ref{detC}), in which 
$\dd{}{\tau}\langle p(\tau ),\, j\rangle$ is expressed in terms 
of the determinant of the vectors $u$, $\omega$, and $I\,\omega$. 
Under a transformation $S\in\Sigma$ as in 
(\ref{Sdef}), the right hand side of (\ref{detC}) gets multiplied by 
$\epsilon _2$. Therefore only the square of (\ref{detC}) is 
a single\--valued function on $M/\Sigma$, whereas 
(\ref{detC}) is single\--valued on the unbranched (= unramified)
double covering $M/\Sigma _0$ of $M/\Sigma$, where 
$\Sigma _0$ denotes the group of $S$ in (\ref{Sdef}) 
such that $\epsilon _2=0$. In view of the isomorphism 
of $M/\Sigma$ with $\op{Jac}(C)/\int_{\infty}^C$ in Proposition 
\ref{mapprop}, 
only $\langle\dd{p}{\tau},\, j\rangle ^2$ can be a single 
valued function of Jacobi's elliptic coordinates.

If $C$ denotes the matrix $\left( u,\,\omega ,\, I\,\omega\right)$, then 
\[
(\op{det}C)^2=\op{det}(C^*\circ C)=\op{det}
\left(\begin{array}{ccc}
\langle u,\, u\rangle &\langle u,\,\omega\rangle &
\langle u,\, I\,\omega\rangle\\
\langle \omega ,\, u\rangle &\langle\omega ,\,\omega\rangle &
\langle\omega ,\, I\,\omega\rangle\\
\langle I\,\omega ,\, u\rangle &\langle I\,\omega ,\,\omega\rangle &
\langle I\,\omega ,\, I\,\omega\rangle  
\end{array}\right) .
\]
Using that $\langle u,\, u\rangle =1$ and 
$\langle I\,\omega ,\, u\rangle =0$,\, cf. (\ref{Iomegau}), 
we therefore obtain that 
\[
(\op{det}C)^2=\left(\langle\omega ,\,\omega\rangle 
-\langle u,\,\omega\rangle ^2\right)
\cdot\langle I\,\omega ,\, I\,\omega\rangle
-\langle\omega ,\, I\,\omega\rangle ^2, 
\]
in which we can substitute 
\[
\langle\omega ,\,\omega\rangle -\langle u,\,\omega\rangle ^2
=\langle u\times\omega ,\, u\times\omega\rangle 
=\langle\dot{u},\,\dot{u}\rangle ,
\]
cf. (\ref{udot}). Substituting (\ref{dotudotu}) in 
combination with (\ref{Xi2ai}) with $X=\dot{x}$, 
(\ref{4IomIom}), and (\ref{Tlambdadot}), we obtain that 
\begin{eqnarray*}
16(\op{det}C)^2&=&\frac{\left(\lambda _2-\lambda _3\right) ^2}
{{\lambda _2}^2\, {\lambda _3}^2}\,\left\{ 
\left(
\frac{{\dot{\lambda}_2}^2}{S\left(\lambda _2\right)}
-\frac{{\dot{\lambda}_3}^2}{S\left(\lambda _3\right)}
\right)\,\left(
\frac{{\lambda _2}^2\,{\dot{\lambda}_2}^2}{S\left(\lambda _2\right)}
-\frac{{\lambda _3}^2\, {\dot{\lambda}_2}^2}{S\left(\lambda _2\right)}
\right)\right.\\
&&\left. -\left(
\frac{\lambda _2\,{\dot{\lambda}_2}^2}{S\left(\lambda _2\right)}
-\frac{\lambda _3\,{\dot{\lambda}_3}^2}{S\left(\lambda _3\right)}
\right) ^2\right\}\\
&=&\frac{\left(\lambda _2-\lambda _3\right) ^2}
{{\lambda _2}^2\, {\lambda _3}^2}\,
\frac{\left( -{\lambda _3}^2-{\lambda _2}^2-2\lambda _2\,\lambda _3\right)
\, {\dot{\lambda}_2}^2\, {\dot{\lambda}_3}^2}
{S\left(\lambda _2\right)\, S\left(\lambda _3\right)}
=\, -\frac{\left(\lambda _2-\lambda _3\right) ^4\, 
{\dot{\lambda}_2}^2\, {\dot{\lambda}_3}^2}
{{\lambda _2}^2\, {\lambda _3}^2\, 
S\left(\lambda _2\right)\, S\left(\lambda _3\right)} .
\end{eqnarray*}
Now we can write (\ref{lambda2dot}) in the form 
\[
-\frac{\left(\lambda _2-\lambda _3\right) ^2\, {\dot{\lambda}_2}^2}
{{\lambda _3}^2\, S\left(\lambda _2\right)}
=4\frac{\left(\| j\| ^2-2\rho\, T\right)\,\left(\gamma -\lambda _2\right)}
{\rho\,\left( -1/\rho -\lambda _3\right)} ,
\]
and we have a similar expression for $\dot{\lambda} _3$, 
obtained by interchanging the indices 2 and 3. 
Substituting these, we obtain from (\ref{detC}) that 
\[
\langle\dd{p}{\tau},\, j\rangle ^2=r^2\, X (u)\, (\op{det}C)^2
=\, -\frac{r^2\, X (u)\, \left(\| j\| ^2-2\rho\, T\right) ^2\, 
\left(\gamma -\lambda _2\right)\,\left(\gamma -\lambda _3\right)}
{\rho ^2\,\left( 1/\rho +\lambda _3\right)
\,\left( 1/\rho +\lambda _2\right)}.
\]
In combination with (\ref{philambda}), this leads to 
\begin{equation}
\left(\dd{\langle p(\tau ),\, j\rangle}{\tau}\right) ^2
=\, -d^2\,\left(\gamma -\lambda _2\right)\,\left(\gamma -\lambda _3\right) ,
\label{dpjdtau}
\end{equation}
with 
\begin{equation}
d^2:=\frac{r^2\,\left(\| j\| ^2-2\rho\, T\right) ^2}
{\rho ^4\,\prod_{i=1}^3\,\left( a_i+1/\rho\right)} .
\label{d2def}
\end{equation}

\begin{remark}
It follows from (\ref{Squadr}) with $\lambda =\gamma$ and 
(\ref{aiIi}) that 
\begin{equation}
-\left(\gamma -\lambda _2\right)\,\left(\gamma -\lambda _3\right) =
\frac{S(\gamma )}{\gamma}\,\left( 1-\sum_{i=1}^3\,
\frac{{x_i}^2}{a_i-\gamma}\right) =
\frac{S(\gamma )}{\gamma}\,\left( 1-\langle u,\, (1-\gamma\, I)^{-1}\, u
\rangle\right) .
\label{gammabranch}
\end{equation}
The right hand side in (\ref{gammabranch}) is equal to zero 
if and only if $u\in U_{\C}$ is in the branch locus for the 
branched covering $(u,\, v )\mapsto u$ from $L$ to $U_{\C}$, 
cf. (\ref{gamma}). 

In fact, we could have concluded this at an earlier stage. 
The investigation of the branch locus 
started with the observation that this is the set of points 
where $u$, $\omega$ and $I\, \omega$ are linearly dependent, 
cf. (\ref{partialf3}), and according to (\ref{detC}), this 
is equal to the condition that $\op{d}\langle p(t),\, j\rangle /\op{d}\! t=0$. 

In terms of the description of the real points $u(t)$ on the unit 
sphere $U$ in Remark \ref{uremark}, the instances when $u(t)$ 
reaches the boundary curves of the annulus in $U$ coincide with 
the instances that $\langle p(t),\, j\rangle$ has a turning point. 
In view of the quasiperiodic 
motion of the rotational motion, this happens infinitely often.  
Keeping in mind that the $j\times e_3$\--component of the velocity of 
$p(t)$ is equal to a positive constant, we obtain that the point 
of contact $p(t)$ performs a swaying motion in the direction 
of $j\times e_3$. According to Corollary \ref{pjqpcor}, the 
function $\tau\mapsto\langle p(\tau ),\, j\rangle$ is 
quasiperiodic if the rotational motion is not periodic and the 
irrational ratio $\nu _1(\epsilon )/\nu _2(\epsilon )$ 
mentioned after (\ref{Fourierp}) 
is sufficiently slowly approximated by rational numbers.  
\label{turnrem}
\end{remark}

As a function on the Jacobi variety $\op{Jac}(C)$ of 
the hyperelliptic curve $C$, the function 
\[
f:=\,-\left(\gamma -\lambda _2\right)\,\left(\gamma -\lambda _3\right)
\] 
is rational, with poles along $\int_{\infty}^C$, 
zeros at $\int_{(\gamma ,\, 0)}^C$ and undetermined 
at the two points of  
\begin{equation}
\int_{\infty}^C\cap\int_{(\gamma ,\, 0)}^C
=\{ 0,\,\int_{\infty}^{(\gamma ,\, 0)}\} . 
\label{undetermined}
\end{equation}
Note that $\int_{\infty}^{(\gamma ,\, 0)}$ is a point of order two 
in $\op{Jac}(C)$, in the sense that $2\int_{\infty}^{(\gamma ,\, 0)}=0$. 

$f$ has double zeros along $\int_{(\gamma ,\, 0)}^C$, 
because when $\lambda$ is near $\gamma$, then 
on the hyperelliptic curve 
\[
z^2=P(\lambda )=\prod_{i=1}^3\,\left( a_i-\lambda \right)\, 
\left( -1/\rho -\lambda\right)\, \left(\gamma -\lambda \right)
\]
the variable $z$ is a local parameter, in terms 
of which $\gamma -\lambda$ is of order $y^2$. 
Also $f$ has double poles along $\int_{\infty}^C$, 
because if $\theta$ is a local parameter for $C$ at infinity, 
then $\lambda =\theta ^{-2}\, u$, $z=\theta ^{-5}\, v$, where $u$ 
and $v$ are units, cf. Shafarevich \cite[III.5.5]{Shaf}. 
Therefore, following a choice of a square root of $f$ 
along curves in $\op{Jac}(C)$, we find that there exists a rational function 
$g$ such that $f=g^2$, either on 
$\op{Jac}(C)$, or on a double covering $\widetilde{\op{Jac}(C)}$, 
unbranched, of $\op{Jac}(C)$. Ben Moonen told me that 
$g$ cannot be single\--valued on $\op{Jac}(C)$. Actually, if in the choice 
of the basis of $\Lambda (C)$ corresponding to 
the closed curves $A_1,\, A_2,\, B_1,\, B_2$ as in Mumford 
\cite[p. 3.76]{TthII}, with his point $a_1$ in the role of 
$(\gamma ,\, 0)$, then $g$ turns into its opposite if we 
follow $A_1$ around, and remains unchanged if we follow 
any of the other curves. Therefore $g$ is a single valued 
rational function on the unbranched double covering 
\[
\widetilde{\op{Jac}(C)}={\cal H}^1(C)^*/\Lambda _0(C),
\]
where $\Lambda _0(C)$ 
is the sublattice of index two in $\Lambda (C)$ which is generated 
by the basis elements corresponding to $2A_1,\, A_2,\, B_1,\, B_2$. 

It follows that $g$ is a rational function on 
$\widetilde{\op{Jac}(C)}$ with simple poles along the preimage of $\int_{\infty}^C$ 
under the double covering $:\widetilde{\op{Jac}(C)}\to \op{Jac}(C)$, 
simple zeros along the preimage of $\int_{(\gamma ,\, 0)}$, and undetermined 
at the preimage of (\ref{undetermined}). 
These properties characterize the function $g$ 
on $\widetilde{\op{Jac}(C)}$ up to a constant factor. 
The function $g$ can be identified with 
a quotient of theta functions as in 
Mumford \cite[p. 3.80, 81]{TthII}, \cite[Ch. II]{TthI}.  

\begin{question}
We know that 
\begin{equation}
\dd{\langle p(\tau ),\, j\rangle}{\tau}=
r\, X(u(\tau ))^{1/2}\,\langle A(\tau )\,\omega (\tau ),\, e_3\times j\rangle 
\label{pjtauA}
\end{equation}
is a quasiperiodic function of $\tau$, because the rotational 
motion is a quasiperiodic function of $\tau$. We now have the 
additional information that the derivative 
of $\langle p(\tau ),\, j\rangle$ is given by $d\, g$ 
along a straight line path in $\widetilde{\op{Jac}(C)}$, where $g$ is the 
quotient of theta functions described above. 
Can this additional information 
be used in order to decide whether 
the problems with the integration (with respect to $\tau$) 
of the quasiperiodic function of $\tau$ in the right 
hand side of (\ref{pjtauA}), which are mentioned 
in Subsection \ref{Fouriersubsec}, really occur? 
\label{translthetaq}
\end{question}
 
\subsection{Chaplygin}
The last part of Chaplygin \cite[\S 5]{chaplsphere}, 
starting with ``Now we discuss the curves traced out on the surface of the sphere \ldots'', corresponds to our 
Remark \ref{uremark}. 

The proof that the elliptic coordinates (\ref{xlambda})
for $u$ lead to a motion on the Jacobi variety of the hyperelliptic curve 
(\ref{Cdef}) with constant velocity is contained in Chaplygin 
\cite[\S 3 up to (30)]{chaplsphere}. 
The notations in Chaplygin \cite[\S 3]{chaplsphere} correspond to ours 
according to the following list, which is a continuation of the list in 
Subsection \ref{chaplremeqmot}. 
\[
\begin{array}{cc}
\mbox{\rm Chaplygin 
\cite[\S 3 up to (30)]{chaplsphere}} & \mbox{\rm our notation}\\
\mbox{\rm (14)} & \mbox{\rm (\ref{TI})}\\
\mbox{\rm (15)} & \mbox{\rm (\ref{SO3eq}), (\ref{Iomegau}) 
and (\ref{omegausigma})}\\
\mbox{\rm (16) and (17)} &\mbox{\rm (\ref{udotIomega})}\\
\mbox{\rm (18)} & \mbox{\rm (\ref{2Tomegau}) and (\ref{udotIomega})}\\
\mbox{\rm (19)} & \mbox{\rm formulas following (\ref{Tlambdadot})}\\
a^2,\, b^2,\, c^2, \partial ^2,\, g,\, k,\,\sigma\;\mbox{\rm in (20)}
& a_1,\, a_2,\, a_3\;\mbox{\rm (cf. (\ref{aiIi}))},\, 1/\rho ,\, 
8T,\, 4\| j\| ^2,\, 1\\ 
x,\, y,\, z\;\mbox{\rm in (21)} & x_1,\, x_2,\, x_3\;\mbox{\rm in (\ref{aiIi})}\\
\mbox{\rm (22)} & \mbox{\rm (\ref{xi2/ai}), (\ref{2Tomegau}), (\ref{Tudot}), 
(\ref{dotudotu}), (\ref{uIu})}\\
\mbox{\rm (23)} & \mbox{\rm formulas preceding (\ref{j2-2rhoT})}\\
\mbox{\rm (24)} & \mbox{\rm (\ref{deltaxsquared}) and (\ref{Xi2ai})}\\
\mbox{\rm (25)} & \mbox{\rm (\ref{S})}\\
\mbox{\rm (26)} & \mbox{\rm (\ref{chaplsum})}\\
\mbox{\rm formulas between (26) and (27)} & \mbox{\rm (\ref{prodlambda}), 
(\ref{Tlambdadot}), (\ref{j2-2rhoT})}\\
\mbox{\rm (27)},\; j,\; -g\,\partial ^2/j & 
\mbox{\rm (\ref{lambda2dot})},\; 8T-4\| j\| ^2/\rho,\;\gamma\\
\mbox{\rm (28)} & \mbox{\rm (\ref{philambda})}\\
\mbox{\rm (29)} & \mbox{\rm (\ref{formCxCint}) }\\
\mbox{\rm (30)} & \mbox{\rm (\ref{Pdef})}\\
\mbox{\rm (33)} & \mbox{\rm (\ref{pe3j})}\\
\mbox{\rm (34)} & \mbox{\rm (\ref{dpjdtau})}
\end{array}
\]
In our Subsection \ref{ellcoordsubsec} 
we have followed the miraculous calculations of Chaplygin 
\cite[\S 3 up to (30)]{chaplsphere} quite closely, adding some more explanations 
in the hope to make these easier to read. In (\ref{philambda}) we identified 
(up to a constant factor) the factor 
\[
\sqrt{\left(\lambda _2+1/\rho\right)\,
\left(\lambda _3+1/\rho\right)}
\]
introduced in Chaplygin \cite[(28)]{chaplsphere} with the 
integrating factor $ X (u)^{1/2}$ of  Lemma \ref{vollem},  
the same as the factor $\sqrt{X}$ at the end of  \cite[\S 2]{chaplsphere}, 
and of  Corollary \ref{commvvcor}. 
In Subsection \ref{jacsubsec} we added a discussion of the relation 
between the phase space of the rotational motion and the 
Jacobi variety of the hyperelliptic curve, about which Chaplygin did not say anything in \cite{chaplsphere}. 

Chaplygin did not tell how he came to the idea of using 
the elliptic coordinates (\ref{xlambda}). The last part of Chaplygin 
\cite[\S 5]{chaplsphere} indicates that 
he had calculated the branch locus of the projection 
$(u,\, v )\mapsto u$ from the $(j,\, T)$\--level surface 
onto the $u$\--sphere. Therefore he might have observed that 
in the complex domain it is equal to the union of two quadrics 
which, together with the sphere, belong to a 
one\--parameter family of confocal quadrics, and 
this might have prompted him to use the elliptic coordinates 
(\ref{xlambda}). 
He might have refrained from mentioning this in his article, 
because in his time the 
use of elliptic coordinates in the presence of families of 
confocal quadrics was standard. 

Although Chaplygin described the motion of $u(t)$ in the 
annulus on the sphere in \cite[end of \S 5]{chaplsphere}, 
he did not observe that $u(t)$ reaches the boundary 
curves of the annulus precisely when 
$\op{d}\langle p(t),\, j\rangle/\op{d}\! t=0$, cf. Remark \ref{turnrem}.

\section{A Geometric Interpretation}
\label{geomsec}
The kinetic energy equation in the form $Y^2-X\, Z=0$, 
cf. (\ref{phipsichi}),  
is equal to the discriminant equation for the quadratic equation 
$X(u)\,\lambda ^2+2Y(u,\, v)\,\lambda +Z(v) =0$ in the variable $\lambda$. 
Let $\alpha$ be an auxiliary parameter. 
Using 
(\ref{SO3eq}), the equation for $\lambda$ can be written in the form 
\begin{equation}
\langle\lambda\, u- v ,\, 
(J+\alpha )\, (\lambda\, u- v )\rangle
=(\rho ^{-1}+\alpha )\,\lambda ^2-2\alpha\, j_3\,\lambda +2T +\alpha\,\| j\| ^2.
\label{alphalambda}
\end{equation}
The discriminant of the right hand side is equal to zero if and only if 
\begin{equation}
\left(\| j\| ^2-{j_3}^2\right)\,\alpha ^2+
\left(2T+\| j\| ^2/\rho\right)\,\alpha +2T/\rho =0.
\label{alphadiscr}
\end{equation}
If (\ref{alphadiscr}) holds, then the right hand side of 
(\ref{alphalambda}) is equal to 
\[
\left(\rho ^{-1}+\alpha\right)\,
\left\{\lambda -\alpha\, j_3/\left(\rho ^{-1}+\alpha\right)\right\} ^2,
\]
and the equation (\ref{alphalambda}) is equivalent to 
\begin{equation}
\langle v,\, \left[\rho\, (I+\rho )^{-1}+\alpha\right]\, v\rangle 
=\rho ^{-1}+\alpha , 
\label{alphaJ}
\end{equation}
with $v=\theta\, u-\eta\, v$ and 
\[ 
\theta =\lambda/\left\{\lambda -\alpha\, j_3/
\left(\rho ^{-1}+\alpha\right)\right\}, \quad 
\eta =1/\left\{\lambda -\alpha\, j_3/
\left(\rho ^{-1}+\alpha\right)\right\} . 
\]
It follows that $\eta = (\theta -1)\, 
\left(\rho ^{-1}+\alpha\right) /\alpha\, j_3$, 
and the conclusion is that the straight line $\widetilde{l}$ 
passing through $u$,  
with the direction vector equal to $\alpha\, j_3\, u-
\left(\rho ^{-1}+\alpha\right)\, v$,  
is tangent to the quadric $Q$ defined by the equation (\ref{alphaJ}). 
Note that $\widetilde{l}=A^{-1}(l)$, where $l$ is equal to 
the straight line passing through $-e_3$, with direction vector 
equal to $\alpha\, j_3\, e_3-\left(\rho ^{-1}+\alpha\right)\, j$. Here the 
rotation $A\in\op{SO}(3)$ varies in the level surface in 
$\op{SO}(3)$, corresponding to the equations 
(\ref{SO3eq}) and (\ref{phipsichi}). Also note that 
$-Q=Q$, and therefore the same properties hold with 
$l$ replaced by $-l$.   

The equation (\ref{alphadiscr}) has two solutions $\alpha _1$ and 
$\alpha _2$, leading to two pairs of straight lines $\pm l_1$ and 
$\pm l_2$ and 
quadrics $Q_1$ and $Q_2$ to which $A^{-1}\left(\pm l_1\right)$ 
and $A^{-1}\left(\pm l_2\right)$ are tangent, respectively. 
The inner product of the two direction vectors is equal to 
\begin{eqnarray*}
&&\langle\alpha _1\, j_3\, u-\left(\rho ^{-1}+\alpha _1\right)\, v ,\, 
\alpha _2\, j_3\, u-\left(\rho ^{-1}+\alpha _2\right)\, v \rangle \\
&&=\left(\alpha _1/\rho +\alpha _2/\rho +\alpha _1\,\alpha _2\right) 
\,\left(\| j\| ^2-{j_3}^2\right) +\| j\| ^2/\rho ^2.
\end{eqnarray*}
On the other hand it follows from (\ref{alphadiscr}) 
that 
\[
\left(\| j\| ^2-{j_3}^2\right)\,\left(\alpha _1+\alpha _2\right) 
=\, -2T -\| j\| ^2/\rho ,\quad 
\left(\| j\| ^2-{j_3}^2\right)\,\alpha _1\,\alpha _2 
=2T/\rho ,
\]
and we conclude that the direction vectors are perpendicular. 
Equivalently, the four straight lines 
$\pm l_1$, $\pm l_2$ form a rectangle. 

\medskip\noindent
The above arguments work under the assumption that the moment $j$ 
is neither vertical nor horizontal. If $j\neq 0$ is horizontal,  
$j_3=0$, then the right hand side of (\ref{alphalambda}) 
is a square if $\alpha =\, -2T\,\rho/\| j\| ^2$, and a constant 
if $\alpha =\, -1/\rho$. In the first case the line passing through 
$u$ with direction vector $ v$ is tangent to the 
quadric defined by (\ref{alphaJ}), with $\alpha =\, -2T\,\rho/\| j\| ^2$. 
In the second case the line passing through $- v$ 
with the direction vector $u$ is tangent to the 
quadric 
$\langle v,\,\left[(I+\rho )^{-1}-\rho ^{-1}\right]\, v\rangle 
=2T-\| j\| ^2/\rho$. 

\begin{question}
What are {\em all} the straight lines $l$ in the 
plane spanned by $e_3$ and $j$, and quadrics $Q$, such that 
for each $A\in\op{SO}(3)$ in the level surface 
corresponding to the equations 
(\ref{SO3eq}) and (\ref{phipsichi}) we have that 
$A^{-1}(l)$ is tangent to $Q$? This question 
may be related to Question \ref{cushmanq}. 
\end{question}

\subsection{Chaplygin}
Section \ref{geomsec} corresponds to
Chaplygin's \cite[\S 5]{chaplsphere}. 
Chaplygin multiplied the figures by the radius $r$ of the sphere, 
in order to have the corner point $-r\, e_3$ of the rectangle 
$\pm l_1$, $\pm l_2$ attached to the 
point of contact of the sphere with the plane.

\end{document}